\let\noarrow = t

\input eplain


\magnification=\mag

\topskip24pt
\def\pagewidth#1{
  \hsize=#1
}

\def\pageheight#1{
  \vsize=#1
}

\pageheight{50pc} \pagewidth{30pc}

\abovedisplayskip=3mm \belowdisplayskip=3mm
\abovedisplayshortskip=0mm \belowdisplayshortskip=2mm
\def\spacing{{\smallskip}}
\frenchspacing
\parindent1pc

\voffset=0pc
\hoffset=5pc


\newdimen\abstractmargin
\abstractmargin=3pc


\newdimen\footnotemargin
\footnotemargin=1pc


\font\eightrm=cmr8 \relax 
\font\sixrm=cmr6 \relax 
\font\eighti=cmmi8 \relax     \skewchar\eighti='177 
\font\sixi=cmmi6 \relax       \skewchar\sixi='177   
\font\eightsy=cmsy8 \relax    \skewchar\eightsy='60 
\font\sixsy=cmsy6 \relax      \skewchar\sixsy='60   
\font\eightbf=cmbx8 \relax 
\font\sixbf=cmbx6 \relax   
\font\eightit=cmti8 \relax 
\font\eightsl=cmsl8 \relax 
\font\eighttt=cmtt8 \relax 

\catcode`\@=11
\newskip\ttglue

\def\eightpoint{\def\rm{\fam0\eightrm}%
 \textfont0=\eightrm \scriptfont0=\sixrm
 \scriptscriptfont0=\fiverm
 \textfont1=\eighti \scriptfont1=\sixi
 \scriptscriptfont0=\fivei
 \textfont2=\eightsy \scriptfont2=\sixsy
 \scriptscriptfont2=\fivesy
 \textfont3=\tenex \scriptfont3=\tenex
 \scriptscriptfont3=\tenex
 \textfont\itfam\eightit \def\it{\fam\itfam\eightit}%
 \textfont\slfam\eightsl \def\sl{\fam\slfam\eightsl}%
 \textfont\ttfam\eighttt \def\tt{\fam\ttfam\eighttt}%
 \textfont\bffam\eightbf \scriptfont\bffam\sixbf
   \scriptscriptfont\bffam\fivebf \def\bf{\fam\bffam\eightit}%
 \tt \ttglue=.5em plus.25em minus.15em
 \normalbaselineskip=9pt
 \setbox\strutbox\hbox{\vrule height7pt depth3pt width0pt}%
 \let\sc=\sixrm \let\big=\eifgtbig \normalbaselines\rm}


 \font\titlefont=cmbx12 scaled\magstep1
 \font\sectionfont=cmbx12
 \font\ssectionfont=cmsl10
 \font\claimfont=cmsl10

 \font\normalfont=cmr10

\catcode`\@=11 \font\teneusm=eusm10 
\font\seveneusm=eusm7  \font\fiveeusm=eusm5
\newfam\eusmfam \textfont\eusmfam=\teneusm
\scriptfont\eusmfam=\seveneusm \scriptscriptfont\eusmfam=\fiveeusm
\def\hexnumber@#1{\ifcase#1
0\or1\or2\or3\or4\or5\or6\or7\or8\or9\or         A\or B\or C\or
D\or E\or F\fi } \edef\eusm@{\hexnumber@\eusmfam}
\def\euscr{\ifmmode\let\next\euscr@\else
\def\next{\errmessage{Use \string\euscr\space only in math mode}}\fi\next}
\def\euscr@#1{{\euscr@@{#1}}} \def\euscr@@#1{\fam\eusmfam#1} \catcode`\@=12

\catcode`\@=11 \font\teneuex=euex10 
 \font\seveneuex=euex7  \newfam\euexfam
\textfont\euexfam=\teneuex  \scriptfont\euexfam=\seveneuex
 \def\hexnumber@#1{\ifcase#1
0\or1\or2\or3\or4\or5\or6\or7\or8\or9\or         A\or B\or C\or
D\or E\or F\fi } \edef\euex@{\hexnumber@\euexfam}
\def\euscrex{\ifmmode\let\next\euscrex@\else
\def\next{\errmessage{Use \string\euscrex\space only in math mode}}\fi\next}
\def\euscrex@#1{{\euscrex@@{#1}}} \def\euscrex@@#1{\fam\euexfam#1}
\catcode`\@=12

\catcode`\@=11 \font\teneufb=eufb10 
\font\seveneufb=eufb7  \font\fiveeufb=eufb5
\newfam\eufbfam \textfont\eufbfam=\teneufb
\scriptfont\eufbfam=\seveneufb \scriptscriptfont\eufbfam=\fiveeufb
\def\hexnumber@#1{\ifcase#1
0\or1\or2\or3\or4\or5\or6\or7\or8\or9\or         A\or B\or C\or
D\or E\or F\fi } \edef\eufb@{\hexnumber@\eufbfam}
\def\euscrfb{\ifmmode\let\next\euscrfb@\else
\def\next{\errmessage{Use \string\euscrfb\space only in math mode}}\fi\next}
\def\euscrfb@#1{{\euscrfb@@{#1}}} \def\euscrfb@@#1{\fam\eufbfam#1}
\catcode`\@=12

\catcode`\@=11 \font\teneufm=eufm10 
\font\seveneufm=eufm7  \font\fiveeufm=eufm5
\newfam\eufmfam \textfont\eufmfam=\teneufm
\scriptfont\eufmfam=\seveneufm \scriptscriptfont\eufmfam=\fiveeufm
\def\hexnumber@#1{\ifcase#1
0\or1\or2\or3\or4\or5\or6\or7\or8\or9\or         A\or B\or C\or
D\or E\or F\fi } \edef\eufm@{\hexnumber@\eufmfam}
\def\euscrfm{\ifmmode\let\next\euscrfm@\else
\def\next{\errmessage{Use \string\euscrfm\space only in math mode}}\fi\next}
\def\euscrfm@#1{{\euscrfm@@{#1}}} \def\euscrfm@@#1{\fam\eufmfam#1}
\catcode`\@=12

\catcode`\@=11 \font\teneusb=eusb10 
\font\seveneusb=eusb7  \font\fiveeusb=eusb5
\newfam\eusbfam \textfont\eusbfam=\teneusb
\scriptfont\eusbfam=\seveneusb \scriptscriptfont\eusbfam=\fiveeusb
\def\hexnumber@#1{\ifcase#1
0\or1\or2\or3\or4\or5\or6\or7\or8\or9\or         A\or B\or C\or
D\or E\or F\fi } \edef\eusb@{\hexnumber@\eusbfam}
\def\euscrsb{\ifmmode\let\next\euscrsb@\else
\def\next{\errmessage{Use \string\euscrsb\space only in math mode}}\fi\next}
\def\euscrsb@#1{{\euscrsb@@{#1}}} \def\euscrsb@@#1{\fam\eusbfam#1}
\catcode`\@=12

\catcode`\@=11 \font\tenmsa=msam10 
\font\sevenmsa=msam7  \font\fivemsa=msam5
\font\tenmsb=msbm10  \font\sevenmsb=msbm7
 \font\fivemsb=msbm5 \newfam\msafam
\newfam\msbfam \textfont\msafam=\tenmsa
\scriptfont\msafam=\sevenmsa
  \scriptscriptfont\msafam=\fivemsa
\textfont\msbfam=\tenmsb  \scriptfont\msbfam=\sevenmsb
  \scriptscriptfont\msbfam=\fivemsb
\def\hexnumber@#1{\ifcase#1 0\or1\or2\or3\or4\or5\or6\or7\or8\or9\or
        A\or B\or C\or D\or E\or F\fi }
\edef\msa@{\hexnumber@\msafam} \edef\msb@{\hexnumber@\msbfam}
\mathchardef\square="0\msa@03 \mathchardef\subsetneq="3\msb@28
\mathchardef\supsetneq="3\msb@29 \mathchardef\ltimes="2\msb@6E
\mathchardef\rtimes="2\msb@6F \mathchardef\dabar="0\msa@39
\mathchardef\daright="0\msa@4B \mathchardef\daleft="0\msa@4C

\def\Bbb{\ifmmode\let\next\Bbb@\else
        \def\next{\errmessage{Use \string\Bbb\space only in math mode}}\fi\next}
\def\Bbb@#1{{\Bbb@@{#1}}}
\def\Bbb@@#1{\fam\msbfam#1}
\catcode`\@=12



\newcount\senu
\def\senum{\number\senu}
\newcount\ssnu
\def\ssnum{\number\ssnu}
\newcount\fonu
\def\fonum{\number\fonu}

\def\num{{\senum.\ssnum}}
\def\numfo{{\senum.\ssnum.\fonum}}


\outer\def\section#1\par{\vskip0pt
  plus.3\vsize\penalty-20\vskip0pt
  plus-.3\vsize\bigskip\vskip\parskip
  \message{#1}\centerline{\sectionfont\senum\enspace#1.}
  \nobreak\smallskip}

\def\endsection{\advance\senu by1\penalty-20\smallskip\ssnu=1}
\outer\def\ssection#1\par{\bigskip
  \message{#1}{\noindent\bf\num\ssectionfont\enspace#1.\thinspace}
  \nobreak\normalfont}

\def\endssection{\advance\ssnu by1\smallskip\ifdim\lastskip<\medskipamount
\removelastskip\penalty55\medskip\fi\fonu=1\normalfont}

\def\proclaim #1\par{\bigskip
  \message{#1}{\noindent\bf\num\enspace#1.\thinspace}
  \nobreak\claimfont}

\def\cor{\proclaim Corollary\par}
\def\defi{\proclaim Definition\par}
\def\lemma{\proclaim Lemma\par}
\def\prop{\proclaim Proposition\par}
\def\rmk{\proclaim Remark\par\normalfont}
\def\thm{\proclaim Theorem\par}

\def\endcor{\endssection}
\def\enddefi{\endssection}
\def\endlemma{\endssection}
\def\endprop{\endssection}
\def\endrmk{\endssection}
\def\endthm{\endssection}

\def\Proof{{\noindent\sl Proof: \/}}


\def\maplefto#1{\ \smash{\mathop{\longleftarrow}\limits^{#1}}\ }

\def\llongrightarrow{\relbar\joinrel\relbar\joinrel\rightarrow}
\def\lllongrightarrow{\hbox to 40pt{\rightarrowfill}}

\def\hooklongrightarrow{\lhook\joinrel\longrightarrow}
\def\twoheadrightarrow{\rightarrow\kern -8pt\rightarrow}

\def\maprighto#1{\smash{\mathop{\longrightarrow}\limits^{#1}}}

\def\mapdownr#1{\Big\downarrow\rlap{$\vcenter{\hbox{$\scriptstyle#1$}}$}}
\def\mapdownl#1{\llap{$\vcenter{\hbox{$\scriptstyle#1$}}$}\Big\downarrow}

\def\llongmaprighto#1{\ \smash{\mathop{\llongrightarrow}\limits^{#1}}\ }

\def\lllongmaprighto#1{\ \smash{\mathop{\lllongrightarrow}\limits^{#1}}\ }

\def\llongleftarrow{\leftarrow\joinrel\relbar\joinrel\relbar}

\def\longleftmapsto{\longleftarrow\kern-2pt\mapstochar\;}

\def\llongmapsto{\,\vert\kern-3.2pt\joinrel\longrightarrow\,}
\def\llongmapsto{\,\vert\kern-3.7pt\joinrel\llongrightarrow\,}
\def\lllongmapsto{\,\vert\kern-5.5pt\joinrel\lllongrightarrow\,}

\def\isomarrow{\maprighto{\lower3pt\hbox{$\scriptstyle\sim$}}}
\def\llongisomarrow{\llongmaprighto{\lower3pt\hbox{$\scriptstyle\sim$}}}
\def\lllongisomarrow{\lllongmaprighto{\lower3pt\hbox{$\scriptstyle\sim$}}}

\def\lisomarrow{\maplefto{\lower3pt\hbox{$\scriptstyle\sim$}}}

\font\labprffont=cmtt8
\def\strutdepth{\dp\strutbox}
\def\labtekst#1{\vtop to \strutdepth{\baselineskip\strutdepth\vss\llap{{\labprffont #1}}\null}}
\def\marglabel#1{\strut\vadjust{\kern-\strutdepth\labtekst{#1\ }}}

\def\label #1. #2\par{{\definexref{#1}{\num}{#2}}}
\def\labelf #1\par{{\definexref{#1}{\numfo}{formula}}}
\def\labelse #1\par{{\definexref{#1}{\num}{section}}}


\def\fibprod{\mathop\times}
\def\tensor{\mathop\otimes}
\def\dirsum{\mathop\oplus}


\def\Def {{\rm Def}}
\def\Gro{{\rm Groups}}

\def\Sch{{\rm Schemes}}


\def\CC{{\bf C}}
\def\E{{\bf E}}
\def\GG {{\bf G}}
\def\NN {{\bf N}}

\def\QQ{{\bf Q}}
\def\RR {{\bf R}}
\def\WW{{\bf W}}
\def\ZZ {{\bf Z}}

\def\A{{\euscrfm A}}
\def\B{{\euscrfm B}}
\def\C{{\euscrfm C}}
\def\J{{\euscrfm D}}
\def\DD{{\euscrfm D}}
\def\e{{\euscrfm e}}

\def\I{{\euscrfm I}}
\def\L{{\euscrfm L}}
\def\l{{\euscrfm l}}
\def\M{{\euscrfm M}}
\def\m{{\euscrfm m}}
\def\P{{\euscrfm P}}
\def\p{{\euscrfm p}}
\def\Q{{\euscrfm Q}}
\def\T{{\euscrfm T}}
\def\cu{{\euscrfm u}}
\def\W{{\euscrfm W}}

\def\HH{{\euscrfb H}}
\def\jj{{\euscrfb j}}
\def\MM{{\euscrfb M}}

\def\D{{\euscr D}}
\def\G{{\euscr G}}

\def\AA{{\rm A}}
\def\E{{\rm E}}
\def\F{{\rm F}}
\def\H{{\rm H}}
\def\j{{\rm j}}
\def\k{{\rm k}}
\def\R{{\rm R}}
\def\V{{\rm V}}


\def\1{{\bf 1}}
\def\AU{{\rm A}^{\rm U}}
\def\Aut{{\rm Aut}}
\def\cand{{\omega^{\rm can}}}
\def\canD{{\KS^{-1}\bigl(\omega^{\rm can}\tensor \omega^{\rm can}\bigr)}}
\def\Cand{{\omega^{\rm can}}}

\def\DL{{\rm D}_L}
\def\End{{\rm End}}
\def\Ext {{\rm Ext}}
\def\EXT {{{\cal E}\kern-.5pt{\it xt} }}
\def\FF{{{\bf F}_p}}

\def\Gr{{\rm GR}}
\def\Hbar{{\overline{H}}}
\def\hslash{{\hbar}}
\def\Hom {{\rm Hom}}
\def\HOM {{{\cal H}\kern-.5pt{\it om} }}
\def\Ker {{\rm Ker}}
\def\Kum{{\rm Kum}}
\def\KS{{\rm KS}}
\def\KSbar{{\rm KS}}
\def\LPj{{\Lambda\bigl(\P,\jj\bigr)}}

\def\MMbar{{\overline{\MM}}}
\def\Norm{{\bf Nm}}
\def\Nm{{\bf Nm}}

\def\phibar{{\overline{\phi}}}

\def\Res{{\rm Res}}

\def\sigmatilde{{\tilde{\sigma}}}
\def\SL{{\rm SL}}
\def\Spec {{\rm Spec}}
\def\Spf {{\rm Spf}}

\def\tame{{\rm t}}
\def\Tate{{\bf Tate}}

\def\Tr{{\rm Tr}}
\def\val{{\rm val}}
\def\wild{{\rm w}}
\def\X{{\Bbb X}}


\def\DeligneRibet{{\rm DeRi}}
\def\DelignePappas{{\rm DePa}}
\def\Chai{{\rm Ch}}

\def\Geer{{\rm vdG}}
\def\Goren{{\rm Go1}}
\def\Gorenn{{\rm Go2}}
\def\Gorennn{{\rm Go3}}
\def\GorenOort{{\rm GoOo}}
\def\Gross{{\rm Gr}}
\def\EGAII{{{\rm EGA II}}}
\def\EGAIVtwo{{{\rm EGA IV$^2$}}}
\def\EGAIVfour{{\rm EGA IV$^4$}}
\def\SGA2{{\rm SGA 2}}
\def\Hida{{\rm Hida}}
\def\Katz{{\rm Ka1}}
\def\Katzz{{\rm Ka2}}
\def\Katzzz{{\rm Ka3}}
\def\Katzzzz{{\rm Ka4}}

\def\Lang{{\rm La}}
\def\Messing{{\rm Me}}
\def\MilnorStasheff{{\rm MiSt}}
\def\vanderPoorten{{\rm vdP}}
\def\Rapoport{{\rm Ra}}
\def\Serre1{{\rm Se}}
\def\Siegel{{\rm Si}}
\def\Washington{{\rm Wa}}
\def\Zink{{\rm Zi}}

\senu=1 \ssnu=1 \fonu=1

\centerline{\titlefont HILBERT MODULAR FORMS:} \spacing
\centerline{\titlefont MOD P AND P-ADIC ASPECTS}
\bigskip
\centerline{ F.~A{\eightpoint NDREATTA} {\eightpoint AND}
E.~Z.~G{\eightpoint OREN}}

\bigskip
\bigskip

{\insert\footins{\leftskip\footnotemargin\rightskip\footnotemargin\noindent\eightpoint
$2000$ {\it Mathematics Subject Classification}. Primary 11G18,
14G35, 11F33, 11F41.
\par\noindent {\it Key words}: Congruences, Hilbert modular forms,
Hilbert modular varieties, zeta functions, filtration.}

\vbox{{\leftskip\abstractmargin \rightskip\abstractmargin
\eightpoint

\noindent A{{\sixrm BSTRACT}}.\enspace We study Hilbert modular
forms in characteristic~$p$ and over $p$-adic rings. In the
characteristic~$p$-theory we describe the kernel and image of the
$q$-expansion map and prove the existence of filtration for
Hilbert modular forms; we define operators $U$,~$V$
and~$\Theta_\chi$ and study the variation of the filtration under
these operators. In particular, we prove that every ordinary
eigenform has filtration in a prescribed box of weights. Our
methods are geometric -- comparing holomorphic Hilbert modular
forms with rational functions on a moduli scheme with level-$p$
structure, whose poles are supported on the non-ordinary locus.

In the $p$-adic theory we study congruences between Hilbert
modular forms. This applies to the study of congruences between
special values of zeta functions of totally real fields. It also
allows us to define $p$-adic Hilbert modular forms ``\`a la Serre"
as $p$-adic uniform limit of classical modular forms, and compare
them with the $p$-adic modular forms ``\`a la Katz" that are
regular functions on a certain formal moduli scheme. We show that
the two notions agree for cusp forms and for a suitable class of
weights containing all the classical ones. We extend the operators
$V$ and~$\Theta_\chi$ to the $p$-adic setting.

}}

\bigskip
\section Introduction\par \noindent This paper is concerned with
developing the theory of Hilbert modular forms along the lines of
the theory of elliptic modular forms. Our main interests in this
paper are: \spacing

\item{{\rm (i)}} to determine the ideal of congruences between Hilbert
modular forms in characteristic~$p$ and to find conditions on the
existence of congruences over artinian local rings. This allows us
to derive explicit congruences between special values of zeta
functions of totally real fields, to establish the existence of
filtration for Hilbert modular forms, to establish the existence
of $p$-adic weight for $p$-adic modular forms (defined as $p$-adic
uniform limit of classical modular forms) and more;

\spacing

\item{{\rm (ii)}} to construct operators $U, V, \Theta_\psi$ (one for
each suitable weight~$\psi$) on modular forms in
characteristic~$p$ and to study the variation of the filtration
under these operators. This allows us to prove that every ordinary
eigenform has filtration in a prescribed box of weights;

\spacing

\item{{\rm (iii)}} to show that there are well defined notions of
a Serre $p$-adic modular form and of a Katz $p$-adic modular form
and to show that the two notions agree for a suitable class of
weights containing all the classical ones. Our argument involves
showing that every $q$-expansion of a  mod~$p$ modular form lifts
to a $q$-expansion of a characteristic zero modular form. We
extend the theta operators~$\Theta_\psi$ to the $p$-adic setting
-- their Galois theoretic interpretation is that of twisting a
representation by a  Hecke character.

\bigskip
\noindent Our approach to modular forms is emphatically geometric.
Our goal is to develop systematically the geometric and arithmetic
aspects of Hilbert modular varieties and to apply them to modular
forms. As to be expected in such a project, we use extensively the
ideas of N.~Katz [\Katz], [\Katzz], [\Katzzz], [\Katzzzz] and
J.-P.~Serre [\Serre1], of the founders of the theory in the case
of elliptic modular forms, and we have benefited much from
B.~Gross' paper~[\Gross]. In regard to previous work on the
subject, we mention that some of the constructions and methods in
this paper were introduced by the second named author, in the
unramified case, in~[\Goren], [\Gorenn], and  the congruences we
list for zeta functions may be derived from the work of P.~Deligne
and K.~Ribet [\DeligneRibet]. For this reason we restrict our
discussion to zeta functions, though the same reasoning applies to
a wide class of $L$-functions.

\bigskip
\noindent We now describe in more detail the main results of this
paper.

\indent Let~$L$ be a totally real field of degree~$g$ over~$\QQ$.
Let~$K$ be a normal closure of~$L$. Let $N\geq 4$ be an integer
prime to $p$. Let~$\MM(S, \mu_N)$ be the fine moduli scheme
parameterizing polarized abelian schemes over~$S$ with RM by~$O_L$
and $\mu_N$-level structure; see~\refn{moduli}. A Hilbert modular
form defined over an $O_K$-scheme $S$ has a weight~$\psi\in \X_S$,
where~$\X_S$ is the group of characters of the algebraic
group~$\G_S={\Res}_{O_L/\ZZ} \GG_{m,O_L} \fibprod_{\Spec(\ZZ)} S$.
We shall mostly be concerned with weights obtained from the
characters~$\X$ of $\G_{O_K} = {\Res}_{O_L/\ZZ} \GG_{m,O_L}
\fibprod_{\Spec(\ZZ)} \Spec(O_K)$; we shall use the
notation~$\X_S^U$ to the denote the group of characters of~$\G_S$
induced from~$\X$ by base change. See~\refn{GGm}.

\indent The group~$\X$ is a free abelian group of rank~$g$ and has
a positive cone~$\X^+$ generated by the characters coming from the
embeddings~$\sigma_1, \dots, \sigma_g\colon L \rightarrow K$.
Indeed, the map $O_L\tensor_\ZZ K \cong \oplus_{i=1}^g K$ induces
a splitting of the torus~$\G_K$, and hence canonical generators
of~$\X$ that we denote accordingly by $\chi_1, \dots, \chi_g$, and
call the {\it fundamental characters}. A complex Hilbert modular
form of weight~$\chi_1^{a_1}\dots \chi_g^{a_g}$ is of
weight~$(a_1, \dots, a_g)$ in classical terminology.

\indent It is important to note that~$\X_S^U$ depends very much
on~$S$. For example, assume that~$L$ is Galois and~$S =
\Spec(O_L/p)$, with~$p$ an inert prime in~$L$, then~$\X \cong
\X_S^U$, while if~$p$ is totally ramified in~$L$, say $p = \P^g$,
then~$\X _S^U\cong \ZZ$; in this case, letting~$\Psi$ denote the
reduction of any fundamental character~$\chi_i$, we obtain that an
$O_L$-integral Hilbert modular form of weight~$\chi_1^{a_1}\dots
\chi_g^{a_g}$ reduces modulo $p$ to a modular form of
weight~$\Psi^{a_1 + \dots + a_g}$.

\indent We denote the Hilbert modular forms defined over~$S$, of
level~$\mu_N$ and weight~$\chi$ by~${\bf M}(S,\mu_N,\chi)$.

\spacing
\noindent Let $p$ be a rational prime. Let~$k$ be a finite field
of characteristic $p$,  which is an $O_K$-algebra. Let~$\X_k(1)$
be the subgroup of~$\X_k$ consisting of characters that are
trivial on~$\bigl(O_L/(p) \bigr)^*$ under the map
$\bigl(O_L/(p)\bigr)^* \hookrightarrow \G(k)= \bigl(
O_L\tensor_\ZZ k\bigr)^*\rightarrow  k^*$. It is proven
in~\refn{GFp} that the map $\X \rightarrow \X_k$ is surjective
and, in particular, $\X_k = \X_k^U$. This allows us to define a
positive cone~$\X_k^+$ in~$\X_k$ as follows. For every~$i$ there
exists $1\leq \tau(i) \leq g$ such that the image
of~$\chi_i^p\chi_{\tau(i)}^{-1}$ in~$\X_k^U$ is in~$\X_k(1)$. The
character~$\chi_i^p\chi_{\tau(i)}^{-1}$ in~$\X_k^U$ does not
depend on the choice of~$\tau(i)$. The positive cone in~$\X_k$ is
the one induced by these generators. The positive cone induces an
order~$\leq_k$ on~$\X_k$; we say that~$\tau_1\leq_k \tau_2$
if~$\tau_1^{-1}\tau_2$ belongs to the positive cone. Note that we
have provided~$\X_k(1)^+:= \X_k(1) \cap \X_k^+$ with a canonical
set of generators.

For every character $\psi\in \X_k(1)^+$ we construct in~\refn{hPi}
a holomorphic modular form~$h_\psi$ over~$k$. By~\refn{qexphPi} it
has the property that its $q$-expansion at any $\FF$-rational cusp
is~$1$. Moreover, the ideal~${\cal I}$ of congruences
$${\cal I}: = \Ker \Bigl\{ \dirsum_{\chi\in \X_k} {\bf M}(k, \mu_N,\chi)\;
\lllongmaprighto{\hbox{{\it q}{\rm -exp}}}\;
k[\![q^\nu]\!]_{\nu\in M}\Bigr\}$$(where $M$ is a suitable
$O_L$-module depending on the cusp used to get the $q$-expansion)
is given by $$\bigl(h_\psi - 1: \psi\in\X_k(1)^+\bigr).$$It is a
finitely generated ideal and a canonical set of generators is
obtained by letting~$\psi$ range over the generators
for~$\X_k(1)^+$ specified above; see~\refn{kerq}; Cf.~[\Gorenn].

\indent Again, it may beneficial to provide two examples. Assume
that~$L$ is Galois. If~$p$ is inert in~$L$, we may
order~$\sigma_1, \dots, \sigma_g$ cyclically with respect to
Frobenius: $\sigma\circ\sigma_i = \sigma_{i+1}$. Let~$k =O_L/(p)$.
Then~$\X_k(1)$ and~$\X_k(1)^+$ are generated by the characters
$\chi_1^p\chi_2^{-1}, \dots, \chi_i^p\chi_{i+1}^{-1}, \dots,
\chi_g^p\chi_1^{-1}$. Note that this positive cone is different
from the one obtained from~$\X^+$ via the reduction map. The
kernel of the $q$-expansion map is generated by~$g$ relations~$h_1
- 1, \dots, h_g-1$, where~$h_i = h_{\chi_{i-1}^p\chi_i^{-1}}$ is a
modular form of weight~$\chi_{i-1}^p\chi_i^{-1}$. On the other
hand, when~$p=\P^g$ is totally ramified,~$k =O_L/\P$, we find
that~$\X_k(1)$ is generated by the characters~$\chi_1^{p-1},
\dots, \chi_g^{p-1}$ that are all the same character in~$\X_k^U$
and the $q$-expansion kernel is generated by a single
relation~$h_{\Psi^{p-1}} - 1$, where~$h_{\Psi^{p-1}}$ is a modular
form of weight~$\Psi^{p-1}$.

\bigskip

\noindent We offer two constructions of the modular
forms~$h_\psi$; see~\refn{hPi} and~\refn{secondhPi}. One
construction allows us to prove in~\refn{redirr} that the divisor
of~$h_\psi$, for~$\psi$ one of the canonical generators
of~$\X_k(1)^+$, is a reduced divisor. The other construction is
related to a compactification of~$\MM(k, \mu_{Np})$.

\indent The proof of the theorem on the ideal of congruences is
based on the isomorphism between the ring $\dirsum_{\chi\in \X_k}
{\bf M}(k, \mu_N,\chi)/{\cal I}$ of {\it modular forms as
$q$-expansions} and the ring of regular {\it functions} on the
quasi-affine scheme~$\MM(k, \mu_{Np})$. The latter scheme can be
compactified by adding suitable roots of certain of the
sections~$h_\psi$. This isomorphism creates a {\it dictionary}
between modular forms of level~$\mu_N$ and meromorphic functions
on the compactification~$\MMbar(k, \mu_{Np})$ that are regular
on~$\MM(k,\mu_{Np})$. Under this dictionary, the weight of a
modular form, a character in~$\X_k$, is mapped to an element
of~$\X_k/\X_k(1)$, the $k^*$-valued characters of the Galois
group~$\bigl(O_L/(p)\bigr)^*$ of the cover~$\MM(k, \mu_{Np})
\rightarrow \MM(k, \mu_{N})^{\rm ord}$. The exact behavior of the
poles is related to a minimal weight (with respect to the
order~$\leq_k$ on~$\X_k$), called {\it filtration}, from which a
$q$-expansion may arise, and here the explicit description of the
compactification is invaluable in studying the properties of the
filtration.

\indent This dictionary also allows us to define operators that
clearly depend only on $q$-expansions, like~$U$,~$ V$ and
suitable~$\Theta_\psi$ operators, first as operators on functions
on~$\MM(k, \mu_{Np})$, and then as operators on modular forms (see
\S\S~15-17). This enables us to read some of the finer properties
of these operators from the corresponding properties on~$\MM(k,
\mu_{Np})$. Our main results are the following:

\spacing
\item{{$\bullet$}} There exists a notion of filtration for Hilbert modular
forms: a $q$-expansion arising from a modular form~$f$ arises from
a modular form of minimal weight~$\Phi(f)$ with respect to~$\leq_k
$; see~\refn{filtrations}. This weight satisfies~$0\leq_k
\Phi(f)$. See~\refn{positivity}.

\spacing
\item{{$\bullet$}} There exists a linear operator $V\colon{\bf M}(k, \chi,
\mu_N)\rightarrow {\bf M}(k, \chi^{(p)}, \mu_N)$, whose effect on
$q$-expansions is $\sum a_\nu q^\nu \mapsto \sum a_\nu q^{p\nu}$.
The character~$\chi^{(p)}$ is the character induced from~$\chi$ by
composing with Frobenius. (Concretely, in the inert
case~$\chi_i^{(p)} = \chi_{i-1}^p$, and in the totally ramified
case~$\Psi^{(p)} = \Psi^p$.) We have~$\Phi(Vf) = \Phi(f)^{(p)}$.
The operator~$V$ comes from the Frobenius morphism of $\MM(\FF,
\mu_{Np})$. See~\refn{V},~\refn{wtV},~\refn{qexpV}
and~\refn{PhiV}.
\spacing

\item{{$\bullet$}} For every fundamental character~$\psi$  there exists a
$k$-derivation $\Theta_\psi$ taking modular forms of weight~$\chi$
to modular forms of weight~$\chi\psi^{(p)}\psi$. Its effect on
$q$-expansions at a suitable cusp is given by $\sum a(\nu)q^\nu
\mapsto \sum \psi(\nu)a(\nu)q^{\nu}$.
See~\refn{ThetaPi}--\refn{qexpThetaPi}. The behavior of filtration
under such operator involves too much notation to include here and
we refer the reader to~\refn{PhiTheta}. In the inert case,
$\Phi(f) = \chi_1^{a_1} \cdots \chi_g^{a_g}$, we have
$\Phi(\Theta_{\chi_i}f) \leq_k \Phi(f)\chi_{i-1}^p \chi_i,$ and
$\Phi(\Theta_{\chi_i}f) \leq_k \Phi(f)\chi_i^2$ iff~$p\vert a_i$
(note: $\chi_i^2 =(\chi_{i-1}^p \chi_i)/(\chi_{i-1}^p
\chi_i^{-1})$). If~$p$ is completely ramified the result resembles
the elliptic case: if~$\Phi(f) = \Psi^a$ then $\Phi(\Theta_\Psi f)
\leq_k \Psi^{a+p+1}$ and~$\Phi(\Theta_\Psi f) \leq_k \Psi^{a+2}$
iff $p\vert a$.

\spacing
\item{{$\bullet$}} There exists a linear operator~$U$ taking holomorphic modular forms of
weight~$\chi$ to  meromorphic modular forms of weight~$\chi$.
See~\refn{U}. We have identity of $q$-expansions $$ VU
f(q)=\prod_{i=1}^g \bigl(I -
\Theta_{\chi_i}\bigr)f(q)$$provided~$p$ is unramified (see
\refn{formulaU} for the general formula) and the effect of~$U$ on
$q$-expansions is
$$\sum a(\nu)q^\nu \mapsto \sum a(p\nu)q^{\nu}.$$If~$\chi$ is
``positive enough" (see \refn{U} for a precise statement) then~$U$
takes holomorphic modular forms of weight~$\chi$ to holomorphic
modular forms of weight~$\chi$.  Every ordinary modular form has
the same $q$-expansion as an ordinary modular form  in an explicit
box of weights (see~\refn{Ueigenform}). For example, this box is
equal to~$[t, p+1]^g$ for~$p$ inert, where~$t=0$ if~$p=2$
and~$t=1$ otherwise. It is equal to~$[2, p+1]$ if~$p$ is totally
ramified.

\spacing

Our main results concerning the $p$-adic theory are the following.
We define a $p$-adic modular form as an equivalence class of
Cauchy sequences of classical modular forms of level~$\mu_N$,
with~$N$ prime to~$p$, with respect to an appropriate $p$-adic
topology; see~\refn{Serrepadic}. A choice of cusp allows one to
identify a $p$-adic modular form with a $p$-adic uniform limit of
$q$-expansions of classical modular forms of a fixed
level~$\mu_N$; see~\refn{qexppadic} and~\refn{padicgeneral}. Thus,
this definition is in the spirit of Serre's original definition in
the one dimensional case~[\Serre1] and we are able to present a
theory very similar to that of loc.~cit. In particular, a $p$-adic
modular form has a well defined weight in a $p$-adic completion~$
\widehat{\X}$ of the group of universal characters~$\X$;
see~\refn{padicgeneral}.

On the other hand, one may define $p$-adic modular forms as
regular functions on a formal scheme along lines established by
Katz~[\Katzz], [\Katzzzz] (who works mostly in the unramified
case); see~\refn{Katzpadicmodfor}. The group~$(O_L\tensor_\ZZ
\ZZ_p)^*$ acts on this ring of regular functions and we
interpret~$\X$ as $p$-adic characters of~$(O_L\tensor_\ZZ
\ZZ_p)^*$. We prove that the notions of a modular form in Serre's
approach  and Katz' approach agree under minor restrictions: one
may then identify a modular form of weight~$\chi \in \widehat{\X}$
with a regular function transforming under the action of the
group~$(O_L\tensor_\ZZ\ZZ_p)^*$ by the character~$\chi$;
see~\refn{CComPPare}.

Using this, we extend the definition of the theta operators
defined modulo~$p$ to $p$-adic operators that agree with those in
Katz~[\Katzzzz] when defined;
see~\refn{padictheta},~\refn{qexppadictheta}
and~\refn{theyareasinKatz}. This allows us to present many
examples of $p$-adic modular forms; see~\refn{padicnonclassic}.

\endsection

\noindent{\sectionfont Contents.}
\spacing
\item{{\rm 1.}} Introduction.
\item{{\rm 2.}} Notations.
\item{{\rm 3.}} Moduli spaces of abelian varieties with real
multiplication.
\item{{\rm 4.}} Properties of $\G$.
\item{{\rm 5.}} Hilbert modular forms.
\item{{\rm 6.}} The $q$-expansion map.
\item{{\rm 7.}} The partial Hasse invariants.
\item{{\rm 8.}} Reduceness of the partial Hasse invariants.
\item{{\rm 9.}} A compactification of $\MM(k,\mu_{pN})$.
\item{{\rm 10.}} Congruences mod $p^n$ and Serre's $p$-adic modular
forms.
\item{{\rm 11.}} Katz's $p$-adic Hilbert modular forms.
\item{{\rm 12.}} Integrality and congruences for values of zeta
functions.
\item{{\rm 13.}} Numerical examples.
\item{{\rm 14.}} Comments regarding values of zeta functions.
\item{{\rm 15.}} The operators $\Theta_{\P,i}$.
\item{{\rm 16.}}  The operator $V$.
\item{{\rm 17.}}  The operator $U$.
\item{{\rm 18.}} Applications to filtrations of modular forms.
\item{{\rm 19.}} Functorialities.
\item{{ }} References.

\section Notations\par\label notAtion. definition\par\defi Let~$L$ be a
totally real number field. Let\/~$p$ be a prime of\/~$\ZZ$. Let
\spacing
\item{{$\bullet$}} $O_L$ be the ring of integers of\/~$L$;
\spacing
\item{{$\bullet$}} $g=[L:\QQ]$ be the degree of\/~$L$
over~$\QQ$;
\spacing
\item{{$\bullet$}} $\DL$ be the different ideal of\/~$L$
over~$\QQ$ and $d_L$ be the discriminant;
\spacing
\item{{$\bullet$}} $\left\{\bigl(\I,\I^+\bigr)\right\}$ be
fractional $O_L$-ideals, with the natural notion of positivity,
forming a set of representatives of the strict class group
of\/~$L$ i.~e., of the isomorphism classes of projective
$O_L$-modules of rank~$1$ with a notion of positivity;
\spacing
\item{{$\bullet$}} $K$ be
a Galois closure of\/~$L$. Fix embeddings $K\hookrightarrow
\overline{\QQ}_p$ and $K \hookrightarrow \CC$;
\spacing
\item{{$\bullet$}} $\{ \gamma\colon L \rightarrow K \}$ be the
set of distinct embeddings of\/~$L$ into~$K$;
\spacing
\item{{$\bullet$}} $\{\P \vert\,\P\,\hbox{{\rm divides }} p\}$ be
the primes of\/~$O_L$ above~$p$. For each prime~$\P$ over~$p$,
let~$\pi_\P\in O_L$ be a generator of\/~$\P\tensor_\ZZ \ZZ_p$ such
that $\pi_\P\not\in\P'$ for any other prime~$\P'$ over~$p$
different from~$\P$;
\spacing
\item{{$\bullet$}} $e_\P$ be the
ramification index of\/~$\P$ over~$p$ for any prime~$\P$
of\/~$O_L$ over~$p$;
\spacing
\item {{$\bullet$}} $\k_\P=O_L/\P$ and $f_\P=[\k_\P:\FF]$
for any prime~$\P$ over~$p$. Denote by $\WW(\k_\P)$ the Witt
vectors of\/~$\k_\P$.
\enddefi

\label basicwt. section\par\ssection More notation\par Let\/~$k$
be a perfect field of characteristic~$p$ such that~$k$ contains
all the residue fields~$\{\k_\P\}_{\P\vert p}$ of the prime ideals
of\/~$O_L$ over~$p$. Let $\sigma\colon k \rightarrow k$ be the
{\it absolute} Frobenius on~$k$ sending $x \mapsto x^p$.
Let~$\sigma\colon \WW(k) \rightarrow \WW(k)$ denote also the
unique lift of~$\sigma$ to the Witt vectors~$\WW(k)$ of~$k$. For
each prime~$\P$ of~$O_L$ over~$p$ let
$$\bigl\{\bar{\sigma}_{\P,i}\colon \k_\P=O_L/\P\llongrightarrow
k\bigl\}_{i=1,\ldots,f_\P}$$
$$(\hbox{{\rm resp. }}\bigl\{\hat{\sigma}_{\P,i}\colon \WW\bigl(\k_\P\bigr)\llongrightarrow
\WW(k)\bigl\}_{i=1,\ldots,f_\P})$$be the set of different
homomorphisms from the residue field~$\k_\P$ of~$O_L$ at~$\P$
to~$k$ (resp. from~$\WW\bigl(\k_\P\bigr)$ to~$\WW(k)$) ordered so
that $$\sigma \circ
\bar{\sigma}_{\P,i}=\bar{\sigma}_{\P,i+1}\qquad(\hbox{{\rm resp.
}}\quad \sigma \circ
\hat{\sigma}_{\P,i}=\hat{\sigma}_{\P,i+1}).$$Then
$$O_L \tensor_{\ZZ_p} \WW(k)=\prod_\P O_{L,\P} \tensor_{\ZZ_p} \WW(k)=
\prod_\P\Bigl(\prod_{i=1}^{f_\P}
O_{L,\P}\tensor_{\WW(\k_\P)}\WW(k) \Bigr),$$where for each
prime~$\P$ the last isomorphism is induced by applying
$O_{L,\P}\tensor_{\WW(\k_\P)}$ to the isomorphism
$$\bigl(\hat{\sigma}_{\P,1},\ldots,\hat{\sigma}_{\P,f_\P}\bigr)\colon
\WW\bigl(\k_\P\bigr)\tensor_{\ZZ_p} \WW(k) \isomarrow
\prod_{i=1}^{f_\P} \WW\bigl(k\bigr).$$For each prime~$\P$ and for
$1\leq i\leq f_\P$, denote by
$$\e_{\P,i}\in O_L \tensor_\ZZ \WW(k)$$the associated idempotent.

\spacing
\noindent Let~$\p$ be a prime of~$O_K$ over~$p$ and fix an
embedding $O_K/\p \hookrightarrow k$. Let~$\sigma$ be a lifting
to~$O_{K,\P}$ of the {\it absolute} Frobenius on~$O_K/\p$. For
every prime~$\P$ of~$O_L$ over~$p$ let
$$\bigl\{\sigma_{\P,i}\colon O_{L,\P}\llongrightarrow
O_{K,\p}\bigl\}_{i=1,\ldots,f_\P}$$be  extensions of the
homomorphisms~$\bigl\{\hat{\sigma}_{\P,i}\bigl\}_{i=1,\ldots,f_\P}$
to homomorphisms from the localization~$O_{L,\P}$ of~$O_L$ to the
localization~$O_{K,\p}$ of~$O_K$.
\endssection

\label grschG. definition\par \defi Let
$$\G:=\Res_{O_L/\ZZ}\bigl(\GG_{m,O_L}\bigr) \colon \Sch
\llongrightarrow \Gro$$be the Weil restriction of\/~$\GG_{m,O_L}$
to~$\ZZ$ i.~e., the functor associating to a scheme~$S$ the
group~$\bigl(\Gamma(S,O_S)\tensor_\ZZ O_L\bigr)^*$. If\/~$T$ is a
scheme, we write $$\G_T:=\G\fibprod_{\Spec(\ZZ)} T.$$If
$T=\Spec(R)$, we write~$\G_R$ for~$\G_T$. For any scheme~$T$
define $$\X_T:=\Hom_\Gr \bigl(\G_T,\GG_{m,T}\bigr)$$as the group
of characters of~$\G_T$. We often write~$\X$ for~$\X_{O_K}$.
\enddefi
\endsection

\section Moduli spaces of abelian varieties
with real multiplication\par

\ssection {\bf Note carefully}\par Throughout this paper we fix a
fractional ideal~$\I$ with its natural positive cone~$\I^+$ among
the ones chosen in~\refn{notAtion}. Below we discuss Hilbert
moduli spaces and Hilbert modular forms, where the polarization
datum is fixed and equal to~$\bigl(\I,\I^+\bigr)$. Our notation,
though, does not reflect that. When we are compelled to consider
the same notions with all the polarization modules chosen
in~\refn{notAtion} at once, this will be explicitly mentioned.
\endssection

\label moduli. definition\par \defi Let $S$ be a scheme. Let\/~$N$
be a positive integer. Denote by
$$\MM\bigl(S,\mu_N\bigr)\rightarrow S$$the moduli stack over~$S$
of\/ $\I$-polarized abelian varieties with real multiplication
by\/~$O_L$ and $\mu_N$-level structure. It is a fibred category
over the category of $S$-schemes. If\/~$T$ is a scheme over~$S$,
the objects of the stack over~$T$ are the $\I$-polarized
Hilbert-Blumenthal abelian schemes over~$T$ relative to~$O_L$ with
$\mu_N $-level structure i.~e., quadruples
$\bigl(A,\iota,\lambda,\varepsilon\bigr)$ consisting of
\item{{a)}} an abelian scheme $A \rightarrow T$ of relative
dimension~$g$;
\item{{b)}} an $O_L$-action i.~e., a ring homomorphism
$$\iota\colon O_L\hookrightarrow \End_T(A);$$
\item{{c)}} a polarization $$\lambda\colon (M_A,M_A^+)
\isomarrow (\I,\I^+)$$i.~e., an $O_L$-linear isomorphism of
sheaves on the  \'etale site of\/~$T$ between the invertible
$O_L$-module $M_A$ of symmetric $O_L$-linear homomorphisms
from~$A$ to its dual~$A^\vee$ and the ideal\/~$\I$ of\/~$L$,
identifying the positive cone of polarizations~$M_A^+$
with~$\I^+$;
\item{{d)}} an $O_L$-linear injective
homomorphism $$\varepsilon\colon\mu_N\tensor_\ZZ\DL^{-1}
\hooklongrightarrow A,$$where for any scheme $S$ over~$T$ we
define
$$\bigl(\mu_N\tensor_\ZZ\DL^{-1}\bigr)(S):=\mu_N(S)\tensor_\ZZ\DL^{-1}.$$
\spacing
\noindent We require that the following condition, called the {\it
Deligne-Pappas condition}, holds:
\spacing
\indent{{\bf (DP)}\enspace the morphism  $A \tensor_{O_L} M_A
\llongrightarrow A^t$ is an isomorphism.}

\enddefi

\label rigidity. remark\par\rmk If $N\geq 4$, the level structure
of type~$\mu_N$ is rigid~[\Gorenn, Lemma 1.1] and hence
$\MM\bigl(S,\mu_N\bigr)$ is represented by a scheme over~$S$.
\endrmk

\label Gamma0p. section\par\ssection
$\Gamma_0\bigl(\J\bigr)$-level structure\par Let $N$ be an integer
and let\/~$\J$ be an ideal of~$O_L$ prime to~$N$. Let\/~$T$ be a
scheme. A $\I$-polarized abelian scheme over~$T$ with real
multiplication by~$O_L$ and~$\mu_N$ and
$\Gamma_0\bigl(\J\bigr)$-level structures is a quintuple
$$\bigl(A,\iota,\lambda,\varepsilon,H\bigr)/T,$$where
$\bigl(A,\iota,\lambda,\varepsilon\bigr)$ is a $\I$-polarized
abelian scheme over~$T$ with real multiplication by\/~$O_L$ and
$\mu_N$-level structure and~$H$ is an $O_L$-invariant closed
subgroup scheme
$$ H \hooklongrightarrow
A,$$such that $H$ is isomorphic to the constant group scheme $
\bigl(O_L/\J\bigr)$ locally \'etale on~$T$.

\spacing
\noindent If\/~$S$ is a scheme, let
$$\MM\bigl(S,\mu_N,\Gamma_0(\J)\bigr)\llongrightarrow S$$be the moduli stack over~$S$ of
$\I$-polarized abelian varieties with real multiplication
by\/~$O_L$ and $\mu_N \fibprod \Gamma_0(\J) $-level structure.
If~$T$ is an $S$-scheme, the objects of the stack over~$T$ consist
of $\I$-polarized abelian schemes over~$T$ with real
multiplication by~$O_L$ and~$\mu_N$ and
$\Gamma_0\bigl(\J\bigr)$-level structures.
\endssection

\label Rapo. definition\par \defi Let\/ $T$ be a scheme.  We say
that an abelian scheme~$A$ over~$T$ with $O_L$-action $\iota\colon
O_L\hookrightarrow \End_T(A)$ satisfies the {\it Rapoport
condition} if
\spacing\indent{{\bf (R)}\enspace
${{\rm \Omega^1}}_{A/T}$ is a locally free $O_T \tensor_\ZZ
O_L$-module.}
\spacing
\noindent Let\/~$S$ be a scheme. Denote by
$$\MM\bigl(S,\mu_N\bigr)^\R$$the open substack
of\/~$\MM\bigl(S,\mu_N\bigr)$ whose objects are the $\I$-polarized
abelian schemes with real multiplication by\/~$O_L$, with
$\mu_N$-level structure and satisfying the Rapoport condition.

\spacing
\noindent Let\/~$p$ be a prime. Suppose that\/~$S$ is a scheme
over~$\ZZ_p$. Denote by
$$\MM\bigl(S,\mu_N\bigr)^{\rm ord}$$the open substack
of\/~$\MM\bigl(S,\mu_N\bigr)$ whose objects are the $\I$-polarized
abelian schemes, with real multiplication by\/~$O_L$ and with
$\mu_N$-level structure, that are {\it geometrically ordinary}.
\enddefi

\label propRapo. remark\par\rmk A quadruple
$\bigl(A,\iota,\lambda,\varepsilon\bigr)$ as in (a)-(d)
of~\refn{moduli} satisfying~(R) {\it automatically} satisfies
(DP). It follows from~[\Rapoport, Cor.~1.13] and the
characteristic~$0$ theory.
\spacing
\noindent \noindent If $d_L$ is invertible in~$S$, then (R) and
(DP) are equivalent. See [\DelignePappas, Cor.~2.9].

\spacing
\noindent Suppose that\/~$S$  is a scheme over~$\ZZ_p$. A
quadruple $\bigl(A,\iota,\lambda,\varepsilon\bigr)$ as in (a)-(d)
of~\refn{moduli}, such that $A\rightarrow S$ is geometrically
ordinary, {\it automatically} satisfies (R). Indeed, by the
Serre-Tate theory of canonical lifts we can lift any ordinary
$\I$-polarized abelian variety with RM by~$O_L$ in
characteristic~$p$ to a $\I$-polarized abelian variety with RM
by~$O_L$ in characteristic~$0$.

\noindent We have open immersions
$$\MM\bigl(S,\mu_N\bigr)^{\rm ord}\hooklongrightarrow
\MM\bigl(S,\mu_N \bigr)^\R\hooklongrightarrow \MM\bigl(S,\mu_N
\bigr).$$Fiberwise over~$S$, the complement
of\/~$\MM\bigl(S,\mu_N\bigr)^\R$ in $\MM\bigl(S,\mu_N\bigr)$ has
codimension at least~$2$. See~[\DelignePappas, Prop.~4.4].

\endrmk

\label Galois. remark\par\rmk Let $N_1$ and $N_2$ be positive
integers such that~$N_2\vert N_1$ and~$N_2\geq 4$. The morphism
$$\MM\bigl(S,\mu_{N_1}\bigr)\llongrightarrow
\MM\bigl(S,\mu_{N_2}\bigr)$$is \'etale. Its  image is open.
Let~$s\in S$ be a geometric point with residue field of
characteristic~$l$. Consider the induced \'etale morphism
$$\MM\bigl(k(s),\mu_{N_1}\bigr)\llongrightarrow
\MM\bigl(k(s),\mu_{N_2}\bigr).$$Then \item{{$\bullet$}} if $l\vert
N_1$, but $l\not\vert N_2$, the map is not surjective and its
image coincides with the ordinary locus;
\item{{$\bullet$}} otherwise the map is surjective.
\smallskip
\noindent Assume that~$R$ is a local artinian ring with residue
field of characteristic~$l$. Suppose that~$N_1 N_2^{-1}=l^k$.
Then, the morphism from $\MM\bigl(S,\mu_{N_1}\bigr)$ to its image
in $\MM\bigl(S,\mu_{N_2}\bigr)$ is Galois with Galois group
isomorphic to $$\Aut_{O_L}\bigl(\mu_{N_1 N_2^{-1}} \tensor_\ZZ
\DL^{-1}\bigr) \simeq \bigl(O_L/N_1N_2^{-1}
O_L\bigr)^*=\bigl(O_L/l^k O_L\bigr)^*.$$
\endrmk
\endsection

\section Properties of\/~$\G$\par\noindent The
notation is as in~\refn{grschG}. The characters of~$\G$ over a
scheme~$S$ are the {\it weights} of Hilbert modular forms
over~$S$. For this reason we study them in some detail in this
section. Special emphasis is reserved to positive characteristic.
An important phenomenon is that if~$p$ ramifies in~$L$, there
exist characters over~$O_K$ (a normal extension of~$O_L$) which
coincide in characteristic~$p$. This explains why the kernel of
the $q$-expansion map in characteristic~$p$ is generated by~$g$
relations if~$p$ is inert in~$L$ and by a single relation if~$p$
is totally ramified. We also investigate the existence of
``exotic" weights (as opposed to the universal ones) over artinian
local rings such as Witt vectors of finite length.

\label GGm. definition\par \defi ({\it The universal
characters})\enspace The homomorphism of $O_K$-algebras $$O_L
\tensor_\ZZ O_K \lllongmaprighto{\prod_\gamma \gamma \tensor {\rm
id}} \prod_{\gamma\colon L \rightarrow K} O_K$$defines a morphism
of $O_K$-schemes
$$\G_{O_K}
\lllongmaprighto{ \prod_\gamma \chi_\gamma}\prod_{\gamma\colon L
\rightarrow K} \GG_{m,O_K}.$$The characters $\{\chi_\gamma\vert
\gamma\colon L \rightarrow K\}$ are called the {\it fundamental
characters} of\/~$\G_{O_K}$. Denote by
$$\X_{O_K}^+ \subset
\X_{O_K}:=\Hom_{O_K}\bigl(\G_{O_K},\GG_{m,O_K}\bigr)$$the cone
spanned by the fundamental characters. Define~$\Norm$ to be the
character which is the product of the fundamental ones.
\spacing

\noindent Let\/~$T$ be a scheme over~$O_K$. The subgroup~$\X^U_T$
of the character group~$\X_T$ of\/~$\G_T$, spanned by the
fundamental characters, is called the subgroup of {\it universal
characters}.
\enddefi

\label geniso. section\par\ssection The generic fiber\par The
group scheme $\G_\QQ$ is a {\it torus}. Let\/~$K$ be the fixed
Galois closure of\/~$L$ over~$\QQ$, then the universal characters
induce an isomorphism
$$\G_K
\lllongmaprighto{\prod_\gamma \chi_\gamma}\prod_{\gamma\colon L
\rightarrow K} \GG_{m,K}.$$In fact, $O_L \tensor_\ZZ
K\lllongmaprighto{\prod_\gamma\gamma \times {\rm id}}
\prod_{\gamma\colon L \rightarrow K} K$ is an isomorphism.
\endssection

\ssection The group scheme $\G_{\ZZ_p}$\par We have
$$\G_{\ZZ_p}\isomarrow \prod_{\P \vert p} \G^\P,$$where
$\G^\P=\Res_{O_{L,\P}/\ZZ_p}\bigl(\GG_{m,O_{L,\P}}\bigr)$ for all
primes~$\P$ of\/~$O_L$ over~$p$. This follows from the isomorphism
$O_L \tensor_\ZZ\ZZ_p \isomarrow \prod_{\P\vert p} O_{L,\P}$. For
each~$\P$  the natural inclusion $\WW(\k_\P) \rightarrow O_{L,\P}$
defines a closed immersion of group schemes
$$\Res_{\WW(\k_\P)/\ZZ_p}\bigl(\GG_{m,\WW(\k_\P)}\bigr)
\hooklongrightarrow\G^\P.$$The group scheme
$\Res_{\WW(\k_\P)/\ZZ_p}\bigl(\GG_{m,\WW(\k_\P)}\bigr)$ is a
torus. Let~$k$ be a field containing the fields~$\k_\P$ for all
primes~$\P$. Then for every\/~$\P$ we have canonical isomorphisms
$$\Res_{\WW(\k_\P)/\ZZ_p}\bigl(\GG_{m,\WW(\k_\P)}\bigr)
\fibprod_{\Spec\bigl(\ZZ_p\bigr)}\Spec\bigl(\WW(k)\bigr)
\isomarrow \prod_{\k_\P \rightarrow k}\GG_{m,\WW(k)} \isomarrow
\GG_{m,\WW(k)}^{f_\P}.$$
\endssection

\label GFp. section\par\ssection The special fiber of\/
$\G_{\ZZ_p}$\par Let $\P$ be a prime of\/~$O_L$ over~$p$. We have
an exact sequence
$$0 \llongrightarrow \G^{\P,u}_\FF
\llongrightarrow \G^\P_\FF \llongrightarrow \G^{\P,ss}_\FF
\llongrightarrow 0,$$where $\G^{\P,u}_\FF$ stands for the
unipotent subgroup of\/~$\G^\P_\FF$ and $\G^{\P,ss}_\FF$ is the
semisimple part. The natural inclusion $\k_\P\hookrightarrow
O_{L,\P}/p O_{L,\P}$ defines a subgroup scheme
$$\Res_{\k_\P/\FF}\bigl(\GG_{m,\k_\P}\bigr)
\hooklongrightarrow \G^\P_\FF$$mapping isomorphically to
$\G^{\P,ss}_\FF$. In particular, the exact sequence above is
split.

\spacing \noindent If~$k$ is a field containing
$\k_\P=O_L/\P$, then $$\Res_{\k_\P/\FF}\bigl(\GG_{m,\k_\P}
\bigr)\fibprod_{\Spec\bigl(\FF\bigr)} \Spec(k)\isomarrow
\prod_{\k_\P \rightarrow k}\GG_{m,k} \isomarrow
\GG_{m,k}^{f_\P}$$is a split torus. In particular, if~$k$ contains
all the residue fields of~$O_K$ at all primes above~$p$, the
reduction map
$$\X_{O_K}   \llongrightarrow
\X_k$$is surjective and, hence,~$\X_k$ is spanned by the universal
characters. As shown in example~\refn{extotram}, the reduction map
is not an isomorphism in general.
\endssection

\label chiPi. remark\par \rmk We have
$$\X_k =
\prod_\P\left(\prod_{i=1,\ldots,f_\P}
\chi_{\P,i}^\ZZ\right),$$where $\chi_{\P,i}$ is the character of
the torus $\G_{\P,\FF}^{ss}\cong \Res_{\k_\P/\FF}(\GG_{m,\k_\P})$
defined over~$k$ by the homomorphism~$\bar{\sigma}_{\P,i}$
of~\refn{basicwt}.
\endrmk

\ssection Example: the inert case\par Suppose that $p$ is inert
in~$L$. Then the morphisms defined in~\refn{GGm} are isomorphisms
after tensoring with~$\ZZ_p$. In particular, $\G_{\ZZ_p}$ is a
torus. The natural morphisms of character groups
$$\X_{\overline{\bf F}_p} \llongleftarrow
\X_{O_K} \llongrightarrow \X_K$$are all isomorphisms.
\endssection

\label extotram. section\par\ssection Example: the totally
ramified case\par Suppose that the prime $p$ is totally ramified
in~$L$ and~$g>1$. Let~$\P$ be the unique prime above it and fix an
isomorphism $O_L/pO_L\cong \k_\P[T]/\bigl(T^g\bigr)$. Then
$$\left(\G^\P_\FF\right)(R)=
\Bigl(R \tensor \k_\P[T]/(T^g)\Bigr)^*$$for any $\k_\P$-algebra
$R$. In particular, the toric part of $\G^\P_\FF$ is one
dimensional. The natural morphism of character groups
$$\X_{O_K}\llongrightarrow \X_{\overline{\bf
F}_p}$$has a non-trivial kernel. In particular, all the
fundamental characters have the same reduction~$\Psi$.
\endssection

\label exoticchar. section\par\ssection Exotic characters of\/
$\G$\par Let~$k$ be a perfect field of positive
characteristic~$p$. Let~$\WW_{m+1}(k) \rightarrow \WW_m(k)$ be the
canonical surjective homomorphism on truncated Witt vectors. It is
defined by a principal ideal~$(\epsilon)$ satisfying~$p
\epsilon=0$. Consider the induced reduction homomorphism
$$\alpha\colon \X_{\WW_{m+1}(k)}
\llongrightarrow\X_{\WW_m(k)}.$$Let $$\chi \colon
\G_{\WW_{m+1}(k)} \llongrightarrow \GG_{m,\WW_{m+1}(k)}$$be a
character. Write
$$\G_{\WW_{m+1}(k)}:=\Spec(A)\qquad\hbox{{\rm
and}}\qquad \GG_{m,\WW_{m+1}(k)}=\Spec\left(\WW_{m+1}(k)\Bigl[t,{1
\over t }\Bigr]\right).$$Then
$$\chi\in\Ker(\alpha)\qquad\hbox{{\rm
iff}}\qquad\chi(t)=1+\epsilon a, \quad a\in A/pA$$with $a$
satisfying
$$\Delta(a)=a \tensor 1 + 1 \tensor a, \qquad\hbox{{\rm
and}}\qquad \cu(a)=0,$$where~$\Delta$ (resp.~$\cu$) is the
comultiplication (resp.~counit) of~$A$. Hence the kernel of
$\alpha$ is $$\Ker\bigl(\alpha
\bigr)=\Hom_{k-\Gr}\bigl(\G_k,\GG_{a,k}\bigr).$$Let $\G_k^u$ be
the unipotent part. As soon as it is non-trivial, the
group~$\Ker(\alpha)$ is not finitely generated. In fact,  there
exists a surjective homomorphism~$\G_k\rightarrow \GG_{a,k}$ and
it is acted upon by
$$\End_{k-\Gr}\bigl(\GG_{a,k} \bigr)\isomarrow \bigl\{f \in
k[X]\vert\, f\in k[X^p],\, f(0)=0\bigr\},$$where the group
structure on the RHS is induced by composition. In particular, as
soon as~$p$ {\it is ramified}, $\X_{\WW_{m+1}(k)}$ is {\it not
finitely generated} and hence {\it does not} consist of universal
characters.

\spacing
\noindent The reader may have noticed that this phenomenon cannot
occur for~$g=1$. The phenomenon of exotic weights, and the exotic
Hilbert modular forms associated to them, may only occur for~$g>1$
in the presence of ramification. However, this ``pathology" is a
peculiarity of artinian bases and disappears in the truly $p$-adic
situation. The situation is different when we consider the
characters that lift to~$\WW_m(k)$ for all integers~$m$.

\endssection

\defi ({\it The formal characters})\enspace Let $\widehat{\G}$ be
the smooth formal group over $\Spf(\ZZ_p)$ associated to the group
scheme $\G$. Let\/~$\widehat{O}_K$ be the completion of\/~$O_K$
with respect to the ideal~$p O_K$. Define
$$\X_{\widehat{O}_K}\left(\widehat{\G}\right)=
\Hom_{\widehat{O}_K}\left(\widehat{\G}
\widehat{\fibprod}_{\Spf\bigl(\ZZ_p\bigr)}
\Spf\bigl(\widehat{O}_K\bigr),\widehat{\GG}_{m,
\widehat{O}_K}\right)$$as the group of formal characters
of\/~$\widehat{\G}$ over $\Spf(\widehat{O}_K)$. For any character
$\chi \in\X_{O_K}$ define $\widehat{\chi}$ to be the induced
element of\/~$\X_{\widehat{O}_K}(\widehat{\G})$.
\enddefi

\prop The natural morphism
$$\prod_{\gamma\colon L \rightarrow K} \widehat{\chi}_{\gamma}^\ZZ
\llongrightarrow  \widehat{\X}_{O_K}\left(\widehat{\G}\right)$$is
an isomorphism.
\endprop
\Proof Consider the following natural diagram over~$\widehat{O}_K$
$$\matrix{ \G_{\widehat{O}_K} &\llongleftarrow & \G_{\widehat{K}}\cr \Big\uparrow &
&\Big\uparrow \cr \G_{\widehat{O}_K}^{\backslash p} & &
\G_{\widehat{K}}^{\backslash\1}\cr \Big\uparrow & & \Big\downarrow
\cr \G_{\widehat{O}_K}^{\backslash p,\1}& \maprighto{\xi} &
\G_{\widehat{O}_K}^{\backslash \1},\cr}$$where
$\G_{\widehat{O}_K}^{\backslash \1}$
(resp.~$\G_{\widehat{K}}^{\backslash \1}$) stands for the
completion of~$\G_{\widehat{O}_K}$ (resp.~$\G_{\widehat{K}}$)
along the identity section, where~$\G_{\widehat{O}_K}^{\backslash
p}$ stands for the completion of~$\G_{\widehat{O}_K} $along the
special fiber~$\G_{\widehat{O}_K/p\widehat{O}_K}$ and
where~$\G_{\widehat{O}_K}^{\backslash p,\1}$ stands for the
completion of~$\G_{\widehat{O}_K}$ at the identity section
of~$\G_{\widehat{O}_K/p \widehat{O}_K}$. Let $\G=\Spec(A)$. Then
\item{{\rm a)}} using that $A$ is a domain and Krull's theorem, one
deduces that all the arrows in the diagram are inclusions at the
level of the algebras of regular functions;
\item{{\rm b)}} since $\G$ is smooth over~$\Spec(\ZZ)$,
we have that~$\G_{\widehat{O}_K}^{\backslash \1}$ is the formal
spectrum of a power series ring over~${\widehat{O}_K}$. Hence, the
associated algebra is $p$-adically complete. In particular, the
natural map $\xi\colon \G_{\widehat{O}_K}^{\backslash
p,\1}\rightarrow \G_{\widehat{O_K}}^{\backslash \1}$ induces an
isomorphism at the level of the algebras of regular functions;
\item{{\rm c)}} under this isomorphism the underlying diagram of
algebras is commutative.

\noindent A similar diagram exists replacing~$\G$
with~$\GG_{m,\ZZ}$. To any element of $
\widehat{\X}_{\widehat{O}_K}\bigl(\,\widehat{\G}\,\bigr)$, one can
associate a formal character of~$\G_{\widehat{O_K}}^{\backslash
p}$ and hence, by the diagram above, an element in
$$\X_{\widehat{K}}\left(\G_{\widehat{K}}^{\backslash\1}\right):=\Hom_{\widehat{K}}
\left( \G_{\widehat{K}}^{\backslash\1},
\GG_{m,\widehat{K}}^{\backslash\1}\right).$$On the other hand to
any character  of~$\G_K$ (resp. of~$\G$), one can associate an
element
in~$\X_{\widehat{K}}\bigl(\,\G_{\widehat{K}}^{\backslash\1}\,\bigr)$
(resp.
in~$\widehat{\X}_{\widehat{O}_K}\bigl(\,\widehat{\G}\,\bigr) $).
In particular, we obtain a commutative  diagram
$$\matrix{\X & \llongisomarrow & \X_K \cr \Big\downarrow &
&\mapdownr{t}\cr
\widehat{\X}_{\widehat{O}_K}\bigl(\,\widehat{\G}\,\bigr) &
\llongmaprighto{s} &
\X_{\widehat{K}}\left(\G_{\widehat{K}}^{\backslash\1}\right).
\cr}$$By~\refn{geniso}, we have $\G_K\cong \prod_{\gamma\colon
K\rightarrow L} \GG_{m,K}$. Hence, $t$ is an isomorphism. By~(a)
the  map~$s$ is injective. Hence, all the arrows in the diagram of
characters are isomorphisms. This proves the claim.

\bigskip In the rest of this section we define a certain filtration~$\{ \X_{B,\p}(n)
\}$ on the group of characters~$\X_B$ of~$\G\fibprod_{\Spec(\ZZ)}
\Spec(B)$. If~$B$ is a $p$-adic ring and an $O_K$-algebra, the
topology of~$\X_B$ induced by this filtration is separated and we
study the resulting completion~$\widehat{\X}_{B,\p}$. The
motivation is that $p$-adic modular forms ``\`a la Serre"
over~$B$, defined in~\refn{Serrepadic}, have a well defined weight
in $\widehat{\X}_{B,\p}$.

\label XB. definition\par\defi Let\/~$B$ be a ring with an
ideal\/~$\p$ such that\/~$B/\p$ is of characteristic~$p$. For any
non-negative integer~$n$ define
$$\X_{B,\p}(n):= \Bigl\{ \chi \in
\X_B \, \vert\, \chi\colon \bigl(O_L/p^nO_L\bigr)^*
\llongrightarrow \bigl(B/\p^n B\bigr)^* \,\hbox{{\rm is
trivial}}\,\Bigr\}.$$Define
$$\widehat{\X}_{B,\p}:= \lim_{\infty\leftarrow n}
\X_B/\X_{B,\p}(n).$$We suppress the index~$\p$ if no confusion is
likely to  arise.

\enddefi

\label structure. section\par\ssection The structure
of\/~$\widehat{\X}_{B,\p}$\par Suppose that~$B$ is an
$O_K$-algebra and that~$\p$ is a maximal ideal. Note that
$\X_B/\X_{B,\p}(1)$ is isomorphic to $\X_{B/\p B}/\X_{B/\p B}(1)$.
It follows from~\refn{GFp} that $\X_B/\X_{B,\p}(1)$  consists of
universal characters and $$\X_B/\X_{B,\p}(1)\isomarrow \Hom\bigl(
(O_L/pO_L)^*,(B/\p)^*\bigr) = \prod_{\P\vert p}
\Hom\bigl(\k_\P^*,(B/\p)^*\bigr).$$Since $$\Ker\Bigl(\bigl(
B/\p^{n+1}\bigr)^* \rightarrow \bigl(
B/\p^n\bigr)^*\Bigr)=1+\p^n\quad\hbox{{\rm (mod
}}\p^{n+1})\cong\p^n/\p^{n+1},$$it is killed by~$p$. In
particular, for any $n\geq 1$ we have that~$p$
kills~$\X_{B,\p}(n)/\X_{B,\p}(n+1)$.

\spacing
\noindent Suppose that\/~$B$ is  flat over~$O_K$.  We aim to prove
that $$\cap_n \X_{B,\p}\bigl(n\bigr)=\bigl\{1\bigr\}.$$Let
$\widehat{B}:=\lim_n B/\p^n$. Let $\sigma_1,\ldots,\sigma_g\colon
L \rightarrow K$ denote the distinct embeddings of~$L$ in~$K$. Let
$\chi_1,\ldots,\chi_g$ be the associated fundamental characters
of~$\Res_{O_K/\ZZ}\bigl(\GG_{m,O_K}\bigr)$ as in~\refn{GGm}. The
flatness of~$B$ over~$\ZZ$ implies that~$\X_B=\chi_1^\ZZ\times
\cdots \times \chi_g^\ZZ$; see~\refn{Bflat}. Suppose $\chi\in
\cap_n \X_{B,\p}(n)$. Write $\chi=\chi_1^{a_1}\cdots
\chi_g^{a_g}$. By assumption the homomorphism $\chi \colon (O_L
\tensor_\ZZ \ZZ_p)^* \rightarrow \widehat{B}^*$ is trivial. Let $U
\subset (O_L \tensor_\ZZ \ZZ_p)^*$ and $V \subset O_L \tensor_\ZZ
\ZZ_p$ be subgroups of finite index where the exponential map is
defined and induces an isomorphism $\exp\colon V \isomarrow U$.
The map $\log_{\widehat{B}}\circ \chi \circ \exp\colon V
\llongrightarrow \widehat{B}$, given by $l \tensor z \longmapsto
\sum_{i=1}^g a_i \sigma_i(l) z$, is the zero map and factors
through the completion of~$K$ at the prime ideal~$\p \cap O_K$.
The independence of the embeddings $\sigma_1,\ldots,\sigma_g$
gives that~$a_i = 0$ for every $i=1,\ldots,g$ i.~e., that $\chi=1$
as claimed.
\endssection

\label Bflat. lemma\par\lemma Suppose that\/~$B$ is a  ring flat
over~$O_K$. Then
\spacing
\item{{\rm 1.}} we have
$$\X_B=\X_{O_K}=\prod_{\gamma\colon L \rightarrow K}
\chi_\gamma^\ZZ;$$
\spacing
\item{{\rm 2.}} the topology on~$\X_B$ induced
by the system of subgroups~$\{\X_{B,\p}(n)\}_{n \in \NN}$ is
separated, i.~e.  $$\cap_{n \in \NN} \X_{B,\p}(n) =\bigl\{1
\bigr\};$$
\spacing
\item{{\rm 3.}} finally, $\widehat{\X}_{B,\p}$ is independent of\/~$\p$
and is the product of
\itemitem{$\bullet$} a finite group of order prime to~$p$ isomorphic to
$\prod_{\P\vert p} \Hom\bigl(\k_\P^*,(B/\p)^*\bigr)$;
\itemitem{$\bullet$} a topological group isomorphic to~$\ZZ_p^g$.
\endlemma
\Proof The flatness implies that $\X_B \hooklongrightarrow
\X_{B\tensor_\ZZ \QQ}$. Claim~(1) follows from~\refn{GGm}.
Claim~(2) follows from~\refn{structure}. Note that $\X_{B,\p}(1)$
is a free abelian group of rank~$g$. Consider
$$\widehat{\X}_{B,\p}(1):=\lim_{\infty \leftarrow n} \X_{B,\p}(1)/\X_{B,\p}(n).$$
By~\refn{structure} the group $\X_{B,\p}(n)/\X_{B,\p}(n+1)$ is
killed by~$p$ for~$n\geq 1$. Hence~$p^n\X_{B,\p}(1)$ is contained
in~$\X_{B,\p}(n)$.  We obtain a continuous surjective homomorphism
$$\lim_{\infty \leftarrow n} \X_{B,\p}(1)/p^n\X_{B,\p}(1) \llongrightarrow
\widehat{\X}_{B,\p}(1).$$The group on the LHS is isomorphic
to~$\ZZ_p^g$ as topological groups. It follows from~(2)
that~$\widehat{\X}_{B,\p}(1)$ has no $p$-torsion. Hence,
$\widehat{\X}_{B,\p}(1)$ is isomorphic to~$\ZZ_p^g$. We conclude
from~\refn{structure} that $\widehat{\X}_{B,\p}=
\widehat{\X}_{B,\p}(1)\times \X_{B,\p}/\X_{B,\p}(1)$. This proves
claim~(3).

\ssection Examples\par In the inert case~$\widehat{\X}_{B,\p}$ is
the direct product of a cyclic subgroup of order~$p^{f_p}-1$ and a
free $\ZZ_p$-module with basis
$\left\{\chi_1^p\chi_2^{-1},\ldots,\chi_g^p\chi_1^{-1}\right\}$.
In the totally ramified case~$\widehat{\X}_{B,\p}$ is the direct
product of a cyclic subgroup of order~$p-1$ and a free
$\ZZ_p$-module with basis
$\left\{\chi_1^{p-1},\chi_2\chi_1^{-1},\ldots,\chi_g\chi_1^{-1}\right\}$
(where one can obtain a similar basis with any~$\chi_i$ in place
of~$\chi_1$).
\endssection

\endsection

\section Hilbert modular forms\par

\label modularforms. definition\par\defi Let\/~$S$ be a scheme.
Let\/~$\chi$ be an element of\/~$\X_S$. A  $\I$-polarized Hilbert
modular form~$f$ over~$S$ of weight\/~$\chi$ and level\/~$\mu_N$
is a rule associating to
\spacing
\item{{i.}} any affine scheme~$\Spec(R)$
over~$S$;

\item{{ii.}} any $\I$-polarized Hilbert-Blumenthal
variety $(A,\iota,\lambda,\varepsilon)$ over~$\Spec(R)$ with
$\mu_N$-level structure;

\item{{iii.}} any generator~$\omega$
of~$\Omega^1_{A/R}$ as $R\tensor_\ZZ O_L$-module
\spacing

\noindent an element
$f\bigl(A,\iota,\lambda,\varepsilon,\omega\bigr)$ of\/~$R$ i.~e.,
$$(A,\iota,\lambda,\varepsilon,\omega)
\longmapsto f(A,\iota,\lambda,\varepsilon,\omega),$$with the
following properties:
\spacing

\item{{I.}} the value $f\bigl(A,\iota,\lambda,\varepsilon,\omega\bigr)$ depends
only on the isomorphism class over~$\Spec(R)$ of the
$\I$-polarized Hilbert-Blumenthal variety
$(A,\iota,\lambda,\varepsilon,\omega)$ with $\mu_N$-level
structure and with section~$\omega$;
\spacing

\item{{II.}} the rule $f$ is compatible with base change i.~e., if
$\phi\colon R \rightarrow B$ is a ring homomorphism
$$f\left(\bigl(A,\iota,\lambda,\varepsilon,\omega\bigr)\fibprod_{\Spec(R)}\Spec(B)\right)=
\phi\bigl(f(A,\iota,\lambda,\varepsilon,\omega)\bigr);$$
\spacing

\item{{III.}} if $\alpha \in \G(R)= (R\tensor_\ZZ O_L)^*$,then
$$f(A,\iota,\lambda,\varepsilon,\alpha^{-1} \omega)
=\chi(\alpha)f(A,\iota,\lambda,\varepsilon,\omega).$$

\noindent Denote by ${\bf M}\bigl(S,\mu_N, \chi\bigr)$ the
$\Gamma(S,O_S)$-module of such functions.
\enddefi

\label modforGamma0p. remark\par\rmk The notation is as
in~\refn{Gamma0p}. One defines $\I$-polarized Hilbert modular
forms of weight\/~$\chi$ and level\/~$\mu_N \times \Gamma_0(\J)$
in the obvious way. Denote by ${\bf M}\bigl(S,\mu_N,\Gamma_0(\J),
\chi\bigr)$ the $\Gamma(S,O_S)$-module of such functions.

\endrmk

\rmk The space of modular forms does not commute with base change.
In fact, it is not even true that, given a morphism of schemes $T
\rightarrow S$, the map $${\bf M}\bigl(S,\mu_N,
\chi\bigr)\tensor_{\Gamma(S,O_S)} \Gamma(T,O_T) \llongrightarrow
{\bf M}\bigl(T,\mu_N, \chi\bigr)$$is surjective. For example, one
can take~$S=\Spec(O_L)$ and~$T=\Spec(\k_\P)$ for some prime
ideal~$\P$ of~$O_L$. By~[\Geer], the weight $\chi=\prod_\gamma
\chi_\gamma^{a_\gamma}$ of non-zero modular forms over~$O_L$ must
satisfy the condition that $a_\gamma\geq 0$. See~\refn{GGm} for
the notation. Examples show that this is not necessary if\/~$T$ is
of characteristic~$p$. See~\refn{hPi}.
\endrmk

\label Lchi. definition\par\defi Suppose that $N\geq 4$. Let\/~$R$
be a ring and let~$\chi\in \X_R$. By~\refn{rigidity} there exists
a universal $\I$-polarized Hilbert-Blumenthal abelian scheme with
$\mu_N$-level structure $$(\AU, \iota^{\rm U}, \lambda^{\rm
U},\varepsilon^{\rm U})$$ over~$\MM\bigl(\ZZ,\mu_N\bigr)$. The
pull-back of the coherent sheaf $ \Omega^1_{\AU/\M(\ZZ,\mu_N)^\R}$
to the Rapoport locus~$\MM\bigl(R,\mu_N\bigr)^\R$, defined
in~\refn{Rapo}, is a locally free $O_{\M(R,\mu_N)^\R}\tensor_\ZZ
O_L$-module of rank~$1$. Hence, it defines a cohomology class
$$c \in \H^1\biggl({\MM\bigl(R,\mu_N\bigr)^\R
,\bigl({O_{\M(R,\mu_N)^\R}\tensor_\ZZ O_L}\bigr)^*}\biggr).$$Let
${\cal L}_\chi$  be the invertible sheaf
over~$\MM\bigl(R,\mu_N\bigr)^\R$ associated to the cohomology
class defined by the push-forward of\/~$c$ via the map~$\chi$
$$\H^1\biggl({\MM\bigl(R,\mu_N\bigr)^\R
,\bigl({O_{\M(R,\mu_N)^\R}\tensor_\ZZ O_L}\bigr)^*}\biggl)
\llongmaprighto{\chi}
\H^1\biggl({\MM\bigl(R,\mu_N\bigr)^\R,O_{\M(R,\mu_N)^\R}^*}\biggl).$$
See {\rm [\Rapoport, \S6.8]}.
\enddefi

\label caniso. proposition\par \prop {\rm 1.} \enspace To define a
$\I$-polarized Hilbert modular form~$f$ of weight\/~$\chi$ and
level~$\mu_N$ over a ring~$R$ is equivalent to define a section
of\/~${\cal L}_\chi$ over\/~$\MM\bigl(R,\mu_N\bigr)^\R$.
\spacing
\noindent {\rm 2.}\enspace For $\chi_1,\chi_2 \in \X_R$ we have a
canonical isomorphism
$${\cal L}_{\chi_1} \tensor{\cal L}_{\chi_2}\isomarrow {\cal L}_{\chi_1\,\chi_2}$$as
invertible sheaves on~$\MM\bigl(R,\mu_N\bigr)^\R$.
\endprop
\endsection

\section The $q$-expansion map\par

\label uniformization. section\par \ssection Uniformization of
semiabelian schemes with $O_L$-action\par Let $R$ be a noetherian
normal excellent domain complete and separated with respect to an
ideal\/~$I$ such that $\sqrt{I}=I$. Let $S:=\Spec(R)$
and\/~$S_0=\Spec(R/I)$. Let $S\backslash S_0$ be the open
subscheme defined by the complement of~$S_0$ in~$S$. Then we have
an equivalence of categories between
\spacing
\noindent{{\bf ($O_L$-DD)}}\enspace the category
of\/~$1$-motives$$\underline{q}\colon \B \longrightarrow \A^{-1}
\DL^{-1}\tensor_\ZZ\GG_m\bigl(S\backslash S_0\bigr),$$where
\itemitem{{$\bullet$}} $\A$ and\/~$\B$ are projective $O_L$-modules of rank\/~$1$;
\itemitem{{$\bullet$}} $\underline{q}$ is $O_L$-linear;
\itemitem{{$\bullet$}} the $O_L$-module $M:=\A\B$ is endowed with a notion of positivity so that
$\bigl(\A\B^{-1},(M\B^{-2})^+\bigr)\cong \bigl(\I,\I^+\bigr)$.

\noindent The motive is subject to the {\it degeneration
condition} that for any $ m=ab\in \A\B$ such that~$m\in M^+$ the
element of the fraction field of\/~$A$ associated to
$a\bigl(\underline{q}(b)\bigr) \in \GG_m(S\backslash S_0)$ belongs
to~$I$. Here, we identify~$\A$ with the character group of the
torus~$\A^{-1} \DL^{-1}\tensor_\ZZ\GG_{m,\ZZ}$.

\noindent The morphisms from an object  $\underline{q}\colon \B
\longrightarrow \A^{-1} \DL^{-1}\tensor_\ZZ\GG_m\bigl(S\backslash
S_0\bigr)$ to a second object $\underline{q}'\colon \B'
\longrightarrow (\A')^{-1}
\DL^{-1}\tensor_\ZZ\GG_m\bigl(S\backslash S_0\bigr)$ are defined
by the $O_L$-linear maps $$\bigl(\A\bigr)^{-1} \llongrightarrow
\bigl(\A'\bigr)^{-1}$$such that in
$$\matrix{
\B & \maprighto{\underline{q}} & \A^{-1}
\DL^{-1}\tensor_\ZZ\GG_m\bigl(S\backslash S_0\bigr)\cr
 & &\Big\downarrow \cr \B' &
\maprighto{\underline{q}'} & (\A')^{-1}
\DL^{-1}\tensor_\ZZ\GG_m\bigl(S\backslash S_0\bigr)\cr}$$there
exists a left vertical arrow making the diagram commute. Note that
the degeneration condition implies that~$\underline{q}$
and~$\underline{q}'$ are injective. Hence, if such map exists it
is unique.

\spacing
\noindent{{\bf ($O_L$-Deg)}}\enspace the category whose objects
consist of semiabelian schemes $G$ over~$S$ endowed with an
$O_L$-action such that\itemitem{{$\bullet$}} $G \fibprod_S
\bigl(S\backslash S_0\bigr)$ is an abelian scheme with
$O_L$-action and with  polarization module isomorphic
to~$(\I,\I^+)$;
\itemitem{{$\bullet$}} $G \fibprod_S S_0 \cong
\A\DL^{-1} \tensor_\ZZ \GG_{m,S_0}$, where~$\A$ is a projective
$O_L$-module of rank\/~$1$.

\noindent The morphisms are the homomorphisms as semiabelian
schemes commuting with the $O_L$-action.

\spacing
\noindent The semiabelian scheme $G$ associated to a $1$-motive
$\underline{q}\colon \B \rightarrow \A^{-1}\DL^{-1}\tensor_\ZZ
\GG_m(S\backslash S_0)$ is usually denoted by
$$\left(\A^{-1}\DL^{-1} \tensor_\ZZ\GG_{m,S}\right)/\underline{q}(\B).$$See~[\Rapoport] for details.
\endssection

\label indiff. remark\par\rmk We gather some properties of this
construction.
\spacing
\noindent {i)\enspace } We have an isomorphism :
$$\biggl(\left(\A^{-1}\DL^{-1}\tensor_\ZZ \GG_{m,S}\right)/\underline{q}(\B)\biggr)\fibprod_S S_0
\cong \A^{-1}\DL^{-1} \tensor_\ZZ\GG_{m,S_0}.$$
\spacing
\noindent {ii)\enspace } For any integral $O_L$-ideal $\J$ we have
an exact sequence of group schemes {\it over}~$S\backslash S_0$ $$
0 \rightarrow \left({\A \over \J \A }\right)(1) \longrightarrow
\biggl(\A^{-1}\DL^{-1} \tensor_\ZZ
\GG_{m,S}/\underline{q}(\B)\biggr)\bigl[\J\bigr] \longrightarrow
\J^{-1} \B/\B\rightarrow 0,$$where $(\A/\J \A)(1)$ is the Cartier
dual of the constant group scheme~$\A/\J\A$.
\spacing
\noindent {iii)\enspace } The module of invariant
differentials~$\omega_{G/S}$ of\/~$G$ relative to~$S$ satisfies
$$\omega_{G/S} \isomarrow \Bigl((\A^{-1} \DL^{-1})^\vee\tensor_\ZZ R\Bigr)
{dt \over t} \isomarrow \Bigl(\A \tensor_\ZZ R\Bigr) {dt \over
t},$$where $R\,dt/t$ is the module of invariant differentials
of\/~$\GG_{m,S}$ relative to~$S$.
\spacing
\noindent {iv)\enspace } The $O_L$-module of symmetric
$O_L$-linear homomorphisms from the abelian scheme $G \fibprod_S
\bigl(S\backslash S_0\bigr)$ to its dual is canonically isomorphic
to
$$\Hom_{O_L}\bigl(\B,(\A^{-1}\DL^{-1})^\vee\bigr)=\B^{-1} \A=M
\B^{-2}= M^{-1} \A^2$$and it is endowed with a notion of
positivity induced by the one on~$M$.
\endrmk

\label tateobjects. definition\par \ssection Tate objects\par Fix
projective $O_L$-modules $\A$ and\/~$\B$ of rank~$1$. Fix a notion
of positivity on~$M:=\A \B$ such that
$\bigl(\A\B^{-1},(M\B^{-2})^+\bigr)\cong\bigl(\I,\I^+\bigr)$. Fix
a rational polyhedral cone decomposition $\{\sigma_\beta\}_\beta$
of the dual cone to $M_\RR^+ \subset M_\RR$, which is invariant
under the action of the totally positive units of~$O_L$ and such
that, modulo this action, the number of polyhedra is finite. Let
$$S:=M^{\vee}\tensor_\ZZ
\GG_{m,\ZZ}.$$For any $\sigma_\beta$ we obtain an affine torus
embedding $S \subset S_{\sigma_\beta}$.
Let\/~$S_{\sigma_\beta}^\wedge$ be the spectrum of the ring
obtained by completing the affine scheme $S_{\sigma_\beta}$ along
the closed subscheme
$S_{\sigma_\beta,0}=S_{\sigma_\beta}\backslash S$ with reduced
structure. Over the base~$S_{\sigma_\beta}$ one has a canonical
$1$-motive giving rise over~$S_{\sigma_\beta}^\wedge\backslash
S_{\sigma_\beta,0}$ to a {\it $\I$-polarized generalized Tate
object}
$$\Tate(\A,\B)_{\sigma_\beta}=\Bigl(\A^{-1}\DL^{-1}
\tensor_\ZZ\GG_{m,S_{\sigma_\beta}^\wedge}/
\underline{q}\bigl(\B\bigr)\Bigr)\fibprod_{S_{\sigma_\beta}^\wedge}\bigl(
S_{\sigma_\beta}^\wedge\backslash S_{\sigma_\beta,0}\bigr).$$See
[\Rapoport] or~[\Katz, \S1.1]. Define
$$\ZZ\bigl(\bigl(\A,\B,\sigma_\beta\bigr)\bigr):=\ZZ((q^\nu))_{\nu
\in \sigma_\beta}.$$It  can be interpreted as
$$\Spec\Bigl(\ZZ\bigl(\bigl(\A,\B,\sigma_\beta\bigr)\bigr)\Bigr)=S_{\sigma_\beta}^\wedge\backslash
S_{\sigma_\beta,0}.$$
\endssection

\label cusp. definition\par \ssection Unramified cusps\par A
$\I$-polarized unramified cusp of level\/~$\mu_N$ over~$\Spec(R)$
is a quadruple $(\A,\B,\varepsilon,\j)$, where
\spacing
\item{{\rm a)}} $\A$ and $\B$ are fractional ideals such that $\A\B^{-1}=\I$;
\item{{\rm b)}} $\varepsilon \colon N^{-1}
O_L/O_L \isomarrow N^{-1} \A^{-1} /\A^{-1}$ is an $O_L$-linear
isomorphism;
\item{{\rm c)}} $\j\colon \A  \tensor_\ZZ R
\isomarrow O_L\tensor_\ZZ R$ is an $O_L\tensor_\ZZ R$-linear
isomorphism.

\endssection

\label jepsilon. remark\par\rmk {i) \enspace } By~\refn{indiff},
the equality $\A\B^{-1}=\I$ identifies~$\I$
with~$M_{\Tate(\A,\B)_{\sigma_\beta}}$, the group of $O_L$-linear
symmetric homomorphisms from~$\Tate(\A,\B)_{\sigma_\beta}$ to its
dual, and the cone~$\I^+$ with the
cone~$M_{\Tate(\A,\B)_{\sigma_\beta}}^+$ of polarizations. Thus,
we obtain a canonical polarization datum
$$\lambda_{\rm can}\colon
\left(M_{\Tate(\A,\B)_{\sigma_\beta}},M_{\Tate(\A,\B)_{\sigma_\beta}}^+
\right)\llongisomarrow \bigl(\I,\I^+\bigr).$$

\noindent {ii) \enspace} By~\refn{indiff}, the isomorphism
$\varepsilon$ defines  a canonical $\mu_N$-level structure on the
abelian scheme~$\Tate(\A,\B)_{\sigma_\beta}$.

\noindent {iii) \enspace} The isomorphism~$\j$ defines a canonical
isomorphism
$$\left(\Omega^1_{\Tate(\A,\B)_{\sigma_\beta}/\bigl(S_{\sigma_\beta}^\wedge\backslash
S_{\sigma_\beta,0}\bigr)}\right)\tensor_\ZZ R \isomarrow  O_L
\tensor_\ZZ \left(R\tensor_\ZZ
O_{S_{\sigma_\beta}^\wedge\backslash
S_{\sigma_\beta,0}}\right).$$In some situations one has a
canonical choice of~$\j$. For example:
\item{{a)}} if $R$ is a $\QQ$-algebra we get $$\j_{\rm can}\colon \A  \tensor_\ZZ R
\isomarrow \bigl(\A \tensor_\ZZ \QQ\bigr) \tensor_\QQ R \isomarrow
L \tensor_\QQ R = \bigl( O_L\tensor_\ZZ \QQ\bigr)\tensor_\QQ R;$$
\item{{b)}} if $\A$ is an integral $O_L$-ideal of norm invertible in~$R$
(e.g.~$\A=O_L$), then the natural inclusion $\A \hookrightarrow
O_L$ induces a canonical isomorphism
$$\j_{\rm can}\colon \A \tensor_\ZZ R \isomarrow O_L\tensor_\ZZ
R.$$
\item{{c)}} if $p^n R=0$ and~$p^n$ divides~$N$, then $\varepsilon$
induces by duality an isomorphism $$\A/N\A \isomarrow O_L/NO_L$$
and, consequently, a canonical isomorphism
$$\j_\varepsilon\colon \A \tensor_\ZZ R
\isomarrow O_L\tensor_\ZZ R.$$
\endrmk

\label qexpansion. definition\par \ssection The $q$-expansion\par
Let\/~$\Spec(R)$ be an affine scheme. Let $f$ be an element of
${\bf M}\bigl(R,\mu_N,\chi\bigr)$ i.~e., a $\I$-polarized Hilbert
modular form over~$\Spec(R)$ of weight $\chi$ and level~$\mu_N$.
Denote
$$f\bigl(\Tate(\A,\B),\varepsilon,\j\bigr):=f\Bigl(\Tate(\A,
\B)_{\sigma_\beta}\fibprod_{\Spec(\ZZ)} \Spec(R),\lambda_{\rm
can},\varepsilon,{dt \over t}\Bigr)$$and call it the $q$-expansion
of\/~$f$ at the unramified cusp $(\A,\B,\varepsilon,\j)$.
\endssection

\label Gamma0pqexp. section\par\ssection A variation\par  Let $\J$
be an integral $O_L$-ideal prime to~$N$. A $\I$-polarized cusp of
level $\mu_N\times \Gamma_0(\J)$ over an affine scheme~$\Spec(R)$
is a quintuple $(\A,\B,\varepsilon,H,\j)$, where
$(\A,\B,\varepsilon,\j)$ satisfy (a)-(c) in~\refn{cusp} and
$$ H \hooklongrightarrow
\Tate(\A,\B)_{\sigma_\beta}\fibprod_{\Spec(\ZZ)} \Spec(R)$$
denotes an $O_L$-invariant closed subgroup scheme isomorphic
\'etale locally to the $O_L$-constant group scheme
$\bigl(O_L/\J\bigr)$.

\noindent Let $f$ be an element of ${\bf
M}\bigl(R,\mu_N,\Gamma_0(\J),\chi\bigr)$ i.~e., a $\I$-polarized
Hilbert modular form over~$\Spec(R)$ of weight $\chi$ and
level~$\mu_N\times \Gamma_0(\J)$ in the sense
of~\refn{modforGamma0p}. Denote
$$f\bigl(\Tate(\A,\B),\varepsilon,H,\j\bigr):=f\Bigl(\Tate(\A,\B)_{\sigma_\beta}\fibprod_{\Spec(\ZZ)}
\Spec(R),\lambda_{\rm can},\varepsilon,H,{dt \over t}\Bigr)$$and
call it the $q$-expansion of\/~$f$ at the  cusp
$(\A,\B,\varepsilon,H, \j)$.

\endssection

\label subring. theorem\par\thm Let\/ $(\A,\B,\varepsilon,\j)$ be
an unramified cusp over~$R$. Then the element
$f\bigl(\Tate(\A,\B),\varepsilon,\j\bigr)$ does not depend on the
cone decomposition $\{\sigma_\beta\}$. Moreover, it is of the form
$$f\bigl(\Tate(\A,\B),\varepsilon,\j\bigr)=\sum_{\nu \in\A\B^+\cup\{0\}} c_\nu q^\nu
\in R[\![q^\nu]\!]_{\nu \in \A\B^+ \cup \{0\}}.$$\endthm\Proof
[\Rapoport, \S4.6] and [\Rapoport, Prop.~4.9].

\label cuspmincomp. remark\par\rmk It is proven in [\Chai,
Thm.~4.3 (X)] that $R[\![q^\nu]\!]_{\nu \in \A\B^+ \cup
\{0\}}^{U^2}$, the invariants under the action of squared
units~$U^2$ of~$O_L$, is the completed local ring of the minimal
compactification of $\MM\bigl(\Spec(R),\mu_N\bigr)$ at the cusp
$(\A,\B,\varepsilon,\j)$.
\endrmk

\label Koecher. theorem\par\thm ({\it $q$-expansion
principle})\enspace Let $R$ be a ring and let
$(\A,\B,\epsilon,\j)$ be an unramified cusp defined over~$R$.
Let\/~$f$ be an element of ${\bf M}\bigl(R,\mu_N,\chi\bigr)$.
\spacing
\item{{\rm i.}} If $f\bigl(\Tate(\A,\B),\varepsilon,\j\bigr)=0$, then
$f=0$;
\spacing
\item{{\rm ii.}} If $(\A,\B,\epsilon,\j)$ is defined
over a subring $R_0$ of\/~$R$ and
$f\bigl(\Tate(\A,\B),\varepsilon,\j\bigr)$ belongs
to~$R_0[\![q^\nu]\!]_{\nu \in \A\B^+ \cup\{0\}}$, then $f$ belongs
to~${\bf M}\bigl(R_0,\mu_N,\chi\bigr)$.
\endthm
\Proof See [\Rapoport, Thm.~6.7].

\label complex. section\par\ssection The comparison with the
complex theory\par Let
$$\sigma_1,\ldots,\sigma_g\colon L \llongrightarrow \RR$$be the real
embeddings of\/~$L$.  Let\/~$\A$ and\/~$\B$ be ideals of\/~$O_L$
such that $\I=\A\B^{-1}$. Fix an $O_L$-linear isomorphism
$$\varepsilon \colon N^{-1} O_L/O_L \isomarrow N^{-1} \A^{-1}
/\A^{-1}.$$Define the group
$$\eqalign{\Gamma_N\bigl(\B\dirsum(\A\DL)^{-1}
\bigr):=& \Bigl\{\left(\matrix{ a&b\cr c&d\cr}\right)\vert\, a\in
1+N O_L,d\in O_L,b\in(\B\A\DL)^{-1}, \cr & c\in N\B\A\DL,ad-bc=1
\Bigr\}.\cr }$$It acts on the $g$-fold product of the Poincar\'e
upper half plane~$\HH^g$ by
$$\left(\matrix{a&b\cr c&d\cr} \right)
(z_1,\ldots,z_g):=\left(\ldots,{\sigma_i(a)z_i+\sigma_i(b)\over
\sigma_i(c)z_i+\sigma_i(d)},\ldots\right)_{i=1,\ldots,g}.$$The
moduli space~$\MM\bigl(\CC,\mu_N)$ classifying $\I$-polarized
abelian varieties~$A$ over~$\CC$ with real multiplication by~$O_L$
and $\mu_N$-level structure, in the sense of~\refn{moduli}, is
isomorphic, as an analytic manifold, to
$$\Gamma_N\bigl(\B\dirsum(\A\DL)^{-1} \bigr)\backslash
\HH^g.$$The abelian variety corresponding to~$\tau\in \HH^g$ is
$$A_\tau:=\CC^g/\bigl(\B \tau+(\A\DL)^{-1} \bigr).$$The $\mu_N$-level
structure on~$A_\tau$ is induced by~$\varepsilon$. The vector
space ${\bf M}\bigl(\CC,\mu_N,\chi \bigr)$ of $\I$-polarized
modular forms of level\/~$\mu_N$ and weight\/~$\chi=\prod_{i=1}^g
\chi_{\sigma_i}^{a_i}$ can be viewed, more classically, as the
vector space of holomorphic functions
$$\matrix{
\HH^g& \maprighto{f}& \CC\cr \tau=(z_1,\ldots,z_g) &\mapsto
&f(\tau),\cr}$$on which the action of the modular
group~$\Gamma_N\bigl(\B\dirsum(\A\DL)^{-1} \bigr)$ is defined by
the automorphic factor
$$j_\chi\Bigl(\mu,\bigl(z_1,\ldots,z_g\bigr)\Bigr)=\prod_{i=1}^g
\bigl(\sigma_i(c)z_i+\sigma_i(d)\bigr)^{a_i}\qquad\hbox{{\rm with
}}\mu=\left(\matrix{ a&b\cr c&d\cr}\right);$$c.f.~[\Geer]
or~[\Gorennn]. Fix a modular form~$f$. Assuming $\B=O_L$, one
deduces that $f$, as a function, is invariant with respect to the
translations on~$\HH^g$
$$\tau \mapsto \tau+\alpha,\qquad \alpha \in(\A\DL)^{-1}.$$In
particular, it has $q$-expansion at the
cusp~$\bigl(i\infty,\ldots,i\infty\bigr)$
$$f\bigl(\,\underline{q}\, \bigr):=a_0+\sum_{\nu\in{\A}^{+}}
a_\nu q^\nu\qquad\hbox{{\rm with }}q^\nu:=\exp^{2\pi i\,{\rm
Tr}_{L/\QQ}(\nu \tau)},$$where $${\rm Tr}_{L/\QQ}(\nu
\tau):=\sigma_1(\nu)z_1+\ldots+\sigma_g(\nu)z_g.$$By~\refn{indiff}
we have a natural $O_L$-linear isomorphism
$$\j_{\rm can}\colon \A \tensor_\ZZ \CC \isomarrow O_L\tensor_\ZZ \CC.$$
Note that under the exponentiation map $z\mapsto \exp^{2\pi
i\,{\rm Tr}_{L/\QQ}(z)}$, where we use the identification
$\CC^g=(\A\DL)^{-1}\tensor_\ZZ \CC$, we have $$\CC^g/\bigl(\B
\tau+(\A\DL)^{-1} \bigr)\cong (\CC^*)^g/\underline{q}(\B),$$with
$\underline{q}(\B)$ equal to the image of~$\B\tau$. By the
discussion in~[\Rapoport, \S\S6.13-6.15] it follows that
$$ f\bigl(\,\underline{q}\,\bigr)=f\bigl(\Tate(\A,\B), \varepsilon,\j_{\rm
can}\bigr).$$
\endssection

\endsection

\section The partial Hasse invariants\par \noindent In this
section we define a canonical set of Hilbert modular forms that we
call ``partial Hasse invariants". The name comes from the fact
that these modular forms factor the Hasse invariant according to
the $O_L$-structure. The partial Hasse invariants are defined  in
characteristic~$p$ and, in general, do not lift to characteristic
zero. However, they play a crucial role in the theory in several
ways: they provide canonical generators for the kernel of the
$q$-expansion map; they allow one to compactify the moduli
scheme~$\MM(k, \mu_{pN})$.

\noindent The results and the methods of this section appear
already in~[\Gorenn] in the unramified case. In that case case,
the divisors of the partial Hasse invariants yield an interesting
stratification of the moduli space~$\MM(k, \mu_N)$, [\Goren]
and~[\GorenOort].

\label k. section\par\label conv. section\par \ssection
Notation\par Let $k$ be a perfect field of characteristic~$p$.
Assume that it contains  the residue fields~$\k_\P=O_L/\P$ for all
primes~$\P$ of~$O_L$ over~$p$. Let~$R$ be a local ring with
maximal ideal~$\m$ and residue field~$R/\m=k$.

\spacing
\noindent For each representative~$(\I,\I^+)$ of the strict class
group of~$L$ fixed in~\refn{notAtion} choose an element of~$\I$
generating~$\I\tensor_\ZZ\ZZ_p$ as $O_L\tensor_\ZZ \ZZ_p$-module.
It provides each $\I$-polarized abelian scheme with real
multiplication by~$O_L$ over a $\ZZ_p$-scheme with a polarization
of degree prime to~$p$.

\endssection

\label phi. definition\par \defi Let $n\geq m\geq 1$ be positive
integers. Let~$N$ be an integer such that $N\geq 4$ and~$N$ is
prime to~$p$. Assume that~$R$ is an $O_K$-algebra. Define

$$ \MM\bigl(R/\m^m,\mu_{p^nN}\bigr) \llongmaprighto{\psi}
\MM\bigl(R/\m^m,\mu_N\bigr)^{\rm ord}$$by
$$\bigl(A,\iota,\lambda,\varepsilon_{p^nN}\bigr) \llongmapsto
\bigl(A,\iota,\lambda,\varepsilon_N\bigr).$$By~\refn{Galois}, it
is a Galois cover of smooth schemes over~$R/\m^m$ with group
$$\Gamma_n:=\Aut_{R/\m^m}\bigl(\mu_{p^n}\tensor_\ZZ
\DL^{-1}\bigr)=\bigl(O_L/p^nO_L\bigr)^*.$$

\noindent Let $$\MM\bigl(R/\m^m,\mu_{p^nN}\bigr)^\Kum
\llongmaprighto{\phi}\MM\bigl(R/\m^m,\mu_N \bigr)^{\rm ord}$$be
the quotient of the above cover by the group  of  elements
of~$\bigl(O_L/p^n O_L)^*$ killed by
$$\bigl\{\chi\colon \bigl(O_L/p^nO_L\bigr)^*\rightarrow
\bigl(R/\m^m\bigr)^*\vert\, \chi\quad\hbox{{\rm is a universal
character}}\,\bigr\}.$$Let~$G_{m,n}$ be its Galois group. By {\rm
[\DeligneRibet, Thm.~4.5]} it acts transitively on the fibers
of~$\phi$.
\enddefi

\label G1. remark\par\rmk If $m=n=1$, then $G:=G_{1,1} \isomarrow
\prod_{\P\vert p}\bigl(O_L/\P\bigr)^*$.
\endrmk

\label universal. definition\par \defi  Let\/~$\chi\in
\X_{R/\m^m}$ be a character in the sense of~\refn{grschG}. Write
$$\mu_{p^n}=\ZZ[t]/\bigl(t^{p^n}-1\bigr).$$Denote by
$${dt\over t} \in \Omega^1_{\DL^{-1}\tensor_\ZZ\mu_{p^n}}
\fibprod_{\Spec(\ZZ)} \Spec\bigl(\ZZ/p^n\ZZ\bigr)$$the canonical
generator, as a free $O_L/p^n O_L$-module of rank\/~$1$, of the
submodule of invariant differentials of\/
$\Omega^1_{\DL^{-1}\tensor_\ZZ \mu_{p^n}} \fibprod_{\Spec(\ZZ)}
\Spec\bigl(\ZZ/p^n\ZZ\bigr)$.
\spacing

\noindent Suppose that $n\geq m$. Let $$\bigl(\AU,\iota^{\rm
U},\lambda^{\rm U},\varepsilon_{p^nN}^{\rm U}\bigr)
\llongrightarrow \MM\bigl(R/\m^m,\mu_{p^nN}\bigr)$$be the
$\I$-polarized universal abelian scheme with real multiplication
by~$O_L$ and $\mu_{p^nN}$-level structure. We use the notation:
$$\Cand:=\varepsilon_{p^nN,*}\biggl({dt \over t}\biggr)\in
\Omega^1_{\AU/\MM\bigl(R/\m^m,\mu_{p^nN}\bigr)}.$$See
\refn{Hodgetrivial}. \spacing \noindent Define $a\bigl(\chi\bigr)$
to be the unique $\I$-polarized modular form of weight\/~$\chi$
and level~$\mu_{p^n N}$ over~$R/\m^m$ satisfying the following.
Let~$T$ be a $R/\m^m$-algebra. Let~$(A,\iota,\lambda,\varepsilon)$
be a $\I$-polarized Hilbert-Blumenthal abelian scheme over~$T$
with $\mu_{p^nN}$-level structure. Let\/~$\omega$ be an
$O_L\tensor_\ZZ T$-generator
of\/~$\H^0\bigl(A,\Omega^1_{A/T}\bigr)$. Then, there is a unique
element\/~$ \gamma$ of\/~$\bigl(O_L\tensor_\ZZ T\bigr)^*$ such
that
$$\omega=\gamma^{-1}\, \varepsilon_*\Bigl({dt \over t}\Bigr).$$Define
$$a\bigl(\chi\bigr)\Bigl(A,\iota,\lambda,
\varepsilon,\omega\Bigr):=\chi\,\left(\gamma\right).$$
\enddefi

\label Hodgetrivial. section\par \ssection Some explanations\par
Decompose $\varepsilon^{\rm U}_{p^nN}:=\varepsilon_{p^n}^{\rm
U}\times \varepsilon_N^{\rm U}$. Since $n\geq m$, the embedding
$$\varepsilon_{p^n }^{\rm U}\colon \bigl(\mu_{p^n}\tensor_\ZZ
\DL^{-1}\bigr)\fibprod \MM\bigl(R/\m^m,\mu_{p^nN}\bigr)
\rightarrow \AU$$defines an isomorphism on tangent spaces at the
origin relative to~$\MM\bigl(R/\m^m,\mu_{p^nN}\bigr)$. By duality,
one gets a canonical isomorphism of the translation invariant
differentials relative to~$\MM\bigl(R/\m^m,\mu_{p^nN}\bigr)$. This
defines a canonical $O_L\tensor_\ZZ
O_{\MM\bigl(R/\m^m,\mu_{p^nN}\bigr)}$-generator of the translation
invariant differentials of~$\AU$
over~$\MM\bigl(R/\m^m,\mu_{p^nN}\bigr)$:
$$\Cand:=\varepsilon_{p^nN,*} \left( {dt\over t}\right)\in
\Omega^1_{\AU/\MM\bigl(R/\m^m,\mu_{p^nN}\bigr)}.$$By~\refn{Lchi},
the pull-back~$\psi^*({\cal L}_\chi)$ of the invertible sheaf
${\cal L}_\chi$, defined on~$\MM\bigl(R/\m^m,\mu_N\bigr)^\R$,
to~$\MM\bigl(R/\m^m,\mu_{p^nN}\bigr)$  is obtained as the push-out
of~$\Omega^1_{\A^{\rm U}/\MM\bigl(R/\m^m,\mu_{p^nN}\bigr)}$. The
section~$\Cand$ of~$\Omega^1_{\A^{\rm
U}/\MM\bigl(R/\m^m,\mu_{p^nN}\bigr)}$ defines by push-out a
section of~$\psi^*({\cal L}_\chi)$. Using the equivalence between
sections of~$\psi^*\bigl({\cal L}_\chi \bigr)$ and modular forms
of weight~$\chi$ and level~$\mu_{p^nN}$ over~$R/\m^m$, we find
that this section  coincides with the modular form~$a(\chi)$.
\endrmk

\prop We have $$ \Bigl(O_L\tensor_\ZZ
O_{\MM\bigl(R/\m^m,\mu_{p^nN}\bigr)}\Bigr)\,\Cand
=\Omega^1_{\AU/\MM\bigl(R/\m^m,\mu_{p^nN}\bigr)}.$$For any~$\alpha
\in \bigl(O_L/p^nO_L\bigr)^*$, the induced action by pull-back
$$\bigl[\alpha\bigr]\colon \Omega^1_{\AU/\MM\bigl(R/\m^m,\mu_{p^nN}\bigr)} \llongrightarrow
\Omega^1_{\AU/\MM\bigl(R/\m^m,\mu_{p^nN}\bigr)}$$sends
$$\Cand \llongmapsto \alpha^{-1}\, \Cand.$$
\endprop
\Proof The first claim follows from~\refn{Hodgetrivial}. Decompose
$\varepsilon^{\rm U}_{p^n N}:=\varepsilon_{p^n}^{\rm U}\times
\varepsilon_N^{\rm U}$. Since the universal $\I$-polarized abelian
scheme~$\AU$ over~$\MM\bigl(R/\m^m\mu_{p^nN} \bigr)$ is the
pull-back of the universal $\I$-polarized abelian scheme
over~$\MM\bigl(R/\m^m,\mu_N \bigr)$, the automorphism~$\alpha$
lifts to an automorphism of~$\AU$, which we denote by~$\alpha$ and
induces the automorphism~$[\alpha]$
on~$\Omega^1_{\AU/\MM\bigl(R/\m^m,\mu_{p^nN}\bigr)} $. By
definition, $\alpha\bigl(\varepsilon_{p^n}^{\rm U}\bigr)=
\bigl(\varepsilon_{p^n}^{\rm U} \circ 1 \tensor \alpha \bigr)$.
Hence,
$$\eqalign{\bigl(\varepsilon_{p^n}^{\rm U} \circ 1\tensor
\alpha\bigr)_*\left({dt\over t}\right)
 &=\varepsilon_{p^n,*}^{\rm U}\left((1\tensor \alpha)\,{dt\over t}\right)\cr
 &=\alpha^{-1}\Cand.\cr}$$

\label actsonachi. corollary\par\cor For any~$\alpha \in
\bigl(O_L/p^nO_L\bigr)^*$, the induced action by pull-back,
$$\bigl[\alpha\bigr]\colon
\Gamma\Bigl(\MM\bigl(R/\m^m,\mu_{p^nN}\bigr),\psi^*\bigl({\cal
L}_\chi\bigr)\Bigr) \llongrightarrow
\Gamma\Bigl(\MM\bigl(R/\m^m,\mu_{p^nN}\bigr),\psi^*\bigl({\cal
L}_\chi\bigr)\Bigr),$$maps
$$a\bigl(\chi\bigr)\llongmapsto\chi^{-1}(\alpha)a\bigl(\chi\bigr).$$
\endcor

\label actiona(chi). corollary\par\cor  The section~$\Cand$
descends to a section over~$\MM\bigl(R/\m^m,\mu_{p^nN}\bigr)^\Kum$
of the relative differentials of the pull-back by~$\phi$ of the
universal $\I$-polarized abelian scheme
over~$\MM\bigl(R/\m^m,\mu_N\bigr)$. We denote it~$\cand$ by abuse
of notation.

\spacing
\noindent The section $a\bigl(\chi\bigr)$ of~$\psi^*\bigl({\cal
L}_\chi\bigr)$ descends to a section, denoted~$a\bigl(\chi\bigr)$
by abuse of notation, of~$\phi^*\bigl({\cal L}_\chi\bigr)$
over~$\MM\bigl(R/\m^m,\mu_{p^nN}\bigr)^\Kum$. It is the push-out
via~$\chi$ of\/~$\cand$.

\spacing
\noindent If~$\alpha$ belongs to the Galois group~$G_{m,n}$
of~$\phi$, then $$\bigl[\alpha\bigr]\bigl(\cand\bigr)=\alpha^{-1}
\cand \qquad\hbox{{\rm
and}}\qquad\bigl[\alpha\bigr]\Bigl(a\bigl(\chi\bigr)\Bigr)=\chi^{-1}(\alpha)
a\bigl(\chi\bigr).$$
\endcor

\label qexpa(chi). proposition\par\prop {\rm 1.}\enspace The
section~$a\bigl(\chi\bigr)$ induces a {\it trivialization} of the
invertible sheaf
$\phi^*\bigl(\cal{L}_\chi\bigr)$on~$\MM\bigl(R/\m^m,\mu_{p^n
N}\bigr)^\Kum$.
\spacing
\noindent {\rm 2.}\enspace For any universal characters $\chi_1$
and~$\chi_2$ the section $a\bigl(\chi_1\bigr) \tensor
a\bigl(\chi_2\bigr)$ of\/~${\cal L}_{\chi_1}\tensor {\cal
L}_{\chi_2}$ is sent by the isomorphism defined in~\refn{caniso}
to the section $a\bigl(\chi_1\,\chi_2\bigr)$ of\/~${\cal
L}_{\chi_1\chi_2}$.
\spacing
\noindent {\rm 3.}\enspace For any universal character~$\chi$ the
$q$-expansion of~$a\bigl(\chi\bigr)$ at {\it any} $\I$-polarized
cusp $(\A,\B,\varepsilon,\j_\varepsilon)$ (see~\refn{jepsilon}
and~\refn{qexpansion}) of level~$\mu_{p^nN}$ over~$R/\m^m$ is~$1$.
\endprop
\Proof (1) and (2) are clear. Let $\Tate(\A,\B)_{\sigma_\beta}$ be
a Tate object with $\mu_{p^nN}$-level structure~$\varepsilon$. The
modular form~$a(\chi)$ takes the value~$1$ on the
section~$\varepsilon_*\bigl(dt/t\bigr)$. The latter coincides with
the section of the translation invariant relative differential
defined by~$\j_\varepsilon$ as in~\refn{jepsilon}. This
proves~(3).

\rmk The proposition justifies the
introduction of the covering~$\phi$ in~\refn{phi}.
\endrmk

\ssection First construction of the Hasse invariants: the
geometric definition\par\noindent Let $R$ be a $k$-algebra. Let
$(A,\iota,\lambda,\varepsilon)$ be a $\I$-polarized
Hilbert-Blumenthal abelian variety over~$R$ satisfying the
condition~(R)  defined in~\refn{Rapo}. Let~$\omega$ be an
$O_L\tensor_\ZZ R$-basis of~$\H^0\bigl(A,\Omega^1_{A/R}\bigr)$.
Using~$\lambda$, and the choices in~\refn{conv}, we get an
isomorphism
$$\H^0\bigl(A,\Omega^1_{A/R}\bigr)
\isomarrow \H^0\bigl(A^\vee,\Omega^1_{A^\vee/R}\bigr).$$Here
$A^\vee$ is the dual abelian scheme. Via the canonical isomorphism
$$ \H^0\bigl(A^\vee,\Omega^1_{A^\vee/R}\bigr)
\isomarrow\Hom_R\bigl(\H^1(A,O_A),R),$$the element~$\omega$
determines a generator $$\eta \in \H^1(A,O_A)$$as $O_L\tensor_\ZZ
R$-module.
\spacing
\noindent The absolute Frobenius on~$A$ induces a $\sigma$-linear
map~$O_A \rightarrow O_A$ and consequently a $\sigma$-linear map
$$\F\colon \H^1(A,O_A) \llongrightarrow
\H^1(A,O_A).$$Let~$\P$ be a prime of~$O_L$ above~$p$ and let
$1\leq i\leq f_\P$. With the notation of~\refn{basicwt}, the
$\sigma$-linearity of~$\F$ implies that for each idempotent
$\e_{\P,i}\in O_L\tensor_\ZZ R$
$$\F\bigl(\e_{\P,i-1}\cdot\eta\bigr) \in
\bigl(O_L\tensor_\ZZ R\bigr) \e_{\P,i}\cdot\eta.$$
\endssection

\label hPi. definition\par\defi Define the modular form of
weight~$\chi_{\P,i-1}^p\chi_{\P,i}^{-1}$,see~\refn{chiPi} for the
notation, over~$k$
$$h_{\P,i} \in {\bf
M}\bigl(k,\chi_{\P,i-1}^p\chi_{\P,i}^{-1}\bigr)$$by the
rule$$\F\bigl(\e_{\P,i-1}\cdot \eta\bigr) \equiv h_{\P,i}
\Bigl(\bigl(A,\iota,\lambda,\omega\bigr)\Bigr) \e_{\P,i}\cdot
\eta\qquad\hbox{{\rm mod }} \P.$$The modular forms $\bigl\{
h_{\P,i} \bigr\}_{\P,i}$ are called the {\it generalized}  or {\it
partial Hasse invariants}. If\/~$\psi \in \X_k$ is a weight of the
form
$\psi=\prod_{\P,i}\bigr(\chi_{\P,i-1}^p\chi_{\P,i}^{-1}\bigr)^{a_{\P,i}}$,
with\/~$a_{\P,i}\in \NN$, define the modular form of
weight\/~$\psi$
$$h_\psi:=\prod_{\P,i} h_{\P,i}^{a_{\P,i}}.$$Define
$h:=h_{\Norm^{p-1}}$ as the {\it Hasse invariant}.
\enddefi

\rmk It is easy to verify that all the
conditions given in the definition~\refn{modularforms} of modular
forms of a given weight are satisfied.
\endrmk

\label qexphPi. proposition\par\prop {\rm 1.}\enspace The
$q$-expansion of the partial Hasse invariant~$h_{\P,i}$ at {\it
any} $\I$-polarized cusp defined over~$\FF$ is~$1$.
\spacing
\noindent {\rm 2.}\enspace The modular form $h$ coincides, up to a
sign, with the determinant of the Hasse-Witt matrix.
\spacing
\noindent {\rm 3.}\enspace The subgroup of\/~$\X_k$ spanned by the
weights of the partial Hasse invariants coincides with the
subgroup\/~$\X_k(1)$ of the elements $\chi\in\X_k$ such that
$\chi\bigl((O_L/p O_L)^*\bigr)\equiv 1$ in~$k$.
\endprop
\Proof Part~(1) is proved via an explicit computation using Tate
objects; c.f.~[\Goren, Thm.~2.1(2)]. Part~(2) is clear. We next
prove part~(3). The set of fundamental characters
$\left\{\chi_{\P,i}:\, \P\vert p,\, 1\leq i \leq f_\P\right\}$
forms a basis for the lattice~$\X_k$. The subgroup~$H$ of~$\X_k$
spanned by the weights $\left\{\chi_{\P,i-1}^p\chi_{\P,i}^{-1}:\,
\P\vert p,\, 1\leq i \leq f_\P\right\}$ of the partial Hasse
invariants is contained in~$\X_k(1)$.  We conclude by remarking
that $\X_k/\X_k(1)\cong\prod_{\P\vert p}\Hom_{\rm
Gr}\bigl(\k_\P^*,k^*\bigr)$ has cardinality~$\prod_{\P\vert
p}\bigl(p^{f_\P}-1\bigr)$, the same as~$\X_k/H$.

\ssection Examples\par Assume that the prime $p$ is inert in~$L$.
Let $\chi_1,\ldots,\chi_g$ be the fundamental characters of~$\G_k$
ordered so that $\sigma\circ \chi_i=\chi_{i+1}$. Then we get
precisely~$g$ partial Hasse invariants $h_1,\ldots, h_g$ of
weights
$\chi_g^p\chi_1^{-1},\chi_1^p\chi_2^{-1},\ldots,\chi_{g-1}^p\chi_g^{-1}$.
In this case $$\bigl(\chi_g^p\chi_1^{-1}\bigr)^\ZZ\times \cdots
\times\bigl(\chi_{g-1}^p\chi_g^{-1}\bigr)^\ZZ=\X_k(1)
\hooklongrightarrow\X_k=\chi_1^\ZZ \times\cdots\times\chi_g^\ZZ.$$

\spacing
\noindent Assume that $p=\P^g$ is totally ramified in~$L$.
Let~$\Psi$ be the unique fundamental character of~$\G_k$. Then we
get a unique partial Hasse invariant~$h_{\Psi^{p-1}}$ which is a
$g$-th root of the Hasse invariant and is of weight~$\Psi^{p-1}$.
We have
$$\bigl(\Psi^{p-1}\bigr)^\ZZ=\X_k(1)\hooklongrightarrow
\X_k=\Psi^\ZZ.$$
\endssection

\ssection Second construction of the Hasse invariants: the
definition by descent theory\par From~\refn{actiona(chi)} we
deduce the following  result:
\endssection

\cor  The section~$a\bigl(\chi\bigr)$ of the
sheaf~$\phi^*\bigl(\cal{L}_\chi\bigr)$
on~$\MM\bigl(R/\m^m,\mu_{p^nN}\bigr)^\Kum$ descends to a section
of\/~${\cal L}_\chi$ on~$\MM\bigl(R/\m^m,\mu_N\bigr)^{\rm ord}$ if
and only if $$\chi(\alpha)=1\in (R/\m^m)^*\qquad\hbox{{\rm for all
}}\alpha\in\bigl(O_L/p^n O_L\bigr)^*.$$\endcor

\label secondhPi. corollary\par \cor Let\/~$\P$ be a prime
of\/~$O_L$ over~$p$ and let $1 \leq i\leq f_\P$. The section
$a\bigl(\chi_{\P,i-1}^p\chi_{\P,i}^{-1}\bigr)$ descends to a
section of the invertible sheaf~${\cal
L}_{\chi_{\P,i-1}^p\chi_{\P,i}^{-1}}$
over~$\MM\bigl(k,\mu_N\bigr)^{\rm ord}$. It coincides with the
restriction of the partial Hasse invariant~$h_{\P,i}$ defined
in~\refn{hPi}.\endcor \Proof The two modular forms have the same
weight and the same $q$-expansion at any $\FF$-rational cusp.
Hence, they must be equal.

\label r(f). definition\par\defi  Let
$$R_{m,n}:=\Gamma\Bigl(\MM\bigl(R/\m^m,\mu_{p^nN}\bigr)^\Kum,O_{\MM\bigl(R/\m^m,\mu_{p^nN}\bigr)^\Kum}
\Bigr).$$Define the map
$$r\colon\dirsum_{\chi
\in \X^U_{R/\m^m}} {\bf M}\bigl(R/\m^m,\mu_N,\chi\bigr)
\llongrightarrow R_{m,n}$$by the formula $$ f=\dirsum_\chi f_\chi
\llongmapsto r(f):=\sum_\chi{\phi^*\bigl(f_\chi\bigr)\over
a\bigl(\chi\bigr)}.$$Let~$\alpha\in\bigl(O_L\tensor_\ZZ
\WW(k)\bigr)^*$. If~$f$ is a $\I$-polarized modular form
on~$\MM\bigl(R/\m^m,\mu_N\bigr)$ define
$$[\alpha] f\bigl(A,\iota,\lambda,\omega\bigr):=
f\bigl(A,\iota,\lambda,\alpha^{-1}\omega\bigr).$$This provides a
graded action of~$\bigl(O_L\tensor_\ZZ R\bigr)^*$
on~$\dirsum_{\chi\in\X^U_{R/\m^m}}{\bf
M}\bigl(R/\m^m,\mu_N,\chi\bigr)$.

\spacing
\noindent  On the other hand the action of~$G_{m,n}$, the Galois
group of~$\MM\bigl(R/\m^m,\mu_{p^nN}\bigr)^\Kum \rightarrow
\MM\bigl(R/\m^m,\mu_N\bigr)^{\rm ord}$, and the canonical
projection $\bigl(O_L\tensor_\ZZ R\bigr)^* \rightarrow G_{m,n}$
defined in~\refn{phi} induce an action of~$\bigl(O_L\tensor_\ZZ
R\bigr)^*$ on~$R_{m,n}$.
\enddefi

\label propr. proposition\par\prop The map~$r$ defined
in~\refn{r(f)} has the following properties:
\spacing
\item{{\rm 1.}} it is $\bigl(O_L\tensor_\ZZ
R\bigr)^*$-equivariant;

\spacing
\item{{\rm 2.}} let $\bigl(\A,\B,\varepsilon,\j_\varepsilon\bigr)$
be a $\I$-polarized unramified cusp over~$R/\m^m$ of
level~$\mu_{p^nN}$; see~\refn{cusp}--\refn{jepsilon}. The
following diagram is commutative
$$\matrix{\dirsum_{\chi\in\X^U_{R/\m^m}}{\bf M}\bigl(R/\m^m,\mu_N,\chi\bigr)
& \maprighto{r} &R_{m,n}  \cr & \searrow & \Big\downarrow \cr & &
R/\m^m\tensor_\ZZ
\ZZ\bigl(\bigl(\A,\B,\sigma_\beta\bigr)\bigr).\cr}$$The notation
is the following. The ring
$\ZZ\bigl(\bigl(\A,\B,\sigma_\beta\bigr)\bigr)$ is the base over
which $\Tate(\A,\B)_{\sigma_\beta}$ lives by~\refn{tateobjects},
the vertical map is the unique one so that the pull-back of the
universal abelian scheme over~$\MM\bigl(R/\m^m,\mu_{p^nN}\bigr)$
is~$\Tate(\A,\B)_{\sigma_\beta}$  and the diagonal map is the
$q$-expansion map  at the
cusp~$\bigl(\A,\B,\varepsilon,\j_\varepsilon\bigr)$  defined
in~\refn{qexpansion}.
\endprop
\Proof Claim~(1) follows from~\refn{actiona(chi)}. Claim~(2)
follows from~\refn{qexpansion} and~\refn{qexpa(chi)}.

\label kerq. proposition\par\prop Assume that $m=n=1$. In
particular, $R=k$. Then
\spacing
\item{{\rm 1)}} the
kernel of~$r$ consists of the ideal~${\cal I}$
$${\cal I}:=\langle{h_{\P,i}-1\,:\,\P\vert p,\,1\leq i\leq
f_\P}\rangle,$$where the $h_{\P,i}$'s are the partial Hasse
invariants defined in~\refn{hPi};

\spacing
\item{{\rm 2)}} the map $r$ is surjective. In particular, the
ring~$R_{1,1}$ is canonically isomorphic to the ring of
$\I$-polarized modular forms $\dirsum_{\chi\in \X_k}  {\bf M}
\bigl(R/\m^m,\mu_N,\chi\bigr)$ modulo the ideal defined by the
kernel of the $q$-expansion map.

\endcor
\Proof We prove~(1). The fact that the $q$-expansion map is zero
on the ideal~${\cal I}$ follows from~\refn{qexphPi}. Let
$$f=\sum_\chi f_\chi \in \dirsum_{\chi \in
\X_k} {\bf M}\bigl(k,\mu_N,\chi\bigr)$$be the sum of non-zero
$\I$-polarized modular forms~$f_\chi$ of weight~$\chi$ such that
$r(f)=0$. Thanks to part~(3) of~\refn{qexphPi}, multiplying the
modular forms~$f_\chi$'s by suitable powers of the partial Hasse
invariants, we can assume that the set of weights~$\{\chi\}$
appearing in the decomposition of~$f$ does not contain two
distinct elements defining the same character $\bigl(O_L/p
O_L\bigr)^* \rightarrow k^*$. Consider the function
$r(f)=\sum_\chi\phi^*\bigl(f_\chi\bigr)/a(\chi)$. It is the
constant function~$0$ by assumption. By~(2) of~\refn{propr} the
group~$\bigl(O_L\tensor_\ZZ k \bigr)^*$ acts
on~$\phi^*\bigl(f_\chi\bigr)/a(\chi)$ via the character~$\chi$.
Hence, each~$\phi^*\bigl(f_\chi\bigr)/a(\chi)$ is zero i.~e.,
$f_\chi=0$. We conclude that~$f=0$. We prove claim~(2). The Galois
group~$G$ of~$\MM\bigl(k,\mu_{pN}\bigr)^\Kum \rightarrow
\MM\bigl(k,\mu_N\bigr)^{\rm ord}$ is isomorphic to~$\prod_{\P\vert
p } \k_\P^*$ and, hence, has order prime to~$p$. For every
$$\chi\in \X_k/\X_k(1)\cong \Hom_{\rm Gr}\bigl((O_L/p O_L)^*,k^*
\bigr),$$let $$R_{1,1}^\chi:=\left\{ b\in R_{1,1}\vert g \cdot
b=\chi(g)\, b\quad \forall g\in G \right\}.
$$By Kummer theory we get a direct decomposition
$$R_{1,1}=\dirsum_{\chi\in \X_k/\X_k(1)} R_{1,1}^\chi$$into
$R_{1,0}$-modules of rank~$1$. If $\chi\in\X_k$ and~$b\in
R_{1,1}^\chi$, then $b\cdot a(\chi)$ is a modular form of
weight~$\chi$ on~$\MM\bigl(k,\mu_{pN}\bigr)^\Kum$.
By~\refn{actiona(chi)} it descends to a $\I$-polarized modular
form of weight~$\chi$ on~$\MM\bigl(k,\mu_N\bigr)^{\rm ord}$. By
multiplying it by a suitable power~$h^s$ of the Hasse invariant,
defined in~\refn{hPi}, we may assume that~$b\cdot a(\chi)h^s$
extends to a modular form~$f$ of weight~$\chi\cdot \Norm^{s(p-1)}$
on~$\MM\bigl(k,\mu_N\bigr)^\R$. By construction, $r(f)=b$.

\ssection Exotic modular forms\par We prove the existence of
$\I$-polarized modular forms of weight~$\chi$ over artinian bases
for non-universal~$\chi$. A fortiori such modular forms can not be
lifted to characteristic~$0$. We first need the following
\endssection

\label condi. lemma\par\lemma Suppose that $p$ is {\it ramified}.
Let~$k$ be a finite field and let~$m>0$ be an integer. There exist
infinitely many characters $\chi\in \X_{\WW_{m+1}(k)}$ such that
\item{{\rm 1)}} $\chi$ is non-trivial and is not a universal character;
\item{{\rm 2)}} $\chi\bigl((O_L/p^{m+1}O_L)^*\bigr)= 1$ in~$\WW_{m+1}(k)$.
\endlemma
\Proof We use the notation and the results of~\refn{exoticchar}.
Let $\chi \in \Ker(\alpha)$ be a non-trivial character
of~$\G_{\WW_{m+1}(k)}$. Let $\tilde{\chi}\colon \G_k \rightarrow
\GG_{a,k}$ be the associated non-trivial homomorphism of group
schemes. Denote by~$\sigma$ the absolute Frobenius on~$k$.
Since~$k$ is a finite field,  there exists a positive integer~$s$
such that~$k$ is killed by~$ \sigma^s-{\rm Id}$. Let $b \in
\G\bigl(\WW_{m+1}(k)\bigr)$
(e.g.~$b\in\G(\ZZ/p^{m+1}\ZZ)=(O_L/p^{m+1}\ZZ)^*$) and let
$\bar{b}\in \G_k(k)$ be the reduction of~$b$ modulo~$p$. Then
$\chi \circ b =1$ if and only if $\tilde{\chi} \circ \bar{b}=0$.
In particular, the character $\chi'\in \Ker(\alpha)$ associated to
$(\sigma^s-{\rm Id})\circ \tilde{\chi}$ is non-trivial and
kills~$\bar{b}$. Therefore, $\chi'$ satisfies the requirements of
the lemma.

\label exoticmodfor. section\par \ssection Construction\par  Let
$\chi$ be as in~\refn{condi}. We have constructed
in~\refn{actiona(chi)} a section~$a(\chi)$ of~$\phi^*\bigl({\cal
L}_\chi\bigr)$. By the properties of~$\chi$ and
using~\refn{actiona(chi)}, it descends to a section of~${\cal
L}_\chi$ over~$\MM\bigl(\WW_{m+1}(k),\mu_N\bigr)^{\rm ord}$.
Let~$h$ be the Hasse invariant defined on~$\MM(k,\mu_N)^\R$
in~\refn{qexphPi}. We shall prove in~\refn{lifting} that there
exists an integer~$t$ such that~$h^t$ lifts to a modular
form~$\hslash$ over~$\MM\bigl(\WW_{m+1}(k),\mu_N\bigr)^\R$. It has
the property that it vanishes exactly on the complement of the
ordinary locus. In particular, there exists an integer~$s$ such
that~$a(\chi) \hslash^s$ extends to a global section of~${\cal
L}_{\chi_1}$ over~$\MM\bigl(\WW_{m+1}(k),\mu_N\bigr)^\R$,
where~$\chi_1=\chi\cdot \Norm^{st(p-1)}$. Note that~$\chi_1$ is
non-trivial and is not a universal character.
\endssection

\endsection

\section Reduceness of the partial Hasse invariants\par
\noindent We prove in this section, using local deformation theory
of displays, that the divisor of a partial Hasse invariant is
reduced. This is an analogue of Igusa's theorem that states that
the zeroes of the supersingular polynomial are simple. This result
is the basis to establishing a notion of filtration for Hilbert
modular forms of any (not necessarily parallel) weight, and for
later computations regarding the operators~$U$,~$V$
and~$\Theta_{\P, i}$.

\label needtofixdisplay. section\par\ssection Deformation theory
of abelian varieties with RM over the Rapoport locus\par Let
$\bigl(A_0,\iota_0,\lambda_0\bigr)$ be a $\I$-polarized abelian
variety with RM by~$O_L$ over a perfect field~$k$ of
characteristic~$p$ satisfying condition~(R) defined
in~\refn{Rapo}. Let~$\bigl(A_0[p^\infty],\iota_0\bigr)$ be the
associated $p$-divisible group with $O_L$-action. Let
$$\bigl({\rm P}_0,{\rm Q}_0,\F_0,\V_0^{-1}\bigr)$$be the
$3n$-display  in the sense of~[\Zink] associated
to~$A_0[p^\infty]$. Note that~$\lambda_0$ and the choices
in~\refn{conv} induce a polarization on~$A_0$, and hence
on~$A_0[p^\infty]$, of degree prime to~$p$. Let
$$\langle\_\,,\_\,\rangle_0\colon {\rm P}_0 \times {\rm P}_0
\llongrightarrow \WW(k)$$be the associated non-degenerate pairing
of $3n$-displays as defined in~[\Zink]. The action of~$O_L$
on~$A_0[p^\infty]$ induces an action of~$O_L$ on~${\rm P}_0$
and~${\rm Q}_0$ such that~$\F_0$ and~$\V_0^{-1}$ are
$O_L$-equivariant. Moreover,
$$\langle l\gamma,\delta\rangle_0 = \langle \gamma, l\delta\rangle_0$$for
all $l\in O_L\tensor_\ZZ \WW(k)$ and for all~$\gamma$ and~$\delta$
in~${\rm P}_0$.

\noindent Note that the $O_L$-action on~$A_0[p^\infty]$ induces a
decomposition over all primes~$\P$ of~$O_L$ over~$p$:
$$A_0[p^\infty]=\prod_{\P\vert p} A_0[\P^\infty].$$Analogously the
$O_L$-action on~$\bigl({\rm P}_0,{\rm Q}_0,\F_0,\V_0^{-1}\bigr)$
and the $O_L$-linearity of~$\F_0$ and~$\V_0^{-1}$ induce a
decomposition
$$\bigl({\rm P}_0,{\rm Q}_0,\F_0,\V_0^{-1}\bigr)= \prod_{\P\vert
p} \bigl(O_{L,\P}\tensor_{O_L}{\rm P}_0,O_{L,\P}\tensor_{O_L}{\rm
Q}_0,\F_0,\V_0^{-1}\bigr).$$For each prime~$\P$ the
display~$\bigl(O_{L,\P}\tensor_{O_L}{\rm
P}_0,O_{L,\P}\tensor_{O_L}{\rm Q}_0,\F_0,\V_0^{-1}\bigr)$ is
associated to the $p$-divisible group~$A_0[\P^\infty]$.
\endssection

\label Canna. proposition\par\prop The $O_L\tensor_\ZZ
\WW(k)$-module ${\rm P}_0$ is free of rank~$2$. There exist
$\alpha$ and\/~$\beta$ in~${\rm P}_0$ such that
$${\rm P}_0=\bigl(O_L\tensor_\ZZ \WW(k)\bigr) \alpha \dirsum
\bigl(O_L\tensor_\ZZ \WW(k)\bigr) \beta$$and
$${\rm Q}_0=\bigl(O_L\tensor_\ZZ \WW(k)\bigr) p\alpha
\dirsum \bigl(O_L\tensor_\ZZ \WW(k)\bigr) \beta.$$Moreover,
$$\T_0:=\bigl(O_L\tensor_\ZZ \WW(k)\bigr) \alpha \qquad
\hbox{{\rm and}}\qquad \L_0=\bigl(O_L\tensor_\ZZ \WW(k)\bigr)
\beta$$are totally isotropic with respect
to~$\langle\_\,,\_\,\rangle_0$.
\endprop

\Proof Noting that ${\rm P}_0 \isomarrow \H_{1,{\rm
crys}}\bigl(A_0/\WW(k)\bigr)$ and using [\Rapoport, Lemma 1.3], we
deduce that~${\rm P}_0$ is a free $O_L\tensor_\ZZ \WW(k)$-module
of rank~$2$. The image~$\bar{{\rm Q}}_0$ of~${\rm Q}_0$
in~$\bar{{\rm P}}_0:={\rm P}_0/p{\rm P}_0$ is isomorphic
to~$\H^0\bigl(A_0^\vee,\Omega^1_{A_0^\vee/k}\bigr)$.
Via~$\lambda_0$ it is isomorphic
to~$\H^0\bigl(A_0,\Omega^1_{A_0/k} \bigr)$. In particular, since
condition~(R) holds, it is a free~$O_L\tensor_\ZZ k$ module of
rank~$1$. Let~$\beta$ be an element of~${\rm Q}_0$
generating~$\bar{{\rm Q}}_0$ as $O_L\tensor_\ZZ k$-module. The
quotient $\bar{P}_0/\bar{{\rm Q}}_0$ is isomorphic
to~$\Hom_k\bigl(\H^0(A_0,\Omega^1_{A_0/k}),k\bigr)$ and, hence, it
is a free $O_L \tensor_\ZZ k$-module of rank~$1$. Let~$\alpha \in
{\rm P}_0$ be an element generating~$\bar{P}_0/\bar{{\rm Q}}_0$ as
$O_L \tensor_\ZZ k$-module. For all~$\gamma\in{\rm P}_0$
$$\langle l \gamma,\gamma\rangle_0=\langle
\gamma,l\gamma\rangle_0=- \langle l\gamma,\gamma\rangle_0.$$We
conclude that $\langle l\gamma,\gamma\rangle_0=0$ for all~$l\in
O_L\tensor_\ZZ \WW(k)$ and all~$\gamma \in {\rm P}_0$. Hence the
conclusion.

\label Not. section\par\label cPi. definition\par\ssection
Notation\par For each prime~$\P$ of~$O_L$ dividing~$p$ and each
integer~$1\leq i\leq f_\P$, let $\e_{\P,i}$ be the associated
idempotent of~$O_L\tensor_\ZZ \WW(k)$ defined in~\refn{basicwt}
and let $\pi_\P\in O_L$ be the generator of the
ideal~$\P\,O_{L,\P}$ chosen in~\refn{notAtion}. The elements
$$\e_{\P,i}^{[1]}:=\e_{\P,i},\quad \e_{\P,i}^{[2]}:=\pi_\P\e_{\P,i},\quad\ldots\quad,
\e_{\P,i}^{[e_\P]}:=\pi_\P^{e_\P-1}\e_{\P,i}$$form a
$\WW(k)$-basis of the module $\bigl(O_L\tensor_\ZZ
\WW(k)\bigr)\cdot\e_{\P,i}$. Let\/~$\alpha\in{\rm P}_0$ be as
in~\refn{Canna}. Define $\alpha_{\P,i}^{[j]}\in {\rm P}_0$ by
$$\alpha_{\P,i}^{[j]}:=\e_{\P,i}^{[j]}
\alpha=\pi_\P^{j-1}\e_{\P,i}\alpha$$for every prime~$\P$ over~$p$
and for all $1\leq i\leq f_\P$ and $1\leq j\leq e_\P$. Denote by
$$\beta_{\P,i}^{[j]}$$the element of\/~$\L_0$ such that
$$\langle\alpha_{\P,i}^{[j]},\beta_{\P,i}^{[j]}\rangle_0=1,
\quad\hbox{{\rm
and}}\quad\langle\alpha_{\Q,s}^{[t]},\beta_{\P,i}^{[j]}\rangle_0=0$$if
$\alpha_{\Q,s}^{[t]}\neq \alpha_{\P,i}^{[j]}$.
\spacing

\noindent For every prime~$\P$ dividing~$p$, define a total
ordering
$$\alpha_{\P,i_1}^{[j_1]}<\alpha_{\P,i_2}^{[j_2]}\qquad\cases
{&if $i_1<i_2$;\cr&if  $i_1=i_2$ and $j_1<j_2$.\cr}$$Analogously
for the $\beta_{\P,i}^{[j]}$'s.  The elements
$\left\{\alpha_{\P,i}^{[j]}\right\}_{i,j}$ form an ordered
$\WW(k)$-basis of the module~$O_{L,\P}\tensor_{O_L}\T_0$.
Analogously the elements $\left\{\beta_{\P,i}^{[j]}\right\}_{i,j}$
form an ordered $\WW(k)$-basis of~$O_{L,\P}\tensor_{O_L}\L_0$. By
construction ${\cal
B}_\P:=\left\{\alpha_{\P,i}^{[j]},\beta_{\P,i}^{[j]}\right\}_{i,j}$
is a symplectic basis for~$\bigl(O_{L,\P}\tensor_{O_L}{\rm
P}_0,\langle\_\,,\_\,\rangle_0\bigr)$ as a $\WW(k)$-module.

\spacing
\noindent Let
$$\left({\matrix{
A_\P &B_\P \cr C_\P & D_\P \cr}}\right)$$be the matrix of
$\F_0\dirsum \V_0^{-1}$
on~$\left(O_{L,\P}\tensor_{O_L}\T_0\right)\dirsum
\left(O_{L,\P}\tensor_{O_L}\L_0\right)$ with respect to the given
basis. Define $c_{\P,i}^{[1]}$ and~$a_{\P,i}^{[1]} $ in~$\WW(k)$
by
$$C_\P\left(\alpha_{\P,i-1}^{[1]}\right)=c_{\P,i}^{[1]}\,\beta_{\P,i}^{[1]}+ \cdots\qquad
\hbox{{\rm and }} \qquad
A_\P\left(\alpha_{\P,i-1}^{[1]}\right)=a_{\P,i}^{[1]}\,\alpha_{\P,i}^{[1]}+\cdots.$$
\endssection

\rmk For any prime $\P$ and~$\Q$ of~$O_L$ over~$p$, we have by
construction
$$\langle\alpha_{\Q,s}^{[t-1]},\pi_\Q\beta_{\P,i}^{[j]}\rangle_0=
\langle\pi_\Q\alpha_{\Q,s}^{[t-1]},\beta_{\P,i}^{[j]}\rangle_0
=\cases{\,\langle\alpha_{\Q,s}^{[t]},\beta_{\P,i}^{[j]}\rangle_0&
if $2\leq t\leq e_\Q$
\cr\,\langle\gamma,\beta_{\P,i}^{[j]}\rangle_0 & $\gamma\in
p\,{\rm P}_0$, if $t=e_\Q+1$.\cr}$$Hence, for any $2\leq j\leq
e_\P$, we have $\pi_\P \beta_{\P,i}^{[j]}=\beta_{\P,i}^{[j-1]}$
modulo~$p\,{\rm P}_0$.
\endrmk

\label A,C. lemma\par\lemma Let\/~$\P$ be a prime over~$p$;
\spacing
\item{{\rm 1.}} the $({f_\P e_\P}) \times 2\,({f_\P e_\P})$ matrix
$$\left(\matrix{ A_\P\cr C_\P}\right)$$has rank~$g$;
\spacing
\item{{\rm 2.}} for any $i,j$ as
above we have
$$A_\P\bigl(\alpha_{\P,i}^{[j]}\bigr) \in \dirsum_{s=1}^{e_\P} \WW(k)\,\alpha_{\P,i+1}^{[s]}
\qquad \hbox{{\rm and}}\qquad C_\P \bigl(\alpha_{\P,i}^{[j]}
\bigr)\in \dirsum_{s=1}^{e_\P}\WW(k)\,
\beta_{\P,i+1}^{[s]}.$$Analogously,
$$B_\P\bigl(\beta_{\P,i}^{[j]}\bigr) \in \dirsum_{s=1}^{e_\P} \WW(k)\,\alpha_{\P,i+1}^{[s]}
\qquad \hbox{{\rm and}}\qquad D_\P \bigl(\beta_{\P,i}^{[j]}
\bigr)\in \dirsum_{s=1}^{e_\P}\WW(k)\, \beta_{\P,i+1}^{[s]};$$
\spacing
\item{{\rm 3.}} the matrix~$A_\P$ is invertible if and only
if~$A_0[\P^\infty]$ is ordinary.

\endlemma
\Proof The first assertion follows since the map~$\F_0\dirsum
\V_0^{-1}$ is an isomorphism. The second assertion follows the
$\sigma$-linearity of\/~$\F_0$ and of~$\V_0^{-1}$ and from the
definition of the elements~$\alpha_{\P,i}^{[j]}$
and~$\beta_{\P,i}^{[j]}$. Note that~$A_0[\P^\infty]$ is ordinary
if and only if the reduction of~$\F_0$
on~$(O_L/\P)\tensor_{O_L}\bigl(\bar{P}_0/\bar{{\rm
Q}}_0\bigr)=\T_0/\P\T_0$ is an isomorphism. This proves the last
assertion.

\label prank0. corollary\par\cor Let\/~$\P$ be a prime over~$p$.
Then, either~$A_0[\P^\infty]$ is ordinary or it is connected.
\endcor \Proof The $3n$-display defined by~$A_0[\P^\infty]$
is~$\bigl(O_{L,\P}\tensor_{O_L}{\rm P}_0,O_{L,\P}\tensor_{O_L}{\rm
Q}_0,\F_0,\V_0^{-1}\bigr)$. By~\refn{Canna}, we have
that~$O_{L,\P}\tensor_{O_L}\T_0$ is free of rank~$1$
over~$O_{L,\P}\tensor_{\ZZ_p} \WW(k)$. Consider the
reduction~$\bar{A}_\P$ on~$\T_0/\P\T_0$ of the matrix~$A_\P$.
By~\refn{A,C}, Part~(2) either $\bar{A}_\P$ is  invertible or it
nilpotent. The $p$-divisible group~$A_0[\P^\infty]$ is connected
if and only if~$A_0[\P]$ is connected. This is equivalent to ask
that~$\bar{A}_\P$ is nilpotent. We conclude by~\refn{A,C},
Part~(3).

\label equi. proposition\par\prop Let\/~$\P$ be a prime over~$p$.
Assume that\/~$A_0[\P^\infty]$ has $p$-rank equal to~$0$
(equiv.~is non-ordinary by~\refn{prank0}). The universal
equi-characteristic deformation space of~$A_0[\P^\infty]$, as a
principally polarized $p$-divisible group, is
$R_\P:=k[\![t_{a,b}]\!]_{1\leq a,b\leq f_\p e_\P}$ with the
relations~$t_{a,b}=t_{b,a}$. The universal display, denoted
by\/~$\bigl({\rm P}_\P,{\rm Q}_\P,\F_\P,\V_\P^{-1} \bigr)$, is
given by
$$ {\rm P}_\P:=\left(O_{L,\P}\tensor_{O_L}{\rm
P}_0\right)\tensor_{\WW(k)} \WW(R_\P), \quad {\rm
Q}_\P:=\left(O_{L,\P}\tensor_{O_L}{\rm Q}_0\right)
\tensor_{\WW(k)} \WW(R_\P)$$and $$
\L_\P:=\left(O_{L,\P}\tensor_{O_L}\L_0\right) \tensor_{\WW(k)}
\WW(R_\P),\quad \T_\P:=\left(O_{L,\P}\tensor_{O_L}\T_0\right)
\tensor_{\WW(k)} \WW(R_\P).$$The matrix of\/~$\F_\P\dirsum
\V_\P^{-1}$ on~$\T_\P \dirsum\L_\P$ with respect to~${\cal B}_P$
is
$$\left({\matrix{ A_\P+ T_\P C_\P & B_\P+ T_\P D_\P \cr C_\P & D_\P \cr}}\right),$$where
$T_\P$ is the symmetric matrix of Teichm\"{u}ller
lifts~$\bigl(w(t_{a,b})\bigr)_{1\leq a,b\leq e_\P f_\P}$. The
pairing
$$\langle\_\,,\_\,\rangle_\P,$$defined
extending $O_{L,\P}\tensor_{O_L}\WW(R_\P)$-linearly the
pairing~$\langle\_\,,\_\,\rangle_0$, is a non-degenerate pairing
of displays such that $\L_\P$ and~$\T_\P$ are maximal isotropic
submodules of~${\rm P}_\P$.
\endprop
\Proof It follows form the assumption on the $p$-rank
of~$A_0[p^\infty]$ that~$({\rm P}_0,{\rm Q}_0,\F_0,\V_0^{-1} )$ is
a display. The theorem follows from~[\Zink].

\label condipasta. section\par\ssection The  universal
equi-characteristic deformation space of~$(A_0,\iota_0)$\par By
the Serre-Tate theorem the universal equi-characteristic
deformation space~$\Spf\bigl(R_\iota\bigr)$ of~$A_0$ with the
$O_L$-action  coincides with the universal equi-characteristic
deformation space~$\Spf\bigl(R_\iota\bigr)$ of~$A_0[p^\infty]$
with the $O_L$-action. Hence,
$$\Spf\bigl(R_\iota\bigr)=\prod_{\P\vert p} \Spf\bigl(R_{\P,\iota}\bigr),$$where for
any prime~$\P$ we define~$\Spf\bigl(R_{\P,\iota}\bigr)$ as the
universal deformation space of~$A_0[\P^\infty]$ with $O_L$-action.

\spacing Fix a prime~$\P$. If~$A_0[\P^\infty]$ is
non-ordinary,~$\Spf\bigl(R_{\P,\iota}\bigr)$  is the closed
subscheme of~$\Spf(R_\P)$ defined by the condition that~$\F_\P$
and~$\V_\P^{-1}$ commute with the $O_L$-action on~${\rm P}_\P$
and~${\rm Q}_\P$ induced by the $O_L$-structures of~${\rm P}_0$
and~${\rm Q}_0$. Since
$\big\langle\F_\P(x),\V_\P^{-1}(y)\big\rangle=\langle
x,y\rangle^\sigma$ for all $x$ in~${\rm P}_\P$ and all~$y$
in~${\rm Q}_\P$, we deduce that~$\F_\P\dirsum \V_\P^{-1}$ is a
symplectic isomorphism. Hence, $\V_\P^{-1}$ is $O_L$-linear if and
only if~$\F_\P$ is. This is equivalent to require that~$\F_\P$
restricted to~$\T_\P$ i.~e.,~$A_\P+T_\P C_\P$, is equivariant with
respect to the $O_L$-structure on~$\T_\P$.
\endssection

\label condiriso. lemma\par\lemma The conditions that the
restriction of\/~$\F_\P$ to~$\T_\P$ is $O_L$-equivariant are
linear in the $t_{a,b}$'s and are equivalent to the conditions
\spacing
\item{{\rm 1.}} $\F_\P\bigl(\alpha_{\P,i}^{[1]} \bigr)  \in
\WW(R) \alpha_{\P,i+1}^{[1]}\dirsum\ldots\dirsum \WW(R)
\alpha_{\P,i+1}^{[e_\P]}$ for all  $1\leq i \leq f_\P$;
\spacing
\item{{\rm 2.}} $\F_\P\bigl(\alpha_{\P,i}^{[j]}
\bigr)=\pi_\P^{j-1}\,\F\bigl(\alpha_{\P,i}^{[1]} \bigr)$ for  all
$1\leq i \leq f_\P$ and all $1\leq j\leq e_\P$.
\endlemma
\Proof Clear.

\defi For any prime~$\P$ and integers
$1\leq i\leq f_\P$ and $1\leq j\leq e_\P$ let
$$t_{\P,i}^{[j]}:=t_{a,b}$$with $a=(i-1)e_\P+1$
and~$b=(i-1)e_\P+1+(j-1)$.
\enddefi

\label unidefRM. theorem\par\thm Let\/~$\P$ be a prime over~$p$.
Assume that~$A_0[\P^\infty]$ is non-ordinary.
Let\/~$\Spf\bigl(R_{\P,\iota}\bigr)$ be the universal
equi-characteristic deformation space of\/~$A_0[\P^\infty]$ as
$p$-divisible group with $O_L$-action. Then, $R_{\P,\iota}$ is a
$(f_\P e_\P) $-dimensional power series ring
$$R_{\P,\iota}=k[\![t_{\P,i}^{[j]}]\!]_{i,j}.$$The
associated universal display~$\bigl({\rm P}_\P,{\rm
Q}_\P,\F_\P,\V_\P^{-1} \bigr)$ is given by

$$ {\rm P}_\P:=\left(O_{L,\P}\tensor_{O_L}{\rm
P}_0\right)\tensor_{\WW(k)} \WW(R_{\P,\iota}), \quad {\rm
Q}_\P:=\left(O_{L,\P}\tensor_{O_L}{\rm Q}_0\right)
\tensor_{\WW(k)} \WW(R_{\P,\iota})$$and $$
\L_\P:=\left(O_{L,\P}\tensor_{O_L}\L_0\right) \tensor_{\WW(k)}
\WW(R_{\P,\iota}),\quad
\T_\P:=\left(O_{L,\P}\tensor_{O_L}\T_0\right) \tensor_{\WW(k)}
\WW(R_{\P,\iota}).$$The matrix of\/~$\F_\P\dirsum \V_\P^{-1}$
on~$\T_\P \dirsum\L_\P$ with respect to~${\cal B}_P$ is
$$\left({\matrix{ A_\P+ T_\P C_\P & B_\P+ T_\P D_\P \cr C_\P & D_\P
\cr}}\right),
$$where $T_\P$ is the  matrix determined by
$$T_\P\bigl(\alpha_{\P,i}^{[j]}\bigr)=\sum_{l=1}^{e_\P}w\bigl(t_{\P,i}^{[l]}\bigr)
\,\pi_\P^{j-1}\,\alpha_{\P,i}^{[l]}$$for all integers $1\leq i\leq
f_\P$ and $1\leq j\leq e_\P$.
\endthm
\Proof By~\refn{condipasta} and~\refn{condiriso}, the formal
scheme~$\Spf\bigl(R_{\P,\iota}\bigr)$ is the universal deformation
space of $\bigl(O_{L,\P}\tensor_{O_L}{\rm
P}_0,O_{L,\P}\tensor_{O_L}{\rm Q}_0,\F_0,\V^{-1}_0\bigr)$ as
display with~$O_L$-action. By~[\Zink], it is also the universal
deformation space of~$A_0[\P^\infty]$ with the $O_L$-action.

\cor The notation is as above. The formal scheme
$\Spf\bigl(R_{\P,\iota}\bigr)$ is formally smooth of
dimension~$f_\P e_\P$.
\endcor

\ssection Example: the inert case\par In this case $p$ remains a
prime ideal in~$O_L$. We omit the subscripts~$\P$ and~$[j]$ in the
formulas above. The matrix $A+TC$ mod~$p$ of~\refn{unidefRM} is
$$\left(\matrix{
0 & 0 &\ldots &0 & \bar{a}_1+ t_1 \bar{c}_1 \cr
  \bar{a}_2+t_2 \bar{c}_2 & 0 & \ldots & 0&0\cr
    0        &\bar{a}_3+t_3 \bar{c}_3 & \ldots & 0&0\cr
 \vdots & \vdots & \vdots & \vdots &\vdots\cr
 0 & 0 & \ldots & \bar{a}_g+ t_g \bar{c}_g& 0\cr}\right).$$

\endssection

\ssection Example: the totally ramified case\par Suppose that
there is only one prime over~$p$ which is totally ramified. In
this case we omit the subscripts~$\P$ and~$i$. We write the
coefficients $\bar{a}^{[j]}_{\P, i}, \bar{c}^{[j]}_{\P, i}$ and
the variables~$t^{[j]}_{\P, i}$ as $\bar{a}_{[j]}, \bar{c}_{[j]}$
and~$t_{[j]}$, respectively. The matrix $A+TC$ mod~$p$
of~\refn{unidefRM} is
$$\left(\matrix{
 \bar{a}_{[1]}+t_{[1]} \bar{c}_{[1]} & 0 & \ldots & 0\cr
 \bar{a}_{[2]}+t_{[2]} \bar{c}_{[1]}+t_{[1]} \bar{c}_{[2]}&\bar{a}_{[1]}+t_{[1]} \bar{c}_{[1]} & \ldots & 0\cr
 \vdots & \vdots & \vdots & \vdots \cr
 \bar{a}_{[g]}+\sum_{i=1}^g t_{[g-i+1]} \bar{c}_{[i]} & \bar{a}_{[g-1]}+ \sum_{i=1}^{g-1}
 t_{[g-i]}
 \bar{c}_{[i]}&\ldots & \bar{a}_{[1]}+t_{[1]} \bar{c}_{[1]}
 \cr}\right).$$

\endssection

\label HASSE. section\par\ssection The Hasse-Witt matrix of the
universal equi-characteristic deformation\par Let
$\bigl(A,\iota,\lambda\bigr)$ be the universal equi-characteristic
object over the universal deformation space
of~$\bigl(A_0,\iota_0,\lambda_0\bigr)$. Let
$$(A^\P,\iota^\P,\lambda^\P):=(A,\iota,\lambda)
\fibprod_{\Spf(R_\iota)}\Spf(R_{\P,\iota});$$the morphism
$\Spf(R_{\P,\iota})\rightarrow \Spf(R_\iota)=\prod_{\Q\vert p}
\Spf(R_{\Q,\iota})$ (see~\refn{condipasta} for the last equality)
is defined to be the identity on the factor~$\Spf(R_{\P,\iota})$
and the map $\Spf(R_{\P,\iota}) \rightarrow \Spf(k) \rightarrow
\Spf(R_{\Q,\iota})$ if~$\Q\neq \P$. Let~$I_{R_{\P,\iota}}$ be the
kernel of the reduction map $\WW\bigl(R_{\P,\iota}\bigr)
\rightarrow R_{\P,\iota}$. By~[\Messing], we get a canonical
identification of~$\H_{1,{\rm dR}}\bigl(A^\P/R_{\P,\iota}\bigr)$
with the Lie algebra of the universal vector extension
of~$A^\P[p^\infty]$. By the definition of~$R_{\P,\iota}$
and~$A^\P$ we have
$$A^\P[p^\infty]\cong A^\P[\P^\infty]\fibprod \prod_{\Q\neq \P}
\left(A_0[\Q^\infty]\fibprod_k \Spf(R_{\P,\iota})\right).$$Hence,
for any prime~$\Q$ over~$p$ different from~$\P$ we get a canonical
isomorphism $(O_L/\Q) \tensor_{O_L} \H_{1,{\rm
dR}}\bigl(A^\P/R_{\P,\iota}\bigr)\cong \bigl({\rm P}_0/\Q{\rm
P}_0\bigr)\tensor_k \R_{\P,\iota}$. By~\refn{unidefRM} we have a
canonical identification
$$(O_L/\P) \tensor_{O_L} \H_{1,{\rm dR}}\bigl(A^\P/R_{\P,\iota}\bigr)\cong
{\rm P}_\P/I_{R_{\P,\iota}}{\rm P}_\P.$$The exact sequence

\bigskip \noindent $0 \llongrightarrow
\Hom\Bigl(\H^1\bigl(A^\P,O_{A^\P}\bigr),R_{\P,\iota}\Bigr)
\llongrightarrow \H_{1,{\rm dR}}\bigl(A^\P/R_{\P,\iota}\bigr)$

\spacing
\rightline{$\llongrightarrow
\Hom\Bigl(\H^0\bigl(A^\P,\Omega^1_{A^\P/R_{\P,\iota}}\bigr),R_{\P,\iota}\Bigr)
\llongrightarrow 0$,}\bigskip \noindent tensored over~$O_L$
with~$O_L/\P$ (resp.~$O_L/\Q$ for~$\Q\neq \P$), and the exact
sequence
$$0 \llongrightarrow \L_\P/I_{R_{\P,\iota}}\L_\P
 \llongrightarrow {\rm P}_\P/I_{R_{\P,\iota}}{\rm P}_\P
\llongrightarrow \T_\P/I_{R_{\P,\iota}}\T \llongrightarrow
0$$(resp. $$0 \llongrightarrow  \bigl(\L_0/\Q\L_0\bigr)\tensor_k
\R_{\P,\iota} \llongrightarrow  \bigl({\rm P}_0/\Q{\rm
P}_0\bigr)\tensor_k \R_{\P,\iota}\llongrightarrow
\bigl(\T_0/\Q\T_0\bigr)\tensor_k \R_{\P,\iota}\llongrightarrow
0)$$are identified. Using~$\lambda^\P$ and the choices
in~\refn{conv}, we obtain a polarization on~$A$ of degree prime
to~$p$. This induces perfect pairings
between~$\H^0\bigl(A^\P,\Omega^1_{A^\P/R_{\P,\iota}}\bigr)$ and~$
\H^1\bigl(A^\P,O_{A^\P}\bigr)$, and between~$\L/I_{R_\iota}\L$
and~$ \T_\P/I_{R_{\P,\iota}}\T_\P$ (resp.~$
\bigl(\T_0/\Q\T_0\bigr)\tensor_k \R_{\P,\iota} $ and~$
\bigl(\L_0/\Q\L_0\bigr)\tensor_k \R_{\P,\iota}$), compatible with
the identifications given above.

\noindent Hence, we get  canonical isomorphisms
$$\H^0\bigl(A^\P,\Omega^1_{A^\P/R_{\P,\iota}}\bigr)
\isomarrow \L_0\tensor_{\WW(k)} \R_{\P,\iota}$$and
$$\H^1\bigl(A^\P,O_{A^\P}\bigr) \isomarrow
\T_0\tensor_{\WW(k)} \R_{\P,\iota}$$so that Frobenius on the LHS,
induced by Frobenius on~$O_{A^\P}$, corresponds to Frobenius on
the RHS.

\noindent By~\refn{Not}, the $k$-vector space~$\T_0/p\T_0$
(resp.~$\L_0/p\L_0 $) is endowed with $k$-generators
$\bar{\alpha}_{\Q,i}^{[j]}$ (resp.~$\bar{\beta}_{\Q,i}^{[j]}$)
defined as the reduction of~$\alpha_{\Q,i}^{[j]}$
(resp.~$\beta_{\Q,i}^{[j]}$). Their images induce canonical
$R_{\P,\iota}$-generators
$$\bigl\{\eta_{\Q,i}^{[j]}\bigr\}_{\Q,i,j}\subset \H^1\bigl(A^\P,O_{A^\P}\bigr))$$ (resp.
$$\bigl\{\omega_{\Q,i}^{[j]}\bigr\}_{\Q,i,j}\subset \H^0\bigl(A^\P,\Omega^1_{A^\P/R_{\P,\iota}}\bigr)).$$
\endssection

\lemma Let\/~$\P$ be a prime over~$p$.
Let\/~$\left\{\Q_1,\ldots,\Q_d\right\}$ be the set of primes
over~$p$ different from~$\P$. Assume that~$A_0[\P^\infty]$ is not
ordinary. The Hasse-Witt matrix
of~$\bigl(A^\P,\iota^\P,\lambda^\P\bigr)$ with respect to the
basis~$\bigl\{\eta_{\Q,i}^{[j]}\bigr\}_{\Q,i,j}$ is canonically
identified with the reduction of the matrix
$$\left(\matrix{
A_{\Q_1} & 0 &  0& 0\cr
 0    &\cdots & A_{\Q_d}&  0\cr
 0 & \cdots &  0 & A_\P+T_\P C_\P}\right)$$via the quotient
map~$\WW\bigl(R_\iota\bigr) \rightarrow R_\iota$.
\endlemma

\label hvst. corollary\par\cor Let~$\P$ be a prime of~$O_L$
over~$p$. Assume that~$A_0[\P^\infty]$ is not ordinary. Let~$i$ be
an integer satisfying $1\leq i\leq f_\P$. With the notation
of\/~\refn{hPi}, we have
$$h_{\P,i}\bigl(A,\iota,\lambda,
\omega_\alpha\bigr)=\bar{a}_{\P,i}^{[1]}+\bar{c}_{\P,i}^{[1]}
t_{\P,i}^{[1]}.$$The elements~$\bar{a}_{\P,i}^{[1]}$
and\/~$\bar{c}_{\P,i}^{[1]}$ are the reduction mod~$p$ of the
element\/~$a_{\P,i}^{[1]}$ and\/~$c_{\P,i}^{[1]}$ defined
in~\refn{cPi}. If\/~$\bar{a}_{\P,i}^{[1]}=0 $,
then~$\bar{c}_{\P,i}^{[1]}$ is invertible by~\refn{A,C}.
\endcor

\label redirr. corollary\par\cor The zero locus~$W_{\P,i}$ of the
partial Hasse invariant\/~$h_{\P,i}$ is a reduced, non-singular
divisor. In particular, it is locally irreducible.

\spacing \noindent The divisor of the Hasse-invariant is equal to~$\sum_{\P,
i} e_\P W_{\P, i}$ and the~$W_{\P, i}$ are normal crossing
divisors.
\endcor

\prop  Let  $f\in {\bf M}\bigl(k,\mu_N,\chi\bigr)$ be a
$\I$-polarized modular form over~$k$ of weight~$\chi$. There is a
{\it unique} $\I$-polarized modular form~$g$ over~$k$ having the
same $q$-expansion as\/~$f$ at a (any) cusp and such that
if\/~$g'$ is a $\I$-polarized modular form over~$k$ with the same
$q$-expansion of\/~$f$, then there exist non-negative
integers~$b_{\P,i}$ for each prime~$\P$ of\/~$O_L$ over~$p$ and
each integer $1\leq i\leq f_\P$ such that
$$g'=g\, \prod_{\P,i} h_{\P,i}^{b_{\P,i}};$$see~\refn{hPi} for the
definition of the partial Hasse invariants~$h_{\P,i}$.
\endprop
\Proof Define $g$ to be a $\I$-polarized modular form satisfying
$$f=g\prod_{\P,i} h_{\P,i}^{a_{\P,i}},$$where the $a_{\P,i}$ are
chosen maximal non-negative so that~$g$ is a holomorphic modular
form. The modular form~$g$ is unique with this property by the two
previous corollaries. By~\refn{qexphPi}, the modular form~$g$ has
the same $q$-expansion of~$f$. It satisfies the requirement of the
proposition by~\refn{kerq}.

\label filtrations. section\par\ssection Filtrations on modular
forms\par The notation is as above. Define the {\it filtration}
of~$f$, denoted by $$ \Phi(f),$$to be the weight of the {\it
unique} $\I$-polarized modular form~$g$  with the properties
described in the proposition. Since the weight of~$h_{\P,i}$ is
$\chi_{\P,i}^p\chi_{\P,i}^{-1}$, we have
$$\chi=\Phi(f)\,\prod_{\P,i}
\Bigl(\chi_{\P,i-1}^p\chi_{\P,i}^{-1}\Bigr)^{a_{\P,i}}$$for
suitable non-negative integers~$a_{\P,i}$.
\endssection
\endsection

\section A compactification of $\MM\bigl(k,\mu_{pN}\bigr)$\par
\noindent  Fix a field $k$ of characteristic~$p$ containing all
the finite fields~$\k_\P$; see~\refn{k}. Let $N\geq 4$ be an
integer prime to~$p$. In this section we construct a
compactification of~$\MM(k, \mu_{pN})$, which is well suited for
the study of the arithmetic of modular forms. First, the
compactification we construct is normal and its only singularities
come from the singularities of~$\MM(k, \mu_N)$. Secondly, it is
explicit, in the sense that it is defined, up to codimension~$2$,
as the scheme resulting from adjoining to~$\MM(k, \mu_N)$ roots of
explicitly given modular forms. The notation is as in~\refn{phi}.

\label phibar. definition\par\defi Let $\MMbar\bigl(k,\mu_N\bigr)$
be the minimal compactification of~$\MM\bigl(k,\mu_N\bigr)$
constructed in~{\rm [\Chai, Thm.~4.3]}. It is a projective normal
scheme over~$k$ obtained by adding finitely many cusps. Its
singular locus consists precisely of the complement of the
Rapoport locus~$\MM\bigl(k,\mu_N\bigr)^\R$ defined in~\refn{Rapo}.
\spacing
\noindent Define
$$\phibar\colon\MMbar\bigl(k,\mu_{pN}\bigr)^\Kum\llongrightarrow\MMbar\bigl(k,
\mu_N\bigr)$$as the normal closure
of~$\MMbar\bigl(k,\mu_{N}\bigr)$
in~$\MM\bigl(k,\mu_{pN}\bigr)^\Kum$ via the Galois cover~$\phi$
with group~$G=G_{1,1}$, defined in~\refn{phi}.
\enddefi

\lemma The following properties hold:
\spacing
\item{{\rm 1.}} the morphism~$\phibar$ is finite;
\spacing
\item{{\rm 2.}} the
scheme~$\MMbar\bigl(k,\mu_{pN}\bigr)^\Kum$ is projective,
irreducible and normal;
\spacing
\item{{\rm 3.}} the
scheme~$\MMbar\bigl(k,\mu_{pN}\bigr)^\Kum$ is endowed with an
action of~$G$ and~$\phibar$ represents the quotient map;
\spacing
\item{{\rm 4.}} the branch locus of~$\phibar$ is a divisor
contained in the complement of~$\MMbar\bigl(k,\mu_{N}\bigr)^{\rm
ord}$. We use the convention that the cusps are in the ordinary
locus.
\endlemma
\Proof Since~$\MMbar(k,\mu_N)$ is of finite type over~$k$, it is
excellent and, hence, universally japanese. By~[\EGAIVtwo,
\S7.8.3] we conclude that~$\phibar$ is finite. By~[\EGAII,
Cor.~6.1.11] we deduce that~$\phibar$ is projective and,
by~[\EGAII, Prop.~5.5.5 (ii)],
that~$\MMbar\bigl(k,\mu_{pN}\bigr)^\Kum$ is projective. The
quotient of~$\MMbar\bigl(k,\mu_{pN}\bigr)^\Kum$ by~$G$ is finite
and birational over~$\MMbar\bigl(k,\mu_{N}\bigr)$. We deduce that
it coincides with~$\MMbar\bigl(k,\mu_{N}\bigr)$. This concludes
the proofs of claims (1)-(3). By purity of branch locus,
see~[\SGA2, X, Thm.~3.4 (i)], the map~$\phibar$ is ramified along
a divisor of~$\MMbar\bigl(k,\mu_{N}\bigr)$. By construction the
pre-image of~$ \MM\bigl(k,\mu_{N}\bigr)^{\rm ord}$
in~$\MMbar\bigl(k,\mu_{pN}\bigr)^\Kum $
is~$\MM\bigl(k,\mu_{pN}\bigr)^\Kum $. Hence, $\phibar$ is \'etale
over~$ \MM\bigl(k,\mu_{N}\bigr)^{\rm ord}$. Since the cusps are
isolated points in the complement of~$
\MM\bigl(k,\mu_{N}\bigr)^{\rm ord}$, the map~$\phibar$ is
unramified also at the cusps. This proves part~(4).

\label Mbar. section\par\ssection Local charts
of~$\MMbar\bigl(k,\mu_{pN}\bigr)^\Kum$\par Fix a prime~$\Q$
of~$O_L$ over~$p$ and an integer $1\leq j\leq f_\Q$. Let
$$\MMbar\bigl(k,\mu_N\bigr)^\R\hookrightarrow\MMbar\bigl(k,\mu_N\bigr),$$be
the locus where condition~(R) holds. The convention is that the
cusps satisfy~(R).  We are going to give an explicit description
of~$\phibar$ in a neighborhood of the generic point of the
divisor~$W_{\Q,j}$ defined by the partial Hasse
invariant~$h_{\Q,j}$; see~\refn{redirr} . Define a scheme
$$\phibar^\R_{\Q,j}\colon \MMbar\bigl(k,\mu_{pN}\bigr)_{\Q,j}^{\Kum,
\R}\llongrightarrow \MMbar\bigl(k,\mu_N\bigr)^\R\backslash
\sum_{(\P,i)\neq (\Q,j)}W_{\P,i}$$over the complement
in~$\MMbar\bigl(k,\mu_N\bigr)^\R$ of the divisor~$\sum_{(\P,i)\neq
(\Q,j)}W_{\P,i}$ by adjoining a $p^{f_\P}-1$-th root of the
modular forms~$h_{\P,i+1}^{p^{f_\P-1}}\,
h_{\P,i+2}^{p^{f_\P-2}}\cdots h_{\P,i}$ for any prime~$\P$
of~$O_L$ over~$p$ and a fixed integer $1\leq i\leq f_\P$ (we
require $i=j$ if~$\P=\Q$).

\spacing
\noindent {1) \enspace} The identity
$$\left(h_{\P,i+1}^{p^{f_\P-1}}\, h_{\P,i+2}^{p^{f_\P-2}}\cdots
h_{\P,i}\right)^p\left(h_{\P,i+2}^{p^{f_\P-1}}\,
h_{\P,i+3}^{p^{f_\P-2}}\cdots
h_{\P,i+1}\right)^{-1}=h_{\P,i}^{p^{f_\P}-1}$$for all $1\leq i\leq
f_\P$ implies that the construction is independent of~$i$
for~$\P\neq\Q$.

\spacing
\noindent {2) \enspace} For every prime~$\P$ of~$O_L$ over~$p$ and
any integer~$1\leq i\leq f_\P$ we have the following equality of
modular forms on~$\MM\bigl(k,\mu_{pN}\bigr)^\Kum$:
$$\eqalign{a\bigl(\chi_{\P,i}\bigr)^{p^{f_\P}-1}&
=a\bigl(\chi_{\P,i}^p\chi_{\P,i+1}^{-1}\bigr)^{p^{f_\P-1}}\,
a\bigl(\chi_{\P,i+1}^p\chi_{\P,i+2}^{-1}\bigr)^{p^{f_\P-2}}\cdots
a\bigl(\chi_{\P,i-1}^p\chi_{\P,i}^{-1}\bigr)^{p^{0}}\cr
 &=\phi^*\bigl(h_{\P,i+1}\bigr)^{p^{f_\P-1}}\,\phi^*\bigl(h_{\P,i+2}\bigr)^{p^{f_\P-2}}
 \cdots\phi^*\bigl(h_{\P,i}\bigr).\cr}$$See~\refn{actiona(chi)} for the notation.
In particular, we have a commutative diagram
$$\matrix{\MM\bigl(k,\mu_{pN}\bigr)^\Kum & \llongrightarrow
&\MMbar\bigl(k,\mu_{pN}\bigr)^{\Kum,\R}_{\Q,j} \cr\mapdownl{\phi}
& &\mapdownr{\phibar^\R_{\Q,j}}\cr\MM\bigl(k,\mu_N\bigr)^{\rm ord}
&\llongrightarrow &\MMbar\bigl(k,\mu_N\bigr)^\R\backslash
\sum_{(\P,i)\neq (\Q,j)}W_{\P,i}.\cr}$$By~\refn{phi} the LHS is a
finite \'etale morphism, Galois under the group~$G$. In
particular, its  degree is equal to~$\prod_{\P\vert
p}\bigl(p^{f_\P}-1\bigr)$.
\endssection

\label casino. proposition\par\prop The
scheme~$\MMbar\bigl(k,\mu_{pN}\bigr)_{\Q,j}^{\Kum,\R}$ has the
following properties
\spacing
\item{{\rm 1)}} it is irreducible and normal;
\spacing
\item{{\rm 2)}} the morphism $\phibar_{\Q,j}^\R$ is finite and
its branch locus is~$(p^{f_\Q}-2)\,W_{\Q,j}$.
\spacing
\item{{\rm 3)}} there is a commutative
diagram $$\matrix{\MM\bigl(k,\mu_{pN}\bigr)^\Kum &
\hooklongrightarrow& \MM\bigl(k,\mu_{pN}\bigr)_{\Q,j}^{\Kum,\R} &
\hooklongrightarrow &\MMbar\bigl(k,\mu_{pN}\bigr)^\Kum
\cr\mapdownl{\phi} & & \mapdownl{\phibar^R_{\Q,j}} &
&\mapdownr{\phibar}\cr\MM\bigl(k,\mu_N\bigr)^{\rm ord}
&\hooklongrightarrow & \MMbar\bigl(k,\mu_N\bigr)^\R\backslash
\sum_{(\P,i)\neq (\Q,j)}W_{\P,i} &\hooklongrightarrow
&\MMbar\bigl(k,\mu_N\bigr),\cr}$$where the squares are cartesian
and the horizontal arrows are open immersions.

\endprop
\Proof  The finiteness in claim~(2) is clear from the
construction. We prove claim~(1), claim~(2) and claim~(3) for the
square on the LHS. The existence of the diagram and the rest of
claim~(3) follow then from the definition
of~$\MMbar\bigl(k,\mu_{pN}\bigr)^\Kum$ and the finiteness in~(2).

\spacing
\noindent Let $x$ be a closed point of~$W_{\Q,j}$ such that
$x\not\in W_{\P,i}$ for~$(\P,i)\neq (\Q,j)$. Let~$U_x=\Spec(A)$ be
an affine  open neighborhood of~$x$
in~$\MMbar\bigl(k,\mu_N\bigr)^\R\backslash \sum_{(\P,i)\neq
(\Q,j)}W_{\P,i} $ over which every invertible sheaf~${\cal
L}_{\chi_{\P,i}}$ is trivial. Choose an ordering~$\Q=\P_1 < \cdots
<\P_s$ of the primes of~$O_L$ over~$p$. Define
$$B:=A\Bigl[X_1,\ldots,X_s\Bigr]
/\Bigl(X_t^{p^{f_{\P_t}}-1}-\bigl(h_{\P_t,i+1}\bigr)^{p^{f_{\P_t}}}\cdots
\bigl(h_{\P_t,i}\bigr)\Bigr)_{t=1,\ldots,s},$$with the abuse of
notation that, via the chosen trivialization, the
elements~$h_{\P,i}$ are now considered as elements of~$A$. As
remarked in~\refn{Mbar}  the definition does not depend on the
choice of~$1\leq i\leq f_{\P_t}$ for~$t>1$. Then:
\item{{i)}} $B$ is
finite and flat over~$A$ of
degree~$\prod_{t=1}^s\left(p^{f_{\P_t}}-1\right)$;
\item{{ii)}} the group $\prod_{t=1}^s \k_{\P_t}^*$ acts $A$-linearly
on~$B$ through roots of unity.

\noindent Note that $ \phi^{-1}\left(U_x\backslash
W_{\Q,j}\right)\rightarrow U_x\backslash W_{\Q,j}$ is endowed with
an action of~$G_{1,1}\cong \prod_{t=1}^s \k_\P^*$ as remarked
in~\refn{Mbar}. By~\refn{actiona(chi)} the morphism
$\phi^{-1}\left(U_x\backslash W_{\Q,j}\right) \rightarrow
\Spec\bigl(B\bigr)$ defined in~(2) of~\refn{Mbar} is equivariant
with respect to the action of~$\prod_{t=1}^s \k_\P^*$. Hence, the
map
$$\phi^{-1}\left(U_x\backslash
W_{\Q,j}\right) \llongrightarrow
\Spec\left(B\bigl[h_{\Q,j}^{-1}\bigr]\right)
$$is an isomorphism. In particular, we conclude that~$B$ is a
domain. Define $$B_1:=A\Bigl[X_1]
/\Bigl(X_1^{p^{f_{\Q}}-1}-\bigl(h_{\Q,i+1}\bigr)^{p^{f_{\Q}}}\cdots
\bigl(h_{\Q,i}\bigr)\Bigr).$$We know it is a domain. Let $$Y=
\sum_{d=0}^{p^{f_\Q}-2}Y_d X_t^d\qquad \hbox{{\rm with }}Y_d\in
{\rm Frac}\left(A\right)$$be an element in the fraction field
of~$B_1$ which is a zero of a monic polynomial~$g(X)\in A[X]$. For
every prime~$Q$ of~$A$, the equation
$$X_1^{p^{f_{\Q}}-1}-\bigl(h_{\Q,i+1}\bigr)^{p^{f_{\Q}}}\cdots
\bigl(h_{\Q,i}\bigr) $$ in the local ring~$A_Q$ is either
Eisenstein, if $h_{\Q,j}\in Q$, or separable otherwise. Hence,
$Y_d\in A_Q$ for every~$d$. Since~$A$ is normal, we conclude
that~$Y_d\in A$ for every~$d$. In particular, $B_1$ is normal.
Moreover, the extension $A\subset B_1$ is ramified only
along~$h_{\Q,j}=0$ with ramification index~$p^{f_\Q}-2$. The
extension $B_1 \subset B$ is \'etale. Hence, $B$ is normal and~$A
\subset B$  is ramified only along~$h_{\Q,j}=0$ with ramification
index~$p^{f_\Q}-2$.

\cor The
scheme~$\cup_{\Q,j}\MM\bigl(k,\mu_{pN}\bigr)_{\Q,j}^{\Kum,\R}$ is
endowed with an action of~$G$ so that the map~$\phibar^\R$ is the
quotient map.  The open
subscheme~$\cup_{\Q,j}\MM\bigl(k,\mu_{pN}\bigr)_{\Q,j}^{\Kum,\R}$
of~$\MMbar\bigl(k,\mu_{pN}\bigr)^\Kum$ has codimension~$2$.
\endcor
\Proof The first claim is clear. The second follows
from~\refn{propRapo}.

\label branch. corollary\par\cor The branch locus of\/~$\phibar$
is exactly the complement of the ordinary locus
in~$\MMbar\bigl(k,\mu_{N}\bigr)$. For each prime~$\P$ and
each~$1\leq i\leq f_\P$, the ramification index of~$W_{\P,i}$
in~$\MMbar\bigl(k,\mu_{pN}\bigr)^\Kum$ is $p^{f_\P}-2$.
\endcor
\Proof It follows from~\refn{casino}.

\rmk The locus $\MM\bigl(k,\mu_N\bigr)^\R$ is the locus where
the modular forms of level~$\mu_N$ are defined;
see~\refn{modularforms} and~\refn{Lchi}. This is why we need an
explicit description of the map~$\phibar$ over such locus, at
least up to codimension~$2$, as given in the proposition.
\endrmk

\endsection

\section Congruences mod~$p^n$ and Serre's $p$-adic modular forms\par
\noindent There is already a notion of $p$-adic Hilbert modular
forms in the literature~[\Katzzzz, \S1.9], [\Hida,\S4]. Although
this notion is important and useful, the authors of this paper are
not aware of a reference that explains how it stands vis-a-vis a
more direct approach. Recall that the theory of $p$-adic modular
forms began when Serre introduced the notion of a $p$-adic modular
form of a given level as a $q$-expansion which is a $p$-adic
uniform limit of $q$-expansions of classical modular forms of that
same level.

In this setting, Katz's approach of defining $p$-adic modular
forms as certain regular functions on a formal scheme obtained
from schemes of the sort~$\MM\bigl(\ZZ_p, \mu_{p^nN}\bigr)$ merged
nicely with Serre's approach. See~[\Katzz, Prop.~A1.6].

In the Hilbert modular case the development did not follow the
same lines. It seems that the interest was mainly devoted to
understanding the phenomenon of analytic families of Hilbert
modular forms and the connection to completed Hecke algebras; the
ensuing theory is now known as Hida's theory. For this purpose
one's interest is more in varying the level than varying the
weight, and thus the mere existence of suitable ambient space of
$p$-adic modular forms in which such families exist is sufficient.

The authors of this paper are interested in following Serre's
original approach. Congruences between modular forms imply
congruences between their weights that suggest defining a $p$-adic
modular form as a $q$-expansion that is a $p$-adic uniform limit
of classical modular forms. Such a limit has a well defined weight
in the completion~$\widehat{\X}$ of~$\X$ with respect to a system
of subgroups depending on~$p$.

We prove that a $p$-adic modular form of weight~$\chi \in
\widehat{\X}$ defined in this fashion is the same thing as a
$p$-adic modular form in Katz's approach, which is an
eigenfunction of character~$\chi$. We remark that $p$-adic modular
forms \`{a} la Katz are certain regular functions on the formal
scheme $$\lim_{\infty \leftarrow m}\Bigl(\lim_{n\rightarrow
\infty} \MM\bigl(W_m(k), \mu_{p^n N}\bigr)\Bigr)$$(see
definition~\refn{Katzpadicmodfor}), and thus correspond to
ordinary $p$-adic modular forms in the case $g=1$, i.~e., to
$p$-adic modular forms of growth condition~$1$.

One virtue of this isomorphism is that the extension of certain
derivation operators $\Theta_{\P,i}$ (see~\refn{ThetaPi}) to
$p$-adic modular forms is easily proven. This yields an ample
supply of examples of $p$-adic modular forms. The results of this
section follow the presentation in Serre, [\Serre1].

\label NNNotation. section\par\ssection Notation\par In this
section we fix a complete discrete valuation ring~$R$ with
fraction field~$F$ of characteristic~$0$ and residue field~$k$  of
characteristic~$p$. Let~$\m=(\pi)$ be the maximal ideal of~$R$.
Suppose that~$R$ is a $O_K$-algebra where~$K$ is a normal closure
of~$L$.
\endssection

\defi  Let\/~$f\in {\bf
M}\bigl(F,\mu_N,\chi\bigr)$ be a $\I$-polarized modular form
over~$F$ with $N$ not necessarily prime to~$p$. Let
$\bigl(\A,\B,\varepsilon,\j\bigr)$ be a $\I$-polarized cusp.
Consider the $q$-expansion
$f\bigl(\Tate(\A,\B),\varepsilon,\j\bigr)=a_0+\sum_{\nu\in
(\A\B)^+} a_\nu q^\nu$ of~$f$ at the given cusp. Define
$${\rm val}(f):={\rm sup}\bigl\{n\in\ZZ\vert\, a_\nu\in \m^n
\,\forall\nu\bigr\}={\rm inf}\bigl\{{\rm val}_\pi(a_\nu)\bigr\}.$$
\enddefi

\label finiteness. proposition\par \prop The notation is as in the
definition. Then ${\rm val}(f) > -\infty$. Moreover,
$$\pi^{-{\rm val}(f)}\, f \in{\bf
M}\bigl(R,\mu_N,\chi\bigr)$$i.~e., is a $\I$-polarized modular
form over~$R$.
\endprop
\Proof  To prove that ${\rm val}(f) > -\infty$ note that, with the
notation of~\refn{tateobjects}, we have
$$f\bigl(\Tate(\A,\B),\varepsilon,\j\bigr)\in F\tensor_\ZZ
\ZZ\bigl(\bigl(\A,\B,\sigma_\beta\bigr)\bigr).$$In particular, the
valuation of the coefficients of the $q$-expansion of~$f$ is
bounded from below. Since the $q$-expansion of~$\pi^{-{\rm
val}(f)}\, f$ at the given cusp has integral coefficients, we
conclude by~\refn{Koecher} that~$\pi^{-{\rm val}(f)} f$ is defined
over~$R$.

\label InDiPcusp. lemma\par\lemma The number ${\rm val}(f)$ is
independent of the chosen cusp.
\endlemma
\Proof By~\refn{finiteness} we may assume that ${\rm val}(f)\geq
0$ at any cusp. By the $q$-expansion principle explained
in~\refn{Koecher} we have that
$f\bigl(\Tate(\A,\B),\varepsilon,\j\bigr)/\pi^n$ has coefficients
in~$R$ if and only if~$f/\pi^n$ is in~${\bf
M}\bigl(R,\mu_N,\chi\bigr)$.

\label cong. proposition\par \prop  Let\/~$N\geq 4$ be an integer
prime to~$p$. Let
$$f_i\in {\bf M}\bigl(R,\mu_N,\chi_i\bigr)$$for $i=1,2$ be two
$\I$-polarized modular forms of weights~$\chi_1$ and~$\chi_2$ and
level~$\mu_N$. Suppose that their  $q$-expansions at a
$\I$-polarized unramified cusp $\bigl(\A,\B,\varepsilon,\j\bigr)$
in the sense of\/~\refn{qexpansion} satisfy
$$f_1\bigl(\Tate(\A,\B),\varepsilon,\j\bigr)\not\equiv 0\quad\hbox{{\rm mod }}\m$$
and $$ f_1\bigl(\Tate(\A,\B),\varepsilon,\j\bigr)\equiv
f_2\bigl(\Tate(\A,\B),\varepsilon,\j\bigr)\quad\hbox{{\rm mod
}}\m^n.$$Then
$$\chi_1\equiv \chi_2\quad\hbox{{\rm mod }}\X_R(n).$$See~\refn{XB}
for the notation~$\X_R(n)$.\endthm

\Proof Let $\bar{f}_i$  be the image of~$f_i$ in~${\bf
M}\bigl(R/\m^n,\mu_N,\chi_i\bigr)$ for $i=1,2$. Consider the
forgetful morphism
$$\psi\colon\MM\bigl(R/\m^n,\mu_{p^n N}\bigr)\llongrightarrow
\MM\bigl(R/\m^n,\mu_N\bigr).$$It is a Galois cover with group
$$\Aut_{O_L}\bigl(\mu_{p^n}\tensor\DL\bigr)=\bigl(O_L/p^n
O_L\bigr)^*.$$Let $$ r_1=r\bigl(\bar{f}_1\bigr)\qquad\hbox{{\rm
and}}\qquad r_2=r\bigl(\bar{f}_2\bigr)$$be the associated regular
functions on~$\MM\bigl(R/\m^n,\mu_{p^n N}\bigr)$ defined
in~\refn{r(f)}. The hypothesis guarantees that~$r_1=r_2$.
Therefore, if $b\in\bigl(O_L/p^n O_L\bigr)^*$ is an element of the
Galois group, then $$\chi_1(b) r_1=[b]\, r_1=[b]\, r_2=\chi_2(b)
r_2.$$Hence, $\bigl(1-\chi_1\chi_2^{-1}(b)\bigr)r_1=0$. This
implies the claim.

\label CoNg. corollary\par\cor Let $f_i\in{\bf
M}\bigl(F,\mu_N,\chi_i\bigr)$ for $i=1,2$ be two $\I$-polarized
modular forms. Assume  their $q$-expansions at a $\I$-polarized
unramified cusp satisfy
$$f_1\bigl(\Tate(\A,\B),\varepsilon,\j \bigr)\equiv
f_2\bigl(\Tate(\A,\B),\varepsilon,\j\bigr) \quad\hbox{{\rm mod
}}\m^n.$$Then $$\chi_1\equiv \chi_2\quad\hbox{{\rm mod
}}\X_R\Bigl(n-{\rm min}\bigl\{{\rm val}(f_1),{\rm val}(f_2)
\bigr\}\Bigr).$$
\endcor
\Proof By~\refn{finiteness} we may assume that~$f_1$ and~$f_2$ are
defined over~$R$. Let $m_i:={\rm val}\bigl(f_i\bigr)$ for $i=1,2$.
Without loss of generality we may assume that $m_1\leq m_2$. Let
$F_i:=\pi^{-m_1} f_i$ for $i=1,2$. By assumption $F_i\in{\bf
M}\bigl(R,\mu_N,\chi_i \bigr)$ and $F_1 \neq 0$ modulo~$\m$ and
$F_1\equiv F_2$ modulo $\m^{n-m_1}$. We conclude
using~\refn{cong}.

\label a0. corollary\par\cor Let $f\in {\bf
M}\bigl(F,\mu_N,\chi\bigr)$ be a $\I$-polarized modular form.
Consider its $q$-expansion
$f\bigl(\Tate(\A,\B),\varepsilon,\j\bigr)=a_0+\sum_\nu a_\nu
q^\nu$  at a $\I$-polarized unramified cusp. Let
$$n(\chi):={\rm min}\bigl\{n\in\NN\vert\,\chi\notin \X_R(n)
\bigr\}.$$Then
$${\rm val}(a_0)\geq -n(\chi)+{\rm val}\bigl(f-a_0\bigr).$$
\endcor
\Proof By~\refn{finiteness} we may assume that~$f$ is defined
over~$R$. Consider~$a_0$ as a modular form of weight~$1$ (the
trivial character). By corollary~\refn{CoNg}, we have $\chi\equiv
1$ modulo $\X_R\bigl({\rm val}\bigl(f-a_0\bigr)-{\rm
val}\bigl(f\bigr) \bigr)$. Hence, ${\rm val}_\pi(a_0)\geq {\rm
val}\bigl(f\bigr)\geq -n(\chi)+{\rm val}\bigl(f-a_0\bigr)$ as
wanted.

\label Serrepadic. definition\par\defi Suppose that $R$ is
$\pi$-adically complete. A {\it $\I$-polarized $p$-adic Hilbert
modular form \`a la Serre over~$F$ of level\/~$\mu_N$} ($N$ prime
to~$p$) is the equivalence class of a Cauchy sequence
$\bigl\{f_i\in {\bf M}\bigl(F,\mu_N,\chi_i \bigr)\bigr\}_{i\in
\NN}$ of classical modular forms. `Cauchy' means Cauchy with
respect to~${\rm val}$ i.~e., that for any $n\in \NN$ there exists
$m\in \NN$ such that
$${\rm val}(f_i-f_j)\geq n\quad\hbox{{\rm for all }} i,j\geq
m.$$We say that two Cauchy sequences $\{f_i\}_{i\in\NN}$
and~$\{g_j\}_{j\in\NN}$ are equivalent iff ${\rm
val}(f_i-g_i)\rightarrow \infty$ for $i\rightarrow \infty$. The
next lemma shows that this is a well defined concept.
\enddefi

\lemma Let\/~$(\A_1,\B_1,\varepsilon_1,\j_1)$ be a $\I$-polarized
unramified cusp. Let\/~$\{f_i\}_i$ be a Cauchy sequence of
$\I$-polarized Hilbert modular forms with respect to the
valuation~${\rm val}_1$ from~$(\A_1,\B_1,\varepsilon_1,\j_1)$.
If~$(\A_2,\B_2,\varepsilon_2,\j_2)$ is another unramified
$\I$-polarized cusp with associated valuation~${\rm val}_2$,
then~$\{f_i\}_i$ is Cauchy also with respect to the~${\rm val}_2$.
\endlemma

\label qexppadic. definition\par\defi ({\it Weight and
$q$-expansions of $p$-adic modular forms})\enspace
 Let $$f=\bigl\{f_i\in {\bf M}\bigl(F,\mu_N,\chi_i
\bigr)\bigr\}_{i\in\NN}$$be a $\I$-polarized $p$-adic Hilbert
modular form \`a la Serre over~$R$ of level\/~$\mu_N$. Define the
weight\/~$\chi \in \widehat{\X}_R$ of\/~$f$ as
$$\chi:=\lim_{i\rightarrow \infty} \chi_i\in \widehat{\X}_R.$$
Fix a $\I$-polarized unramified cusp $\bigl(\A,\B,\varepsilon,\j
\bigr)$ over~$F$. Define the $q$-expansion of~$f$ at the given
cusp by
$$f\bigl(\Tate(\A,\B),\varepsilon,\j\bigr):=\lim_{i\rightarrow
\infty}f_i\bigl(\Tate(\A,\B),\varepsilon,\j\bigr).$$Finally,
define
$${\rm val}(f):={\rm sup}\bigl\{n\in\ZZ\vert\, a_\nu\in \m^n
\,\forall\nu\bigr\}={\rm inf}\bigl\{{\rm val}_\pi(a_\nu)\bigr\}.$$
\enddefi

\label padicgeneral. proposition\par\prop The notation is as in
the definition.

\item{{\rm 1)}} The weight and the $q$-expansion at the $\I$-polarized cusp
$\bigl(\A,\B,\varepsilon,\j \bigr)$ of a $\I$-polarized $p$-adic
modular form~$f$ \`a la Serre are well defined  i.~e., the limits
exist and do not depend on the choice of Cauchy sequence of
classical modular forms~$f_i$ defining it;
\item{{\rm 2)}} the map
$$
\hbox{{\rm $p$-adic Hilbert modular forms of wt~$\chi$}}
\lllongrightarrow\, \hbox{{\rm $q$-expansions at }}
\bigl(\A,\B,\varepsilon,\j \bigr),$$associating to a $p$-adic
modular form $f$ its $q$-expansion
$f\bigl(\Tate(\A,\B),\varepsilon,\j\bigr)$, is injective;
\item{{\rm 3)}} The assertions in~\refn{finiteness}-\refn{a0} hold if one replaces
$\I$-polarized Hilbert modular forms of level\/~$\mu_N$ with
$\I$-polarized $p$-adic Hilbert modular forms \`a la Serre of
level\/~$\mu_N$ and\/~$\X_R$ with\/~$\widehat{\X}_R$.

\endlemma
\Proof Assertions (1) and~(2) follow from the previous
corollaries. The last assertion follows as in~[\Serre1, \S1].

\rmk  See~\refn{Gamma0paspadic} for  examples of how
Hilbert modular forms of level~$\mu_{Np^n}$ and trivial nebentypus
at~$p$ define $\I$-polarized $p$-adic Hilbert modular forms.
See~\refn{padicEisenstein} for examples of $p$-adic, but not
classical, $\I$-polarized Hilbert modular forms arising from
Eisenstein series. Other examples are given
in~\refn{padicnonclassic} by applying suitable $p$-adic theta
operators to classical Hilbert modular forms.

\endrmk

\defi We say
that a $\I$-polarized $p$-adic modular form~$\{f_i\}_i$ of
level~$\mu_N$ over~$F$ is a cusp form if the constant coefficient
of its $q$-expansion at any cusp is~$0$.
\enddefi

\endsection

\section Katz's $p$-adic Hilbert modular forms\par
\noindent The notation is as in~\refn{NNNotation}.

\label MM(n,n). definition\par\defi Let $m\geq 1$, $n\geq 0$
and~$N\geq 4$ be integers. Let\/~$p$ be a prime not dividing~$N$.
Consider the affine schemes $$\MM(m,n) =\MM\bigl(R/\m^m,
\mu_{p^nN}\bigr).$$For $n=0$ we use the convention
$\MM(m,0)=\MM(R/\m^m,\mu_N)^{\rm ord}$. For $m' \leq m$ we have a
closed immersion $$\MM(m', n) \hooklongrightarrow \MM(m,n).$$For
$n' \geq n$ we have a Galois covering $$\MM(m, n')
\llongrightarrow \MM(m, n)$$with Galois group $$\Gamma(n',n)
=\bigl(O_L/p^{n'}O_L\bigr)^*/\bigl(O_L/p^n O_L\bigr)^*,$$where the
group $\bigl(O_L/p^0\,O_L)^*$ is understood as the trivial group.
Define
$$\Gamma_{n'}:=\Gamma(n',0)=\Aut_{O_L}(\mu_{p^{n'}}\tensor_\ZZ
\DL^{-1})=\bigl(O_L/p^{n'}\,O_L\bigr)^*.$$

\noindent Define $\MM(\infty, \infty)$ as the formal scheme
$$\MM(\infty, \infty) = \lim_{m\rightarrow\infty} \MM(m, \infty) =
\lim_{m\rightarrow\infty}\Bigl(\lim_{\infty\leftarrow
n}\MM(m,n)\Bigr).$$The group $$\Gamma  =\lim_{\infty\leftarrow n}
\Gamma_n=\lim_{\infty\leftarrow n}\bigl(O_L/p^n\,O_L\bigr)^*=
\bigl(\ZZ_p\tensor_\ZZ O_L \bigr)^*$$acts as Galois automorphisms
on~$\MM(m, \infty)$ for every~$m$ with quotient~$\MM(m, 0)$.
\enddefi

\label MM(n,n)Kummer. definition\par\ssection The Kummer
description\par Let $n\geq m$. Let
$$\MM(m,n)^\Kum \llongrightarrow \MM(m,0)$$be the Kummer part of the
cover $\MM(m,n)\rightarrow \MM(m,0)$; see~\refn{phi}. Its Galois
group is
$$G_{m,n}\cong \Gamma(n,0)/ \bigl(\cap_\chi
\Ker(\chi)\bigr);$$the intersection is taken over all the
universal characters~$\chi$ restricted to $$\chi\colon
\bigl(O_L/p^n O_L\bigr)^*\llongrightarrow
\bigl(R/\m^m\bigr)^*.$$Since the natural map
$$\lim_{\infty\leftarrow n}\Gamma_n \llongrightarrow
\lim_{\infty\leftarrow n} G_{m,n}$$is an isomorphism, we deduce
that
$$\MM(\infty, \infty) = \lim_{m\rightarrow\infty} \MM(m, \infty)^\Kum =
\lim_{m\rightarrow\infty}\Bigl(\lim_{\infty\leftarrow
n}\MM(m,n)\Bigr)^\Kum.$$

\endssection

\label padicsequence. lemma\par\lemma To give a regular
function~$f$ on~$\MM(\infty, \infty)$ is equivalent to give a
sequence~$\{ f_n \}_n$ of regular functions~$f_n$ on~$\MM(n, n)$
(equivalently on~$\MM(n,n)^\Kum$) such that~$f_{n+1}\vert_{\MM(n,
n+1)}$ (resp.~$f_{n+1}\vert_{\MM(n, n+1)^\Kum}$) is the pull back
of~$f_n$ under the forgetful morphism $\MM(n, n+1) \rightarrow
\MM(n, n)$ (resp.~$\MM(n, n+1)^\Kum \rightarrow \MM(n, n)^\Kum$).
We shall call such a sequence a {\it compatible sequence}.
\endlemma
\Proof A regular function~$f$ on~$\MM(\infty, \infty)$ may be
identified with a sequence~$f_n$ with~$f_n$ a regular function
on~$\MM(n, \infty)$ satisfying $f_{n+1}\equiv f_n$ mod~$\m^n$. The
function~$f_n$ is thus a function on~$\MM(n, s(n))$ for a suitable
integer~$s(n)$. We may assume w.l.o.g.~that the sequence~$s(n)$ is
monotone increasing and~$s(n) \geq n$ for all~$n$. We then have
the compatibility that~$f_{n+1}\vert_{\MM(n, s(n+1))}$ is the pull
back of~$f_n$ under the morphism $\MM(n, s(n+1)) \rightarrow
\MM(n, s(n))$. Since the morphism $\MM(m, n) \rightarrow \MM(m+1,
n)$ is a closed immersion of affine schemes, we may find for every
$n$ a regular function~$f_n'$ on~$\MM(s(n), s(n))$ equal to~$f_n$
under restriction. The function~$f$ is identified with the
sequence~$\{f_n'\}_n$. The same proof applies to the Kummer parts.

\ssection Weights and characters of $\Gamma$\par Let $\chi \in
\X_R$ be a character. We may apply~$\chi$ to~$\Gamma$ by forming
the composition $$\Gamma = (\ZZ_p\tensor_\ZZ O_L)^*
\hooklongrightarrow (R \tensor_\ZZ O_L)^* = \G(R)
\llongmaprighto{\chi} R^*.$$Suppose that $\chi \in \X_R(n)$
(see~\refn{XB}), then $\chi(\Gamma) \equiv 1$ modulo~$\m^n$. It
follows that every element of $\chi \in \widehat{\X}_R$ defines a
well defined homomorphism $$\Gamma \llongmaprighto{\chi} R^*.$$
\endssection

\label Katzpadicmodfor. definition\par\defi ({\it $p$-adic Hilbert
modular forms \`a la Katz; c.f.~{\rm [\Katzzzz, \S1.9]}})\enspace
Let $\chi$ be a character in $\widehat{\X}_R$. A {\it
$\I$-polarized Katz modular form} of weight\/~$\chi$ and
level~$\mu_N$ defined over~$R$ is a regular function~$f$
on~$\MM(\infty, \infty)$ such that for every~$\alpha \in \Gamma$
we have
$$\alpha^*\bigl(f\bigr)=\chi(\alpha)\, f,$$where
$\alpha^*\bigl(f\bigr):=f\circ \alpha$. Denote the $R$-module of
such functions by
$${\bf M}(R,\mu_N,\chi)^{p-{\rm adic}}.$$The $F$-module of $\I$-polarized Katz
modular forms of weight\/~$\chi$ and level~$\mu_N$ defined
over~$F$ is defined by $${\bf M}(F,\mu_N,\chi)^{p-{\rm
adic}}:=F\tensor_R {\bf M}(R,\mu_N,\chi)^{p-{\rm adic}}.$$
\enddefi

\rmk The notion of $\I$-polarized Katz $p$-adic Hilbert modular forms commutes
with base change. More precisely, let $R\rightarrow R'$ be an
extension of $O_K$-algebras, which are discrete valuation rings of
unequal characteristic~$p$ and~$0$. Let~$F$ and~$F'$ be the
associated fraction fields. If $\chi\in \widehat{\X}_R\subset
\widehat{\X}_{R'}$, then
$${\bf M}(R',\mu_N,\chi)^{p-{\rm adic}}\cong {\bf
M}(R,\mu_N,\chi)^{p-{\rm adic}}\tensor_R R'$$and $${\bf
M}(F',\mu_N,\chi)^{p-{\rm adic}}\cong {\bf
M}(F,\mu_N,\chi)^{p-{\rm adic}}\tensor_F F'.$$This follows from
the fact that the formation of the moduli spaces~$\MM(m,n)$
commutes with base change. \endrmk

\label qexpKatz. definition\par\defi ({\it $q$-expansions of Katz
modular forms})\enspace Define a $\I$-polarized unramified cusp
of~$\MM(\infty,\infty)$ to be equivalently \item{{\rm a)}} an
unramified cusp $(\A,\B,\varepsilon_{p^\infty N
},\j_\varepsilon\bigr)$ over~$R$ with $\mu_{p^\infty N}$-level
structure; \item{{\rm b)}} a compatible system of unramified cusps
$(\A,\B,\varepsilon_{p^n N },\j_\varepsilon\bigr)$ over~$R/\m^n$.

\noindent Note that we write~$\j_\varepsilon$ in place
of\/~$\j_{\varepsilon_{p^\infty}}$
(resp.~$\j_{\varepsilon_{p^n}}$). We refer to~\refn{jepsilon} for
the latter notation.  Given such a cusp, one has a $q$-expansion
map
$$\dirsum_\chi {\bf
M}(F,\mu_N,\chi)^{p-{\rm adic}} \llongrightarrow
F[\![q^\nu]\!]_{\nu\in (\A\B)^+\cup\{0\}},$$which is injective;
see~{\rm [\Katzzzz, Thm.~1.10.15]}. With the notation
of~\refn{tateobjects}, it is defined by evaluating a
$\I$-polarized Katz modular form~$f$ defined over~$R$ via
$$\Spec\left( \ZZ((\A,\B,\sigma_\beta ))\tensor_\ZZ R\right)
\llongrightarrow \MM(\infty,\infty).$$We say that a $\I$-polarized
Katz modular form is a cusp form if the constant coefficient of
its $q$-expansion at any cusp is~$0$.
\enddefi

\label totalrecall. section\par\ssection Recall\par If the
inequality $m \leq n$ holds and $\chi \in \X_R$ there is a
canonically defined modular form
$$a(\chi)\in {\bf M}\bigl(R/\m^m,\mu_{p^n
N},\chi\bigr);$$see~\refn{universal}. Let
$(\A,\B,\varepsilon_{p^\infty N },\j_\varepsilon\bigr)$ be a
$\I$-polarized  unramified cusp of~$\MM(m, n)$; see~\refn{cusp}.
Then

\item{{\rm i)}} the modular form $a(\chi)$
descends to a modular form
on~$\MM\bigl(R/\m^m,\mu_{p^nN}\bigr)^\Kum$.
See~\refn{actiona(chi)};
\item{{\rm ii)}}
the modular form~$a(\chi)$ transforms under~$\Gamma_n$
(equivalently under $G_{m,n}$) according to the
character~$\chi^{-1}$. See~\refn{actsonachi};
\item{{\rm iii)}} $a(\chi)$  has $q$-expansion $1$ at the
cusp~$(\A,\B,\varepsilon_{p^\infty N },\j_\varepsilon\bigr)$;
\item{{\rm iv)}} for $m'\leq m\leq n$ and $n'\geq n\geq m$ the
modular forms~$a(\chi)$ defined on~$\MM(m, n)$, on~$\MM(m', n)$
and on~$\MM(m, n')$ (resp. on~$\MM(m, n)^\Kum$, on~$\MM(m',
n)^\Kum$ and on~$\MM(m, n')^\Kum$) agree.

\endssection

\label compaRe. section\par\ssection The comparison between Serre
and Katz modular forms\par Let $g=\{g_i\}_i$ be a $\I$-polarized
Serre $p$-adic modular form of weight~$\chi$, level~$\mu_N$ and
defined over~$R$ as in~\refn{Serrepadic}. We may assume
w.l.o.g.~that
$$g_{n+1} \equiv g_n\qquad \hbox{{\rm mod }}\m^n.$$Let $\chi_n$ be the
weight of~$g_n$ so that $\chi=\lim_n \chi_n$ by~\refn{qexppadic}.
For every $n\in\NN$ define $$f_n := g_n/a(\chi_n)$$as a regular
function on~$\MM(n, n)$. By~\refn{padicsequence}  the
sequence~$\{f_n\}_n$ defines a $\I$-polarized Katz $p$-adic
modular form of weight~$\chi$ and level~$\mu_N$ over~$R$.
Moreover, the $q$-expansion at a cusp of~$\MM(n,n)$ of~$f_n$ is
the $q$-expansion of~$g_n$. Hence, the $q$-expansion of~$f$ at a
cusp of~$\MM(\infty,\infty)$ coincides with that of~$g$. Thus, we
obtain a map $$r\colon \hbox{{\rm  $\I$-pol. Serre $p$-adic HMF
over }}F \llongrightarrow  \hbox{{\rm $\I$-pol. Katz $p$-adic HMF
over }}F,$$which preserves  weights. We also conclude that for any
$\I$-polarized unramified cusp $(\A,\B,\varepsilon_{p^\infty
N},\j_\varepsilon)$ of~$\MM(\infty,\infty)$ the following diagram
is commutative
$$\matrix{
\hbox{{\rm  Serre $p$-adic HMF over }}F  & \llongmaprighto{r} &
\hbox{{\rm Katz $p$-adic HMF over }}F\cr
 \hbox{{\rm $q$-exp}}\searrow & &\swarrow\hbox{{\rm $q$-exp}}\cr
 & F[\![q^\nu]\!]_{\nu\in (\A\B)^+\cup\{0\}}.\cr}$$We deduce
from~\refn{qexpKatz}, and the fact that~$r$ preserves weights,
that the $q$-expansion map is injective also on the graded ring of
$\I$-polarized Serre $p$-adic Hilbert modular forms.
\endssection

\label lifting. lemma\par\lemma Let $N\geq 4$ be an integer.
\spacing
\item{{\rm 1)}} For a suitable
integer $n_0> 0$ the modular form $h^{n_0}$ (the $n_0$-th power of
the Hasse invariant defined in~\refn{hPi}) admits a lift to a
modular form~$\hslash$ over~$\ZZ_p$ of
weight\/~$\Norm^{(p-1)n_0}$. We may choose~$\hslash$ so that the
leading coefficient of its $q$-expansion at a given cusp is~$1$.
\spacing
\item{{\rm 2)}} For any  $\I$-polarized modular form $f\in {\bf
M}\bigl(R/\m^n,\mu_N,\Nm^s\bigr)$ there exists a $\I$-polarized
modular form $g_n\in{\bf M}\bigl(R,\mu_N,\Nm^{s'}\bigr)$ such that
\itemitem{{\rm
2.a)}} $g_n$ mod~$\m^n$ and\/~$f$ have the same $q$-expansion at
one (any) cusp; \itemitem{{\rm 2.b)}} $\Nm^s\equiv \Nm^{s'}$
mod~$\X_R(n)$.

\spacing
\item{{\rm 3)}} For any character $\chi\in\X_R$, any~$n\in\NN$
and any $\I$-polarized cusp form $f\in {\bf
M}\bigl(R/\m^n,\mu_N,\chi\bigr)$ there exists a $\I$-polarized
modular form $g_n\in{\bf M}\bigl(R,\mu_N,\chi'\bigr)$ such that

\itemitem{{\rm
3.a)}} $g_n$ mod~$\m^n$ and\/~$f$ have the same $q$-expansion at
one (any) cusp; \itemitem{{\rm 3.b)}} $\chi\equiv \chi'$
mod~$\X_R(n)$.
\endlemma
\Proof We fix some notation. Let $$\delta\colon
\MMbar(R,\mu_N)_{\sigma_\beta}\llongrightarrow \MMbar(R,\mu_N)$$be
the morphism from a toroidal compactification to the minimal
compactification of~$\MM(R,\mu_N)$. See~[\Rapoport, \S5]
and~[\Chai, \S4]. Let $$
\pi\colon\AA\longrightarrow\MM(R,\mu_N)$$be the universal
$\I$-polarized abelian scheme with real multiplication by~$O_L$.
By~[\Rapoport, \S5.4], the abelian scheme~$\AA$ extends to a
semiabelian scheme over~$\MMbar(R,\mu_N)_{\sigma_\beta}$. Hence,
the sheaf $\Omega^1_{\AA/\MM(R,\mu_N)}$ extends to an invertible
$O_L\tensor_\ZZ O_{\MMbar(R,\mu_N)_{\sigma_\beta}}$-module
$\Omega$ on~$\MMbar(R,\mu_N)_{\sigma_\beta}$. The line bundle
$$\wedge^g\pi_*\Bigl(\Omega^1_{\AA/\MM(R;\mu_N)}\Bigr)={\cal L}_\Norm$$is
the Hodge bundle; see~\refn{Lchi} for the notation. It extends to
a line bundle~$\wedge^g\pi_*\bigl(\Omega\bigr)$
on~$\MMbar(R,\mu_N)_{\sigma_\beta}$. By~[\Chai, Thm.~4.3 (IX)] the
latter descends to an ample line bundle on~$\MMbar(R,\mu_N)$,
which we denote in the same way. This implies that there exists an
integer~$n_0\gg 0$ such that~${\cal L}_{\Norm^{n_0(p-1)}}$ is very
ample. Note that~$h$ is a section of~${\cal L}_\Norm$
over~$\MM(k,\mu_N)^\R$ which has codimension at least two~$2$
in~$\MMbar(k,\mu_N)$. By~[\Chai, Thm.~4.3 (V)], the
scheme~$\MMbar(k,\mu_N)$ is  normal. In particular, $h$ extends to
a section over~$\MMbar(k,\mu_N)$. Hence, $h^{n_0}$ lifts to a
modular form~$\hslash$ of~${\cal L}_{\Norm^{n_0(p-1)}}$
over~$\MMbar(R,\mu_N)$.

\noindent Fix a $\I$-polarized unramified cusp
$(\A,\B,\varepsilon,\j_\varepsilon)$ of~$\MM(R,\mu_N)$. The
$q$-expansion principle implies that the $q$-expansion
of~$\hslash$ is congruent to that of~$h^{n_0}$ modulo~$\m$.
By~\refn{qexphPi}, letting
$$a_0 + \sum_{\nu\in (\A\B)^+} a_\nu q^\nu$$be the $q$-expansion
of~$\hslash$ at the given  cusp, we have that~$a_0\equiv 1$
mod~$\m$ and~$a_\nu\in\m$ for all~$\nu\neq 0$. In particular,
$a_0$ is a unit. Replacing~$\hslash$ by~$a_0^{-1}\hslash$ we may
assume~$\hslash$ to have a leading Fourier coefficient equal to
one at the given cusp. This proves~(1). The proof of~(2) is
analogous: one replaces~$f$ by~$f':=f\hslash^n$ for~$n\gg 0$,
argues that~$f'$ extends to the minimal compactification and lifts
to~$R$ provided~$n$ is suitably chosen.

\spacing
\noindent Define $$\delta^o\colon
\MMbar(R,\mu_N)_{\sigma_\beta}^o\llongrightarrow
\MMbar(R,\mu_N)^o$$as the restriction of~$\delta$ to the
complement of the divisor~$\hslash$ in~$\MMbar(R,\mu_N)$. Let
$\chi\in \X_R$. Since~$\Omega^1_{\AA/\MM(R,\mu_N)}$ extends on~$
\MMbar(R,\mu_N)_{\sigma_\beta}^o$ as an invertible $O_L\tensor_\ZZ
O_{\MMbar(R,\mu_N)_{\sigma_\beta}^o}$-module, proceeding as
in~\refn{Lchi} we see that the invertible sheaf~${\cal L}_\chi$
on~$\MM(R,\mu_N)$ extends to an invertible sheaf~$\bar{{\cal
L}}_\chi$ on~$ \MMbar(R,\mu_N)_{\sigma_\beta}^o$. Since~$\delta^o$
is proper, $\delta^o_*\bigl(\bar{{\cal L}}_\chi \bigr)$ is a
coherent sheaf. By construction its restriction to~$\MM(R,\mu_N)$
coincides with~${\cal L}_\chi$. Since~$\hslash$ is very ample,
$\MMbar(R,\mu_N)^o$ is affine. In particular, for any $n\in \NN$
the map
$$\Gamma\left(\MMbar(R,\mu_N)^o,\delta^o_*\bigl(\bar{{\cal L}}_\chi \bigr)\right)
\llongrightarrow
\Gamma\left(\MMbar(R/\m^n,\mu_N)^o,\delta^o_*\bigl(\bar{{\cal
L}}_\chi \bigr)\right)$$is surjective. Note that
$${\bf
M}\bigl(R/\m^n,\mu_N,\chi\bigr)=\Gamma\left(\MM(R/\m^n,\mu_N),{\cal
L}_\chi\fibprod_{\Spec(R)} \Spec(R/\m^n)\right).$$In particular,
if $f\in{\bf M }\bigl(R/\m^n,\mu_N,\chi\bigr)$ is a cusp form,  we
can extend~$f$ by~$0$ to a global section~$f'$
of~$\delta^o_*\bigl(\bar{{\cal L}}_\chi \bigr)$
over~$\MMbar\bigl(R/\m^n,\mu_N\bigr)^o$. There exists a global
section~$g'$ of~$\delta^o_*\bigl(\bar{{\cal L}}_\chi \bigr)$
over~$\MMbar\bigl(R,\mu_N\bigr)^o$ lifting~$f'$. Hence, there
exists~$r\gg 0$ such that~$g:=g'\hslash^r $ extends to a section
of~${\cal L}_\chi$ over~$\MM(R,\mu_N)$ and $\Norm^{(p-1)n_0
r}\equiv 1$ mod~$\X_R(n)$. Since
$${\bf M}\bigl(R,\mu_N,\chi\bigr)=\Gamma\bigl(\MM(R,\mu_N),{\cal
L}_\chi\bigr),$$this proves~(3).

\label U1(N)isNorm. lemma\par\lemma Let\/~$N$ be an integer.
Let~$U_1(N)$ denote the elements of~$O_L^*$ congruent to~$1$
modulo~$N$. Let $\chi \in \X_R$ such that we
have~$\chi\bigl(U_1(N)\bigr) = 1$. Then~$\chi$ is  a power
of\/~$\Norm$. \endlemma \Proof The character $\chi$ belongs
to~$\X$ by~\refn{geniso}. Let us write the complex embeddings
of~$L$ as~$\sigma_1, \dots, \sigma_g$. Then we may write
$\chi\tensor \CC = \sigma_1^{a_1} \dots \sigma_g^{a_g}$.
Replacing~$\chi$ by~$\chi^2$ if necessary, we may assume that the
$a_1,\ldots,a_g$ are even. By multiplying~$\chi$ by a suitable
power of~$\Norm^2$ we may assume~$a_i \geq 0$ for all $1\leq i\leq
g$ and w.l.o.g.~$a_1 = 0$. By Dirichlet's units theorem there
exists a unit~$u\in O_L^*$ such that~$\sigma_1(u)
>1$ and~$0<\sigma_i(u) <1$ for~$i = 2, \ldots, g$. Since~$U_1(N)$ is
of finite index in~$O_L^*$, there exists a power~$u^n$ of~$u$ such
that~$u^n \in U_1(N)$. But then $1= \chi\bigl(u^n\bigr) = \prod_{i
= 2}^g \sigma_i^{a_i} \bigl(u^n)\bigr) \leq 1$. We have equality
if and only if~$a_i = 0$ for $2\leq i\leq g$ i.~e., $\chi$ is a
multiple of~$\Norm$.

\label CComPPare. theorem\par\thm The notions of a $\I$-polarized
$p$-adic Hilbert modular form over~$F$ in the sense of Serre and
in the sense of Katz are the same i.~e., the map~$r$ is an
isomorphism, in the following cases:
\spacing
\item{{\rm i.}} cusp forms;
\spacing
\item{{\rm ii.}} forms of weight~$\chi \in \X$;
\spacing
\item{{\rm iii.}} forms of weight~$\Nm^z$ with~$z\in\ZZ_p$.
\spacing

\noindent Moreover, a modular form of non-parallel
weight~$\chi\in\X$ i.~e., whose weight is not of the
form~$\Norm^{s}$ for a suitable integer~$s$, is a cusp form.
\endthm
\Proof  The injectivity of~$r$ follows from the injectivity of the
$q$-expansion map on Serre $p$-adic modular forms proven
in~\refn{padicgeneral}. The fact that~$r$ preserves the notion of
cusp form in the sense of Katz and of Serre follows
from~\refn{compaRe}. We are left with the proof of the
surjectivity of~$r$. By~\refn{finiteness} it suffices to prove it
for modular forms defined over~$R$.

\noindent Let $\{f_n\}_n$ be a sequence as in~\refn{padicsequence}
giving a Katz modular form of weight $\chi\in~\widehat{\X}$. For
every~$n$, the regular function~$f_n$ on~$\MM(n, n)$ is of
weight~$\chi_n:= \chi$ mod~$\m^n$, a character of~$\Gamma_n$. In
particular, $f_n$ is invariant under $H:=\cap_\psi \Ker(\psi)$,
where the intersection is taken over all characters $\psi\colon
\Gamma_n\rightarrow (R/\m^n)^*$. Let~$\MM(n,n)\rightarrow
\MM(n,n)^\Kum$ be the Galois cover with group~$H$; see~\refn{phi}.
Note that
$$f_n\in O_{\MM(n,n)^\Kum}^{\chi_n}\cong {\cal L}_{\chi_n}$$the
isomorphism being as $O_L\tensor_\ZZ O_{\MM(n,0)}$-modules.  It
follows from~\refn{actiona(chi)} that the modular
form~$g_n':=a(\chi_n)\, f_n$ on~$\MM(n,n)^\Kum$ descends to a
modular form on~$\MM(n, 0)$ of weight~$\chi_n$ mod.~$\X_R(n)$ and
with the same $q$-expansion at any cusp as that of~$f_n$. By
multiplying it by a high enough power of the modular
form~$\hslash$, constructed in~\refn{lifting} as a lifting
to~$\ZZ_p$ a power of the Hasse invariant, we may assume
that~$g_n'$ extends to a modular form defined
over~$\MM(R/\m^m,\mu_N)$.

\noindent Assume first that~$f$ is a cusp form. Then so is
each~$f_n$ and~$g_n'$. It follows from~\refn{lifting} that we may
find a modular form~$g_n$ over~$R$ such that $$ g_n \equiv g_n'
\qquad \hbox{{\rm mod }}\m^n.$$Hence, the sequence of modular
forms~$g_n$ of weight~$\psi_n \equiv \chi_n$ mod~$\X_R(n)$ of
modular forms over~$R$ converges to a Serre $p$-adic modular form
with the same $q$-expansion as that of~$f$. It follows
from~\refn{compaRe} that~$f$ is the image of the Serre $p$-adic
modular form~$\{g_n\}_n$. This proves~(i).

\noindent Suppose  that~$f$ is not a cusp form, but has
weight~$\chi\in\X$. The $q$-expansion of~$g_n'$ lies in the ring
$$ (R/\m^n)[\![q^\nu]\!]_{\nu \in O_L^+}^{U_1(N)},$$where the
action of $U_1(N)$ (the units of~$O_L$ congruent to~$1$ mod~$N$)
is given by the ``factor of automorphy" $$ q^\nu \mapsto
\chi(\epsilon)q^{\epsilon^2 \nu}.$$Looking at the
coefficient~$q^0$, this implies that $\chi(\epsilon) \equiv 1$
mod~$\m^n$ for all~$n$ and~$\epsilon \in U_1(N)$ and, therefore,
that~$\chi$ induces the trivial character $\chi\colon U_1(N)
\rightarrow R^*$. We conclude from~\refn{U1(N)isNorm} that~$\chi$
is of the form~$\Norm^r$ for a suitable $r\in \ZZ$.

In particular, to prove~(ii) and~(iii) we may assume that $\chi_n
= \Norm^{s(n)}$ mod~$\X_R(n)$ for a suitable integer~$s(n)$. There
exists a positive integer~$t$, depending on~$g_n'$, such that the
modular form~$g_n' \hslash^t$ extends to a modular form
on~$\MM(R/\m^n,\mu_N)$. By~\refn{lifting}, we may take~$t$ such
that~$g_n' \hslash^t$  can be lifted to a modular form~$g_n$ of
parallel weight, defined over~$R$ and whose $q$-expansion
mod~$\m^n$ is the $q$-expansion of~$g_n'$. The sequence
$\{g_n\}_n$ defines a Serre $p$-adic modular form whose associated
Katz $p$-adic modular form is~$f$ by construction.

\label Gamma0paspadic. section\par\ssection Modular forms of level
$\Gamma_0(p^n)$ as $p$-adic modular forms\par Let ${\bf
M}(R,\mu_{p^nN},\chi)^\1$ be the $R$-module of $\I$-polarized
modular forms of level~$\mu_{p^nN}$, character~$\chi$ and trivial
nebentypus at~$p$ i.~e., invariant
under~$\Gamma_n=\Aut\bigl(\DL^{-1}\tensor_\ZZ \mu_{p^n} \bigr)$.
Define an $R$-linear map
$$\tau_F\colon {\bf
M}(F,\mu_{p^nN},\chi)^\1\llongrightarrow {\bf
M}(F,\mu_N,\chi)^{p-{\rm adic}}$$as follows. For any $m\in\NN$,
consider the Galois morphism with group~$\Gamma_n$
$$\psi\colon \MM(R/\m^m,\mu_{p^nN}) \llongrightarrow
\MM(R/\m^m,\mu_N).$$Let~$f\in {\bf M}(R,\mu_{p^nN},\chi)^\1$. For
any~$m\in \NN$ the reduction~$g_m$ of~$f$ defines a modular form
on~$\MM(R/\m^m,\mu_{p^nN})$ of weight~$\chi$ invariant under the
action of~$\Gamma_n$. Hence, $g_m$ descends to a modular form of
weight~$\chi$ on~$\MM(m,0)=\MM(R/\m^m,\mu_N)^{\rm ord}$ (we freely
use the interpretation of modular forms of given weight as
sections of line bundles; see~\refn{Lchi}). In particular, for any
$s\geq m$ we obtain a regular function
$${g_m\over a(\chi)}\in
\Gamma\left(\MM(m,s),O_{\MM(m,s)}\right);$$see~\refn{totalrecall}
for the notation. By~\refn{padicsequence} this defines a
$\I$-polarized $p$-adic Hilbert modular form~$g$ \`a la Katz of
weight~$\chi$. Define $$\tau_F(f):=g.$$If~$f$ is a $\I$-polarized
modular form of weight~$\chi$ and level~$\mu_{p^nN}$ defined
over~$F$, as argued in~\refn{finiteness}, there exists  an
integer~$t\gg 0$ such that $\pi^t\, f$ is defined over~$R$. Define
$$\tau_F\bigl(f\bigr):=\pi^{-t}\tau_F\bigl(\pi^t
f\bigr).$$

\spacing
\noindent We now explain the connection to $\I$-polarized modular
forms of level~$\mu_N\times\Gamma_0(p^n)$. The notation is as
in~\refn{modforGamma0p}. The morphism
$$\MM\bigl(\bar{F},\mu_{p^nN}\bigr)\llongrightarrow
\MM\bigl(\bar{F},\mu_N,\Gamma_0(p^n)\bigr)$$is finite, \'etale and
Galois with group~$\Aut\bigl(\DL^{-1}\tensor_\ZZ \mu_{p^n}
\bigr)$. In particular, we obtain the isomorphism $${\bf
M}\bigl(\bar{F},\mu_N,\Gamma_0(p^n),\chi\bigr) \isomarrow {\bf
M}\bigl(\bar{F},\mu_{p^nN},\chi\bigr)^\1.$$Hence, we get
$\bar{F}$-linear maps
$${\bf M}\bigl(\bar{F},\mu_N,\Gamma_0(p^n),\chi\bigr) \isomarrow {\bf
M}\bigl(\bar{F},\mu_{p^nN},\chi\bigr)^\1\llongmaprighto{\tau_{\bar{F}}}
{\bf M}\bigl(\bar{F},\mu_N,\chi\bigr)^{p-{\rm adic}}.$$By
construction we obtain the following important property
$$\tau_{\bar{F}}(f)\Bigl(\Tate(\A,\B),\varepsilon_{p^\infty
N},\j_\varepsilon\Bigr)= f\Bigl(\Tate(\A,\B),\varepsilon_{p^n
N},\j_\varepsilon \Bigr);$$see~\refn{qexpansion}
and~\refn{qexpKatz} for the notation.

\endsection

\section Integrality and congruences for values of zeta functions\par
\noindent In this section we apply the results of Section~10 to
derive congruences between values of Dedekind zeta functions and
bounds on the denominators of these values. In fact, the method is
applicable for a wide range of $L$-functions; [\DeligneRibet].

\defi Let $$\zeta_L(s):=
\sum_{I} \Norm(I)^{-s}\qquad \bigl({\rm Re}(s)>1\bigr)$$be the
Dedekind $\zeta$-function associated to~$L$.
\enddefi

\label Siegel. theorem\par \thm The function $\zeta_L(s)$ can be
continued to a meromorphic function on~$\CC$, holomorphic for $s
\neq 1$. Moreover, $\zeta_L(1-k)$ is in~$\QQ$ for every integer $k
\geq 1$.\endthm \Proof See [\Siegel].

\label Eisenstein. theorem\par\thm Let $k \geq 2$ be an even
integer. There exists a $\I$-polarized modular form $$ \E_k \in
{\bf M}\bigl(\CC,\mu_N,\Norm^{k-1}\bigr)$$such that the
$q$-expansion of\/~$\E_{k}$ at a $\I$-polarized unramified cusp
$(\A,\B,\j_{\rm can})$, as in~\refn{cusp} and\/~\refn{complex}, is
$$\Norm^{k-1}(\A) \biggl( 2^{-g} \zeta_L(1-k) +
\sum_{\nu\in (\A\B)^+} \Bigl(\sum_{\nu\in \C\subset \A\B}
\Norm(\nu\C^{-1})^{k-1}\Bigr) q^{\nu}\biggr).$$
\endthm
\Proof See [\Geer, Chap.~I, \S6] or [\DeligneRibet, Thm.~6.1].

\label maps. definition\par\defi  Let $\J$ be an ideal of\/~$O_L$
dividing~$p$ and prime to~$N$. See~\refn{Gamma0p} for the notion
of\/ $\Gamma_0(p)$-level structures. Let
$$\pi_\J\colon \MM\bigl(\CC,\mu_{N}, \Gamma_0(p)\bigr)
\lllongrightarrow\MM\bigl(\CC,\mu_{N}\bigr)$$be the map
associating to a Hilbert-Blumenthal abelian scheme $A$ over a
scheme~$S$  with $O_L$-action, $\mu_{N}$-level structure and
$\Gamma_0(p)$-level structure $ H \hooklongrightarrow A$ the
Hilbert-Blumenthal abelian scheme $A/H[\J]$ over~$S$ with induced
$O_L$-action and $\mu_N$-level structure.\enddefi

\rmk Here we are forced to work with several polarization modules.
Indeed, if the abelian scheme~$A$ in~\refn{maps} is
$\I$-polarized, the quotient~$A/H[\J]$ is $\J\I$-polarized.
\endrmk

\label modified. theorem\par \thm The notation is as
in~\refn{maps}. Consider the $\J\I$-polarized Eisenstein series
$$\E_k \in {\bf M}\bigl(\CC,\mu_N,\Norm^{k-1}\bigr).$$Let $$H=
\left({O_L \over \J }\right)(1) \hooklongrightarrow
\Tate\bigl(O_L,\I^{-1}\bigr)_{\sigma_\beta}$$be the subgroup
defined in~\refn{indiff}. The $q$-expansion of
$\pi_\J^*\bigl(\E_k\bigr)$ at the $\I$-polarized cusp
$\bigl(O_L,\I^{-1},\varepsilon,H,\j_{\rm can}\bigr)$ is
$$\Norm(\J)^{k-1}\biggl(2^{-g} \zeta_L(1-k) + \sum_{\nu\in (\I^{-1}\J)^+}
\Bigl(\sum_{\nu\in\C \subset \J\I^{-1}}
\Norm(\nu\C^{-1})^{k-1}\Bigr)
q^{\nu}\biggr).$$See~\refn{Gamma0pqexp} for the notation.
\endthm

\Proof The $q$-expansion of $\pi_\J^*\bigl(\E_k\bigr)$ is defined
by
$$\pi_\J^*\bigl(\E_{k}\bigr)
\biggl(\Tate(O_L,\I^{-1})_{\sigma_\beta}\tensor_\ZZ
R,\varepsilon,H,{dt \over t}\biggr);$$the notation is as
in~\refn{cusp} and in~\refn{Gamma0pqexp}.  This is equal to
$$\E_{k}\Bigl(\pi_\J\bigl(\Tate(O_L,\I^{-1})_{\sigma_\beta}\tensor_\ZZ
R,\varepsilon,H,{dt \over t}\bigr)\Bigr).$$As explained
in~\refn{tateobjects}, the abelian scheme
$\Tate(O_L,\I^{-1})_{\sigma_\beta}$ is defined by restricting the
semiabelian scheme  $\DL^{-1}\tensor_\ZZ
\GG_{m,S_{\sigma_{\beta}}}/ \underline{q}\bigl(\I^{-1}\bigr)$ to
the open $S_{\sigma_\beta}\backslash S_{\sigma_\beta,0}$. Using
the dictionary of\/~\refn{uniformization} we get that
$\pi_\J\Bigl(\Tate(O_L,\I^{-1})_{\sigma_\beta}\Bigr)$ coincides
with the restriction to $S_{\sigma_\beta}\backslash
S_{\sigma_\beta,0}$ of $\J^{-1} \DL^{-1}\tensor_\ZZ\GG_{m,
S_{\sigma_{\beta}}}/\underline{q}\bigl(\I\bigr)$. Observe that
$\I^{-1} \J \subset\I^{-1}$, where the inclusion is as
$O_L$-modules of rank~$1$ with a notion of positivity. Any
rational polyhedron $\{\sigma_\beta\}_\beta$ in the given rational
polyhedral cone decomposition of the dual cone to
$\bigl(\I^{-1}\bigr)^+_\RR \subset \bigl(\I^{-1}\bigr)_\RR$
induces a rational polyhedron  of the dual cone to
$\bigl(\I^{-1}\J\bigr)^+_\RR\subset \bigl(\I^{-1}\J\bigr)_\RR$.
Hence $\pi_\J\Bigl(\Tate(O_L,\I^{-1})_{\sigma_\beta}\Bigr)$ is the
pullback of $\Tate(\J,\I^{-1})_{\sigma_\beta}$ via the morphism
induced by completing along the boundaries the affine torus
embeddings associated to $\sigma_\beta$ to the isogeny of the tori
$$(\I^{-1}\bigr)^{\vee} \tensor_\ZZ \GG_{m,\ZZ} \subset(\I^{-1}\J\bigr)^{\vee}
\tensor_\ZZ \GG_{m,\ZZ}.$$Since the differential~$dt/t$ descends
to~$\Tate(\J,\I^{-1})$ and corresponds to the differential~$dt/t$
on~$\Tate(\J,\I^{-1})$ we conclude.

\label otherEisenstein. corollary\par\cor  There exists a
$\I$-polarized modular form $\E^\dagger_{k}$ of
weight\/~$\Norm^{k-1}$ and level\/~$\mu_{N}\times \Gamma_0(p)$
i.~e., $\E^\dagger_k\in {\bf M}\bigl(\CC,\mu_{N},\Gamma_0(p),
\Norm^{k-1}\bigr)$, whose $q$-expansion at the
cusp~$(O_L,\I^{-1},\varepsilon,(O_L/p),\j_{\rm can})$ is
$$\left(\prod_{\P\vert p} \bigl(1 - \Norm(\P)^{k-1}\bigr)
\bigl(2^{-g}\zeta_L(1-k)\bigr)\right) + \sum_{\nu\in (\I^{-1})^+ }
\left(\sum_{\C \subset \I^{-1}} \Norm'(\nu
\C^{-1})^{k-1}\right)q^{\nu},$$where
$$\Norm'(\nu
\C^{-1})=\cases{\, \Norm(\nu \C^{-1}) & if $\Norm(\nu \C^{-1})$ is
prime to~$p$\cr\,  0 & otherwise.\cr}$$
\endcor

\Proof Let $Z$ be the set of primes of~$L$ over~$p$. If~$I$ is a
subset of~$Z$, define~$\J_I:=\prod_{\P\in I} \P$. Define
$$\E^\dagger_{k}:= \sum_{I \subset Z} (-1)^{\vert
I\vert} \pi_{\J_I}^*\bigl(\E_{k}\bigr).$$Fix $\nu\in (\I^{-1})^+$.
Let~$I$ be a subset of~$Z$. It follows from~\refn{modified} that
the coefficient of~$q^\nu$ in the $q$-expansion at the given cusp
of~$\pi_{\J_I}^*(\E_k)$ is~$0$, if~$\nu\not\in (\I^{-1}\J_I)^+$,
and is
$$\eqalign{\Norm^{k-1}(\J_I) \left(\sum_{\nu\in\C \subset \I^{-1}\J_I }
\Norm(\nu\C^{-1})^{k-1} \right)& =\sum_{\nu\in\C \subset
\I^{-1}\J_I }
\Norm\left(\nu\bigl(\C\J_I^{-1}\bigr)^{-1}\right)^{k-1}\cr
 &= \sum_{\W \subset \I^{-1}\vert \nu\in \W\J_I}
\Norm\left(\nu\W^{-1}\right)^{k-1} \cr}$$otherwise. Fix~$\W\subset
\I^{-1}$ such that~$\nu\in\W$. Let~$I$ be the maximal subset
of~$Z$ such that $\J_I\vert (\nu \W^{-1})$. The contribution
of~$\Norm(\nu\W^{-1})$ to the coefficient of~$q^\nu$ in the
$q$-expansion of~$\E^\dagger_{k}$ is
$$\left(\,\sum_{I'\subset I} (-1)^{\vert I' \vert}\right)
\Norm(\nu\W^{-1}).$$Since $\sum_{I'\subset I} (-1)^{\vert I'
\vert}$ is~$0$ if~$I\neq \emptyset$ and is~$1$ if~$I=\emptyset$,
we conclude.

\label padicEisenstein. section\par\ssection $p$-adic Eisenstein
series\par  Let $p$ be a prime not dividing~$N$. Let $k_1< k_2<
\ldots<k_n< \ldots$ be a sequence of even integers $\geq 2$
converging $p$-adically to~$k\in \ZZ_p$. Let
$\bigl(O_L,\I^{-1},\j\bigr)$ be a $\I$-polarized unramified cusp
defined over~$\ZZ_{(p)}$. It follows from~\refn{Eisenstein} that
for every $\nu\in \bigl(\I \bigr)^+$ the coefficient~$a_{k_i,\nu}$
of~$q^\nu$ in the $q$-expansion of~$E_{k_i}$ at the given cusp is
$$a_{k_i,\nu}=\left(\,\sum_{\nu\in \C\subset \I^{-1}}
\Norm(\nu\C^{-1})^{k_i-1}\right)\in \QQ.$$It follows from the
$q$-expansion  principle, see~\refn{Koecher}, that
$$E_{k_i} \in {\bf M}\bigl(\QQ,\mu_N,\Norm^{k_i-1}\bigr).$$For $i
\rightarrow \infty$ we have
$$\lim_{i\rightarrow \infty}a_{k_i,\nu}= \Bigl(\sum_{\nu\in\C \subset \I^{-1}}
\Norm'(\nu \C^{-1})^{k-1}\Bigr)$$and the convergence is uniform
in~$k_i$. See~\refn{otherEisenstein} for the definition
of~$\Norm'$. It follows from~\refn{a0} that the sequence
$$a_{k_i,0}=2^{-g}\,\zeta_L(1-k_i)\in\QQ_p$$is bounded. As in~[\Serre1,
Cor.~2, \S1.5], one concludes  that it converges $p$-adically.
Define $$\zeta_L^*(1-k):=\lim_{i\rightarrow \infty}
\zeta_L(1-k_i)\in\QQ_p.$$One may interpret this formula as the
value of the $p$-adic zeta function~$\zeta_L^*$ associated to~$L$
at~$1-k$; see [\Serre1, Thm.~3, \S1.6, Thm.~20, \S5.3]. We also
get that
$$\E^*_k:={\zeta_L^*(1-k)\over
2^g}\,+\,\sum_{\nu\in\in\I^+}\Bigl(\sum_{\nu\C \subset \I}
\Norm'(\nu \C^{-1})^{k-1}\Bigr)q^\nu$$is a $p$-adic modular form
\`a la Serre; see~\refn{Serrepadic}.

\endssection

\label Eisensteinpadic. section\par\ssection $\E_k^\dagger$ as a
$p$-adic Eisenstein series\par Let $k$ be an even, positive
integer. By~\refn{Gamma0paspadic} the $\I$-polarized modular
form~$\E_k^\dagger$ defines a $\I$-polarized $p$-adic modular form
\`a la Katz, and hence \`a la Serre by~\refn{CComPPare}, of
level~$\mu_N$ over~$\QQ_p$. It has the property that its
$q$-expansion at a $\I$-polarized  cusp
$(\A,\B,\varepsilon_{p^\infty N},\j_\varepsilon)$, in the sense
of~\refn{qexppadic}, is the $q$-expansion of~$\E_k^\dagger$ at the
$\I$-polarized cusp $(\A,\B,\varepsilon_N,H,\j_\varepsilon)$. The
latter is given in~\refn{otherEisenstein}. In particular, it has
the same coefficients for $\nu\neq 0$ as the $q$-expansion of the
Serre $\I$-polarized $p$-adic Hilbert modular form~$\E_k^*$ of
level~$\mu_N$ at the cusp~$(\A,\B,\varepsilon,\j_\varepsilon)$. It
follows from~\refn{padicgeneral}, more precisely from the
generalization of~\refn{a0} to the case of $p$-adic modular forms
of level~$\mu_N$, that~$\E_k^\dagger$ has {\it the same
$q$-expansion} as~$\E_k^*$. Hence,
$${\zeta_L^*(1-k)} =\prod_{\P\vert p} \left(\bigl(1 -
\Norm(\P)^{k-1}\bigr)\right)\cdot{\zeta_L(1-k)}\qquad \hbox{{\rm
for }}k\in 2\ZZ,\, k\geq 2.$$Compare with [\Serre1, Rmk 1, \S1.6].
\endssection

\label CFT1. section\par\ssection Notation\par Let~$\P$ be a prime
dividing~$p$ in~$O_L$. Let~$B_\P$ be the maximal abelian
subextension over~$\QQ_p$ of the completion~$L_\P$ of~$L$ at~$\P$
and let~$A_\P$ be its ring of integers. Let~$e'(\P/p)$ be the
ramification index of~$B_\P$ over~$\QQ_p$. Note that~$e'(\P/p)$
divides~$e_\P$. Write~$e'(\P/p)=e'(\P/p)^\tame\cdot
e'(\P/p)^\wild$, where~$e'(\P/p)^\tame$ is the prime to~$p$ part
of~$e'(\P/p)$. Local class field theory gives that
$$H(\P):=\Norm_{B_\P/\QQ_p}(A_\P^*)
$$is equal to~$\Norm_{L_\P/\QQ_p}(O_{L_\P}^*)$. Moreover, for~$p
\neq 2$, $$H(\P)=H_1(\P)\, H_2(\P),$$where $H_1(\P)\subset
\mu_{p-1}$ is the unique subgroup of~$\mu_{p-1}$ of
index~$e'(\P/p)^\tame$ and~$H_2(\P)\subset 1+p\ZZ_p$ is equal
to~$1+p e'(\P/p)^\wild \ZZ_p$, is the unique subgroup
of~$1+p\ZZ_p$ of index~$e'(\P/p)^\wild$. For~$p=2$, $H(\P)$ is the
unique subgroup of index~$e'(\P/p)$ of~$\ZZ_2^*\cong \{\pm
1\}\times \ZZ_2$.

Let $$e_p^\tame:={\rm min}\{e'(\P/p)^\tame:\, \P\vert p\}, \qquad
e_p^\wild:={\rm min}\{e'(\P/p)^\wild:\, \P\vert p\}.$$For~$p=2$ we
use
$$e_2:=e_2^\wild.$$We find that for~$p\neq 2$,
$$H:=\Nm\left((O_L\tensor\ZZ_p)^*\right)=\prod_{\P\vert p} H_1(\P)
\times \prod_{\P\vert p} H_2(\P)=H_1 \times H_2,$$where~$H_1$ is
the unique subgroup of index~$e_p^\tame$ of~$\mu_{p-1}$ and~$H_2$
is the unique subgroup of~$e_p^\wild$ of~$1+p\ZZ_p$.

\noindent For~$p=2$, we find that
$$H:=\Nm\left((O_L\tensor\ZZ_p)^*\right)=\prod_{\P\vert p} H(\P)$$is a subgroup of
index~$e_2$ of~$\ZZ_2^*$.
\endssection

\label CFT2. section\par\ssection Calculation of exponents\par Let
$n\geq 1$ be an integer. We compute the exponent of the abelian
group $$\Hbar:={\rm Image}(H)\quad\hbox{{\rm via the map }}\ZZ_p^*
\rightarrow \bigl(\ZZ_p/p^n\ZZ_p\bigr)^*.$$If $x\in\RR$, we use
the notation
$$\lceil{ x}\rceil:={\rm min}\left\{n\in\ZZ\vert x\leq n\right\}.$$

\noindent If~$p\neq 2$, the exponent is equal to $${p-1 \over
e_p^\tame} \times \left\lceil{{p^{n-1}\over
e_p^\wild}}\right\rceil.$$If~$p=2$, the exponent is equal to
$$2^\epsilon \cdot \left\lceil{{2^{\l(n)}\over
e_2}}\right\rceil,$$where~$\epsilon=0,1$ depending on the case and
$$\l(n):=\cases{ n-1 & if $n\leq 2$,\cr
                 n-2 & if $n \geq 3$;\cr}$$(note that~$2^{\l(n)}$
is the exponent of~$\bigl(\ZZ_2/2^n\ZZ_2\bigr)^*$).
\endssection

\label integrality. theorem\par\thm  Let\/~$k> 1$ be an even
integer. Suppose that\/~$2^{-g}\zeta_L(1-k)$ is not $p$-integral
and let\/~$n=-{\rm val}_p\bigl(2^{-g}\zeta_L(1-k)\bigr)$. Then,
\spacing
\item{{\rm i.}} if $p\neq 2$, $$ k \equiv 0
\qquad\hbox{{\rm mod}} \quad {p-1 \over e_p^\tame}\cdot
\left\lceil{{p^{n-1}\over e_p^\wild}}\right\rceil;$$

\item{{\rm ii.}} if $p=2$, $$ k \equiv 0 \qquad\hbox{{\rm mod}} \quad
\left\lceil{{2^{\l(n)}\over e_2}}\right\rceil.$$\endthm

\Proof Consider the modular form~$f_1:=p^n\E_{k}$. It is a modular
form of weight~$\Norm^k$ over~$\ZZ_p$ and the reduction
modulo~$p^n$ of its $q$-expansion is equivalent
to~$2^{-g}p^n\zeta_L(1-k)$. Let~$f_2$ be the modular
form~$2^{-g}p^n\zeta_L(1-k)$ over~$\ZZ_p$ of weight~$0$.
By~\refn{cong} we conclude that~$\Norm^k\in\X_{\ZZ/p^n\ZZ}(n)$
i.~e., that for any $b\in \bigl(O_L\tensor_\ZZ \ZZ_p)^*$ we
have~$\Norm^k(b)\equiv 1$ modulo~$p^n$. It follows that the
exponent of~$\Hbar$ (see~\refn{CFT2}) divides~$k$ and the theorem
follows.

\thm Let $k,k'\geq 2$ be even integers such that $k \equiv k'$
modulo $(p-1)p^m$ for some non-negative integer~$m$. Then
\spacing
\item{{\rm i. }} if $k \not\equiv 0$
mod~$(p-1)/e_p^\tame$ (and hence~$p\neq 2$), then
$$\kern-5pt\val_p\left\{\biggl(\,\prod_{\P\vert p} (1 - \Norm(\P)^{k-1})\biggr){
\zeta_L(1-k)\over 2^g}- \biggl(\,\prod_{\P\vert p} (1
-\Norm(\P)^{k'-1})\biggr) {\zeta_L(1-k')\over 2^g}\right\}$$ $\geq
m+1;\hfill$
\spacing
\item{{\rm ii. }} if $k \equiv 0$ mod
$(p-1)/e_p^\tame$ and~$p\neq 2$,  then
$$\kern-5pt\val_p\left\{\biggl(\,\prod_{\P\vert p} (1 - \Norm(\P)^{k-1})\biggr){
\zeta_L(1-k)\over 2^g}-  \biggl(\,\prod_{\P\vert p} (1
-\Norm(\P)^{k'-1})\biggr) {\zeta_L(1-k')\over 2^g}\right\}$$ $\geq
m-\val_p(kk')-1-2\val_p\bigl(e_p^\wild\bigr);\hfill$
\spacing
\item{{\rm iii. }} if $p=2$, then  $$\kern-5pt\val_2\left\{\biggl(\,\prod_{\P\vert p} (1 -
\Norm(\P)^{k-1})\biggr){ \zeta_L(1-k)\over 2^g}-
\biggl(\,\prod_{\P\vert p} (1 -\Norm(\P)^{k'-1})\biggr) {
\zeta_L(1-k')\over 2^g}\right\}$$ $\geq
m-2-\val_2(kk')-2\val_2(e_2).\hfill$\endthm

\Proof Let $$\ell:={\rm max}\Bigl\{-{\rm
val}_p\bigl(2^{-g}\zeta_L(1-k)\bigr),-{\rm val}_p
\bigl(2^{-g}\zeta_L(1-k')\bigr),0\Bigr\}$$and let $$\beta=p^\ell
\biggl( 2^{-g}\zeta_L(1-k)\prod_{\P\vert p} (1 - \Norm(\P)^{k-1})
- 2^{-g}\zeta_L(1-k')\prod_{\P\vert p} (1
-\Norm(\P)^{k'-1})\biggr).$$Let $i:=1$ if~$p\neq 2$ and let~$i:=2$
if~$p=2$. Note that if~$x$ is an integer prime to~$p$, then $$
x^k-x^{k'}\equiv 0\quad\hbox{{\rm mod }} p^{m+i}.$$It follows
that$$f:=p^\ell \E^\dagger_{k}-p^\ell \E^\dagger_{k'} -\beta\equiv
0 \qquad\hbox{{\rm mod}}\quad
p^{m+i+\ell}.$$Using~\refn{Eisensteinpadic}, we interpret~$f$ as a
$p$-adic modular form \`a la Katz. It reduces to  function~$0$
on~$\MM\bigl(\ZZ/p^{m+i+\ell}, \mu_{N p^{m+i+\ell}}\bigr)$,
invariant under~$\Gamma_{m+i+\ell}$. (Here~$N$ is any auxiliary
integer~$\geq 4$ and prime to~$p$). It follows that for
all~$\alpha \in \Gamma_{m+i+\ell}=(O_L/p^{m+i+\ell})^*$ we have
$$\alpha^*f-f=\bigl(\Nm^k(\alpha)-1\bigr)p^\ell
E_k^\dagger-\bigl(\Nm^{k'}(\alpha)-1\bigr)p^\ell
E_{k'}^\dagger\equiv 0\qquad\hbox{{\rm mod}}\quad
p^{m+i+\ell}.$$Consider this equation modulo~$p^{m+i}$. Using that
for~$\alpha\in\Gamma_{m+i+\ell}$ we
have~$\Nm^k(\alpha)\equiv\Nm^{k'}(\alpha)$ modulo~$p^{m+i}$, we
find $$\alpha^*f-f=\bigl(\Nm^k(\alpha)-1\bigr)\left(p^\ell
E_k^\dagger-p^\ell E_{k'}^\dagger\right)\equiv 0\qquad\hbox{{\rm
mod}}\quad p^{m+i},$$and from here that
$$\bigl(\Norm^k(\alpha)-1\bigr)\beta\equiv 0\qquad\hbox{{\rm
mod}}\quad p^{m+i}.$$Let~$t$ be the $p$-adic valuation of $\beta$.
Then~$\Nm^k(\alpha)-1\equiv 0$ modulo~$\left\lceil{{p^{m+i-t}}
}\right\rceil$. Let~$n:=m+i-t$. In the notation of~\refn{CFT1}
and~\refn{CFT2}, $H=\Norm\bigl((O_L\tensor \ZZ_p)^*\bigr)$ and we
have:

\item{{1)}} If $p\neq 2$, then~$H=H_1 \times H_2$ where
$H_1$,~$H_2$ are as in loc.~cit.
\itemitem{{1.a)}} If $k\not\equiv  0$ mod~${p-1\over e_p^\tame}$, then
we must have~$t\geq m+i$. In this case also~$\ell=0$
by~\refn{integrality}.  Part~(i) follows.
\itemitem{{1.b)}} If $k\equiv  0$ mod~${p-1\over e_p^\tame}$,
then~\refn{CFT1} implies that~$k\equiv 0$
mod~$\left\lceil{{{p^{n-1}\over e_p^\wild}}}\right\rceil$.
Therefore, $n-1-\val_p(e_p^\wild)\leq \val_p(k)$ and we get
that~$t\geq m+i-1-\val_p(e_p^\wild)-\val_p(k)$. The same holds
for~$k'$ and we conclude that $$\val_p(p^{-\ell}\beta)\geq
m+i-1-\val_p(e_p^\wild)-\ell-{\rm min}\bigl\{{\rm
val}_p(k),\val_p(k')\bigr\}.$$However, by~\refn{integrality},
$\ell\leq{\rm max}\bigl\{\val_p(k),
\val_p(k')\bigr\}+\val_p(e_p^\wild)+1$. Put together this yields
$$\val_p(p^{-\ell}\beta)\geq m-1-2\val_p(e_p^\wild)-\val_p(kk').$$This
proves Part~(ii).

\item{{2)}} If $p=2$, $H$ is a subgroup of index~$e_2$
of~$\ZZ_2^*$. By~\refn{CFT2} its image~$\Hbar$
in~$\bigl(\ZZ_2^*/2^n\ZZ_2\bigr)^*$ is killed
by~$\left\lceil{{2^{\l(n)}\over e_2}}\right\rceil $, but not by
any smaller power of~$2$. Thus,
$\val_2(k)\geq\l(n)-\val_2(e_2)=\l(m+2-t)-\val_2(e_2)$. Recall
that~$t\geq 0$ and~$m\geq 1$. Thus, $t\geq
m-\val_2(e_2)-\val_2(k)$. Therefore,
$$\val_2\bigl(2^{-\ell}\beta\bigr)\geq m-\val_2(e_2)-\ell-{\rm min}\bigl\{{\rm val}_2(k),{\rm
val}_2(k')\bigr\}.$$Using~\refn{integrality}, we find that
$\l(\ell)-\val_2(e_2)\leq \val_2(k)$ and therefore, $$\ell\leq
2+{\rm max}\bigl\{{\rm val}_2(k),{\rm
val}_2(k')\bigr\}+\val_2(e_2).$$Part~(iii) follows from the last
two displayed formulas.

\rmk For $p=2$, one can improve the results given the precise
structure of the completions of~$L$ at primes above~$2$ e.~g.,
when~$2$ is inert; c.f.~[\Gorenn].

\endrmk

\endsection

\section Numerical examples\par

\ssection Example 1\par Consider the field~$L = \QQ(\sqrt{3})$. It
is a real quadratic field of discriminant~$12$, equal to the
totally real subfield of the cyclotomic field obtained by
adjoining to~$\QQ$ the roots of unity of order~$12$. The following
table provides some information on the decomposition of rational
primes in~$L$.

\bigskip
\centerline{\vbox{\offinterlineskip
 \hrule
 \halign{&\vrule#&
  \strut\quad\hfil#\quad\cr
  height2pt&\omit&&\omit&&\omit&&\omit&\cr
 & $p$\hfil&& decomposition\hfil &&
 $e_p^{\rm tame}$\hfil &&
$e_p^{\rm Wild}$\hfil &\cr
 height2pt&\omit&&\omit&&\omit&&\omit&\cr
 \noalign{\hrule}
 height2pt&\omit&&\omit&&\omit&\cr
 & $2$\hfil&& ramified \hfil && $1$\hfil && $2$ \hfil &\cr
 & $3$\hfil&& ramified \hfil && $2$\hfil && $1$ \hfil &\cr
 & $5, 7, 17, 19, 29, 31$\hfil&& inert \hfil && $1$\hfil && $1$ \hfil &\cr
 & $11, 13, 23, 37$\hfil&& split \hfil && $1$\hfil && $1$ \hfil &\cr
 height2pt&\omit&&\omit&&\omit&&\omit&\cr}
 \hrule}}
\bigskip

\noindent The results of the Section~11 imply that  the only odd
primes at which~$\zeta_L(1-k)^{-1}$ can have positive $p$-adic
valuation~$n$ are the primes~$p$ such that~$(p-1)\vert k$, and
then~$n-1$ is at most the power of~$p$ dividing~$k$. For the
prime~$p=2$, we find that if~$n$ is the valuation at~$2$
of~$2^g\zeta_L(1-k)^{-1}$ , then $\val_2(k)\geq \l(n)-1$. In this
case, the $\epsilon$ in~\refn{CFT2} is~$1$ since~$-1$ is not a
norm from~$\QQ[\sqrt{3}]$. Hence, the bound may be improved
to~$\val_2(k)\geq \l(n)$. For example, taking~$k = 18$ the
prediction is that the only odd primes at
which~$\zeta_L(-17)^{-1}$ may have positive valuation~$n$
are~$3,7$ and~$19$ and that valuation can be at most~$3, 1$
and~$1$, respectively. At~$2$ the valuation can be at most~$1$.
Indeed:
$$\eqalign{ \zeta_L(-17) = &  514802473837215246476827/7182 \cr
       = & 2^{-1} \cdot 3^{-3} \cdot 7^{-1} \cdot 11\cdot 19^{-1} \cdot 43867 \cdot 1066866320794499171.\cr}$$
\noindent Another interesting value is the denominator
of~$\zeta_L(-35)$, which is $${\rm
denominator}\bigl(\zeta_L(-35)\bigr)=2^2 \cdot 3^3\cdot 5 \cdot 7
\cdot 13 \cdot 19 \cdot 37.$$ We consider the congruences
involving~$2^{-2}\zeta_L(1-2)$ and~$2^{-2}\zeta_L(1-26)$.
Congruences are predicted for the primes~$2, 3, 5, 7, 13$.  The
prediction is $$\val_2\left\{ (1 - 2^{2-1})2^{-2}\zeta_L(1-2) -
(1-2^{26-1}) {2}^{-2}\zeta_L(1-26)\right\} \geq 3-2-2-2=-3.$$Using
$\zeta_L(-1) = 1/6$ and $$ \zeta_L(-25) =
59603426243912408678663547473670548011/6$$one verifies the
congruence since the valuation is~$-1$.

\noindent For the prime~$3$ the prediction is
$$\val_3\left\{(1-3)2^{-2}\zeta_L(-1) -(1 -
3^{25})2^{-2}\zeta_L(- 25)\right\}\geq 1-0-1-2\cdot 0=0.$$This
indeed holds, since the valuation is~$0$. For the prime~$7$, we
expect
$$\val_7\left\{(1 - 49)2^{-2}\zeta_L(-1) -
(1-49^{25})2^{-2}\zeta(-25)\right\}\geq 0+1=1,$$which holds, since
the valuation is~$1$. For the prime~$13$ the predicted congruence
is
$$\val_{13}\left\{(1-13)^2 2^{-2}\zeta_L(-1) - (1-13^{25})^2
2^{-2}\zeta_L(-25)\right\}\geq 0+1=1.$$This is verified, since the
valuation is~$1$.
\endssection

\ssection Example 2\par Consider the cyclic cubic totally real
field~$L$ of discriminant~$49$, equal to the totally real subfield
of the cyclotomic field of roots of unity of order~$7$. The
following table provides some information on the decomposition of
rational primes in~$L$.

\bigskip
\centerline{\vbox{\offinterlineskip
 \hrule
 \halign{&\vrule#&
  \strut\quad\hfil#\quad\cr
  height2pt&\omit&&\omit&&\omit&&\omit&\cr
 & $p$\hfil&& decomposition\hfil &&
 $e_p^{\rm tame}$\hfil &&
 $e_p^{\rm Wild}$\hfil &\cr
 height2pt&\omit&&\omit&&\omit&&\omit&\cr
 \noalign{\hrule}
 height2pt&\omit&&\omit&&\omit&\cr
 & $7$\hfil&& ramified \hfil && $3$\hfil && $1$ \hfil &\cr
 & $2, 3, 5, 11, 17,
19, 23, 31, 37$\hfil&& inert \hfil && $1$\hfil && $1$ \hfil &\cr
 & $13, 29$\hfil&& split \hfil && $1$\hfil && $1$ \hfil &\cr
 height2pt&\omit&&\omit&&\omit&&\omit&\cr}
 \hrule}}
\bigskip

\noindent  The results of Section~11 imply that if~$p$ is odd, not
equal to~$7$, then the only odd primes at
which~$\zeta_L(1-k)^{-1}$ can have positive $p$-adic valuation~$n$
are the primes~$p$ such that~$(p-1)\vert k$ and then~$n-1$ is at
most the power of~$p$ dividing~$k$. If~$p=7$,
then~$\zeta_L(1-k)^{-1}$ can indeed have positive $7$-adic
valuation~$n$ ($k$ is even) and then~$n-1$ is at most the power
of~$7$ dividing~$k$. For the prime~$p=2$, letting
$n:=\val_2(2^3\zeta_L(1-k)^{-1})$, we find that $\val_2(k)\geq
\l(n)$. This implies that $\val_2({\rm denom.}\zeta_L(1-k))\leq
\val_2(k)-1$. For example, taking~$k = 10$ the prediction is that
the only odd primes at which~$\zeta_L(-9)^{-1}$ may have positive
valuation~$n$ are~$3,7$ and~$11$, and that valuation can be at
most~$1$ in each case. At~$2$, the valuation cannot be positive.
Indeed:
$$\zeta_L(-9) = -1141452324871/231 = -3^{-1}\cdot 7^{-1} \cdot
11^{-1}\cdot 1141452324871.$$We consider congruence for~$7$
for~$2^{-3}\zeta_L(1-2)$ and~$2^{-3}\zeta_L(1-14)$. The expected
congruence is
$$\val_7\left\{(1-7^{3(2-1)})2^{-3}\zeta_L(1-2) -
(1-7^{3(16-1)})2^{-3}\zeta_L(1-14)\right\}\geq -1-1= -2.$$It holds
because the denominators of both values,
$$\zeta_L(-1) = -1/21, \zeta_L(-13) =
-5589087133015782866737/147$$are not divisible by~$7^3$. For~$2$
the expected congruence is $$\val_2\left\{ (1-8)2^{-3}\zeta_L(-1)-
(1-8^{15})2^{-3}\zeta_L(-13)\right\} \geq 2-2-2-2\cdot 0=-2.$$This
is visibly true, because both zeta values have odd numerator. For
$p=13$, the expected congruence is
$$\val_{13}\left\{(1-13)^3 2^{-3}
\zeta_L(-1)-(1-13^{25})2^{-3}\zeta_L(-13) \right\}\geq 0+1=1.
$$This holds, since the valuation is~$1$.

\endssection

\ssection Example 3\par Take the non-Galois totally real cubic
field~$L = \QQ[x]/(x^3-9x -6)$. It has discriminant~$2^{3} \cdot
3^{5}$. The prime~$2$ decomposes as~$\P_1^2\P_2$ and
therefore~$e_2=1$. The prime~$3$ decomposes as~$\P^3$, the
field~$L_3$ is cubic non-Galois and, therefore,
$e_3^\tame=e_3^\wild=1$. We conclude that if an odd prime divides
the denominator of~$\zeta_L(1-k)$ with valuation~$n$, then
$k\equiv 0$ modulo~$(p-1)p^{n-1} $. Analogously, if~$n$ is the
valuation at~$2$ of~$2^3\zeta_L(1-k)^{-1}$, then $\val_2(k)\geq
\l(n)$. For example, if $k=6$, we find $$\zeta_L(-5)=-2\cdot
3^{-2}\cdot 5^2\cdot 7^{-1}\cdot 184669\cdot 512249.$$The
prime~$7$ decomposes as a product of two prime ideals in~$L$ and
$e_7^\tame=e_7^\wild=1$. The expected congruence
for~$2^{-3}\zeta_L(1-2)$ and~$2^{-3}\zeta_L(1-14)$ is
$$(1-7)(1-7^2)2^{-3}\zeta_L(-1)\equiv (1-7^{15})(1-7^{30})2^{-3}\zeta_L(-13) \pmod{7}.$$
Using that $\zeta_L(-1) = -70/3$ and $$\zeta_L(-13)=
-433461315504312280903563360244187028747610/3$$are both divisible
by~$7$, the congruence follows trivially. For the
values~$2^{-3}\zeta_L(1-4)$ and~$2^{-3}\zeta_L(1-16)$ we again
predict
$$\val_7\left\{(1-7^3)(1-7^6)2^{-3}\zeta_L(-3)-
(1-7^{15})(1-7^{30})2^{-3}\zeta_L(-15)\right\}\geq 1.$$The zeta
values are $\zeta(-3) = 2556221/15 $ and $$\zeta_L(-15) =
83822500848624173596590790551322515127580563498549957/1020,$$which
are both $7$-adic units. The congruence holds, though, since the
valuation is~$1$. For the prime~$2$, we predict that
$$\val_2\left\{(1-2^3)(1-4^3)2^{-3}\zeta_L(-3) -
(1-2^{15})(1-4^{15})\zeta_L(-15)\right\}\geq 2-2-6=-6.$$This
holds, since the valuation is~$-5$. Again for the prime~$2$ we
predict
$$\val_2\left\{(1-2)^2 2^{-3}\zeta_L(-1)-(1-2^{17})^2 2^{-3}
\zeta_L(-17)\right\} \geq 4-2-2-2\cdot 0=0.$$Using
$$\eqalign{\zeta_L(-17)=&-\bigl(3647421225841578953319613809666454838832065018732125543\cr&
06326430\bigr)\cdot 3591^{-1}\cr}$$one verifies that the valuation
is~$1$. In particular, the congruence holds.

\endssection

\endsection

\section Comments regarding values of zeta functions\par
\noindent  We make a few remarks on values of zeta functions. We
have no real theorem to offer here, but rather we would like to
draw the reader's attention to the amount of mathematics hidden
behind the innocuous Bernoulli numbers, or even just their
denominators. The connection to zeta functions rests on the
identity $$\zeta(1 - 2k) = -{B_{2k}\over 2k}, \quad k \geq 1.$$To
begin with, consider the question of which normalized Eisenstein
series~$\E_{2k}$ for the modular group ${\rm SL}_2(\ZZ)$ are
congruent to one modulo~$p$. Here~$p$ is a fixed odd prime. Since
the Eisenstein series has the form $$ 1 + {2 \over \zeta(1 -
2k)}\sum_{n=1}^\infty \sigma_{2k-1}(n)q^n,$$the question is for
which~$k$ does~$p$ divide the denominator of~$\zeta(1-2k)$? (As an
aside we mention that this means that there is an element of
order~$p$ is the~$4k-1$ stable homotopy groups of the spheres.
Cf.~[\MilnorStasheff,App.~B]; note that~$B_n$ in their notation is
our~$B_{2n}$). The Kummer congruences imply that this is the case
iff $2k\equiv 0 \pmod{p-1}$ and then $$\val_p(\zeta(1 - 2k)) =
-1-\val_p(2k);$$see~[\Serre1, \S1.1]. Assume, for argument's sake,
that~$k$ is prime. Then~$p-1$ is either~$1$, $2$, $k$ or~$2k$.
Apart from $p = 2$, $3$, assuming $k>3$, we are left only with the
possibility that $p = 2k+1$, i.~e., that~$k$ is a Germain prime!
There are infinitely many primes that are {\it not} Germain primes
(take~$k$ to be a prime congruent to~$1$ mod.~$6$), hence the
order of the denominator of~$\zeta(1-2k)$ does not grow to
infinity with~$k$. In fact, one can prove that if all the prime
factors of~$k$ are congruent to~$1$ modulo~$6$ then the
denominator of~$\zeta(1 - 2k)$ is equal to~$12$ ([\MilnorStasheff,
App.~B, Pb~B-1]).

\spacing
\indent Let $p\equiv 1 \pmod{4}$ and let~$K= \QQ(\sqrt{p})$. Let
$\chi = \left( {\cdot\over p} \right)$ be the corresponding
Dirichlet character. Let~$\omega$ be the Teichm\"uller character,
then~$\chi = \omega^{(p-1)/2}$. We note also that the discriminant
of~$L$ is~$p$ and is equal to the conductor~$f_\chi$ of~$\chi$.

\prop Let\/~$p$ be a prime congruent to~$1$ modulo~$4$ and let\/
$K := \QQ(\sqrt{p})$. Assume Vandiver's conjecture ([\Washington,
p.~159]). Then for any integer~$r \geq 1$ we have $$
\val_p(\zeta_K(1 - r(p-1))) = -1-\val_p(r).$$ Hence, in this case,
the Hasse invariant in characteristic~$p$ lifts to the Eisenstein
series~$\E_{K, p-1}$.
\endprop
\Proof We have the identity $$ \zeta_K(1 - r(p-1)) = \zeta_\QQ(1 -
r(p-1))\cdot L(1 - r(p-1),\chi).$$Therefore, the claim would
follow if we prove that $\val_p(L(1 - r(p-1), \chi)) = 0$.
Using~[\Washington, Thms 4.2, 5.11] we get that for any~$n\geq 1$
$$ L_p(1 - n, \chi) = (1 - \chi\omega^{-n}(p)p^{n-1})\cdot L(1
-n, \chi\omega^{-n}).$$ Applying this formula in our situation we
get $$\val_p( L(1 - r(p-1),\chi)) = \val_p( L_p(1 -
r(p-1),\chi)).$$By [\Washington, Cor.~5.13], for any integers~$m,
n$ we have $L_p(m, \chi) \equiv L_p(n, \chi)$ mod.~$p$ and both
numbers are $p$-integral. It is thus enough to prove for a {\it
single} integer~$m$ that~$L_p(m, \chi)$ is a $p$-adic unit. Using
the results cited above we find
$$\eqalign{ \val_p(L_p(1 - (p-1)/2, \chi))
 & = \val_p(L(1-(p-1)/2, \chi\omega^{(p-1)/2})) \cr
 & = \val_p(\zeta(1 - (p-1)/2)) \cr
 & =\val_p(B_{(p-1)/2}).\cr}$$The assertion $\val_p(B_{(p-1)/2}) = 0$
is known as the Ankeny-Artin-Chowla conjecture and  was verified
for $p < 10^{11}$ by A.~J.~van der Poorten, H.~te Riele and
H.~Williams in~[\vanderPoorten]. This conjecture is a consequence
of Vandiver's conjecture. See~[\Washington, Thm 5.34].

\prop Let\/~$p$ be a prime congruent to~$1$ modulo~$4$ and let $K
= \QQ(\sqrt{p})$. For every odd integer~$r \geq 1$ $$
\val_p(\zeta_K(1 - r(p-1)/2)) = -1-\val_p(r).$$Hence the partial
Hasse invariant defined in~\refn{hPi} admits a lift to char.~$0$,
the Eisenstein series $\E_{K, (p-1)/2}$.
\endprop
\Proof We have $$\zeta_K(1 - r(p-1)/2) = \zeta_\QQ(1 -
r(p-1)/2)\cdot L(1 - r(p-1)/2, \chi).$$Note that $r(p-1)/2
\not\equiv 0$ mod.~$p-1$ and therefore $$\val_p(\zeta_K(1 -
r(p-1)/2)) = -1-\val_p(r)\; \Leftrightarrow \; \val_p(L(1 -
r(p-1)/2, \chi)) = -1-\val_p(r).$$However, $$\eqalign{
 L_p(1 -r(p-1)/2, {\bf 1}) & = (1 -
\omega^{r(p-1)/2}(p)p^{r(p-1)/2-1})L(1
-r(p-1)/2,\omega^{r(p-1)/2}) \cr
 & = L(1 - r(p-1)/2, \chi),\cr}$$and therefore,
$$\val_p(L(1 - r(p-1)/2, \chi))  = \val_p(L_p(1 - r(p-1)/2,
{\bf 1})).$$Using [\Washington, Ex.~5.11], we write $$ L_p(s, {\bf
1}) = {p-1\over p}\cdot {1\over s-1} + a_0 + a_1(s-1) + a_2(s-1)^2
+ \dots,$$where $a_i \in \ZZ_p$ for every~$i$. In particular,
$\val_p(L_p(1 - r(p-1)/2, {\bf 1})) = -1-\val_p(r)$.

\endsection

\section The operators $\Theta_{\P,i}$\par
\noindent This section is devoted to the construction of certain
derivation operators on $p$-adic Hilbert modular forms and on
Hilbert modular forms in characteristic~$p$. For~$p$ unramified,
these operators were defined by Katz~[\Katzzzz, \S2.5]. However,
our construction is independent; more importantly, in
characteristic~$p$ the operators defined by Katz are defined only
on the ordinary locus, while our operators are defined on the
whole space.

\noindent Let $R$ be a discrete valuation ring of unequal
characteristic~$p$-$0$ and with maximal ideal~$\m$. Our method is
first to define these operators as derivation operators on
functions on~$\MM\bigl(R/(\m^n), \mu_{p^nN}\bigr)$ and then to use
the map~$r$ of~\refn{r(f)}, which relates modular forms
on~$\MM\bigl(R/\m^n, \mu_{N})$ to functions on~$\MM\bigl(R/\m^n,
\mu_{p^nN}\bigr)$, to transport the operators to modular forms
on~$\MM\bigl(R/\m^n, \mu_{N}\bigr)^{\rm ord}$.

\noindent In characteristic~$p$, we succeed in defining  theta
operators on modular forms defined on~$\MM\bigl(R/\m,
\mu_{N}\bigr)$. After establishing the existence of these
operators and some basic properties we examine how  the divisor of
a modular form changes under such a derivation operator. This is
applied to study how the filtration of a $q$-expansion changes
under these derivation operators.

\spacing
\noindent This section is technically demanding, yet forms the
core of the rest of the paper. To orient the reader we explain its
structure.

\noindent {\it Sections 15.1-15.6} recall the definition of the
Kodaira-Spencer isomorphism in the setting of the schemes~$\MM(m,
n)$ ($n\geq m$). The main use we make of this isomorphism is to
construct a canonical basis for the holomorphic $1$-forms
on~$\MM(m, n)$ from the modular forms~$a(\chi_{\P, i})$
of~\refn{universal}.

\noindent {\it Sections 15.7-15.10} are devoted to constructing
machinery to be used in the following definition of the theta
operators. The complications arise from ramification. One of our
goals is to have a certain operator~$\Lambda$ -- constructed out
of theta operators -- on modular forms, such that~$\Lambda(f)$ is
in the image of the operator~$V$. The operator~$V$ is essentially
raising to a $p$-power (see Section~16). Hence, we
need~$\Lambda(f)$ to have $q$-expansion of the form~$\sum_{\nu\in
O_L^+} a(p\nu)q^{p\nu}$. On the other hand, our theta operator
associated to a character~$\chi$ yield $q$-expansion of the
form~$\sum_{\nu\in O_L^+}  \chi(\nu)a(\nu)q^{\nu}$ in which ``too
many" coefficients are zero; we find that in the case of
ramification the theta operators need an extra modification for
which \S\S~15.7-15.9 provide technical background.

\noindent {\it Sections 15.11-15.13} provide the definition of the
theta operators in the mod $p$ and $p$-adic setting. In essence,
the theta operators are coming from the operation~$f \mapsto df$
and the canonical trivialization of the holomorphic one forms on
the schemes~$\MM(m, n)$. The propagation of this definition to
modular forms of level~$\mu_N$ is carried out later.

\noindent {\it Sections 15.14-15.25} are devoted to the
calculation of the effect of the theta operators on $q$-expansions
and various corollaries.

\noindent {\it Sections 15.26-15.33} are devoted to examining how
the poles of a rational function on~$\MMbar(m,n)^\Kum$, whose
poles are supported on the complement of the ordinary locus,
change under a theta operator. For this we use the local charts
of~$\MMbar(m, n)^\Kum$ constructed in~\refn{Mbar}. Ultimately this
will be applied to determining how weights of modular forms change
under a theta operator.

\noindent {\it Sections 15.34-15.38} apply the previous results to
define and derive the properties of theta operators on modular
forms. We remark that on the level of Galois representations, an
application of a  theta operator corresponds to a twist by a Hecke
character. The change of filtration under a theta operator is
examined in Section~18; that corresponds to the question of the
minimal weight from which a Galois representation arises, up to a
twist.

\ssection Notation\par Let $R$ be a complete discrete valuation
ring with maximal ideal~$\m$, residue field~$k$ of
characteristic~$p$ and fraction field~$F$ of characteristic~$0$.
Suppose that~$R$ is a $O_K$-algebra where~$K$ is a normal closure
of~$L$.
\endssection

\defi Let $n \geq n'\geq  m\geq 1$ be integers. Let $$\MM(m,n)
\llongrightarrow \MM(m,n')$$be the Galois cover
over~$\Spec(R/\m^m)$ defined in~\refn{phi}. We also write
$$\MM(m,0):=\MM(R/\m^m\mu_N)^{\rm ord}.
$$Let
$$d\colon O_{\MM(m,n)} \llongrightarrow
\Omega^1_{\MM(m,n)/(R/\m^m)}\qquad\hbox{{\rm and}}\qquad d'\colon
O_{\MM(m,n')} \llongrightarrow \Omega^1_{\MM(m,n')/(R/\m^m)}$$be
the sheaves of differentials of\/~$\MM(m,n)$
(resp.~of\/~$\MM(m,n')$) over~$(R/\m^m)$. Denote by
$$\pi\colon\AA\longrightarrow\MM(m,n)\qquad\hbox{{\rm and}}\qquad
\pi'\colon\AA'\llongrightarrow\MM(m,n')$$the universal
$\I$-polarized abelian schemes.

\noindent Let
$$\omega_{\AA/\MM(m,n)}:=\pi_*\Bigl(\Omega^1_{\AA/\MM(m,n)}
\Bigr)\qquad\hbox{{\rm and}}\qquad
\omega_{\AA'/\MM(m,n')}:=\pi'_*\Bigl(\Omega^1_{\AA'/\MM(m,n')}
\Bigr).$$
\enddefi

\rmk By~\refn{phi} we have  canonical isomorphisms $$\AA \isomarrow
\AA'\fibprod_{\MM(m,0)}\MM(m,n)\qquad\hbox{{\rm and}}\qquad
\omega_{\AA/\MM(m,n)} \isomarrow
\phi^*\bigl(\omega_{\AA'/\MM(m,0)}\bigr).$$We deduce that
$O_{\MM(m,n)}$,~$\Omega^1_{\MM(m,n)/(R/\m^m)}$
and~$\Omega^1_{\AA/\MM(m,n)}$ are canonically endowed with an
action of\/~$\Gamma_n$ (see~\refn{phi} for the notation). Note
that~$d$ is~$\Gamma_n$-equivariant. It follows
from~\refn{propRapo} that\/~$\Omega^1_{\AA/\MM(m,n)}$ is a locally
free $ O_{\MM(m,n)}\tensor_\ZZ O_L$-module of rank~$1$ and hence
it is endowed with an action of~$\bigl( (R/\m^m)\tensor_\ZZ O_L
\bigr)^*$. Let $\chi$ be a universal character. Let~${\cal
L}_\chi$ the line bundle associated to~$\Omega^1_{\AA/\MM(m,n)}$
and the character~$\chi$ as in~\refn{Lchi}. Then~${\cal L}_\chi$
is endowed with two actions of\/~$\bigl(O_L/p^nO_L\bigr)^* $:

\spacing
\item{{\rm a.}} the first is induced by the Galois action of
$\Gamma_n$ on~$\Omega^1_{\AA/\MM(m,n)}$;

\spacing
\item{{\rm b.}} the second is induced via the natural map
$$\bigl(O_L/p^nO_L\bigr)^* \llongrightarrow \bigl(
(R/\m^m)\tensor_\ZZ O_L\bigr)^* $$by the action of~$\bigl(
(R/\m^m)\tensor_\ZZ O_L\bigr)^*$. The latter is identified with
the automorphism group of~$\Omega^1_{\AA/\MM(m,n)}$ as
$O_{\MM(m,n)} \tensor_\ZZ O_L$-module.
\spacing

\noindent By~\refn{actiona(chi)} the two actions are {\it one the
inverse of the other}.

\endrmk

\label TAU. section\par\ssection Notation\par Fix an
element~$l\in\I^+$ generating~$\ZZ_p\tensor_\ZZ \I$ as
$\ZZ_p\tensor_\ZZ O_L$-module as in~\refn{k}. For any
$\I$-polarized abelian scheme with real multiplication by~$O_L$
over a $\ZZ_p$-scheme $S$ such a choice induces a prime to~$p$ and
$O_L$-linear polarization.

\endssection

\rmk The sequence of differentials \labelf exactdifferentials\par
$$0 \llongrightarrow
\pi^*\bigl(\Omega^1_{\MM(m,n)/(R/\m^m)}\bigr)\llongrightarrow
\Omega^1_{\AA/(R/\m^m)} \llongrightarrow \Omega^1_{\AA/\MM(m,n)}
\llongrightarrow 0\eqno{(\numfo)}$$\advance\fonu by1

\noindent is exact. The exactness on the right and in the middle
is obvious; see~[\EGAIVfour, Cor.~16.4.19]. The exactness on the
left follows from the smoothness of~$\AA$ over~$\MM(m,n)$;
see~[\EGAIVfour, Prop.~17.2.3]. The sheaves appearing in the
sequence above are locally free $O_{\MM(m,n)}$-modules. The
functor from the category of sheaves on~$\AA$, invariant under
translation on~$\AA$, to the category of sheaves on~$\MM(m,n)$
defined by pulling back along the identity section of\/~$\pi$
defines an equivalence of categories. Hence, there is a canonical
isomorphism of $O_\AA$-modules
$$\pi^*\bigl(\Omega^1_{\MM(m,n)/(R/\m^m)}\bigr) \isomarrow
\Omega^1_{\MM(m,n)/(R/\m^m)} \tensor_{O_{\MM(m,n)}}
O_\AA.$$Applying the functor~$\pi_*$ to the
sequence~(\refn{exactdifferentials}), we obtain an exact sequence
of $O_{\MM(m,n)}$-modules

\bigskip
\noindent $0\longrightarrow\Omega^1_{\MM(m,n)/(R/\m^m)}
\longrightarrow \pi_*\Bigl(\Omega^1_{\AA/(R/\m^m)}\Bigr)$
\spacing

\rightline{$\longrightarrow \omega_{\AA/\MM(m,n)} \longrightarrow
\kern-3.0ptR^1\pi_*\Bigl(\pi^*\bigl(\Omega^1_{\MM(m,n)/(R/\m^m)}
\bigr)\Bigr)$.} \bigskip\noindent We have canonical isomorphisms
$$R^1\pi_*\Bigl(\pi^*\bigl(\Omega^1_{\MM(m,n)/(R/\m^m)}
\bigr)\Bigr)\isomarrow \Omega^1_{\MM(m,n)/(R/\m^m)}
\tensor_{O_{\MM(m,n)}} R^1\pi_*\bigl( O_\AA\bigr) $$
$$\isomarrow \Omega^1_{\MM(m,n)/(R/\m^m)} \tensor_{O_{\MM(m,n)}}
\Hom_{O_{\MM(m,n)}}\bigl(\omega_{\AA^\vee/\MM(m,n)},O_{\MM(m,n)}\bigr),$$where
$\AA^\vee$ is the dual abelian scheme of\/~$\AA$. By~\refn{TAU} we
get  a prime to~$p$ and $O_L$-linear polarization $A \rightarrow
A^\vee$. This induces an isomorphism
$$\omega_{\AA/\MM(m,n)}\tensor_{O_{\MM(m,n)}} \omega_{\AA/\MM(m,n)}
\cong \omega_{\AA/\MM(m,n)}\tensor_{O_{\MM(m,n)}} \omega_{\AA^\vee
/\MM(m,n)}$$and hence a
morphism$$\omega_{\AA/\MM(m,n)}\tensor_{O_{\MM(m,n)}}
\omega_{\AA/\MM(m,n)}
 \llongrightarrow \Omega^1_{\MM(m,n)/(R/\m^m)}.$$

\noindent A similar remark holds for $\AA'$ in place of\/~$\AA$
and~$\MM(m,n')$ in place of\/~$\MM(m,n)$.
\endrmk

\label KS. proposition\par \prop The maps described above induce
canonical isomorphisms, called Kodaira-Spencer isomorphisms,
$$\KS'\colon \Omega^1_{\MM(m,n')/(R/\m^m)} \llongrightarrow
\omega_{\AA'/\MM(m,n')}^{\tensor_{O_L}^2}$$and
$$\KS\colon
\Omega^1_{\MM(m,n)/(R/\m^m)} \llongrightarrow
\omega_{\AA/\MM(m,n)}^{\tensor_{O_L}^2},$$where $ \tensor_{O_L}^2$
means tensor product {\it as $O_L$-modules}. In particular,
$\Omega^1_{\MM(m,n')/(R/\m^m)}$
(resp.~$\Omega^1_{\MM(m,n)/(R/\m^m)}$) is endowed with the
structure of free $O_{\MM(m,n')} \tensor_\ZZ O_L$-module (resp.~of
$O_{\MM(m,n)}\tensor_\ZZ O_L$-module) of rank\/~$1$. Moreover, the
following diagram is commutative:

$$\matrix{
\Omega^1_{\MM(m,n)/(R/\m^m)} & \lllongmaprighto{\KS}
&\omega_{\AA/\MM(m,n)}^{\tensor_{O_L}^2} \cr \mapdownr{\wr} & &
\mapdownr{\wr} \cr \phi^*\bigl(
\Omega^1_{\MM(m,n')/(R/\m^m)}\bigr) &
\lllongmaprighto{\phi^*(\KS)}
&\omega_{\AA'/\MM(m,n')}^{\tensor_{O_L}^2}.\cr}$$In particular,
the map~$\KS$ is $\Gamma_n$-equivariantly.

\endprop
\Proof The fact that $\KS$ and~$\KS'$ are isomorphisms is well
known. See~[\Katzzzz, \S1.0.21]. The rest is clear.

\label omegaPi. section\par\ssection Recall\par The notation is as
in~\refn{universal}. Define
$$\Cand \tensor \Cand \in \Gamma\bigl(\MM(m,n),
\omega_{\AA/\MM(m,n)}^{\tensor_{O_L}^2}\bigr)$$ $$\hbox{{\rm
(resp.
}}\canD\in\Gamma\bigl(\MM(m,n),\Omega^1_{\MM(m,n)/(R/\m^m)}\bigr)\hbox{{\rm
)}}$$the canonical generators of
$\omega_{\AA/\MM(m,n)}^{\tensor_{O_L}^2}$ (resp. $
\Omega^1_{\MM(m,n)/(R/\m^m)}$) as $ O_{\MM(m,n)}\tensor_\ZZ
O_L$-module.

\spacing
\noindent For each prime~$\P$ of~$O_L$ dividing~$p$ let~$\pi_\P\in
O_L$ be a generator of the ideal~$\P\,O_{L,\P}$ as
in~\refn{notAtion}. Let~$\p:=\m \cap O_K$. For each integer~$1\leq
i\leq f_\P$, let
$$\sigma_{\P,i}\colon R\tensor_\ZZ O_L \llongrightarrow
R\qquad\hbox{{\rm and}}\qquad e_{\P,i}\in R\tensor_\ZZ O_L$$be as
in~\refn{basicwt}. Recall that $R$ is assumed to be a
$\m$-adically complete $O_K$-algebra.

\endssection

\label sigmatilde. lemma\par\lemma Let $m\geq 1$. Let\/~$\P$ be a
prime of\/~$O_L$ over~$p$ and let $0\leq \jj\leq e_\P-1$. There
exists a maximal non-negative integer~$1 \leq t_\P^{[\jj]}(m)\leq
m$ satisfying the following. For each $1\leq i\leq f_\P$  there
exists a unique morphism
$$\sigmatilde_{\P,i}^{[\jj]}\colon
 (R/\m^m)\tensor_\ZZ \P^\jj \llongrightarrow (R/\m^{t_\P^{[\jj]}(m)})$$of
$O_L\tensor_\ZZ (R/\m^m)$-modules such that the diagram
$$\matrix{  (R/\m^m) \tensor_\ZZ  O_L& \lllongmaprighto{\cdot (1\tensor \pi_\P^\jj)}&
 (R/\m^m)\tensor_\ZZ \P^\jj \cr \mapdownl{\sigma_{\P,i}} & \swarrow
\sigmatilde_{\P,i}^{[\jj]} \cr
R/\m^{t_\P^{[\jj]}(m)}\cr}$$commutes. Moreover, the sequence
$\{t_\P^{[\jj]}(m)\}_{m\in \NN}$ is non-decreasing and $$\lim_{m
\rightarrow \infty} t_\P^{[\jj]}(m)=\infty.$$
\endlemma
\Proof Fix $0\leq \j\leq e_\P-1$. The morphism
$\sigmatilde_{\P,i}^{[\jj]}$, whose existence is claimed in the
lemma exists if and only if $\Ker\bigl(\cdot (1 \tensor
\pi_\P^\jj)\bigr)\subset \Ker\bigl(\sigma_{\P,i}\bigr)$. Consider
first the case~$m=1$. Then~$R/\m^1=k$. Moreover, $a\in
\Ker\bigl(\cdot (1\tensor \pi_\P^\jj)\bigr)$ if and only if
$a=(1\tensor\pi_\P^{e_\P-\jj})\cdot  b$ for some $b\in
k\tensor_\ZZ O_L$. Since $e_\P-\jj>0$, it follows that
$\sigma_{\P,i}(a)=0$. This proves that $t_\P^{[\jj]}(1)=1$ and
that $t_\P^{[\jj]}(m)\geq t_\P^{[\jj]}(1)=1$ for any~$m$. Clearly
$t_\P^{[\jj]}(m+1)\geq t_\P^{[\jj]}(m)$. Since $\cdot
(1\tensor\pi_\P^\jj)\colon R\tensor_\ZZ O_L \rightarrow
R\tensor_\ZZ O_L$ is injective and~$R$ is $\m$-adically complete,
it is clear that
$$\lim_{\infty\leftarrow m} \Ker\Bigl(\cdot (1\tensor\pi_\P^\jj
)\vert_{ (R/\m^m)\tensor_\ZZ O_L}\Bigr)=0.$$Hence, $\lim_{m
\rightarrow \infty} t_\P^{[\jj]}(m)=\infty$ as claimed.

\label Sigmatilde. definition\par \defi Let $\P$ be a prime
of\/~$O_L$ over~$p$ and let $1\leq i\leq f_\P$. Define
$$\sigma_{\P,i}\colon \omega_{\AA/\MM(m,n)}^{\tensor_{O_L}^2} \llongrightarrow {\cal L}_{\chi_{\P,i}^2}
$$to be the unique morphism of $O_{\MM(m,n)}$-modules
such that

\item{{\rm 1)}}
$\sigma_{\P,i}(\Cand\tensor\Cand)=a\bigl(\chi_{\P,i}^2\bigr)$.
See~\refn{actiona(chi)} for the definition
of\/~$a\bigl(\chi_{\P,i}\bigr) $;

\item{{\rm 2)}} for any $l\in O_L$ and any local section~$\omega$ of
$\omega_{\AA/\MM(m,n)}^{\tensor_{O_L}^2}$ we have
$$\sigma_{\P,i}(\alpha\cdot \omega)=\sigma_{\P,i}(\alpha)
\sigma_{\P,i}(\omega).$$

\spacing
\noindent Let $0\leq \jj\leq e_\P-1$ be an integer. Let
$t_\P^{[\jj]}(m)$ be as in~\refn{sigmatilde}. Let
$${\cal L}_{\chi_{\P,i}^2} \llongrightarrow
\MM\bigl(t_\P^{[\jj]}(m),n\bigr)$$be the line bundle defined
in~\refn{Lchi}. Define
$$\sigmatilde_{\P,i}^{[\jj]}\colon
\omega_{\AA/\MM(m,n)}^{\tensor_{O_L}^2}\tensor_{O_L} \P^\jj
\llongrightarrow {\cal L}_{\chi_{\P,i}^2}$$as the unique morphism
of $O_L\tensor_\ZZ O_{\MM(m,n)}$-modules such that the following
diagram commutes
$$\matrix{ \omega_{\AA/\MM(m,n)}^{\tensor_{O_L}^2} & \lllongmaprighto{\cdot (1\tensor \pi_\P^\jj)}&
\omega_{\AA/\MM(m,n)}^{\tensor_{O_L}^2}\tensor_{O_L} \P^\jj\cr
\mapdownl{\sigma_{\P,i}} &  \swarrow \sigmatilde_{\P,i}^{[\jj]}\cr
{\cal L}_{\chi_{\P,i}^2} & &.\cr}$$Its existence is guaranteed
by~\refn{sigmatilde}.
\enddefi

\label invarianceSigmatilde. lemma\par\lemma  The morphisms
$\sigmatilde_{\P,i}^{[\jj]} $ satisfy the following properties
\spacing
\item{{\rm 1)}} they are compatible for different~$m$ and~$n$
i.~e., if $n'\geq n$ and $m'\geq m$ are integers such that $n'\geq
m'$ and $n\geq m$, then
\itemitem{{\rm 1.a)}} the morphism $\sigmatilde_{\P,i}^{[\jj]} $, defined
on~$\omega_{\AA/\MM(m',n')}^{\tensor_{O_L}^2}\tensor_{O_L} \P^\jj
$, restricts mod~$\m^m$ to the
morphism~$\sigmatilde_{\P,i}^{[\jj]} $, defined on
$\omega_{\AA/\MM(m,n')}^{\tensor_{O_L}^2}\tensor_{O_L} \P^\jj $;
\itemitem{{\rm 1.b)}} the morphism $\sigmatilde_{\P,i}^{[\jj]} $,
defined on~$\omega_{\AA/\MM(m,n')}^{\tensor_{O_L}^2}\tensor_{O_L}
\P^\jj $, coincides with the pull-back of the
morphism~$\sigmatilde_{\P,i}^{[\jj]} $ defined
on~$\omega_{\AA/\MM(m,n)}^{\tensor_{O_L}^2}\tensor_{O_L} \P^\jj $.
\spacing
\item{{\rm 2)}} Let~$\Gamma_n$ be the Galois group of $\MM(m,n)\rightarrow \MM(m,0)$.
For any $\alpha\in \Gamma_n$, we have
$$\alpha^*\bigl(\sigmatilde_{\P,i}^{[\jj]}\bigr)=\sigmatilde_{\P,i}^{[\jj]}.$$
\endlemma
\Proof By the definition of~$\sigmatilde_{\P,i}^{[\jj]}$, it
suffices to prove the lemma for $\jj=0$ i.~e., for
$\sigma_{\P,i}$. Both $\omega_{\AA/\MM(m,n)}^{\tensor_{O_L}^2}$
and~${\cal L}_{\chi_{\P,i}^2}$ are defined over~$\MM(m,0)$ and are
compatible for different~$m$'s. In particular, they are endowed
with a canonical action of~$\Gamma_n$ proving that statement~(2)
makes sense. By~\refn{Lchi} the line bundle~${\cal
L}_{\chi_{\P,i}^2}$ over~$\MM(m,0)$ is defined by push-out
of~$\omega_{\AA/\MM(m,0)}^{\tensor_{O_L}^2}$ by the
map~$\sigma_{\P,i}$ of~\refn{sigmatilde}. By~\refn{actiona(chi)}
this defines the map~$\sigmatilde_{\P,i}^{[\jj]}$ over~$\MM(m,0)$.
The conclusions follow.

\label theta. definition\par \defi  For each prime~$\P$ of~$O_L$
dividing~$p$ and each integer~$1\leq i\leq f_\P$, define the
$(R/\m^m)$-derivation $$\Theta_{\P,i}\colon O_{\MM(m,n)}
\,\lllongrightarrow\,  O_{\MM(m,n)}$$by the formula $$ f
\llongmapsto \sigma_{\P,i}\Bigl(\KS\bigl(df\bigr)\Bigr)\cdot
a\bigl(\chi_{\P,i}^2\bigr)^{-1}.$$See~\refn{omegaPi} for the
definition of~$\sigma_{\P,i}$. For~$0 \leq \jj \leq e_\P-1$ define
the subsheaf of~$O_{\MM(m,n)}$
$$O_{\MM(m,n)}^{\P,[\jj]}(U):=\left\{f
\in \Gamma\bigl(U,O_{\MM(m,n)})\vert\,\bigl(df \bigr)\equiv 0\,
\hbox{{\rm mod }}\P^\jj\right\}.$$Define the $(R/\m^m)$-derivation
$$\Theta_{\P,i}^{[\jj]}\colon O_{\MM(m,n)}^{\P,[\jj]}
\lllongrightarrow O_{\MM(t_\P^{[\jj]}(m),n)},$$considered as a
sheaf on~$\MM(m,n)$, by the formula
$$ f  \llongmapsto \sigmatilde_{\P,i}^{[\jj]}\Bigl(\KS\bigl(df\bigr)\Bigr)\cdot
a\bigl(\chi_{\P,i}^2\bigr)^{-1}.$$See~\refn{Sigmatilde} for the
definition of~$\sigmatilde_{\P,i}^{[\jj]}$ and~\refn{sigmatilde}
for the definition of~$t_\P^{[\jj]}(m)$. Note that
$$O_{\MM(m,n)}^{\P,[0]}
=O_{\MM(m,n)}$$and $\Theta_{\P,i}^{[0]}$ coincides
with~$\Theta_{\P,i}$.
\enddefi
\label wttheta. proposition\par\prop Let $f$ be a regular function
on~$\MM(m,n)$. Suppose that\/~$f$ is an eigenfunction for the
action of\/~$\Gamma_n$ of weight $\psi\colon \Gamma_n \rightarrow
(R/\m^m)^*$ i.~e., $\alpha^*\bigl(f\bigr)=\psi\bigl(\alpha\bigr)
f$. Then, $\Theta_{\P,i}^{[\jj]}\bigl(f\bigr)$, whenever defined,
is an eigenfunction for~$\Gamma_n$ of weight $\Psi\cdot
\chi_{\P,i}^2$.
\endprop

\Proof  For any $\alpha\in \Gamma_n$ denote by $\alpha^*$ the
induced action on functions, differentials, etc. Using~\refn{KS}
and\/~\refn{actiona(chi)}, we obtain:

$$\eqalign{\psi(\alpha) \KS\bigl(df\bigr) & = \KS\bigl(\psi(\alpha) df\bigr)
= \KS\left(d\bigl(\psi(\alpha) f\bigr)\right) \cr
 & = \KS\left(d\bigl(\alpha^* f\bigr)\right)
  = \KS\left(\alpha^*\bigl(df\bigr)\right)\cr
 & = \alpha^*\left(\KS\bigl(df\bigr)\right).\cr}$$Applying
$\sigmatilde_{\P,i}^{[\jj]} $ on the first and last terms of these
inequalities and using~\refn{invarianceSigmatilde}, we have
$$\eqalign{\psi(\alpha) \sigmatilde_{\P,i}^{[\jj]}  \KS\bigl(df\bigr) & =
\sigmatilde_{\P,i}^{[\jj]} \bigl( \psi(\alpha) \KS(df)\bigr) \cr
 & =\sigmatilde_{\P,i}^{[\jj]}\Bigl(\alpha^*\left(\KS\bigl(df\bigr)\right)\Bigr)
\cr
 & =\alpha^*\left(\sigmatilde_{\P,i}^{[\jj]}\bigl(\KS(df)\bigr)\right).\cr}
$$By\/~\refn{actiona(chi)} we have that
$$\alpha^*\bigl(a(\chi_{\P,i}^2)\bigr)=\chi_{\P,i}^{-2}(\alpha) \cdot
a(\chi_{\P,i}^2).$$Hence, the conclusion.

\label padictheta. section\par\ssection Theta operators on Katz
$p$-adic modular forms\par The notation is as
in~\refn{Katzpadicmodfor}. For every prime ideal~$\P$ of~$O_L$
over~$p$  and any $1\leq i\leq f_\P$ one can define the $R$-linear
operator
$$\Theta_{\P,i}\colon{\bf M}(R,\mu_N,\chi)^{p-{\rm adic}} \llongrightarrow
{\bf M}(R,\mu_N,\chi\chi_{\P,i}^2)^{p-{\rm adic}}$$as follows. Let
$\{ f_n \}_n$ be a compatible sequence of functions, $f_n$ a
regular function on~$\MM(n, n)$, defining a $\I$-polarized (Katz)
$p$-adic Hilbert modular form~$f$. Let~$\chi_n = \chi$
mod~$\X_R(n)$ be the weight of~$f_n$. Define $$\Theta_{\P,i}
(f_n)\in \Gamma\bigl(\MM(n,n),O_{\MM(n,n)}\bigr)$$as
in~\refn{theta}. By the compatibility of the Kodaira-Spencer
morphisms proven in~\refn{KS}, by~\refn{totalrecall}
and~\refn{invarianceSigmatilde}, it follows that $\{
\Theta_{\P,i}(f_n)\}_n$ is a compatible sequence defining a
$p$-adic modular form~$\Theta_{\P,i}(f)$ of
weight~$\chi\cdot\chi_{\P,i}^2$.

\noindent For any $0\leq \jj\leq e_\P-1$ let $$f\in{\bf
M}(R,\mu_N,\chi)^{{p-{\rm adic}},\P,[\jj]}$$be the $p$-adic
modular form associated to  the compatible sequence of functions
$\{f_n\}_n\in \Gamma(\MM(n,n),O_{\MM(n,n)}^{\P,[\jj]})$. One
defines the $R$-linear operator
$$\Theta_{\P,i}^{[\jj]}\colon{\bf M}(R,\mu_N,\chi)^{{p-{\rm adic}},\P,[\jj]}
\llongrightarrow {\bf M}(R,\mu_N,\chi\chi^2_{\P,i})^{{p-{\rm
adic}},\P,[\jj]}$$as follows. With the notation
of~\refn{sigmatilde}, we have a compatible sequence
$$\Theta_{\P,i}^{[\jj]} (f_n)\in
\Gamma\left(\MM(t_\P^{[\jj]}(n),n),O_{\MM(t_\P^{[\jj]}(n),n)}\right).$$Since
$t_\P^{[\jj]}(n)\rightarrow \infty$ this defines a $p$-adic
modular form by~\refn{padicsequence}.
\endssection

\label qexpTheta. section\par\ssection The behavior of
$\Theta_{\P,i}^{[\jj]}$ on $q$-expansions\par Fix a prime~$\P$
of\/~$O_L$ over~$p$ and integers $1\leq i\leq f_\P$ and $0\leq
\jj\leq e_\P-1$. To compute the effect of~$\Theta_{\P,i}^{[\jj]}$
on $q$-expansions we need to use some definitions and properties
of Tate objects; see~\refn{tateobjects}. Fix a
cusp~$\bigl(\B,\A,\varepsilon,\j\bigr)$ as in~\refn{cusp}. Fix a
rational polyhedral cone decomposition $\{\sigma_\beta\}_\beta$ of
the dual cone to $M_\RR^+ \subset M_\RR$, where $M:=\A\B$, as
in~\refn{tateobjects}. Then
$$R\bigl(\bigl(\A,\B,\sigma_\beta\bigr)\bigr):=R\tensor_\ZZ
\ZZ\bigl(\bigl(\A,\B,\sigma_\beta\bigr)\bigr)$$is a formally
smooth $R$-algebra of dimension~$g$, which can be interpreted as
$$\Spec\Bigl(R\bigl(\bigl(\A,\B,\sigma_\beta\bigr)\bigr)\Bigr)=
\Bigl(S_{\sigma_\beta}^\wedge\backslash
S_{\sigma_\beta,0}\Bigr)\fibprod_\ZZ R.$$Moreover,
$$\Omega^1_{R\bigl(\bigl(\A,\B,\sigma_\beta\bigr)\bigr)/R}=
\Big\langle{dq^\nu\over q^\nu}\Big\rangle_{\nu\in M}.$$The span is
taken as a $R\bigl(\bigl(\A,\B,\sigma_\beta\bigr)\bigr)$-module.
In particular, we conclude that
\endssection

\lemma There is a canonical isomorphism
$$\Omega^1_{R\bigl(\bigl(\A,\B,\sigma_\beta\bigr)\bigr)/R} \cong
R\bigl(\bigl(\A,\B,\sigma_\beta\bigr)\bigr)\tensor_\ZZ M.$$In
particular, the module of relative differentials
$\Omega^1_{R\bigl(\bigl(\A,\B,\sigma_\beta\bigr)\bigr)/R}$ is
endowed with the structure of free $
R\bigl(\bigl(\A,\B,\sigma_\beta\bigr)\bigr)\tensor_\ZZ O_L$-module
of rank\/~$1$.
\endlemma

\bigskip

\rmk By~\refn{indiff} the translation invariant relative
differentials
$$\omega_{\Tate(\A,\B)_{\sigma_\beta}/
R\bigl(\bigl(\A,\B,\sigma_\beta\bigr)\bigr)}$$ of the universal
object~$\Tate(\A,\B)_{\sigma_\beta}$
over~$R\bigl(\bigl(\A,\B,\sigma_\beta\bigr)\bigr)$ are canonically
isomorphic to~$
R\bigl(\bigl(\A,\B,\sigma_\beta\bigr)\bigr)\tensor_\ZZ \A$ as $
R\bigl(\bigl(\A,\B,\sigma_\beta\bigr)\bigr)\tensor_\ZZ O_L
$-module. The $O_L$-structure is of course given by the real
multiplication by\/~$O_L$ on~$\Tate(\A,\B)_{\sigma_\beta}$ over~$
R\bigl(\bigl(\A,\B,\sigma_\beta\bigr)\bigr)$.

\noindent In~\refn{TAU} we have fixed an element $l\in
\I=\Hom_{O_L}(\B,\A)$ inducing an $O_L$-linear isomorphism
$\ZZ_p\tensor_\ZZ\B \isomarrow \ZZ_p\tensor_\ZZ\A$. We get an
$O_L$-linear isomorphism $$\tau\colon \ZZ_p\tensor_\ZZ(\A\B)
\isomarrow \ZZ_p\tensor_\ZZ\A^2.$$
\endrmk

\label tau. proposition\par\prop  The Kodaira-Spencer map is
defined on~$R\bigl(\bigl(\A,\B,\sigma_\beta\bigr)\bigr)$.  The
$O_L$-linear isomorphism~$\tau$ makes the following diagram
commutative
$$\matrix{\Omega^1_{R\bigl(\bigl(\A,\B,\sigma_\beta\bigr)\bigr)/R}
&
\lllongmaprighto{\KS}&R\tensor_\ZZ\left(\omega_{\Tate(\A,\B)_{\sigma_\beta}/
T\bigl(\bigl(\A,\B,\sigma_\beta\bigr)\bigr)}^{\tensor_{O_L}^2}\right)
\cr
 \mapdownr{\wr} & &\mapdownr{\wr}\cr
 R\bigl(\bigl(\A,\B,\sigma_\beta\bigr)\bigr) \tensor_\ZZ M &
\lllongmaprighto{1 \tensor \tau}
 & R\bigl(\bigl(\A,\B,\sigma_\beta\bigr)\bigr)\tensor_\ZZ\A^{\tensor^2_{O_L}},\cr}$$
\endprop
\Proof It follows from~[\Katzzzz, \S1.1.18-\S1.1.20].

\rmk The isomorphisms~$\KS$ and~$\tau$ depend on the choice of the
element~$l\in\I^+$ in~\refn{TAU}. Different choices of~$l$ change
the isomorphisms by multiplication by a unit of\/~$R\tensor_\ZZ
O_L$. In~[\Katzzzz, \S1.0.21] one finds an expression for~$\KS$
which is independent of such choice.
\endrmk

\label qexptheta. corollary\par\cor The notation is as
in~\refn{theta}. Let $f\in
\Gamma\bigl(\MM(m,n),O_{\MM(m,n)}\bigr)$ be a regular function
on~$\MM(m,n)$. Fix a $\I$-polarized unramified
cusp~$\bigl(\A,\B,\varepsilon,\j_\varepsilon\bigr)$
of\/~$\MM(m,n)$. Suppose that the value of\/~$f$
in~$R/\m^m\bigl(\bigl(\A,\B,\sigma_\beta\bigr)\bigr)$ is
$f=\sum_\nu a_\nu q^\nu$. For any $\nu\in \P^\jj M$ define
$$\tilde{\chi}_{\P,i}^{[\jj]}(\nu)=
\sigmatilde_{\P,i}^{[\jj]} \left(\j_\varepsilon
\bigl(\tau(\nu)\bigr)\right) \hbox{{\rm mod
}}\m^{t_\P^{[\jj]}(m)},$$where $\sigmatilde_{\P,i}^{[\jj]}\colon
R\tensor_\ZZ\P^\jj \rightarrow R$ is the $R\tensor_\ZZ O_L$-linear
homomorphism defined in~\refn{sigmatilde}. Then

\spacing
\item{{\rm 1.}} $f\in
\Gamma\bigl(\MM(m,n),O_{\MM(m,n)}^{\P,[\jj]}\bigr)$ in the sense
of\/~\refn{theta} if and only if $\nu\in \P^\jj\, M$ for all~$\nu$
such that~$a_\nu\neq 0$;
\spacing
\item{{\rm 2.}} if $f$ is in
$\Gamma\bigl(\MM(m,n),O_{\MM(m,n)}^{\P,[\jj]}\bigr)$, the value
of\/~$\Theta_{\P,i}^{[\jj]}(f)$
in~$R/\m^{t_\P^{[\jj](m)}}\bigl(\bigl(\A,\B,\sigma_\beta\bigr)\bigr)$
is
$$\Theta_{\P,i}^{[\jj]}(f)=\sum_\nu \tilde{\chi}_{\P,i}^{[\jj]}(\nu)
a_\nu q^\nu.$$

\endcor
\Proof Note that $$\eqalign{d\left(a_0+\sum_{\nu \in M^+} a_\nu
q^\nu\right)& =\sum_{\nu \in M^+} a_\nu q^\nu
{d\bigl(q^\nu\bigr)\over q^\nu}\cr & =\sum_{\nu \in M^+} a_\nu
q^\nu \tensor \nu,\cr}$$where the last element is in~$
R\bigl(\bigl(\A,\B,\sigma_\beta\bigr)\bigr) \tensor_\ZZ M$.
By~\refn{tau} we have that~$f$ belongs to
$\Gamma\bigl(\MM(m,n),O_{\MM(m,n)}^{\P,[\jj]}\bigr)$ if and only
if $\sum_{\nu \in M^+} a_\nu q^\nu \tensor \nu $ is in $
R\bigl(\bigl(\A,\B,\sigma_\beta\bigr)\bigr) \tensor_\ZZ \P^\jj M$.
Since condition~$f\in
\Gamma\bigl(\MM(m,n),O_{\MM(m,n)}^{\P,[\jj]}\bigr)$ is a closed
condition, Part~(1) is proven.  By the definition
of~$\Theta_{\P,i}^{[\jj]}$ in~\refn{theta} and the fact that the
$q$-expansion of\/~$a\bigl(\chi_{\P,i}\bigr)$ at any unramified
cusp is~$1$, we have that the $q$-expansion
of~$\Theta_{\P,i}^{[\jj]}(f)$ at the given cusp is
$$\eqalign{\Theta_{\P,i}^{[\jj]}(f)\bigl(\Tate(\A,\B),\varepsilon,\j_\varepsilon)&
= \sigmatilde_{\P,i}^{[\jj]}\left(\j_\varepsilon
\Bigl(\KS\bigl(d(a_0+\sum_{\nu \in M^+} a_\nu
q^\nu)\bigr)\Bigr)\right)\cr & =\sum_{\nu \in M^+} a_\nu q^\nu
\cdot
\sigmatilde_{\P,i}^{[\jj]}\Bigl(\j_\varepsilon\bigl(\tau(\nu)\bigr)\Bigr)\cr}$$as
claimed.

\label thetacommute. corollary\par\cor The operators
$\Theta_{\P,i}^{[\jj]}$ commute for different primes~$\P$ and
different~$1\leq i\leq f_\P$.\endcor

\label qexppadictheta. corollary\par\cor The notation is as
in~\refn{padictheta} and in~\refn{qexptheta}. Let $f$ be a
$\I$-polarized $p$-adic Hilbert modular form over~$R$ of
level\/~$\mu_N$ and weight\/~$\chi$; see~\refn{Katzpadicmodfor}.
Suppose that~$f$ has $q$-expansion at the $\I$-polarized $p$-adic
cusp $(\A,\B,\varepsilon_{p^\infty N },\j_\varepsilon\bigr)$ equal
to $a_0 + \sum_{\nu\in M^+} a_\nu q^\nu$; see~\refn{qexpKatz} for
the notation. Here $M=\A\B$. Then

\spacing
\item{{\rm 1.}} $f\in{\bf M}(R,\mu_N,\chi)^{{p-{\rm adic}},\P,[\jj]}$  if
and only if $\nu\in \P^\jj\, M$ for all $\nu\in M$;
\spacing
\item{{\rm 2.}} if $f\in{\bf M}(R,\mu_N,\chi)^{{p-{\rm adic}},\P,[\jj]}$, the $q$-expansion
of\/~$\Theta_{\P,i}^{[\jj]}(f)$ at the same cusp is
$$\Theta_{\P,i}^{[\jj]}(f)\bigl(\Tate(\A,\B),\varepsilon,\j_\varepsilon\bigr)=
\sum_\nu \tilde{\chi}_{\P,i}^{[\jj]}(\nu) a_\nu q^\nu.$$
\endcor
\Proof It follows from the definition of~$\Theta_{\P,i}^{[\jj]}$
given in~\refn{padictheta} and the previous corollary.

\label padicthetacommute. corollary\par\cor The $p$-adic theta
operators  $\Theta_{\P,i}^{[\jj]}$ commute for different
primes~$\P$ and different~$1\leq i\leq f_\P$.\endcor

\ssection The comparison with the complex theory\par We use the
notation of~\refn{complex}. Given a class in the strict class
group of~$L$, we may choose a representative~$\I$ which is
fractional ideal prime to~$p$. Take~$\B=O_L$ and~$\A=\I$. Let~$f$
be a $\I$-polarized Hilbert modular forms in~${\bf
M}(\bar{\QQ},\mu_N,\chi)$. Choose embeddings $\bar{\QQ}\subset\CC$
and $\bar{\QQ}\subset\bar{\QQ}_p $ and view~$f$ as a complex or
$p$-adic Hilbert modular form. Assume that the $q$-expansion
of~$f$ at the cusp~$(i\infty,\ldots,i\infty)$ is
$a_0+\sum_{\nu\in\A^+}a_\nu q^\nu$. Then,
$f\in\bar{\QQ},\mu_N,\chi)^{\P,[\jj]}$ if and only if~$a_\nu=0$
for all~$\nu\not\in(\P^\j\A)^+$ and
$$\Theta_{\P,i}^{[\jj]}(f)\bigl(\Tate(\A,\B),\varepsilon,\j_\varepsilon\bigr)=
\sum_{\nu\in \A^+}\chi_{\P,i}\bigl(\nu\pi_\P^{-\jj}\bigr) a_\nu
q^\nu;$$where $\pi_\P$ is the chosen uniformizer of~$O_L$ at~$p$.
\endssection

\label theyareasinKatz. section\par\ssection Katz's $p$-adic theta
operators\par In~[\Katzzzz, \S2.6] one finds a definition of
$p$-adic theta operators on $p$-adic Hilbert modular forms \`a la
Katz (see~\refn{Katzpadicmodfor}) if~$p$ {\it is unramified}
in~$L$. In this case it follows from~[\Katzzzz, Cor.~2.6.25]
and~\refn{qexppadictheta} that Katz's theta operators coincide on
$q$-expansions with the $p$-adic theta operators defined above.
\endssection

\label padicnonclassic. section\par\ssection Other examples of
$p$-adic modular forms\par One way to produce examples of $p$-adic
modular forms (see~\refn{Katzpadicmodfor}) is by applying the
$p$-adic theta operators to classical Hilbert modular forms. In
general, the image of a classical modular form, i.~e., an element
of~$M(F, \mu_N, \chi)$ (see~\refn{modularforms}), under a theta
operator is not a classical modular form. To illustrate that we
consider the case of~$g=2$ and the $\I$-polarized classical
Eisenstein series~$\E_2$ over~$\QQ$ of weight~$2$. Recall
from~\refn{Eisenstein} that its $q$-expansion at a $\I$-polarized
unramified $\QQ$-rational cusp $(\A,\B,\j_{\rm can})$ is
$$\Norm^{k-1}(\A) \biggl( 2^{-g} \zeta_L(1-k) +
\sum_{\nu\in (\A\B)^+} \Bigl(\sum_{\nu\in \C\subset \A\B}
\Norm(\nu\C^{-1})^{k-1}\Bigr) q^{\nu}\biggr).$$Let~$p$ be a split
prime and let~$\chi$, $\chi'$ be the fundamental characters
over~$p$. Then~$\Theta_{\chi}E_2$ is a $p$-adic modular form,
whose $q$-expansion is
$$\Norm^{k-1}(\A) \biggl( 2^{-g} \zeta_L(1-k) +
\sum_{\nu\in (\A\B)^+} \Bigl(\sum_{\nu\in \C\subset \A\B}
\Norm(\nu\C^{-1})^{k-1}\Bigr)\tilde{\chi}(\nu) q^{\nu}\biggr).
,$$which is not a classical modular form. To see that we remark
that there exists a pull-back map from classical modular forms for
our quadratic field to modular forms on~$\QQ$, taking a modular
form of level~$1$ (say) and weight~$\chi_1^{a_1}\dots
\chi_g^{a_g}$ to a modular form on~$\SL_2(\ZZ)$ of weight~$a_1 +
\dots + a_g$. See~\S~19. The pull-back of~$\E_2$ is a multiple of
the Eisenstein series~$\E^{\QQ}_4$ on~$\SL_2(\ZZ)$, where to avoid
confusion, we write~$\E_k^{\QQ}$ for Eisenstein series of
weight~$k$ for the field~$\QQ$.

\noindent  If~$\Theta_{\chi}E_2$ is classical of some weight and
level, so is its Galois conjugate~$\Theta_{\chi'}E_2$ over~$L$.
Hence, the sum~$\Theta_{\chi}E_2+\Theta_{\chi'}E_2$ is in the
graded ring of classical modular forms. It follows
from~\refn{Hilbertvsellitpic} that the pull-back
of~$\Theta_{\chi}E_2+\Theta_{\chi'}E_2$ is proportional
to~$\Theta\E^{\QQ}_4$ and is a $p$-adic cusp form on~$\SL_2(\QQ)$
of weight~$6$ and integral $q$-expansion. But then, reducing
modulo~$p$,  we would have a mod~$p$ cusp form of weight~$4+p+1$,
hence divisible, as a holomorphic modular form, by~$\Delta$.
Take~$p=2,3$ or~$5$ to obtain a contradiction.

\endssection

\bigskip

\noindent From now on we work in characteristic~$p$. The goal of
the rest of this section is to define theta operators on Hilbert
modular forms in characteristic~$p$. See the introduction of the
section for a more detailed discussion.

\ssection The poles of\/~$\Theta_{\P,i}$ in char~$p$\par Denote by
$$\MM:=\MM(1,1)^\Kum \qquad\hbox{{\rm and}}\qquad
\MM':=\MM\bigl(1,0\bigr).$$Then $\phi\colon \MM\rightarrow \MM'$
is a Galois cover with group
$$G\cong \prod_{\P\vert (p)}\bigl(O_L/\P\bigr)^*.$$
By~\refn{phibar} and~\refn{Mbar} we can complete~$\phi$ to a
morphism of projective normal schemes:
$$\matrix{
 \MM & \hooklongrightarrow & \MMbar \cr
 \mapdownr{\phi} & &\mapdownr{\phibar}\cr
 \MM' & \hooklongrightarrow & \MMbar'.\cr}$$Note that
$\MMbar'$ is the minimal compactification of\/~$\MM(k,\mu_N)$ and
the map~$\phibar$ is ramified along the complement of the ordinary
locus of\/~$\MMbar'$. The Kodaira-Spencer isomorphisms extends to
the boundary as follows. Over the open subscheme~$\MM(k,\mu_N)$
of\/~$\MMbar'$ one has a universal abelian scheme denoted, by
abuse of notation, by~$\pi\colon\AA'\rightarrow \MM(k,\mu_N)$.
Proceeding as in~\refn{KS} one obtains an isomorphism
$$ \KSbar'\colon
\Omega^1_{\MM(k,\mu_N)/k} \isomarrow
\omega_{\AA'/\MM(k,\mu_N)}^{\tensor_{O_L}^2}.$$Then we have
$$\KSbar:=\phibar^*\Bigl(\KSbar'\Bigr)\colon \phibar^*\Bigl(\Omega^1_{\MM(k,\mu_N)/k}
\Bigr) \isomarrow \phibar^*\Bigl(
\omega_{\AA'/\MM(k,\mu_N)}^{\tensor_{O_L}^2}\Bigr)$$extending~$\KS$.
Let
$$\MMbar^\R:=\phibar^{-1}\Bigl(\MMbar\bigl(k,\mu_N\bigr)^{\R}
\Bigr)$$the inverse image of the Rapoport locus defined
in~\refn{Rapo} including the cusps. Its complement in~$\MMbar$ has
{\it codimension at least}~$2$.

\endssection

\label WPi. definition\par\defi For every prime~$\P$ over~$p$ and
any integer~$1\leq i\leq f_\P$, define the Weil divisor
of\/~$\MMbar$:
$$W_{\P,i}:=\hbox{{\rm support of the effective divisor }} \phibar^{-1}\bigl(h_{\P,i}
\bigr).$$
\enddefi

\rmk By construction $W_{\P,i}$ is reduced. Moreover $$\Bigl(\phibar^{-1}\bigl(h_{\P,i}
\bigr)\Bigr)=\phibar^{-1}\Bigl(\bigl(h_{\P,i} \bigr)\Bigr)=\bigl(
p^{f_\P}-1\bigr) W_{\P,i}.$$
\endrmk

\noindent If $f$ is a regular function of~$\MM$ lying
in~$\Gamma\bigl(\MM,O_\MM^{\P,[\jj]}\bigr)$, we are interested in
computing the poles of $\Theta_{\P,i}^{[\jj]}(f)$ on~$\MMbar$.
This is achieved as follows:

\label polef. section\par \ssection Some notation\par Let~$\Q$
and~$\P$ be prime ideals of\/~$O_L$ over~$p$ and let $1\leq j\leq
f_\Q$ and $1\leq i\leq f_\P$. Define
$$\bigl(\Omega^1_{\MM/k}\bigr)_{\P,i}:=e_{\P,i}\cdot
\Omega^1_{\MM/k}$$and if $f\in O_{\MM}$ define
$$ df_{\P,i}:=e_{\P,i}\cdot df \in
\bigl(\Omega^1_{\MM/k}\bigr)_{\P,i}.$$See~\refn{omegaPi} for the
definition of the idempotent~$e_{\P,i}$.  Let\/~$\delta$ be a
local uniformizer of an irreducible component of the
divisor~$W_{\Q,j}$. Let $v_\delta$ be the discrete valuation on
the meromorphic functions on~$\MM$ defined by~$\delta$. Express
$$f:={u \over \delta^n},$$where $u$ is a function such
that~$v_\delta(u)=0$. Then
$$\Bigl(\Omega^1_{\MM/k}\Bigr)_{\P,i} \ni
\bigl(df\bigr)_{\P,i}={\bigl(du\bigr)_{\P,i} \over \delta^n}-{n
u\bigl(d\delta \bigr)_{\P,i}\over \delta^{n+1}}.$$
\endssection

\label poledelta. section\par \ssection The poles
of\/~$d\delta$\par To compute $(d\delta)_{\P,i}$ we need to use
the detailed analysis of the deformation theory using displays
defined in~\refn{unidefRM}. By~\refn{Mbar}-\refn{casino} we may
choose~$\delta$ such that
$$\delta^{p^{f_\Q}-1}=h_{\Q,j+1}^{p^{f_\Q-1}}\cdots h_{\Q,j-1}^p
h_{\Q,j},$$where the elements $h_{\Q,j}$ are now viewed not as
sections of line bundles but as functions via the choice of local
generators for the sheaves themselves.  Since~$k$ is of
characteristic~$p$, we conclude that $$\bigl(p^{f_\Q}-1 \bigr)
\delta^{p^{f_\Q}-2} d\delta=\Bigl(h_{\Q,j+1}^{p^{f_\Q-1}}\cdots
h_{\Q,j-1}^p \Bigr) dh_{\Q,j}.$$In~\refn{unidefRM} we have
computed the completion
$$R_\iota=k[\![t_{\P,i}^{[j]}]\!]_{\P,i,j}$$of~$\MM\bigl(k,\mu_N\bigr)^\R$ at a
geometric point. With the notation of~\refn{HASSE}, it follows
from~\refn{hvst} that
$$dh_{\Q,j}=\bar{c}_{\Q,j}^{[1]} dt_{\Q,j}^{[1]}$$for an invertible element~$\bar{c}_{\Q,j}^{[1]}$
of\/~$k$.
\endssection

\lemma We have: $$d\left(t_{\Q,j}^{[1]}\right) \in
\Bigl(\Omega^1_{\MM(k,\mu_N)^\R/k}\Bigr)_{\Q,j}.$$
\endlemma
\Proof The proof relies on understanding the connection between
two deformation theories of abelian varieties: one based on
Grothendieck and Mumford's infinitesimal theory, the other based
on the theory of displays.  Consider the display $\bigl({{\rm
P}}_0,{{\rm Q}}_0,\F_0,\V_0^{-1}\bigr)$ over~$k$ introduced
in~\refn{needtofixdisplay}--\refn{Not}.  The following results
follow from~[\Zink].

\spacing
\indent {\rm a.}  \enspace Let $\bigl({{\rm P}}_1,{{\rm
Q}}_1,\F_1,\V_1^{-1}\bigr)$ and~$\bigl({{\rm P}}_2,{{\rm
Q}}_2,\F_2,\V_2^{-1}\bigr)$ two displays over an artinian
$k$-algebra $D$ deforming~$\bigl({{\rm P}}_0,{{\rm
Q}}_0,\F_0,\V_0^{-1}\bigr)$. For $i=1,2$ let~$\widehat{{{\rm
Q}}}_i$ be the inverse image of~${{\rm Q}}_i$ via ${{\rm
P}}\rightarrow {{\rm P}}_0$. Then there exists a unique
isomorphism
$$\alpha_{{{\rm P}}_1,{{\rm P}}_2}\colon \bigl({{\rm P}}_1,\widehat{{{\rm Q}}}_1,\F_1,\V_1^{-1}\bigr)
\isomarrow \bigl({{\rm P}}_2,\widehat{{{\rm
Q}}}_2,\F_1,\V_2^{-1}\bigr)$$ inducing the identity
on~$\bigl({{\rm P}}_0,{{\rm Q}}_0,\F_0,\V_0^{-1}\bigr)$.
\spacing

\noindent In particular, the functor associating to a nilpotent
pd-thickening $k \subset D$ the $D$-module $\D_{{{\rm
P}}_0}(D):=\bar{{{\rm P}}}:={{\rm P}}/I_D{{\rm P}}$, where $({{\rm
P}},{{\rm Q}},\F,\V^{-1})$ is a display over~$D$ deforming~$({{\rm
P}}_0,{{\rm Q}}_0,\F_0,\V_0^{-1}\bigr)$, defines a crystal
over~$k$.

\indent {\rm b.} \enspace Let~$\DD_{A_0[p^\infty]}$ be the crystal
over~$k$ associated to the formal group~$A_0[\infty]$ by
Grothendieck and Messing [\Messing]. The morphism $\D_{{{\rm
P}}_0}\rightarrow \DD_{A_0[p^\infty]}$, which associates to a
nilpotent pd-thickening $k\subset D$ and a display $({{\rm
P}},{{\rm Q}},\F,\V^{-1})$ over~$D$ deforming~$({{\rm P}}_0,{{\rm
Q}}_0,\F_0,\V_0^{-1}\bigr)$ the Lie algebra of the universal
extension over~$D$ of the formal group associated to~${{\rm P}}$,
defines an isomorphism of crystals.
\spacing

\noindent  Let $D$ be an artinian $k$-algebra  with a section~$s$.
Let~$\bigl({{\rm P}}_{\rm can}, {{\rm Q}}_{\rm can},\F_{\rm
can},\V_{\rm can}^{-1}\bigr)$ be the trivial deformation
of~$\bigl({{\rm P}}_0,{{\rm Q}}_0,\F_0,\V_0^{-1}\bigr)$ to~$D$
defined by pulling-back via~$s$. Let~$\Def_{{{\rm
P}}_0}\bigl(D\bigr)$ be the category of deformations of the
display~$\bigl({{\rm P}}_0,{{\rm Q}}_0,\F_0,\V_0^{-1}\bigr)$ to a
display over~$D$. Let~$\Def_{{\bar{\rm Q}} \subset {\bar{\rm
P}}_{\rm can} }\bigl(D\bigr)$ be the category of liftings of the
Hodge filtration~$\bar{{{\rm Q}}}_0 \subset \bar{{{\rm P}}}_0 $
in~${\bar{\rm P}}_{\rm can}$.

\indent {\rm c.} \enspace  The functor
$$\Def_{{{\rm P}}_0}\bigl(D\bigr) \llongrightarrow
\Def_{{\bar{\rm Q}} \subset {\bar{\rm P}}_{\rm can}
}\bigl(D\bigr),$$associating to~$\bigl({{\rm P}},{{\rm
Q}},\F,\V^{-1}\bigr)\in \Def_{{{\rm P}}_0}\bigl(D\bigr)$ the flag
$\alpha_{{{\rm P}}_{\rm can},{{\rm P}}}(\bar{{{\rm Q}}}) \subset
\bar{{{\rm P}}}_{\rm can}$, defines an equivalence of categories.

\spacing
\noindent Let~$D:=k[\![t_{a,b}]\!]_{1\leq a,b\leq g}$ with the
relations~$t_{a,b} t_{a',b'}=0$. Let~$T$ (resp.~$\bar{T}$) be the
matrix of Teichm\"{u}ller lifts~$\bigl(w(t_{a,b})\bigr)_{1\leq
a,b\leq g}$ (resp.~the matrix~$\bigl(t_{a,b}\bigr)_{1\leq a,b\leq
g}$).

\indent {\rm d.} \enspace The flag $\bar{{{\rm Q}}}_{\rm can}+
\bar{T} \bar{\T}_{\rm can} \subset \bar{{{\rm P}}}_{\rm can}$
comes from the  display~${{\rm P}}:={{\rm P}}_{\rm can}$
and~${{\rm Q}}:={{\rm Q}}_{\rm can}$ via the map
$$\alpha_{{{\rm P}}_{\rm can},{{\rm P}}}:=\left(\matrix{ {\rm Id}_g
& T \cr 0 & {\rm Id}_g \cr}\right);$$the matrix is given on~${{\rm
P}}_{\rm can}=\T_{\rm can} \dirsum\L_{\rm can}$ with respect
to~${\cal B}$. In particular, the matrix of\/~$\F\dirsum \V^{-1}$
with respect to the same basis is given by
$$\left({\matrix{ A+ T C & B+ T D \cr C & D
\cr}}\right).$$Imposing the condition that~$\bigl({{\rm P}},{{\rm
Q}},\F,\V^{-1}\bigr)$ is polarized is equivalent to restrict to
the quotient of~$D$ defined by the relations~$t_{a,b}=t_{b,a}$.
Let~$R_\iota$ be the maximal $R$-algebra over which the
$O_L$-action extends; see~\refn{unidefRM}. Denote
by~$\m_\iota=(t_{a,b})$ its maximal ideal.  Define~$\bigl({{\rm
P}}, {{\rm Q}},\F,\V^{-1}\bigr)$ to be the associated display
over~$R_\iota$.

\spacing \noindent Let $A \rightarrow \Spec(R_\iota)$ be the
abelian scheme with real multiplication by~$O_L$ associated to the
display~$\bigl({{\rm P}}, {{\rm Q}},\F,\V^{-1}\bigr)$. The
isomorphism~$\alpha_{{{\rm P}}_{\rm can},{{\rm P}}}$ induces an
$O_L$-linear isomorphism

\advance\ssnu by-1\labelf GaussManin\par
$$\H_{1,{\rm dR}}(A_0/k)\fibprod_k \Spec(R_\iota)=\bar{P}_{\rm can} \isomarrow \bar{P}=\H_{1,{\rm
dR}}(A/R_\iota)\eqno{(\numfo)}.$$\advance\ssnu by1\advance\fonu
by1

\noindent By~(b) it is compatible with the isomorphism
$$\DD_{A_0[p^\infty]}(k) \fibprod_k \Spec(R_\iota)\isomarrow
\DD_{A_0[p^\infty]}(R_\iota).$$By the comparison theorem between
the crystals~$\DD_{A_0[p^\infty]}$ and~$R^1\pi_{{\rm
crys},*}\bigl(O_{A_0,{\rm crys}}\bigr)$, we get that the
isomorphism~(\refn{GaussManin}) identifies~$\H_{1,{\rm
dR}}(A_0/k)$ with the horizontal sections of the Gauss-Manin
connection on~$\H_{1,{\rm dR}}(A/R_\iota)$.

\spacing
\noindent Via the identifications of~\refn{HASSE}, we get that the
$O_L$-linear map deduced from~(\refn{GaussManin})
$$\Hom\Bigl(\H^1\bigl(A_0,O_{A_0}\bigr),k\Bigr)=\bar{Q}_0 \llongrightarrow
\m_\iota\cdot \bigl( \bar{P}/\bar{Q}\bigr)
=\Hom\Bigl(\H^0\bigl(A_0,\Omega^1_{A_0/k}\bigr),\m_\iota\Bigr)$$is
defined by the matrix~$\bar{T}$. By~[\Katzzzz, \S1.0.11--\S1.0.21]
the induced $k$-linear map
$$\H^0\bigl(A_0,\Omega^1_{A_0/k}\bigr)\tensor_{O_L} \H^0\bigl(A_0,\Omega^1_{A_0/k}\bigr)
\isomarrow \m_\iota=\Omega^1_{R_\iota/k}$$coincides with the
Kodaira-Spencer map. In particular, for any prime~$\Q$, any $1\leq
j\leq f_\Q$ and any $1\leq l\leq e_\Q$, such map is uniquely
defined by the formula
$$\bigl(\omega_{\Q,j}^{[1]} \bigr) \tensor
\bigl(\omega_{\Q,j}^{[l]} \bigr)\llongmapsto d\left(
t_{\Q,j}^{[l]} \right).$$This concludes the proof of the lemma.

\bigskip
\noindent Gathering the results of\/~\refn{polef}
and\/~\refn{poledelta} and using the results of the lemma, we
conclude that for the particular choice of~$\delta$ made
in~\refn{poledelta}
$$v_\delta\Bigl(\bigr(d\delta\bigl)_{\P,i}\Bigr)=\cases{
 \, \infty & if $(\P,i)\neq (\Q,j)$; \cr
 \, 2-p^{f_\P} & if $(\P,i)=(\Q,j)$.\cr}$$Hence,

\prop With the notation of\/~\refn{polef}, we have
$$v_\delta\bigl((df)_{\P,i}\bigr)\cases{
 \,\geq v_\delta\bigl(f\bigr) & if $p\vert v_\delta(f)$ or $(\P,i) \neq
(\Q,j)$;\cr
 \,= v_\delta\bigl(f\bigr)-\bigl( p^{f_\P}-1\bigr) &
 otherwise.\cr}$$
\endprop

\label polethetaf. proposition\par\prop The notation is as
in~\refn{theta}. Fix primes~$\P$ and~$\Q$ over~$p$ and integers
$1\leq i\leq f_\P$ and $1\leq j\leq f_\Q$. Let $0 \leq\jj \leq
e_\P-1$. Let $f$ be a regular function on~$\MM$ such that $f\in
\Gamma(\MM,O_\MM^{\P,[\jj]})$. Let\/~$C$ be an irreducible
component of the divisor~$W_{\Q,j}$ defined in~\refn{WPi}.
Let\/~$v_C$ be the corresponding valuation. Then
$$v_C\left(\Theta_{\P,i}^{[\jj]}\bigl(f\bigr)\right)\cases{
 \,\geq  v_C\bigl(f\bigr) & if $\P\neq\Q$;\cr
 \,\geq v_C\bigl(f\bigr)-2 p^{f_\P-r}  & if $\P=\Q$ and  $i\neq j$;\cr
 \,\geq v_C\bigl(f\bigr)-2  & if $\P=\Q$, $i=j$ and
 $p\vert v_C(f)$;\cr
 \,=v_C\bigl(f\bigr)-2-\bigl( p^{f_\P}-1\bigr) & if $\P=\Q$,
 $i=j$ and $p\,\not\vert \,v_C(f)$;
\cr}$$where~$r=j-i$ if~$j>i$ and~$r=f_\P+j-i$ if~$j<i$.
\endprop
\Proof By~\refn{theta} we have
$$\Theta_{\P,i}^{[\jj]}(f)=\sigmatilde_{\P,i}^{[\jj]}\Bigl(\KS\bigl(df\bigr)\Bigr)\cdot
a\bigl(\chi_{\P,i}^2\bigr)^{-1}$$on~$\MM$. The
map~$\sigmatilde_{\P,i}^{[\jj]}$ defined in~\refn{Sigmatilde}
over~$\MM(k,\mu_N)^{\rm ord}$ extends to an $O_L\tensor_\ZZ
k$-linear map
$$\sigmatilde_{\P,i}^{[\jj]}\colon
\omega_{\AA/\MM\bigl(k,\mu_N\bigr)^\R}^{\tensor_{O_L}^2}\tensor_{O_L}\P^\jj
\llongrightarrow {\cal
L}_{\chi_{\P,i}^2}$$over~$\MM\bigl(k,\mu_N\bigr)^\R$. Hence we
obtain an $O_L\tensor_\ZZ k$-linear map
$$\Omega^1_{\MM(k,\mu_N)^\R/k} \llongmaprighto{\KS}
\omega_{\AA/\MM\bigl(k,\mu_N\bigr)^\R}^{\tensor_{O_L}^2}\tensor_{O_L}\P^\jj
\lllongmaprighto{\sigmatilde_{\P,i}^{[\jj]}} {\cal
L}_{\chi_{\P,i}^2}.$$Pulling-back we get a map
$$\phibar^*\bigl(\Omega^1_{\MM(k,\mu_N)^\R/k}\tensor_{O_L}
\P^\jj\bigr)\llongrightarrow \phibar^*\bigl({\cal
L}_{\chi_{\P,i}^2}\bigr).$$We also have a natural inclusion of
sheaves of $O_{\MM^\R}$-modules
$$0 \llongrightarrow \phibar^*\bigl(\Omega^1_{\MM(k,\mu_N)^\R/k}
\bigr) \llongrightarrow \Omega^1_{ \MM^\R/k}.$$The cokernel is the
branch locus of~$\phibar$; see~\refn{branch}. Note that $df\in
\Omega^1_{ \MM^\R/k}$! The modular form
$a\bigl(\chi_{\P,i}^2\bigr)$ defined in~\refn{actiona(chi)}
extends by~\refn{Mbar}  to a section of\/~$\phibar^*\bigl({\cal
L}_{\chi_{\P,i}^2}\bigr)$ over~$\MM^\RR$, which we denote in the
same way. It is non-vanishing on~$\MM$ and locally
on~$\MM^\R\backslash\MM$ it satisfies
$$a\bigl(\chi_{\P,i}\bigr)^{p^{f_\P}-1}=h_{\P,i+1}^{p^{f_\P-1}}\,
h_{\P,i+2}^{p^{f_\P-2}}\cdots h_{\P,i}.$$Hence, $v_C\bigl(
a(\chi_{\P,i}^2)\bigr)$ is equal to~$0$ if~$\P\neq\Q$, it is equal
to~$2 p^{f_\P-r}$ if~$\P=\Q$ and~$i\neq j$, and it is equal
to~$2$, if~$\P=\Q$ and~$i=j$. We conclude using the previous
proposition.

\label ThetaPi. definition \par \ssection {\bf
Definition.}\enspace (The operator $\Theta_{\P,i}^{[\jj]}$ on
modular forms)\par  Let $\P$ be a prime over~$p$, let $1\leq i\leq
f_\P$ and let $ 0\leq\jj\leq e_\P-1$. Define the subspace of
$\I$-polarized modular forms of weight\/~$\psi$ $${\bf
M}\bigl(k,\mu_N,\psi\bigr)^{\P,[\jj]}:=\left\{f\in {\bf
M}\bigl(k,\mu_N,\psi\bigr)\vert r(f)\in\Gamma\bigl(\MM,
O_\MM^{\P,[\jj]}\bigr)\right\},$$where $r(f)$ is the regular
function on~$\MM$ defined in~\refn{r(f)} and~$ \Gamma\bigl(\MM,
O_\MM^{\P,[\jj]}\bigr)$ is defined in~\refn{theta}. For $f\in{\bf
M}\bigl(k,\mu_N,\psi\bigr)_\P^{[\jj]}$ define
$$\Theta_{\P,i}^{[\jj]}(f):=\Theta_{\P,i}^{[\jj]}\bigl(r(f)\bigr)\cdot a\bigl(\psi\bigr)
\cdot a\bigl(\chi_{\P,i}^2\bigr) \cdot h_{\P,i}.$$We have
$$\Theta_{\P,i}^{[\jj]}(f)\in\Gamma\Bigl(\MM, {\cal L}_{\psi
\chi_{\P,i-1}^p\chi_{\P,i}}\Bigr).$$
\enddefi

\label ThetaPihasnopoles. theorem\par\thm The notation is as
in~\refn{ThetaPi}.
\spacing
\item{{\rm 1.}} The section~$\Theta_{\P,i}^{[\jj]}(f)$
over~$\MM$ descends to a section of the line bundle ${\cal
L}_{\psi \chi_{\P,i-1}^p\chi_{\P,i}} $
over~$\MM\bigl(k,\mu_N\bigr)^{\rm ord}$;
\spacing
\item{{\rm 2.}}  the section $\Theta_{\P,i}^{[\jj]}(f)$
extends to a section of\/ ${\cal L}_{\psi
\chi_{\P,i-1}^p\chi_{\P,i}} $ over~$\MM(k,\mu_N)$.

\spacing
\noindent Hence, we obtain a $k$-{\it derivation}:
$$\Theta_{\P,i}^{[\jj]}\colon \dirsum_\psi{\bf M}\bigl(k,\mu_N,\psi\bigr)^{\P,[\jj]}
\llongrightarrow \dirsum_\psi{\bf
M}\bigl(k,\mu_N,\psi\chi_{\P,i-1}^p\chi_{\P,i} \bigr).$$
\endthm
\Proof Part (1) follows from the description of action of the
Galois group~$G$ on functions and on modular forms given
in~\refn{actiona(chi)}. Part~(2) follows from~\refn{polethetaf}.

\label qexpThetaPi. corollary\par\cor The notation is as
in~\refn{ThetaPi} and in~\refn{qexpTheta}.  Let\/~$f\in{\bf
M}\bigl(k,\mu_N,\psi\bigr)$ be a $\I$-polarized modular form of
level\/~$N$ and weight~$\psi$. Suppose that the $q$-expansion
of\/$f$ at a $\I$-polarized unramified
cusp\/~$\bigl(\B,\A,\varepsilon,\j \bigr)$ is
$$f\bigl(\Tate(\A,\B),\varepsilon,\j\bigr)=\sum_\nu a_\nu
q^\nu.$$Then

\spacing
\item{{\rm 1.}} $f\in{\bf M}\bigl(k,\mu_N,\psi\bigr)^{\P,[\jj]}$  if
and only if $a_\nu=0$ for all $\nu \notin \P^\jj M$;

\spacing
\item{{\rm 2.}} if $f\in{\bf M}\bigl(k,\mu_N,\psi\bigr)^{\P,[\jj]}$, the $q$-expansion
of\/~$\Theta_{\P,i}^{[\jj]}(f)$ at the same cusp is
$$\Theta_{\P,i}^{[\jj]}(f)\bigl(\Tate(\A,\B),\varepsilon,\j\bigr)=
\sum_\nu \tilde{\chi}_{\P,i}^{[\jj]}(\nu) a_\nu
q^\nu.$$See~\refn{qexptheta} for the definition
of\/~$\tilde{\chi}_{\P,i}^{[\jj]}$.
\endcor
\Proof It follows from the definition of~$\Theta_{\P,i}$ given
in~\refn{ThetaPi} and~\refn{qexptheta}.

\label ThetaPicommute. corollary\par\cor The  theta operators
$\Theta_{\P,i}^{[\jj]}$ commute for different primes~$\P$ and
different~$1\leq i\leq f_\P$.\endcor

\ssection A comparison with Katz's definition\par Katz's
definition of theta operators extends to Hilbert modular form in
characteristic~$p$; see~[\Katzzzz, \S2.6]. Note however that the
our operators $\Theta_{\P,i}$ in characteristic~$p$  and Katz's
theta operators change the weights in a different way;
compare~\refn{ThetaPi} with~[\Katzzzz, Cor.~2.6.25]. Indeed, our
theta operators in characteristic~$p$ are Katz's operators
multiplied by suitable partial Hasse invariants in order to ensure
that they send holomorphic modular forms to holomorphic modular
forms. Instead, Katz is interested only in modular forms defined
over the {\it ordinary locus} where the holomorphicity is not an
issue.
\endssection

\endsection

\section The operator $V$\par \noindent In this section we suppose that
$k\subset \bar{{\bf F}}_p$.

\label V. section\par \ssection The definition\par Let $g\in {\bf
M}\bigl(k,\mu_N,\psi\bigr)$ be a $\I$-polarized modular form
over~$k$ of level\/~$\mu_N$ and of weight~$\psi$. Let $R$ be a
$k$-algebra. Let
$$\underline{A}:=\bigl(A,\iota,\lambda,\varepsilon\bigr)\qquad\hbox{{\rm
and}}\qquad \omega\in {\rm H}^0\bigl(R,\Omega^1_{A/R}\bigr)$$be a
$\I$-polarized abelian scheme over~$R$ with $O_L$-action
and\/~$\mu_N$-level structure, as defined in~\refn{moduli}, and a
generator~$\omega$ of\/~${\rm H}^0\bigl(R,\Omega^1_{A/R}\bigr)$ as
a free $O_L\tensor_\ZZ R$-module. Define
$$V\bigl(g\bigr)\bigl(\underline{A},\omega\bigr):=g\bigl(\underline{A}^{(p)},\omega^{(p)}
\bigr),$$where the superscript~$(p)$ stands for the base change by
the {\it absolute Frobenius}.

\endssection

\label Frobchar. definition\par \defi The notation is as
in~\refn{grschG}. Let\/~$\F^{\rm abs}$ be the absolute Frobenius
morphism on~$\G_\FF$.  Let $\F^{\rm abs}_k:=\F^{\rm
abs}\fibprod_{\Spec(\FF)} \Spec(k)$. Define the following
endomorphism of the group of characters~$\X_k$ by
$$\Bigl(\chi\colon \G_k \rightarrow \GG_{m,k}\Bigr) \llongmapsto
\Bigl( \chi^{(p)}:=\chi\circ\F^{\rm abs}_k\Bigr).$$It  preserves
the universal characters defined in~\refn{GGm}.
\enddefi

\label sigma*psi. remark\par\rmk Suppose that $k=\bar{{\bf F}}_p$.
Then the group of characters~$\X_k$ of~$\G_k$ is endowed with a
${\rm Gal}(k/\FF)$-action. In particular, the absolute
Frobenius~$\sigma$ defines an action on~$\X_k$. Explicitly:
$$\sigma^*\colon \chi \mapsto \sigma^{-1}\circ \chi \circ
\sigma.$$Such action preserves the universal characters since $$
\sigma^*\bigl(\chi_{\P,i} \bigr)=\chi_{\P,i-1}.$$Indeed, if
$a\tensor 1 \in  O_L\tensor_\ZZ\bar{{\bf F}}_p$, then
$$\bigl(\sigma^{-1}\circ \chi_{\P,i}\circ \sigma\bigr) (a\tensor
1)=\bigl(\sigma^{-1}\circ \chi_{\P,i}\bigr) (a\tensor
1)=\chi_{\P,i-1}(a).$$For any $\chi\in \X_k$, we have
$$\chi^{(p)}=\Bigl(\sigma^*(\chi)\Bigr)^p.$$
\endrmk

\ssection Examples\par For $p$ inert, $\X_k$ is freely generated
by~$\{\chi_1,\ldots,\chi_g\}$ and~$\chi_i^{(p)}=\chi_{i-1}^{p}$,
where~$\chi_1^{(p)}=\chi_g^p$. If~$p$ is totally ramified,
then~$\X_k$ is freely generated by one character~$\Psi$
and~$\Psi^{(p)}=\Psi^p$.
\endssection

\label wtV. proposition\par\prop Let\/~$g\in{\bf
M}(k,\mu_N,\psi)$, then $V\bigl(g\bigr)\in{\bf
M}(k,\mu_N,\psi^{(p)})$.
\endprop
\Proof We need to verify properties (I)-(III)
of~\refn{modularforms}. The first two clearly hold. For (III)
consider any $\underline{A}$ and any~$\omega$ as in~\refn{V} and
any~$\gamma\in \G_k(R)=\bigl(O_L\tensor_\ZZ R\bigr)^*$. Then, with
the notation $\gamma^{(p)}:=\F^{\rm abs}_k(\gamma)$, we have
$$\eqalign{ V\bigl(g\bigr)\bigl(\underline{A},\gamma^{-1}\omega\bigr)
 &=g\bigl(\underline{A}^{(p)},(\gamma^{(p)})^{-1}\omega^{(p)}\bigr)\cr
 &=\psi\bigl(\gamma^{(p)}\bigr)
 g\bigl(\underline{A}^{(p)},\omega^{(p)}\bigr)\cr
 &=\psi^{(p)}(\gamma)\Bigl(V\bigl(g\bigr)
 \bigl(\underline{A},\omega\bigr)\Bigr).\cr}$$

\noindent We want to compute the effect of~$V$ on $q$-expansions.
To do this we introduce an auxiliary operator.

\label AbSoFrOb. definition\par\defi Let $\sigma$ be the absolute
Frobenius on~$k$. Consider the induced self-equivalence on the
category of schemes over~$k$ given by $$S \llongmapsto
S^\sigma:=S\fibprod_{k,\sigma} k.$$
\spacing
\noindent Let\/~$g$ be a $\I$-polarized modular form of
weight~$\psi$ and level\/~$\mu_N$ over~$k$. Define the rule
$\sigma^*(g)$ by
$$\sigma^*(g)\bigl(\underline{A},\omega\bigr):=
\Bigl(g\bigl(\underline{A}^\sigma,\omega^\sigma\bigr)\Bigr)^{\sigma^{-1}},$$where
$\underline{A}$ and\/~$\omega$ are as in~\refn{V}.
\enddefi

\label doweneed. lemma\par\lemma The rule $\sigma^*(g)$ defines a
$\I$-polarized modular form of weight\/~$\sigma^*(\psi)$ as
defined in~\refn{sigma*psi}. Moreover, if the $q$-expansion
of\/~$g$ at a $\I$-polarized $\FF$-rational unramified
cusp~$\bigl(\A,\B,\varepsilon,\j\bigr)$ is
$$g\bigl(\Tate(\A,\B),\varepsilon, \j\bigr)=\sum_\nu a_\nu q^\nu,$$then the
$q$-expansion of\/~$\sigma^*(g)$ at the same cusp is
$$\sigma^*(g)\bigl(\Tate(\A,\B),\varepsilon, \j\bigr)=\sum_\nu a_\nu^{\sigma^{-1}}\,
q^\nu.$$See~\refn{cusp} and~\refn{qexpansion} for the notions of
cusps and $q$-expansions.
\endlemma
\Proof To prove that $\sigma^*(g)$ defines a modular form, it
suffices to prove (III) of~\refn{modularforms}. For any
$\I$-polarized Hilbert-Blumenthal abelian scheme $\underline{A}$
over~$R$ with $\mu_N$-level structure and for any
generator~$\omega$ of the relative differentials and for
any~$\gamma\in \bigl(O_L\tensor_\ZZ R\bigr)^*$, we have:
$$\eqalign{ \sigma^*(g)\bigl(\underline{A},\gamma^{-1}\omega\bigr)
 &=\Bigl(g\bigl(\underline{A}^\sigma,(\gamma^\sigma)^{-1}\omega^\sigma\bigr)
 \Bigr)^{\sigma^{-1}}\cr
 &=\Bigl(\psi\bigl(\gamma^\sigma\bigr)\Bigr)^{\sigma^{-1}}
 \Bigl(g\bigl(\underline{A}^\sigma,\omega^\sigma\bigr)\Bigr)^{\sigma^{-1}}\cr
 &=\Bigl(\sigma^*\bigl(\psi\bigr)(\gamma)\Bigr)\Bigl(\sigma^*(g)
 \bigl(\underline{A},\omega\bigr)\Bigr).\cr}$$For the assertion on
$q$-expansions note that
$$\eqalign{\sigma^*(g)\bigl(\Tate(\A,\B),\varepsilon, \j\bigr)
 & =\Bigl(g\bigl( \Tate(\A,\B)^\sigma,\varepsilon^\sigma,
 (\j)^\sigma\bigr)\Bigr)^{\sigma^{-1}}\cr
 & =\Bigl(g\bigl(\Tate(\A,\B),\varepsilon,\j
 \bigr)\Bigr)^{\sigma^{-1}},\cr}$$since Tate objects are defined
over~$\FF$.

\label qexpV. proposition\par\prop The notation is as in~\refn{V}.
Fix a $\I$-polarized unramified $\FF$-rational
cusp~$\bigl(\B,\A,\varepsilon,\j\bigr)$. Suppose that\/~$g$ has
$q$-expansion
$$g\bigl(\Tate(\A,\B),\varepsilon, \j\bigr)=\sum_\nu a_\nu q^\nu,$$then the
$q$-expansion of\/~$V(g)$ at the same cusp is
$$V(g)\bigl(\Tate(\A,\B),\varepsilon, \j\bigr)=\sum_\nu a_\nu\,
q^{p\nu}.$$
\endprop
\Proof We use the notations of~\refn{tateobjects}. Let
$k\bigl(\bigl(\A,\B,\sigma_\beta\bigr)\bigr):=\ZZ\bigl(
\bigl(\A,\B,\sigma_\beta\bigr)\bigr)\tensor_\ZZ k$. We can factor
the absolute Frobenius~$\F^{\rm abs}$ on~$
k\bigl(\bigl(\A,\B,\sigma_\beta\bigr)\bigr)$ as follows:
$$\matrix{
k[\![q^\nu]\!]& \maprighto{\sigma} & k[\![q^\nu]\!]&
\maprighto{\xi}&  k[\![q^\nu]\!]\cr
 \sum_\nu a_\nu q^\nu & \mapsto & \sum_\nu a_\nu^p q^\nu \cr
     & & \sum_\nu b_\nu q^\nu & \mapsto & \sum_\nu b_\nu
    q^{p\nu}.\cr}$$Note that $\xi$ is a homomorphism of~$k$-algebras. In
particular,

$$\eqalignno{ V(g)\bigl(\Tate(\A,\B),\varepsilon,\j\bigr)
 &= g\bigl(\Tate(\A,\B)^{(p)},\varepsilon^{(p)},(\j)^{(p)} \bigr)\cr
 &=g\left(\bigl((\Tate(\A,\B),\varepsilon)^\sigma\bigr)^\xi,\bigl((\j)^{\sigma}
 \bigr)^\xi\right)\cr
 &=\xi\left(g\bigl(\Tate(\A,\B)^\sigma,\varepsilon^\sigma,(\j)^\sigma\bigr)\right)&(\refn{modularforms}(II))\cr
 &=\xi\circ
 \sigma\left(\sigma^*\bigl(g\bigr)\bigl(\Tate(\A,\B),\varepsilon,\j\bigr)\right)&(\refn{AbSoFrOb})\cr
 &=\xi \circ \sigma\left( \sum_\nu a_\nu^{\sigma^{-1}}\, q^\nu\right)&(\refn{doweneed})\cr
 &=\sum_\nu a_\nu\, q^{p\nu}.\cr}$$

\label KatzV. section\par\ssection Katz's $V$ operator\par
In~[\Katzzzz, \S1.11.21] one finds a more general notion of
Frobenius operator on Katz's $\I$-polarized $p$-adic Hilbert
modular forms in the sense of~\refn{Katzpadicmodfor}
$$V\colon {\bf M}(R,\mu_N,\chi)^{p-{\rm adic}}\llongrightarrow {\bf M}(R,\mu_N,\chi)^{p-{\rm adic}},$$where
$R$ is as in~\refn{NNNotation}. It is defined as follows.

\noindent Let $f=\{f_n\}_n$ be a $\I$-polarized  $p$-adic Hilbert
modular form of level~$\mu_N$ and weight~$\chi$ over~$R$. The
sequence of functions $\{f_n\in \Gamma(\MM(n,n),O_{\MM(n,n)})\}_n$
is a compatible sequence in the sense of~\refn{padicsequence}. For
$n,m\in \NN$ and let $$ \bigl(\AU,\iota^{\rm U},\lambda^{\rm
U},\varepsilon_{Np^n}^{\rm U}\bigr) \llongrightarrow \MM(m,n)$$be
the universal $\I$-polarized abelian scheme with real
multiplication by~$O_L$ and $\mu_{Np^n}$-level structure. Fix
$n\geq m$. Define
$$\pi\colon \AU \llongrightarrow (\AU)':=\AU/\bigl(\DL^{-1}
\tensor_\ZZ \mu_p\bigr).$$Then $(\AU)'$ inherits a  canonical
$O_L$-action $(\iota^{\rm U})'$ and a canonical $\mu_N\times
\mu_{p^{n-1}}$-level structure~$(\varepsilon_{Np^n}^{\rm U})'$.
By~[\Katzzzz, Lemma 1.11.6] it inherits a canonical polarization
data~$\lambda'$ and it is characterized by the property that for
$m=1$ we have $$ \bigl((\AU)',(\iota^{\rm U})',(\lambda^{\rm
U})',(\varepsilon_{Np^n}^{\rm
U})'\bigr)\cong\left((\AU)^{(p)},(\iota^{\rm
U})^{(p)},(\lambda^{\rm U})^{(p)},(\varepsilon_{Np^n}^{\rm
U})^{(p)}\right),$$where the latter is the abelian scheme obtained
by pull-back via the absolute Frobenius
on~$\MM\bigl(k,\mu_{Np^n}\bigr)$. In particular, for any $n\geq 1$
there exists a unique morphism of schemes over~$R/\m^n$
$$F_n\colon \MM(n,n) \llongrightarrow \MM(n,n-1)$$such that
$$\left((\AU)',(\iota^{\rm U})',(\lambda^{\rm
U})',(\varepsilon_{Np^n}^{\rm U})'\right)\cong\left(\AU
,\iota^{\rm U},\lambda^{\rm U},\varepsilon^{\rm
U}_{Np^{n-1}}\right)\fibprod_{\MM(n,n-1),\F_n} \MM(n,n).$$With the
notation of~\refn{MM(n,n)}, let $\alpha\in \Gamma_n$. By
construction the diagram
$$\matrix{ \MM(n,n) & \llongmaprighto{\alpha} & \MM(n,n)\cr
\mapdownl{F_n} & & \mapdownr{F_n}\cr \MM(n,n-1)
&\llongmaprighto{\alpha} &\MM(n,n-1)\cr}$$is commutative. It
follows that the sequence
$$f_{n-1} \circ
F_n\in\Gamma\bigl(\MM(n,n),O_{\MM(n,n)}\bigr)$$defines a unique
$p$-adic modular form of the {\it same} weight~$\chi$
$$V(f)\in {\bf M}(R,\mu_N,\chi)^{p-{\rm adic}}.$$Fix a $\I$-polarized unramified cusp
$(\A,\B,\varepsilon_{p^\infty N },\j_\varepsilon)$. By~[\Katzzzz,
\S1.11.23] the operator~$V$ changes the $q$-expansions at this
cusp according to the rule
$$a_0+\sum_\nu a_\nu q^\nu \llongrightarrow  a_0+\sum_\nu a_\nu
q^{p\nu}.$$
\endssection

\rmk Let~$f$ be a $p$-adic modular form \`a la Serre over~$F$ of level\/~$\mu_N$
(\refn{Serrepadic}). Then, $\Theta_{\P,i}^{[\jj]}(f)$ and~$V(f)$
are $p$-adic modular form \`a la Serre over~$F$ if~$f$ is either a
cusp form, or of weight~$\chi\in\X$, or of weight~$\Nm^z$
with~$z\in\ZZ_p$. This follows from~\refn{CComPPare}.
\endrmk

\lemma Let $f\in{\bf M}(k,\mu_N,\chi)$. We have
$$r\bigl(V(f)\bigr)=V\bigl(r(f)\bigr),$$where the $V$ on the LHS
is the one defined in~\refn{V}, while the one defined on the RHS
is Katz's~$V$ operator on regular functions on~$\MM(1,\infty)$ and
$$r\colon \dirsum_\psi{\bf M}\bigl(k,\mu_N,\psi\bigr) \llongrightarrow
\Gamma\bigl(\MM(1,1),O_{\MM(1,1)}\bigr)$$is the map $\sum g_\psi
\rightarrow \sum_\psi g_\psi/a(\psi)$ defined in~\refn{r(f)}.
\endlemma
\Proof  By~\refn{qexpa(chi)} the $q$-expansion of the~$a(\psi)$'s
can be assumed to be~$1$. The lemma follows by comparing the
effect on $q$-expansions of Katz's~$V$ operator and our~$V$
operator; see~\refn{qexpV}.

\rmk Note that if $f\in{\bf
M}(k,\mu_N,\chi)$, then $V(f)\in {\bf M}(k,\mu_N,\chi^{(p)})$.
Note also that $\chi^{(p)}\chi^{-1}\in \X_k(1)$ in the notation
of~\refn{XB}. In particular, $r\bigl(V(f)\bigr)$ and
$V\bigl(r(f)\bigr)$ have the same weight as functions
on~$\MM(1,1)$. This also justifies  our approach to the
operator~$V$; it allows to control how the weight changes
(see~\refn{wtV}).
\endrmk

\label VabsFrob. remark\par\rmk Recall that $\MM(k,\mu_N)$ is
defined over~$\FF$. Let~$\F^{\rm abs}$ be the absolute Frobenius
morphism on~$\MM(\FF,\mu_N)$ and let~$\F_k^{\rm abs}$ be its
base-change to~$k$. Using the moduli property of~$\MM(k,\mu_N)$,
one sees that if~$g$ is a modular form of level~$\mu_N$ and
weight~$1$ over~$k$ i.~e., $g \in
\Gamma\bigl(\MM(k,\mu_N),O_{\MM(k,\mu_N)} \bigr)$, then
$V(g)=\F_k^{\rm abs,*}(g)$.

\endrmk
\endsection

\section The operator $U$\par \noindent  We use previous notation and denote
$$\MM:=\MM(1,1)^\Kum \qquad\hbox{{\rm and}}\qquad
\MM':=\MM\bigl(1,0\bigr).$$Then $\phi\colon \MM\rightarrow \MM'$
is a Galois cover with group
$$G\cong \prod_{\P\vert (p)}\bigl(O_L/\P\bigr)^*.$$

\label Thetapsi. definition\par \defi Let $\P$ be a prime of~$O_L$
over~$p$. Let~$0 \leq \jj\leq e_\P-1$ be an integer. Define an
operator $$\LPj\colon \Gamma\bigl(\MM,O_\MM^{\P,[\jj]}\bigr)
\llongrightarrow
\Gamma\bigl(\MM,O_\MM^{\P,[\jj+1]}\bigr)$$(see~\refn{theta} for
the notation, extended by the same formula for~$\jj+1=e_\P$) as
follows. Choose
$\psi=\prod_{\P,i}\chi_{\P,i}^{a_{\P,i}}\in\X_k(1)$, in the
notation of\/~\refn{XB}, with~$a_{\P,i}\geq 0$. Let
$$\LPj:={\rm
Id}-\prod_{i=1}^{f_\P}\bigl(\Theta_{\P,i}^{[\jj]}\bigr)^{a_{\P,i}}.
$$
\enddefi

\label LPj. lemma\par\lemma The operator $\LPj$ is well defined.
\spacing
\item{{\rm 1.}} It has the following effect on
$q$-expansions. If $f=a_0+\sum_{\nu\in(\P^\jj\A\B)^+} a_\nu q^\nu$
at the $\I$-polarized unramified
cusp~$\bigl(\A,\B,\varepsilon_{pN},\j_\varepsilon\bigr)$, then
$$\LPj(f)=\sum_{\nu\in(\P^{\jj+1}\A\B)^+} a_{\nu} q^{\nu}$$at this cusp;
\spacing
\item{{\rm 2.}} let $\chi\colon G \rightarrow k^*$ be a character and
let $f$ be an eigenfunction for~$G$ with character~$\chi$. Then
$\LPj(f)$ is also an eigenfunction for~$G$ with character~$\chi$.

\endlemma

\Proof Choose $\psi$ as in~\refn{Thetapsi}. We calculate the
effect of~$\LPj$ on $q$-expansions using that the operators
$\Theta_{\P,i}$, for different~$i$'s, commute
by~\refn{ThetaPicommute}. Suppose $f \in
\Gamma(\MM,O_\MM^{\P,[\jj]})$ i.~e., by~\refn{qexptheta},
that~$a_\nu=0$ if~$\nu \notin \P^\jj \A\B$. It follows
from~\refn{qexptheta}  that the $q$-expansion of $$\Bigl({\rm
Id}-\prod_{i=1}^{f_\P}\bigl(\Theta_{\P,i}^{[\jj]}\bigr)^{a_{\P,i}}\Bigr)(f)
$$at the given cusp is $$\sum_\nu \Bigl(1-\bigl(\prod_{i=1}^{f_\P}(\tilde{\chi}_{\P,i}^{[\jj]}
)^{a_{\P,i}}\bigr)(\nu)\Bigr) a_\nu q^{\nu}=\sum_{\nu\in
(\P^{\jj+1}\A\B)^+} a_\nu q^\nu.$$The last equality follows since
by~\refn{sigmatilde}
$$\prod_{i=1}^{f_\P}\bigl(\tilde{\chi}_{\P,i}^{[\jj]}(\nu)\bigr)^{a_{\P,i}}=\cases{ 0 & if $\nu\in
\P^{\jj+1}\A\B$, \cr 1 & otherwise.\cr}$$In particular,
by~\refn{qexptheta}
$$\LPj(f)\in
\Gamma\left(\MM,O_\MM^{\P,[\jj+1]}\right)$$and the definition
of~$\LPj$ does not depend on the choice of~$\psi$. Part~(2)
follows from~\refn{wttheta}.

\label L. definition\par\defi Define an operator $$\Lambda\colon
\Gamma\bigl(\MM,O_\MM\bigr) \llongrightarrow
\Gamma\bigl(\MM,O_\MM\bigr)$$by
$$\Lambda:=\circ_{\P\vert p}
\biggl(\Lambda\bigl(\P,e_\P-1\bigr)\circ\cdots\circ\Lambda\bigl(\P,0\bigr)
\biggr).$$
\enddefi

\label qexpLambda. proposition\par\prop The operator $\Lambda$ is
well defined.
\spacing
\item{{\rm 1.}} It has the following effect on
$q$-expansions. If $f=\sum_\nu a_\nu q^\nu$ at the $\I$-polarized
unramified
cusp~$\bigl(\A,\B,\varepsilon_{pN},\j_\varepsilon\bigr)$, then
$$\Lambda(f)=\sum_\nu a_{p\nu} q^{p\nu}$$at this cusp;
\spacing
\item{{\rm 2.}} let $\chi\colon G \rightarrow k^*$ be a character and
let $f$ be an eigenfunction for~$G$ with character~$\chi$. Then
$\Lambda(f)$ is also an eigenfunction for~$G$ with
character~$\chi$.

\endprop
\Proof It follows from~\refn{LPj}.

\label VeqLambda. proposition\par\prop Let $f$ be a regular
function on~$\MM$. There exists a {\it unique} regular
function~$a$ on~$\MM$ such that
$$V(a)=\Lambda(f).$$Moreover, if $f$ is an eigenfunction for~$G$
with character~$\chi$, then~$a$ is also an eigenfunction for~$G$
with character~$\chi$.

\endprop
\Proof  Let~$K$ be the function field of~$\MM$. We have a
commutative diagram:
$$\matrix{
K & \maprighto{d} & \Omega^1_{K/k}\cr
 \Big\downarrow & &\Big\downarrow \cr
 k\bigl(\bigl(\A,\B,\sigma_\beta \bigr)\bigr) & \maprighto{d}
 & \Omega^1_{k\bigl(\bigl(\A,\B,\sigma_\beta \bigr)\bigr)/k}.\cr}$$

\noindent The left vertical arrow is injective being a morphism of
fields. By~[\Rapoport, Thm.~5.1] the
schemes~$S_{\sigma_\beta}^\wedge$, defined in~\refn{tateobjects},
for a suitable choice of the cone
decomposition~$\{\sigma_\beta\}_\beta$  appear as formal
completions of boundary components of smooth toroidal
compactifications of~$\MM$. In particular,
$$\Omega^1_{K/k} \hooklongrightarrow
\Omega^1_{k\bigl(\bigl(\A,\B,\sigma_\beta \bigr)\bigr)/k}$$is
injective. By~\refn{qexpLambda} we conclude that
$$\Omega^1_{K/k}\ni d\left(\Lambda(f)\right)=0.$$We conclude that
$$\Lambda(f)\in
\Gamma\bigl(\MM,O_\MM\bigr)\cap K^pk;$$see~[\Lang, Ch.~X,
Thm.~7.4]. Since~$\MM$ is normal and affine
$$\Gamma\bigl(\MM,O_\MM\bigr)\cap
K^pk:=\Gamma\bigl(\MM,O_\MM\bigr)^pk.$$We conclude
by~\refn{VabsFrob}.

\label Upsi. definition\par \defi  Let\/~$f$ be a~regular function
on~$\MM$. Define~$U(f)$ to be the unique~regular function on~$\MM$
such that
$$V\Bigl(U\bigl(f\bigr)\Bigr)=\Lambda(f).$$Its existence
and uniqueness is guaranteed by~\refn{VeqLambda}.
\enddefi

\label U. theorem\par\label qexpU. corollary\par\thm  Let ${\bf
M}\bigl(k,\mu_N,\chi\bigr)^{\rm ord}$ be the space of\/
$\I$-polarized modular forms over~$k$ of weight\/~$\chi$  defined
on\/~$\MM(k,\mu_N)^{\rm ord}$. There exists a (unique) $k$-linear
operator $$U\colon {\bf M}\bigl(k,\mu_N,\chi\bigr)^{\rm
ord}\llongrightarrow {\bf M}\bigl(k,\mu_N,\chi\bigr)^{\rm
ord}$$with the following effect on $q$-expansions. Let $f\in{\bf
M}\bigl(k,\mu_N,\chi \bigr)^{\rm ord}$. Suppose that its
$q$-expansion at a $\I$-polarized unramified
cusp\/~$\bigl(\B,\A,\varepsilon_{pN},\j_\varepsilon\bigr)$ is
$$ f\bigl(\Tate(\A,\B),\varepsilon_{pN},\j_\varepsilon \bigr)=\sum_\nu a_\nu
q^\nu.$$Then the $q$-expansion of\/~$U(f)$ at the same cusp is
$$U(f)\bigl(\Tate(\A,\B),\varepsilon_{pN},\j_\varepsilon \bigr)=\sum_\nu a_{p\nu}\,
q^{\nu}.$$Moreover, if $\chi=\prod_{\P,i}\chi_{\P,i}^{a_{\P,i}}$
and $a_{\P,i}\geq 2$ for all primes~$\P$  and all integers $1\leq
i\leq f_\P$, then
$$f\in {\bf M}\bigl(k,\mu_N,\chi\bigr) \quad\Longrightarrow\quad
U(f)\in {\bf M}\bigl(k,\mu_N,\chi\bigr).$$
\endthm
\Proof Let~$U$ be the operator on functions on~$\MM$ defined
in~\refn{Upsi}. Let
$$U(f):=U\bigl(r(f) \bigr) a(\chi),$$where $r(f)$ is the
regular function on~$\MM$ associated to~$f$ as in~\refn{r(f)}
and~$a(\chi)$ is the modular form on~$\MM$ defined
in~\refn{actiona(chi)}. In particular, $U(f)$ is a modular form
on~$\MM'$. By~\refn{actiona(chi)} it descends to a modular form on
the ordinary locus~$\MM'=\MM(k,\mu_N)^{\rm ord}$
of~$\MM(k,\mu_N)$. It has weight~$\chi$.  By~\refn{qexpLambda},
\refn{Upsi} and~\refn{qexpV}, we conclude that the effect on
$q$-expansions is as claimed in the theorem.

\indent By comparing weights and $q$-expansions we conclude the
following equality of modular forms

\advance\ssnu by-1 \labelf formulaU\par$$h^{p+1}
V\bigl(U(f)\bigr)=h_\chi\cdot \prod_\P\left(\prod_{\jj=0
}^{e_\P-1}\Bigl(\prod_{i=1}^{f_\P}
h_{\P,i}^{p+1}-\prod_{i=1}^{f_\P}\bigl(\Theta_{\P,i}^{[\jj]}\bigr)^{p-1}\Bigr)\right)
\bigl(f\bigr),\eqno{(\numfo)}$$\advance\fonu by1\advance\ssnu by1

\noindent where $h_{\P,i}$ (resp.~$h_\chi$) are the partial Hasse
invariants defined in~\refn{hPi}.

\noindent Assume that $f\in{\bf M}(k,\mu_N,\chi)$ with
$\chi=\prod_{\P,i}\chi_{\P,i}^{a_{\P,i}}$ and $a_{\P,i}\geq 2$ for
all primes~$\P$  and all integers $1\leq i\leq f_\P$. Consider the
equality of meromorphic modular forms:
$$\prod_\P\Bigl(\prod_{i=1}^{f_\P} h_{\P,i}^{p+1-a_{\P,i}}\Bigr)
\cdot V\bigl(U(f)\bigr)= \prod_\P\left(\prod_{\jj=0
}^{e_\P-1}\Bigl(\prod_{i=1}^{f_\P}
h_{\P,i}^{p+1}-\prod_{i=1}^{f_\P}\bigl(\Theta_{\P,i}^{[\jj]}\bigr)^{p-1}\Bigr)\right)
\bigl(f\bigr).$$The poles of~$U(f)$ are supported on the
complement of~$\MM(k,\mu_N)^{\rm ord}$ in~$\MM(k,\mu_N)$, which is
the union  of the distinct reduced divisors~$W_{\P,i}$ defined
by~$h_{\P,i}$. It follows from~\refn{VabsFrob} that~$V$ increases
the poles by~$p$. On the other hand the order of vanishing
of~$\prod_\P\Bigl(\prod_{i=1}^{f_\P}
h_{\P,i}^{p+1-a_{\P,i}}\Bigr)$  along any~$W_{\P,i}$ is at
most~$p-1$. Hence, if~$U(f)$ is not holomorphic, so is the LHS.
But, the modular form  on the RHS of the equality has no poles
by~\refn{ThetaPihasnopoles}.

\endsection

\section Applications to filtrations of modular forms\par \noindent Let
$f\in{\bf M}\bigl(k,\mu_N,\chi\bigr)$ be a $\I$-polarized modular
form. Let $\Phi(f)=\prod_{\P,i}\chi_{\P,i}^{a_{\P,i}}$
with~$a_{\P,i}\in \ZZ$ be its filtration as defined
in~\refn{filtrations}.

\ssection A summary\par Let~$\bigl(\A,\B,\varepsilon,\j \bigr)$ be
an unramified cusp in the sense of\/~\refn{cusp}. Suppose that the
$q$-expansion of\/~$f$ at the given cusp is
$$f\Bigl(\Tate(\A,\B),\varepsilon,\j\Bigr)=a_0 +\sum_{\nu\in \A\B^+}
a_\nu q^\nu\in k[\![q^\nu]\!]_{\nu \in \{0\}\cup(\A\B)^+}.$$In the
following table we summarize the effect of various operators on
the weight\/~$\chi$ and the $q$-expansion of\/~$f$; the lower part
of the table should be understood on the level of $q$-expansions
only.
\bigskip
\centerline{\vbox{\offinterlineskip
 \hrule
 \halign{&\vrule#&
  \strut\quad\hfil#\quad\cr
  height2pt&\omit&&\omit&&\omit&\cr
 & Operator\hfil&&Weight\hfil &&
 q-expansion\hfil &\cr
 height2pt&\omit&&\omit&&\omit&\cr
 \noalign{\hrule}
 height2pt&\omit&&\omit&&\omit&\cr
 & $\left[h_{\P,i} \right](f)$ \hfil &&
 $\chi \chi_{\P,i-1}^p\chi_{\P,i}^{-1}$\hfil&& $a_0+\sum_\nu a_\nu
 q^\nu$\hfil&\cr
 height2pt&\omit&&\omit&&\omit&\cr
 & $\Theta_{\P,i}^{[\jj]}(f)$\hfil &&
 $\chi \chi_{\P,i-1}^p\chi_{\P,i}$\hfil &&  $\sum_\nu
 \tilde{\chi}_{\P,i}^{[\jj]}(\nu) a_\nu q^\nu\hfil$&\cr
 height2pt&\omit&&\omit&&\omit&\cr
 &$V(f)$ \hfil&&  $\chi^{(p)}$\hfil &&  $a_0+\sum a_\nu
 q^{p\nu}$\hfil&\cr
 height2pt&\omit&&\omit&&\omit&\cr
 & $U(f)$ \hfil&& $\chi$ \hfil &&  $a_0
 +\sum_\nu a_{p\nu} q^\nu$ \hfil&\cr
 height2pt&\omit&&\omit&&\omit&\cr
 \noalign{\hrule}
 height2pt&\omit&&\omit&&\omit&\cr
 & $\LPj(f)$ && $\chi$ mod~$\X_k(1)$ && $a_0+\sum_{\nu\in(\P^{\jj+1}\A\B)^+} a_{\nu}
 q^{\nu}$\hfil&\cr
 height2pt&\omit&&\omit&&\omit&\cr
 & $\Lambda(f)$\hfil && $\chi$ mod~$\X_k(1)$ && $a_0+\sum a_{p\nu}
 q^{p\nu}$\hfil&\cr
 height2pt&\omit&&\omit&&\omit&\cr}
 \hrule}}
\bigskip
\noindent See~\refn{ThetaPi} and~\refn{qexpThetaPi} for the
operator~$\Theta_{\P,i}^{[\jj]}$,~\refn{V},~\refn{wtV}
and~\refn{qexpV} for the operator~$V$,~\refn{U}  for the
operator~$U$. Finally, $[h_{\P,i}]$ means multiplication by the
partial Hasse invariant\/~$h_{\P,i}$ introduced in~\refn{hPi}.
\endssection

\prop Let\/~$\bigl(\A,\B,\varepsilon,\j \bigr)$ be a
$\I$-polarized unramified cusp of level~$\mu_N$.
Let\/~$\widetilde{{\bf M}} (k,\mu_N)$ be the subring
of~$k\bigl(\bigl(\A,\B,\varepsilon,\j \bigr)\bigr)$ generated by
the $q$-expansions at this cusp of modular forms of
level\/~$\mu_N$ and any weight. Define the following differential
operator of degree~$(p-1)\,g$ on~$\widetilde{{\bf M}} (k,\mu_N)$
$$\Theta:= \Lambda-{\rm
Id};$$see~\refn{Thetapsi} for the notation. The following sequence
of $k$-vector spaces is exact:
$$0 \longrightarrow \widetilde{{\bf M}} (k,\mu_N)
\maprighto{V} \widetilde{{\bf M}} (k,\mu_N)
\llongmaprighto{\Theta}\widetilde{{\bf M}} (k,\mu_N)
\maprighto{U}\widetilde{{\bf M}} (k,\mu_N) \longrightarrow 0.$$
\endprop
\Proof By the table above one sees that $U \circ\Theta=0$ and
that~$\Theta \circ V=0$. One also concludes also that~$V$ is
injective and~$U\circ V=1$. In particular, $U$ is surjective. By
the formula~(\refn{formulaU}), we have $V \circ U=1-\Theta$. We
conclude that if $f\in \Ker\bigl(\Theta\bigr)$, then
$f=(1-\Theta)(f)=V\bigl(U(f)\bigr)$. Hence, $f\in{\rm Im}(V)$.
Finally, if~$f\in\Ker(U)$, then~$f=\sum_{(n,p)=1} a_\nu q^\nu$.
Hence, $f=\Theta(f)$ i.~e.,~$f\in{\rm Im}\bigl(\Theta\bigr)$.

\label vPi. definition\par\defi Let $g$ be a rational function
on~$\MM$. For every prime~$\P$ over~$p$ and any $1\leq i\leq
f_\P$, define
$$v_{\P,i}(g):={\rm min}\left\{v_C(g)\vert\, C\,\hbox{{\rm
an irreducible component of }}W_{\P,i}\right\}.$$
\enddefi

\bigskip
\noindent {\bf N.B.:} Below we use additive notation for character
groups in order to ease notation. Thus, for example,
$\X_k=\sum_{\P,i} \ZZ \chi_{\P,i}$ and not $\X_k=\prod_{\P,i}
\chi_{\P,i}^\ZZ$.

\label variouscones. definition\par\defi For any prime $\P$
above~$p$ and any $1\leq i\leq f_\P$, let
$$\psi_{\P,i} := \chi_{\P,i-1}^p\chi_{\P,i}^{-1}$$be the weight of
the partial Hasse invariant~$h_{\P,i}$ defined in~\refn{hPi}.
Define two positive cones in~$\X_k$ by
$$\X_k^+:=\sum_{\P,i}\QQ_{\geq
0}\,\psi_{\P,i}\qquad\hbox{{\rm and}}\qquad
\left(\X_{O_K}^+\right)_k=\sum_{\P,i}\QQ_{\geq
0}\,\chi_{\P,i},$$($\left(\X_{O_K}^+\right)_k$ is the image of the
cone~$\X_{O_K}^+$, defined in~\refn{GGm}, under the reduction
map). Finally, we define a partial ordering~$\leq_k$ on~$\X_k$ by
requiring that if\/~$\chi$,~$\psi$ are in~$\X_k$, $$\chi \leq_k
\psi \quad \Longleftrightarrow\quad\psi\, \chi^{-1}\in \X_k^+.$$
\enddefi

\label comparecones. lemma\par\lemma The elements
$\bigl\{\psi_{\P,i}\bigr\}_{\P,i}$ form a basis
of\/~$\X_k\tensor_\ZZ \QQ$. Moreover,
$$\left(\X_{O_K}^+\right)_k \subset \X_k^+\qquad\hbox{{\rm
and}}\qquad \sum_{\P,i}\ZZ \psi_{\P,i}=\X_k(1).$$See~\refn{XB} for
the definition of\/~$\X_k(1)$.
\endlemma
\Proof For any prime $\P$ the fundamental characters
$\bigl\{\chi_{\P,i}\bigr\}_i$ defined in~\refn{GGm} are expressed
in terms of $\bigl\{\psi_{\P,i}\bigr\}_i$ by the positive rational
$f_\P \times f_\P$-matrix
$${1\over p^{f_\P}-1}\left({\matrix{
  1 & p & p^2  &  \ldots & p^{f_\P-1} \cr
  p^{f_\P-1} & 1 &  p & \ldots & p^{f_\P-2} \cr
 \vdots &  &  &  & \vdots \cr
 p & p^2& p^3 & \ldots & 1
\cr}}\right).$$This proves the inclusion. The other statement is
Part~(3) of~\refn{qexphPi}.

\label powers. lemma\par\lemma For any positive integer~$k$ one
has, $\Phi(f^k) = \Phi(f)^k$.
\endlemma
\Proof We assume that $f$ has weight~$\Phi(f)$. We certainly have
$\Phi(f^k) \leq \Phi(f)^k$. Suppose that the inequality is strict.
Then, for some~$\P$ and~$i$, we have that~$f^k$, therefore~$f$,
vanishes on every component of~$W_{\P,i}$. Hence, $f/h_{\P,i}$ is
also holomorphic. This implies that~$\Phi(f)$ is strictly less
than the weight of~$f$. Contradiction.

\label positivity. proposition\par\prop Let $f$ be a non-zero
$\I$-polarized Hilbert modular form of level~$\mu_N$ and
weight~$\chi$. Then, $0\leq_k \chi$ i.~e.,
$\chi=\prod_{\P,i}\,\psi_{\P,i}^{b_{\P,i}}$ with\/~$b_{\P,i}
\in\QQ$ and\/~$b_{\P,i} \geq 0$.
\endprop
\Proof By the definition of filtration given in~\refn{filtrations}
the weight~$\chi$ is the filtration~$\Phi(f)$ multiplied by
positive multiples of the elements~$\psi_{\P,i}$. Hence, we may
assume without loss of generality that~$\Phi(f)=\chi$.
Let~$\chi=\prod_{\P,i}\,\psi_{\P,i}^{b_{\P,i}}$. By~\refn{powers}
we may assume that $b_{\P,i} \in \ZZ$ for every~$\P$ and any~$i$.
Assume that there exist~$\Q$ and~$j$ such that~$b_{\Q,j}<0$.
Consider the function $r(f) := f/\prod_{\P,i} h_{\P,i}^{b_{\P,i}}$
on~$\MM$. Since~$f$ and the~$h_{\P,i}$'s are holomorphic modular
forms
$$v_{\Q,j}\bigl(r(f)\bigr) = v_{\Q,j}\bigl(f\bigr)
- v_{\Q,j}\bigl(h_{\Q,j}^{b_{\Q,j}}) \geq -b_{\Q,j}> 0.$$This
implies that the modular form $f/h_{\Q,j}=r(f)\cdot
\prod_{(\P,i)\neq (\Q,j)} h_{\P,i}^{b_{\P,i}} \cdot
h_{\Q,j}^{b_{\Q,j}-1}$ is holomorphic, contradicting the
assumption that the weight of~$f$ is~$\Phi(f)$.

\ssection {\bf Question}\par Over the complex numbers, a non-zero
Hilbert modular form  has weight in~$\X_{O_K}^+$. In
characteristic~$p$ this is no longer true as the example of the
partial Hasse invariant shows. However, a notion of positivity is
retained; the weight belongs to~$\X_k^+$, which is a positive cone
depending on~$p$ and containing~$\left(\X_{O_K}^+\right)_k$. Note
though that the filtration of a partial Hasse invariant is the
trivial character~$1$ that lies in~$\left(\X_{O_K}^+\right)_k$. We
therefore ask: {\it is the filtration of any modular form
in~$\left(\X_{O_K}^+\right)_k$}?
\endssection

\label fpoler(f). proposition\par\prop Let $f$ be a $\I$-polarized
modular form of level~$\mu_N$ over~$k$. The filtration of~$f$ and
the poles of~$r(f)$ determine each other by the following
relation. If\/~$\Phi(f)=\prod_{\P,i} \chi_{\P,i}^{a_{\P,i}}$, then
$$v_{\P,i}\bigl(r(f)\bigr)=-v_{\P,i}\bigl(a(\chi)\bigr)=-\sum_{s=0}^{f_\P-1} a_{\P,i+s} p^s.$$
\endprop
\Proof Without loss of generality we may assume that the weight
of~$f$ is equal to its filtration~$\Phi(f)$. In particular, $f$
does not vanish identically along any~$W_{\P,i}$. Else,
by~\refn{redirr} the modular form~$f/h_{\P,i}$ is holomorphic of
weight strictly smaller than~$\Phi(f)$. Since for any irreducible
component~$C$ of~$W_{\P,i}$ we have
that~$v_{\P,i}\bigl(a(\chi)\bigr)=v_C\bigl(a(\chi)\bigr)$
and~$v_C\bigl(a(\chi)\bigr)\geq 0$ by~\refn{positivity}, we
conclude that
$$v_{\P,i}\bigl(r(f)\bigr)=-v_{\P,i}
\bigl(a(\chi)\bigr).$$By~\refn{Mbar} $$\eqalign{v_{\P,i}
\bigl(a(\chi)\bigr)& =\sum_{j=1}^{f_\P} a_{\P,j}
v_{\P,i}\bigl(a(\chi_{\P,j}) \bigr)\cr & =\sum_{s=0}^{f_\P-1}
a_{\P,i+s} p^s.\cr} $$

\label PhiTheta. proposition\par\prop Let\/~$\P$ be a prime
over~$p$ and let $1\leq i\leq f_\P$ and let $0\leq \jj\leq
e_\P-1$. Let\/~$\Theta_{\P,i}^{[\jj]}$ be the operator introduced
in~\refn{ThetaPi}. Suppose that $f\in{\bf
M}\bigl(k,\mu_N,\chi\bigr)_\P^{[\jj]}$. Then
$$\Phi\left(\Theta_{\P,i}^{[\jj]}(f)\right)\leq_k\Phi(f)
\chi_{\P,i-1}^p\chi_{\P,i},$$ with {\it equality} in the
direction~$\chi_{\P,i-1}^p\chi_{\P,i}$ if and only
if\/~$p\not\vert a_{\P,i}$.
\endprop

\Proof We may assume that $\chi=\Phi(f)$. The filtration
$\Phi\bigl(\Theta_{\P,i-1}^{[\jj]}(f)\bigr)$
of~$\Theta_{\P,i}^{[\jj]}(f)$ is less or equal to its weight
$\Phi(f)\chi_{\P,i-1}^p\chi_{\P,i}$. By definition
$$\Theta_{\P,i}^{[\jj]}(f)=\Theta_{\P,i}^{[\jj]}\bigl(r(f)\bigr) a\bigl(\chi\bigr)
a\bigl(\chi_{\P,i}^2\bigr) h_{\P,i};$$see~\refn{ThetaPi}.
By~\refn{fpoler(f)} we have
$$v_{\Q,j}\bigl(r(f)\bigr)=\sum_{s=0}^{f_\Q-1} a_{\Q,j+s} p^s.$$In
particular, $p\vert v_{\Q,j}\bigl(r(f)\bigr)$ if and only
if~$p\vert a_{\Q,j}$.  Let~$C$ be an irreducible component
of~$W_{\Q,j}$. If~$f$ vanishes along~$C$, then
$v_C\bigl(r(f)\bigr) \geq v_{\Q,j}\bigl(r(f)\bigr)+(p^{f_\Q}-1)$
by~\refn{casino}. Using~\refn{polethetaf}, we conclude that
$v_C\Bigl(\Theta_{\Q,j}^{[\jj]}\bigl(r(f)\bigr)\Bigr)\geq
v_{\Q,j}\bigl(r(f)\bigr)-2$, whether $p\vert v_C\bigl(r(f)\bigr)$
or not. Let~$C$ be a component of~$W_{\Q,j}$ along which~$f$ does
not vanish identically. Such exists by the assumption
that~$\chi=\Phi(f)$. At this component
$v_C\bigl(r(f)\bigr)=v_{\Q,j}\bigl(r(f)\bigr)$. Using loc.~cit.,
we find that
$v_C\Bigl(\Theta_{\Q,j}^{[\jj]}\bigl(r(f)\bigr)\Bigr)$ is greater
or equal to~$ v_{\Q,j}\bigl(r(f)\bigr)-2$ if~$p\vert a_{\Q,j}$ and
is equal to~$v_{\Q,j}\bigl(r(f)\bigr)-2-(p^{f_\Q}-1)$ otherwise.
Taking the minimum over all components~$C$ of~$W_{\Q,j}$ and using
loc.~cit., we obtain
$$v_{\Q,j}\left(\Theta_{\P,i}^{[\jj]}\bigl(f\bigr)\right)\cases{
 \,\geq  v_{\Q,j}\bigl(f\bigr) & if $\P\neq\Q$;\cr
 \,\geq v_{\Q,j}\bigl(f\bigr)-2 p^{f_\P-r}  & if $\P=\Q$ and  $i\neq j$;\cr
 \,\geq v_{\Q,j}\bigl(f\bigr)-2  & if $\P=\Q$, $i=j$ and
 $p\vert a_{\Q,j}$;\cr
 \,=v_{\Q,j}\bigl(f\bigr)-2-\bigl( p^{f_\P}-1\bigr) & if $\P=\Q$,
 $i=j$ and $p\,\not\vert \,a_{\Q,j}$.
\cr}$$By~\refn{fpoler(f)} we have
$v_{\Q,j}\bigl(r(f)\bigr)=-v_{\Q,j}\bigl(a(\chi)\bigr)$. One
concludes the proof computing $v_{\Q,j}\bigl(a(\chi_{\P,i})\bigr)$
and~$v_{\Q,j}\bigl(h_{\P,i}\bigr)$.

\rmk Since
$\bigl(\Theta_{\P,i}^{[\jj]}\bigr)^p=\Theta_{\P,i+1}^{[\jj]}$ one
can not hope to strengthen~\ref{PhiTheta} i.~e., that  equality
holds if~$p\not\vert a_{\P,i}$.

\endrmk

\label PhiV. proposition\par\prop See~\refn{V} and\/~\refn{wtV}
for the definition of the operator~$V$ and its effect on weights.
We have $\Phi\bigl( V(f)\bigr)=\Phi(f)^{(p)}$.
\endprop
\Proof By definition $V$ is induced by Frobenius
on~$\MM(k,\mu_N)$. In particular, it increases the poles or the
zeroes of~$f$ by~$p$. Hence, if some partial Hasse invariant
divides~$V(f)$, it must divide~$f$ itself.

\defi We say that\/~$f$ is an {\it ordinary eigenform} if there exists~$\lambda\in k^*$
such that\/~$U(f)=\lambda f$.
\enddefi

\label PhiU. proposition\par\prop  Let $f$ be a $\I$-polarized
Hilbert modular form of level\/~$\mu_N$ and
filtration~$\Phi(f)=\prod_{\P,i}\chi_{\P,i}^{a_{\P,i}}$ such
that\/~$U(f)$ is also holomorphic e.~g., $a_{\P,i}\geq 2$ for
every~$\P$ and every~$i$. Then, $$\Phi\bigl(U(f)\bigr)^{(p)}\leq_k
\Phi(f) \Norm^{p^2-1},$$with strict inequality  if\/~$p$ is
ramified or~$ a_{\P,i}\equiv 1$ modulo~$p$ for some\/~$\P$
and\/~$i$.
\endprop
\Proof By~(\refn{formulaU}), we have
$$h^{p+1} V\bigl(U(f)\bigr)=\prod_\P\prod_{i=1}^{f_\P}
h_{\P,i}^{a_{\P,i}}\cdot \prod_\P\prod_{\jj=0
}^{e_\P-1}\Bigl(\,\prod_{i=1}^{f_\P}
h_{\P,i}^{p+1}-\prod_{i=1}^{f_\P}\bigl(\Theta_{\P,i}^{[\jj]}\bigr)^{p-1}\Bigr)
\bigl(f\bigr).$$By~\refn{PhiV} $V\bigl(U(f)\bigr)$ has filtration
$\Phi\bigl(U(f)\bigr)^{(p)}$. For each $\jj\in\NN$, each
prime~$\P$ and each $g\in{\bf M}\bigl(k,\mu_N,\chi\bigr)$ we have
$$\Phi\left(\Bigl(\prod_{i=1}^{f_\P}
h_{\P,i}^{p+1}-\prod_{i=1}^{f_\P}\bigl(
\Theta_{\P,i}^{[\jj]}\bigr)^{p-1}\Bigr)\bigl(g\bigr)\right)={\rm
max}\left\{ \Phi(g), \Phi\biggl(\,\prod_{i=1}^{f_\P}\bigl(
(\Theta_{\P,i}^{[\jj]})^{p-1}\bigr)\bigl(g\bigr)\biggr)\right\}.$$
By~\refn{PhiTheta} we have $$\Phi\Bigl(\,\prod_{i=1}^{f_\P}
\bigr(\Theta_{\P,i}^{[\jj]}\bigl)^{p-1}\bigl(g\bigr)\Bigr) \leq_k
\Phi(g) \prod_{i=1}^{f_\P}\chi_{\P,i}^{p^2-1}.$$This proves the
first part the proposition. By~\refn{PhiTheta} we have strict
inequality if there exist~$\P$,~$i$, $0\leq j < e_\P(p-1)$ such
that $a_{\P,i}+j\equiv 0$ mod~$p$ i.~e., either~$p$ is ramified
or~$a_{\P,i}\not\equiv 1$ mod~$p$.

\label Ueigenform. corollary\par\cor If $f$ is an ordinary
eigenform, there exists an ordinary eigenform~$g$ of
weight~$\chi=\prod_{\P,i}\chi_{\P,i}^{a_{\P,i}} $ such that~$f$
and~$g$ have the same $q$-expansions and $ a_{\P,i}\in \bigl[t,
\ldots, p+1\bigr]$, where~$t=1$ if\/~$p\neq 2$ and\/~$t=0$
if\/~$p=2$.
\endcor
\Proof Note that
$$\Phi\bigl(VU(f)\bigr)=\Phi\bigl(V(f)\bigr)=\Phi(f)^{(p)}.$$Without
loss of generality we may assume that~$f$ has
weight~$\chi=\Phi(f)$. From~\refn{PhiU} we get the inequality
$$\prod_{\P,i}\bigl(\chi_{\P,i-1}^p\chi_{\P,i}^{-1}\bigr)^{a_{\P,i}}=\Phi(f)^{(p)}\Phi(f)^{-1}\leq_k
\Norm^{p^2-1}=\prod_{\P,i}\bigl(\chi_{\P,i-1}^p\chi_{\P,i}^{-1}\bigr)^{p+1}.$$It
follows that each~$a_{\P,i}$ is less or equal to~$p+1$.

\noindent For every prime~$\P$ we can ensure, by multiplying by
suitable partial Hasse invariants, that $2\leq a_{\P,i}\leq p+1$
for $2\leq i\leq f_\P$. We claim that $a_{\P,1}> -\bigl(1+{2\over
p-1}\bigr)$. Indeed, using the matrix appearing in the proof
of~\refn{comparecones} and~\refn{positivity}, we have $$p^{f_\P-1}
a_{\P,1}+ a_{\P,2}+\ldots+p^{f_\P-2} a_{\P,f_\P}\geq 0.$$This
gives $0\leq p^{f_\P-1} a_{\P,1}+(p+1)(p^{f_\P-1}-1)(p-1)^{-1}$,
which implies the claim. Hence, if~$p\neq 2$, multiplying further
by~$h_{\P,2}$ (if needed), we get the result. If~$p=2$,
multiplying further by at most a second power of~$h_{\P,2}$, gives
this case.

\cor The notation is as in~\refn{PhiU}. If $a_{\P,i}> p+1$ for
some\/~$\P$ and\/~$i$, then $\Phi\bigl(U(f)\bigr)<_k\Phi(f)$.
\endcor
\Proof  By hypothesis we have $\Phi(f)^{(p)-1}=\prod_{\P,i}
\psi_{\P,i}^{a_{\P,i}}
 >_k \prod_{\P,i} \psi_{\P,i}^{p+1}=\Nm^{p^2-1}$;
see~\refn{variouscones} for the notation~$\psi_{\P,i}$. Therefore,
$\Phi(f)^{(p)}>_k\Phi(f)\,\Nm^{p^2-1}$. Using~\refn{PhiU} we get
that $\Phi(f)^{(p)}>_k\Phi\bigl(U(f)\bigr)^{(p)}$. Since the
operation $\chi \mapsto \chi^{(p)}$ induces a bijection on the
positive cone~$\X_k^+$ we conclude
that~$\Phi\bigl(U(f)\bigr)<_k\Phi(f)$.

\rmk The assumption $a_{\P,i}> p+1$ for all\/~$\P$ and\/~$i$, implies that
$\Phi\bigl(U(f)\bigr)<_k\Phi(f)$ with respect to
any~$\psi_{\P,i}$.
\endrmk

\endsection

\section Functorialities\par \noindent Let
$L_1 \subset L_2$ be an extension of totally real number fields.
Let $\I_i\subset L_i$ be fixed fractional ideals such
that~$\I_2=\I_1\tensor_{O_{L_1}}{\rm D}_{L_2/L_1}$. We add
subscripts~$1$ or~$2$ to the usual notations for the degrees
over~$\QQ$, the primes, the ramification indices, the degrees of
the residue fields, the associated moduli spaces and spaces of
modular forms.

\label Upsilon. defi\par\defi Let $N$ be an integer and let\/~$S$
be a scheme. We define a morphism $$\Upsilon\colon
\MM_1\bigl(S,\mu_N\bigr) \llongrightarrow
\MM_2\bigl(S,\mu_N\bigr),$$where $\MM_i\bigl(S,\mu_N\bigr)$ is the
moduli space with respect to $\I_i$-polarization. Consider an
abelian scheme~$\bigl(A_1,\iota_1,\lambda_1,\varepsilon_1\bigr)$
over a $S$-scheme $T$, with real multiplication~$\iota_1$
by~$O_{L_1}$, polarization type $\lambda_1\colon
(M_{A_1},M_{A_1}^+) \isomarrow (\I_1,\I_1^+)$ and $\mu_N$-level
structure $\varepsilon_1$; see~\refn{moduli}. Then
$$\Upsilon\bigl(A_1,\iota_1,\lambda_1,\varepsilon_1\bigr):=
\bigl(A_2,\iota_2,\lambda_2,\varepsilon_2 \bigr),$$where
$$A_2:=A_1\tensor_{O_{L_1}}{\rm D}_{L_2/L_1}^{-1} \rightarrow T$$is an abelian
scheme over~$T$ with real multiplication $$\iota_2:=\iota_1
\tensor {\rm id}\colon O_{L_2} \hooklongrightarrow \End_T\bigl(
A\tensor_{O_{L_1}}{\rm
D}_{L_2/L_1}^{-1}\bigr)=\End_T\bigl(A_2\bigr)$$and $\mu_N$-level
structure $$\varepsilon_2:=\varepsilon_1 \tensor {\rm id}\colon
\mu_N \tensor_\ZZ {\rm D}_{L_1}^{-1}\tensor_{O_{L_1}}{\rm
D}_{L_2/L_1}^{-1} \hooklongrightarrow A_1\tensor_{O_{L_1}}{\rm
D}_{L_2/L_1}^{-1}=A_2.$$Note that
$$\bigl( A_2\bigr)^\vee:=\Ext^1_T\bigl(
A_2,\GG_{m,T}\bigr)\lisomarrow
\Ext^1_T\bigl(A_1,\GG_{m,T}\bigr)\tensor_{O_{L_1}} O_{L_2}
\isomarrow A_1^\vee \tensor_{O_{L_1}}O_{L_2},$$where $\vee$
denotes the dual abelian scheme. Hence, we get the polarization
type
$$\lambda_2:=\lambda_1\tensor {\rm id}\colon
\bigl(M_{A_2},M_{A_2}^+\bigr) \isomarrow
\bigl(\I_2,\I_2^+\bigr)$$with $\I_2:=\I_1 \tensor_{O_{L_1}}{\rm
D}_{L_2/L_1}$.
\enddefi

\label TateUpsilon. section\par\ssection The extension of
$\Upsilon$ to the cusps\par Let $i=1$ or~$2$. Let
$\bigl(\A_i,\B_i\bigr)$ be two fractional ideals of~$L_i$ such
that $\A_i\B_i^{-1}=\I_i$. Fix a rational polyhedral cone
decomposition $\{\sigma_{i,\beta}\}_\beta$ of the dual cone to
$(\A_i\B_i)_{\RR}^+ \subset (\A_i\B_i)_{\RR}$ which is invariant
under the action of the totally positive units of~$O_{L_i}$ and
such that, modulo this action, the number of polyhedra is finite.
Let
$$S_i:=(\A_i\B_i)^{\vee}\tensor_\ZZ
\GG_{m,\ZZ}.$$We have constructed in~\refn{tateobjects}  Tate
objects
$$\Tate\bigl(\A_i,\B_i\bigr)_{\sigma_{i,\beta}}=\Bigl(\A_i^{-1}{\rm D}_{L_i}^{-1}
\tensor_\ZZ\GG_{m,S_{i,\sigma_{i,\beta}}^\wedge}/
\underline{q}(\B_i)\Bigr)\fibprod_{S_{i,\sigma_{i,\beta}}^\wedge}\bigl(
S_{i,\sigma_{i,\beta}}^\wedge\backslash
S_{i,\sigma_{i,\beta},0}\bigr).$$This is an abelian scheme with
real multiplication by~$O_{L_i}$ over the open subscheme  $
S_{i,\sigma_{i,\beta}}^\wedge\backslash S_{i,\sigma_{i,\beta},0}$
of\/~$S_{i,\sigma_{i,\beta}}^\wedge$. The latter is defined as the
spectrum of the ring obtained completing the affine scheme
$S_{i,\sigma_{i,\beta}}$ along the closed subscheme
$S_{i,\sigma_{i,\beta},0}=S_{\sigma_{i,\beta}}\backslash S_i$ with
reduced structure.

\spacing
\noindent  By~\refn{indiff} we have that \bigskip

$\Upsilon\Bigl(
\Tate\bigl(\A_1,\B_1\bigr)_{\sigma_{1,\beta}}\Bigr)\isomarrow$
$$\qquad\qquad \Bigl(\A_1^{-1}{\rm D}_{L_2}^{-1}
\tensor_\ZZ\GG_{m,S_{1,\sigma_{1,\beta}}^\wedge}/
\underline{q}\bigl(\B_1{\rm
D}_{L_2/L_1}^{-1}\bigr)\Bigr)\fibprod_{S_{1,\sigma_{1,\beta}}^\wedge}\bigl(
S_{1,\sigma_{1,\beta}}^\wedge\backslash
S_{1,\sigma_{1,\beta},0}\bigr).$$Assume that $\A_2=\A_1 O_{L_2}$
and~$\B_2=\B_1 {\rm D}_{L_2/L_1}^{-1}$.  The trace map defines a
$O_{L_1}$-linear, surjective homomorphism $1 \tensor
\Tr_{L_2/L_1}\colon \A_2\B_2 \llongrightarrow \A_1\B_1$ and,
hence, a closed immersion $S_1 \hooklongrightarrow S_2$. Choose
cone decompositions $\{\sigma_{i,\beta}\}_\beta$ for $i=1,2$ as
above such that~$\{\sigma_{1,\beta}\}_\beta$ is induced
by~$\{\sigma_{2,\beta}\}_\beta$ via~$1\tensor \Tr_{L_2/L_1}$ for
each index~$\beta$. We get an induced a closed immersion
$\rho_\beta\colon S_{1,\sigma_{1,\beta}}^\wedge
\hooklongrightarrow S_{2,\sigma_{2,\gamma}}^\wedge$ such that
$\rho_\beta^{-1}\bigl(S_{2,\sigma_{2,\beta},0}
\bigr)=S_{1,\sigma_{1,\beta},0} $ and
$$\Upsilon\Bigl(
\Tate\bigl(\A_1,\B_1\bigr)_{\sigma_{1,\beta}}\Bigr)=
\Tate\bigl(\A_2,\B_2\bigr)_{\sigma_{2,\beta}}\,
\fibprod_{S_{2,\sigma_{2,\beta}}}S_{1,\sigma_{1,\beta}}.$$

\endssection

\label Psi. section\par\ssection The canonical homomorphism on
weights\par There is a canonical injective homomorphism of
$\ZZ$-group schemes $\G_1 \llongrightarrow \G_2$ defined on
$R$-valued points by the inclusion $\bigl(O_{L_1}\tensor_\ZZ
R\bigr)^* \hooklongrightarrow \bigl(O_{L_2}\tensor_\ZZ R\bigr)^*$.
This induces for any scheme~$T$ a map of characters defined
over~$T$:
$$\Psi \colon \X_T\bigl(\G_2\bigr) \llongrightarrow
\X_T(\G_1\bigr).$$It follows from~\refn{universal} that it sends
fundamental (resp.~universal) characters to fundamental
(resp.~universal) characters.
\endssection

\label modformUpsilon. section\par\ssection The effect of
$\Upsilon$ on modular forms\par Let $f\in {\bf
M}_2\bigl(S,\mu_N,\chi\bigr)$ be a modular form over~$S$ of
level~$\mu_N$ and weight~$\chi$. Define
$$\Upsilon^*(f)\in{\bf M}_1\bigl(S,\mu_N,\Psi(\chi)\bigr)$$by
requiring that for any affine $S$-scheme $\Spec(R)$, any abelian
scheme~$A_1$ with RM by~$O_{L_1}$ and level~$\mu_N$ and any
generator~$\omega$ of~$\H^0\bigl(A_1,\Omega^1_{A_1/R}\bigr)$ as
$R\tensor_\ZZ O_{L_1}$-module
$$\Upsilon^*(f)\bigl(A_1,\iota_1,\lambda_1,\varepsilon_1,\omega\bigr):=
f\left(A_1\tensor_{O_{L_1}}{\rm D}_{L_2/L_1}^{-1},
\iota_2,\lambda_2,\varepsilon_2,\omega \tensor
1\right).$$See~\refn{Upsilon} for the notation.  Note that
$$\Omega^1_{A_1\tensor_{O_{L_1}}{\rm D}_{L_2/L_1}^{-1}/R}\isomarrow
\Omega^1_{A_1/R}\tensor_{O_{L_1}} O_{L_2}$$and, hence,
$\omega\tensor 1$ is a generator
of~$\Omega^1_{A_1\tensor_{O_{L_1}}{\rm D}_{L_2/L_1}^{-1}/R} $. One
checks that the definition is well posed.
\endssection

\label qexpUpsilon. lemma\par\lemma ({\it The effect of $\Upsilon$
on $q$-expansions})\enspace Let $f\in {\bf
M}_2\bigl(S,\mu_N,\chi\bigr)$ and let
$\bigl(\A,\B,\varepsilon,\j)$ be a $\I_1$-polarized unramified
cusp of~$\MM_1(S,\mu_N\bigr)$. Suppose that the $q$-expansion
of~$f$ at the cusp $\bigl(\A\,O_{L_2},\B\,{\rm
D}_{L_2/L_1}^{-1},\varepsilon,\j\bigr)$, in the sense
of~\refn{qexpansion}, is
$$f\bigl(\Tate(\A\,O_{L_2},\B\,{\rm D}_{L_2/L_1}^{-1}),\varepsilon,\j
\bigr)=a_0+\sum_{\nu\in (\A\B{\rm D}_{L_2/L_1}^{-1})^+}a_\nu
q^\nu.$$Then the $q$-expansion of~$\Upsilon^*(f)$ at the
cusp~$\bigl(\A,\B,\j,\varepsilon\bigr)$ is $$a_0+\sum_{\delta\in
(\A\B)^+} \Bigl(\sum_{\nu\vert
 \Tr_{L_2/L_1}(\nu)=\delta}a_\nu\Bigr) q^\delta.$$
\endlemma
\Proof We calculate from the definitions
$$\eqalign{\Upsilon^*(f)\bigl(\Tate(\A,\B),\varepsilon,\j\bigr)&
:=\Upsilon^*(f)\Bigl(\Tate(\A,\B)_{\sigma_{1,\beta}},
\varepsilon,{dt \over t} \Bigr)\cr
 &=f\Bigl(\Tate(\A,\B)_{\sigma_{1,\beta}}\tensor_{O_{L_1}}{\rm D}_{L_2/L_1}^{-1},
\varepsilon,{dt \over t} \tensor 1\Bigr)\cr
 &=f\Bigl(\Tate(\A O_{L_2},\B{\rm D}_{L_2/L_1}^{-1})_{\sigma_{2,\beta}}
\fibprod_{S_{2,\sigma_{2,\beta}},\rho_\beta}S_{1,\sigma_{1,\beta}},
\varepsilon,{dt \over t}\Bigr)\cr
 &=\rho_\beta\biggl(f\Bigl(\Tate(\A O_{L_2},\B{\rm D}_{L_2/L_1}^{-1}), \varepsilon,{dt \over
 t}\Bigr)\biggr) \cr
 &=\rho_\beta\Bigl( a_0+\sum_{\nu\in (\A\B{\rm D}_{L_2/L_1}^{-1})^+}a_\nu
 q^\nu\Bigr)\cr
 &=a_0+\sum_{\delta\in (\A\B)^+} \left(\sum_{\nu\vert
 \Tr_{L_2/L_1}(\nu)=\delta}a_\nu\right) q^\delta.\cr}$$

\ssection Compatibilities of $U$ and $V$ operators\par
For~$i=1$,~$2$ let~$U_i$ and~$V_i$ be the~$U$ and~$V$ operators on
the space of modular forms $\dirsum_{\chi\in\X_k(\G_i)}{\bf
M}_i\bigl(k,\mu_N,\chi)$ in characteristic~$p$ introduced
in~\refn{U} and~\refn{V}. From the behavior of these operators on
weights and $q$-expansions described in~\refn{qexpU},
in~\refn{wtV} and~\refn{qexpV}, we conclude that
$$U_1\circ\Upsilon^*=\Upsilon^*\circ U_2\qquad \hbox{{\rm and}}\qquad
V_1\circ\Upsilon^*=\Upsilon^*\circ V_2.$$
\endssection

\prop ({\it Compatibilities of $\Theta$ operators})\enspace  Let
$\P_1$ be a prime of~$O_{L_1}$ over~$p$ and let $1\leq i\leq
f_{\P_1}$. We have the following identity of differential
operators on the algebra~$\dirsum_\chi{\bf
M}_2\bigl(k,\mu_N,\chi\bigr)$:
$$ \Theta_{\P_1,i} \circ \Upsilon^* = \Upsilon^*
\circ \left(\sum_{\P_2\vert \P_1,j\vert i
}e_{\P_2/\P_1}\Theta_{\P_2,j} \right).$$Here $\sum_{\P_2\vert
\P_1,j\vert i }$ means summing over all primes~$\P_2$ of~$O_{L_2}$
over~$\P_1$ and all $1\leq j\leq f_{\P_2}$ such that the embedding
$\bar{\sigma}_{\P_2,j}\colon O_{L_2}/\P_2 \rightarrow k$ induces
the embedding~$\bar{\sigma}_{\P_1,i}$ on~$O_{L_1}/\P_1$. As
customary, $e_{\P_2/\P_1}$ denotes the ramification index
of~$\P_2$ relatively to~$\P_1$.

\endprop
\Proof To prove this identity it is enough to show that both sides
change the weight and the $q$-expansion of a modular form~$f\in
{\bf M}_2\bigl(k,\mu_N,\chi\bigr)$ in the same way.

\noindent For weights we argue as follows: the operator
$\Upsilon^* \circ \left(\sum_{\P_2\vert \P_1,j\vert i
}e_{\P_2/\P_1}\Theta_{\P_2,j} \right)$ is equal to
$\sum_{\P_2\vert \P_1,j\vert i }e_{\P_2/\P_1}\Upsilon^*\circ
\Theta_{\P_2,j}$. Therefore, it is enough to calculate the weight
of the modular form~$\bigl(\Upsilon^*\circ
\Theta_{\P_2,j}\bigr)(f)$. By~\refn{ThetaPi}, the modular
form~$\Theta_{\P_2,j}(f)$ has weight $\chi\cdot \chi_{\P_2,j-1}^p
\cdot \chi_{\P_2,j}$. By~\refn{modformUpsilon}, we conclude that
$\Upsilon\bigl(\Theta_{\P_2,j}(f)\bigr)$ has weight
$\Psi(\chi)\cdot \chi_{\P_1,i-1}^p\cdot \chi_{\P_1,i}$, which is
equal to the weight of~$\bigl(\Theta_{\P_1,i} \circ \Upsilon^*
\bigr)(f)$.

\noindent We compute the effect on  $q$-expansions. Let
$\bigl(\A,\B,\varepsilon,\j)$ be a $\I_1$-polarized unramified
cusp of~$\MM_1(k,\mu_N\bigr)$. Denote the $q$-expansion of~$f$ at
the $\I_2$-polarized unramified cusp $\bigl(\A O_{L_2},\B {\rm
D}_{L_2/L_1}^{-1},\varepsilon,\j)$ by~$a_0+\sum_\nu a_\nu q^\nu$.
By~\refn{qexpUpsilon} and~\refn{qexpThetaPi} the effect
of~$\Theta_{\P_1,i}$ on the $q$-expansion of~$\Upsilon^*(f)$ at
the cusp $\bigl(\A,\B,\varepsilon,\j\bigr)$ is
$$\eqalign{\Theta_{\P_1,i}\bigl(\Upsilon^*(f)\bigr)\bigl(\Tate(\A,\B),\varepsilon,\j\bigr)
 &=a_0+\sum_{\delta\in (\A\B)^+} \tilde{\chi}_{\P_1,i}(\delta)\biggl(\sum_{\nu\vert
 \Tr_{L_2/L_1}(\nu)=\delta}a_\nu\biggr) q^\delta \cr
 &= a_0+\sum_{\delta,\Tr(\nu)=\delta} \biggl(
\sum_{\P_2\vert \P_1,j\vert i }e_{\P_2/\P_1}
\tilde{\chi}_{\P_2,j}(\nu) a_\nu\biggr) q^\delta \cr
 &=\Upsilon^*\Bigl(\sum_{\P_2\vert \P_1,j\vert i
}e_{\P_2/\P_1}\Theta_{\P_2,j}(f)\Bigr).\cr}$$

\label Hilbertvsellitpic. corollary\par\cor Let\/~$\Upsilon^*$ be
the homomorphism from Hilbert modular forms over~$k$ of
level\/~$\mu_N$ w.r.t.~$O_L$ to elliptic modular forms over~$k$.
Then,
$$\Theta \circ \Upsilon^* = \Upsilon^*
\circ \left(\,\sum_{\P\vert p,1\leq i\leq f_\P
}e_{\P}\Theta_{\P,i} \right),$$where $\Theta$ is the classical
theta operator of Serre and Swinnerton-Dyer.
\endcor

\endsection

\bigskip
\noindent {\bf References.}

\spacing

\item{[\Chai]} Chai, C.-L.: Arithmetic minimal compactification of
the Hilbert-Blumenthal moduli spaces. Appendix to ``The Iwasawa
conjecture for totally real fields" by A.~Wiles. {\it Ann.~of
Math.} (2) 131 (1990), no.~3, 541--554.

\spacing

\item{[\DelignePappas]} Deligne, P., Pappas, G.:
Singularit\'es des espaces de modules de Hilbert, en les
caract\'eristiques divisant le discriminant. {\it Compositio
Math.} 90 (1994), no.~1, 59--79.

\spacing

\item{[\DeligneRibet]} Deligne, P., Ribet, K.~A.: Values of abelian
$L$-functions at negative integers over totally real fields. {\it
Invent.~Math.} 59 (1980), no.~3, 227--286.

\spacing
\item{[\Geer]} Van der Geer, G.: {\it Hilbert modular surfaces}.
Ergebnisse der Mathematik und ihrer Grenzgebiete (3), 16,
Springer-Verlag, 1988.

\spacing

\item{[\Goren]} Goren, E.~Z.: Hasse invariants for Hilbert modular varieties. {\it Israel
J.~Math.},  122 (2001), 157--174.

\spacing

\item{[\Gorenn]} Goren, E.~Z.: Hilbert modular forms modulo $p^m$ - the unramified case. CICMA
pre-print 1998-10. {\it J.~Number Theory}, to appear.

\spacing

\item{[\Gorennn]} Goren, E.~Z.: {\it Lectures on Hilbert modular varieties and modular forms}.
CRM Monograph Series, to appear.

\spacing

\item{[\GorenOort]} Goren, E.~Z., Oort, F.: Stratifications of
Hilbert modular varieties, {\it J.~Algebraic Geometry} 9 (2000),
111-154.

\spacing

\item{[\Gross]} Gross, B.~H.: A tameness criterion for Galois
representations associated to modular forms (mod $p$). {\it Duke
Math.~J.} 61 (1990), no.~2, 445--517.

\spacing
\item{[\SGA2]} Grothendieck, A.: Cohomologie locale des faisceaux
coh\'erents et th\'eor\`emes de Lefschetz locaux et globaux
(\SGA2). Advanced studies in pure mathematics. Masson \& CIE,
Paris. North-Holland, Amsterdam, 1962.

\spacing
\item{[\EGAII]} Grothendieck, A.: \'Etude globale \'el\'ementaire de quelques
classes de morphismes. \EGAII. {\it Publ. Math. IHES} 8 (1961).

\spacing
\item{[\EGAIVtwo]} Grothendieck, A.: \'Etude locale des sch\'emas et
des morphismes de sch\'emas. \break \EGAIVtwo. {\it Publ. Math.
IHES} 24 (1965).

\spacing
\item{[\EGAIVfour]} Grothendieck, A.: \'Etude locale des sch\'emas et
des morphismes de sch\'emas.\break \EGAIVfour. {\it Publ. Math.
IHES} 32 (1967).

\spacing

\item{[\Hida]}  Hida, H.: On $p$-adic Hecke algebras for ${\rm
GL}\sb 2$ over totally real fields. {\it Ann.~of Math.} (2) 128
(1988), no.~2, 295--384.

\spacing
\item{[\Katz]} Katz, N.~M.: $p$-adic properties of modular schemes
and modular forms. In {\it Modular functions of one variable, III}
(Proc. Internat. Summer School, Univ. Antwerp, Antwerp, 1972),
pp.~69--190. Lecture Notes in Mathematics 350. Springer-Verlag,
1973.

\spacing

\item{[\Katzz]} Katz, N.~M.: Higher congruences between modular
forms. {\it Ann.~of Math.} (2) 101 (1975), 332--367.

\spacing

\item{[\Katzzz]} Katz, N.~M.: A result on modular forms in
characteristic $p$. In {\it Modular functions of one variable, V}
(Proc.~Second Internat.~Conf., Univ.~Bonn, Bonn, 1976),
pp.~53--61. Lecture Notes in Math., Vol.~601, Springer, Berlin,
1977.

\spacing

\item{[\Katzzzz]} Katz, N.~M.: $p$-adic $L$-functions for CM fields.
{\it Invent.~Math.} 49 (1978), no. ~3, 199--297.

\spacing

\item{[\Messing]} Messing, W.: {\it The crystals associated to
Barsotti-Tate groups}, Lecture Notes in Mathematics 264,
Springer-Verlag, 1972.

\item{[\MilnorStasheff]} Milnor, J.~W.; Stasheff, J.~D.:
{\it Characteristic classes.} Annals of Mathematics Studies, No.
76. Princeton University Press, Princeton, N. J.; University of
Tokyo Press, Tokyo, 1974.

\spacing

\item{[\vanderPoorten]} Van der Poorten, A.~J., te Riele, H.~J.~J.,
Williams, H.~C.: Computer verification of the Ankeny-Artin-Chowla
conjecture for all primes less than $100\,000\,000\,000$.
Math.~Comp., to appear.

\spacing

\item{[\Serre1]} Serre, J.-P.: Formes modulaires et fonctions z\^eta $p$-adiques.
{\it Modular functions of one variable, III} (Proc. Internat.
Summer School, Univ. Antwerp, 1972), pp. 191--268. Lecture Notes
in Math., Vol. 350, Springer, Berlin, 1973.

\spacing

\item{[\Siegel]} Siegel, C.~L.: Lectures on the analytic theory of
quadratic forms (1934-35). Institute for Advanced Studies Lecture
Notes. Reprinted~1955.
\spacing

\item{[\Washington]} Washington, L.~C.: {\it Introduction to
Cyclotomic Fields}. Graduate Text in Mathematics 83.
Springer-Verlag, 1982.

\spacing

\item{[\Zink]} Zink, T.:  The display of a formal p-divisible group, to appear in
Asterisque. Available from
http://www.mathematik.uni-bielefeld.de/$\,_{\widetilde{
}}\;$zink/.

\bigskip
\bigskip

\indent D{\eightpoint IPARTIMENTO} {\eightpoint DI} M{\eightpoint
ATEMATICA} P{\eightpoint URA} {\eightpoint E} A{\eightpoint
PPLICATA}, U{\eightpoint NIVERSIT\'A} {\eightpoint DEGLI}\break
S{\eightpoint TUDI} {\eightpoint DI} P{\eightpoint ADOVA},
{\eightpoint VIA} B{\eightpoint ELZONI} 7, P{\eightpoint ADOVA}
35123, I{\eightpoint TALIA}

\indent {\it E-mail address}: {fandreat@math.unipd.it}

\bigskip
\indent D{\eightpoint EPARTMENT} {\eightpoint OF} M{\eightpoint
ATHEMATICS} {\eightpoint AND} S{\eightpoint TATISTICS},
M{\eightpoint C}G{\eightpoint ILL} U{\eightpoint NIVERSITY}, 805
S{\eightpoint HERBROOKE} S{\eightpoint TR.}, M{\eightpoint
ONTREAL} H3A 2K6, C{\eightpoint ANADA}

\indent {\it E-mail address}: {goren@math.mcgill.ca}

\bye